%% file: pdf.tex
\documentclass[12pt]{book}
\usepackage{amssymb}
\usepackage{fancyheadings}
\input{localmacros}
\input{diagrams}
\input{veraformat}

\makeindex

\begin{document}
\input{title}
%\input{zusammenfassung}
\input{vorwort}

%%%%%%%%%%%%%%%%%%%%%%%%%%%%%%%%%%%%%%%
%%%%%%%%% Inhaltsverzeichnis %%%%%%%%%%
%%%%%%%%%%%%%%%%%%%%%%%%%%%%%%%%%%%%%%%
\cleardoublepage
\tableofcontents

%%%%%%%%%%%%%%%%%%%%%%%%%%%%%%%%%%%%%%%%%%%%%%%%%%%%
%%%%%%%%% Introduction %%%%%%%%%%%%%%%%%%%%%%%%%%%%%
%%%%%%%%%%%%%%%%%%%%%%%%%%%%%%%%%%%%%%%%%%%%%%%%%%%%
\cleardoublepage
\pagestyle{fancyplain}
\headrulewidth0.2pt
\lhead[\fancyplain{}{\small\bf\sf\thepage}]%
{\fancyplain{}{\scriptsize\sf\rightmark}}
\rhead[\fancyplain{}{\scriptsize\sf\leftmark}]%
{\fancyplain{}{\small\bf\sf\thepage}}
\cfoot{}
\setcounter{page}{1}
\input{introduction}

%%%%%%%%%%%%%%%%%%%%%%%%%%%%%%%%%%%%%%%%%%%%%%%%%%%%
%%%%%%%%%%%%% Basics %%%%%%%%%%%%%%%%%%%%%%%%%%%%%%%
%%%%%%%%%%%%%%%%%%%%%%%%%%%%%%%%%%%%%%%%%%%%%%%%%%%%
\cleardoublepage
\chapter{Basic Definitions}\label{chapter.definitions}
\input{algebra}
\input{rewriting}
\input{groebner}
%\input{problems}
%%%%%%%%%%%%%%%%%%%%%%%%%%%%%%%%%%%%%%%%%%%%%%%%%%%%
%%%%%%%%%%%% Reduction Rings %%%%%%%%%%%%%%%%%%%%%%%
%%%%%%%%%%%%%%%%%%%%%%%%%%%%%%%%%%%%%%%%%%%%%%%%%%%%
\cleardoublepage
\chapter{Reduction Rings}\label{chapter.reduction.rings}
\input{reductionrings}
\cleardoublepage
\chapter{Function Rings}\label{chapter.function.rings}
\input{generalization_k}
\input{generalization_rr}
\input{generalization_z}

\input{rmodule}
\input{ideals}

%%%%%%%%%%%%%%%%%%%%%%%%%%%%%%%%%%%%%%%%%%%%%%%%%%%%
%%%%%%%%% Applications %%%%%%%%%%%%%%%%%%%%%%%%%%%%%
%%%%%%%%%%%%%%%%%%%%%%%%%%%%%%%%%%%%%%%%%%%%%%%%%%%%
\cleardoublepage
\chapter{Applications of Gr\"obner Bases}\label{chapter.applications}
\input{applications}
%\input{fglm}

%%%%%%%%%%%%%%%%%%%%%%%%%%%%%%%%%%%%%%%%%%%%%%%%%%%%
%%%%%%%%% Examples %%%%%%%%%%%%%%%%%%%%%%%%%%%%%%%%%
%%%%%%%%%%%%%%%%%%%%%%%%%%%%%%%%%%%%%%%%%%%%%%%%%%%%
%\cleardoublepage
%\chapter{Examples of Function Rings}\label{chapter.examples}
%\input{examples}

%%%%%%%%%%%%%%%%%%%%%%%%%%%%%%%%%%%%%%%%%%%%%%%%%%%%
%%%%%%%%% Conclusions %%%%%%%%%%%%%%%%%%%%%%%%%%%%%%
%%%%%%%%%%%%%%%%%%%%%%%%%%%%%%%%%%%%%%%%%%%%%%%%%%%%
\cleardoublepage
\input{conclusions}

%%%%%%%%%%%%%%%%%%%%%%%%%%%%%%%%%%%%%%%%%%%%%%%%%%%%
%%%%%%%%%%%%%%%%%% Index %%%%%%%%%%%%%%%%%%%%%%%%%%%
%%%%%%%%%%%%%%%%%%%%%%%%%%%%%%%%%%%%%%%%%%%%%%%%%%%%

%%%%%%%%%%%%%%%%%%%%%%%%%%%%%%%%%%%%%%%%%%%%%%%%%%%%
%%%%%%%%% Bibliography %%%%%%%%%%%%%%%%%%%%%%%%%%%%%
%%%%%%%%%%%%%%%%%%%%%%%%%%%%%%%%%%%%%%%%%%%%%%%%%%%%
\cleardoublepage
\bibliographystyle{alpha}

{\small
\bibliography{references}
}
\end{document}

%% file: localmacros.tex
\newcommand{\auskommentieren}[1]{}
\newcommand{\betonen}[1]{{\bf #1}}

\newcommand{\rnf}[1]{\!\!\Downarrow_{#1}}

\newcommand{\f}{{\cal F}} %Funktionenring
\newcommand{\myt}{{\cal T}} %Traegermenge
\newcommand{\zero}{o}   %Nullfunktion
\newcommand{\one}{{\bf 1}} %Eins
\newcommand{\supp}{{\sf supp}} %Traeger
\newcommand{\radd}{\oplus} %Ringaddition
\newcommand{\rmult}{\star} %Ringmultiplikation
\newcommand{\monoms}{{\sf M}} % Monome
\newcommand{\terms}{{\sf T}} %Terme

\newcommand{\m}{{\cal M}}      %Monoid M
\newcommand{\mm}{\circ}        %Multiplikation im Monoid
\newcommand{\id}{\equiv} %Zeichenweise Gleichheit

\newcommand{\irr}{{\rm IRR}\/}
\newcommand{\pred}{\prec}
\newcommand{\predeq}{\preceq}

\newcommand{\mswab}[1]{\frak{#1}}

\newcommand{\g}{{\cal G}}      %Gruppe
\newcommand{\h}{{\cal H}}      %Normale Untergruppe
\newcommand{\inv}[1]{{\sf inv}\/(#1)}

\newcommand{\n}{{\mathbb{N}}} 
\newcommand{\z}{{\mathbb{Z}}} 

\newcommand{\myk}{{\mathbb{K}}}
\newcommand{\q}{{\mathbb{Q}}}

\newcommand{\rr}{{\sf R}}

\newcommand{\mrm}{\ast} %Multiplikation im Monoid Ring 
\newcommand{\mmult}{\ast} %Skalarmultiplikation im Modul
\newcommand{\skm}{\cdot} %Ringmultiplication 
\newcommand{\cd}{\cdot}

\newcommand{\ideal}[2]{{\sf ideal}_{#1}^{#2}}
\newcommand{\kernel}{{\sf ker}}
\newcommand{\image}{{\sf im}}
\newcommand{\spol}[1]{{\sf spol}_{#1}}
\newcommand{\spols}{{\sf SPOL}}

\newcommand{\sat}{{\sf SAT}}

\newcommand{\hm}{{\sf HM}}
\newcommand{\hc}{{\sf HC}}
\newcommand{\hterm}{{\sf HT}}
\newcommand{\reductum}{{\sf RED}}
\newcommand{\lcm}{{\sf LCM}}

\newcommand{\R}{\Longrightarrow}
\newcommand{\myr}{\longrightarrow}

\newcommand{\lr}{\longleftrightarrow}

\newcommand{\red}[4]{\mbox{$\,\stackrel{#1}{#2}\!\!\mbox{}^{{\rm #3}}_{#4}\,$}}
\newcommand{\nred}[4]%
{\mbox{$\,\,\,\,{\not\!\!\!\stackrel{#1}{#2}\!\!\mbox{}^{{\rm #3}}_{#4}}\,$}}

\newcommand{\Ba}[1]{{\bf Proof #1}: 
                      \renewcommand{\baselinestretch}{1.1}\small\normalsize}   %Beweisanfang
\newcommand{\qed}{\hspace*{\fill} q.e.d. \par\medskip
        \renewcommand{\baselinestretch}{1}\small\normalsize}    %Beweisende

\newcommand{\problem}[4]%
{\noindent{{\sc #1}}
\hrule height 0.0pt \hfill
\\
\begin{tabular}{p{2.2cm}lp{11cm}}
{\bf Given:} && #2 
\\
{\bf Problem:} && #3
\\
{\bf Proceeding:} #4
\end{tabular}
\\
}

\newtheorem{lemma}{Lemma}[section]
\newtheorem{example}[lemma]{Example}
\newtheorem{remark}[lemma]{Remark}
\newtheorem{definition}[lemma]{Definition}
\newtheorem{theorem}[lemma]{Theorem}
\newtheorem{corollary}[lemma]{Corollary}
\newcommand{\exaend}{\hfill $\diamond$}
\newcommand{\dend}{\hfill $\diamond$}
\newcommand{\remend}{\hfill $\diamond$}
\newcommand{\lemend}{}
\newcommand{\theoend}{}
\newcommand{\corend}{}
\newcommand{\ohnebeweis}{\hfill $\diamond$}

\newcommand{\gb}{\mbox{\sc Gb}}

\newcommand{\procedure}[2]%
{\noindent{{\bf Procedure}: {\sc #1}}\vspace{1mm}
\hrule height 0.5pt \hfill%\\[-1.3ex] 
#2% \\[-1.3ex] 
\hrule height 0.5pt \hfill}
\newcommand{\kommentar}{\small}

%% file: veraformat.tex
\newcommand{\PreserveBackslash}[1]{\let\temp=\\#1\let\\=\temp}

\newcommand{\randnotiz}[1]%
{\marginpar{\small{\begin{flushleft}{\bf #1}\end{flushleft}}}{}}

\renewcommand{\chaptermark}[1]{\markboth{\chaptername\ \thechapter\ -\ #1}{1}}

%% file: title.tex
%%%%%%%%%%%%%%%%%%%%%%%%%%%%%%%%%%%%%%%
%%%%%%%% Titelseite %%%%%%%%%%%%%%%%%%%
%%%%%%%%%%%%%%%%%%%%%%%%%%%%%%%%%%%%%%%
%\maketitle
\begin{titlepage}

\begin{center}
\vspace{10cm}
{\huge\bf \mbox{       }    \\[5ex]
       A Systematic Study   \\[2ex]
              of            \\[3ex]
    Gr\"obner Basis Methods} 
\\[7ex]
Vom Fachbereich Informatik \\[1ex]
der Technischen Universit{\"a}t Kaiserslautern \\[1ex]
genehmigte Habilitationsschrift\\[1ex]
von\\[5ex]
Dr. Birgit Reinert
\end{center}

\vspace*{1cm}

Datum der Einreichung:  6. Januar 2003\\[1ex]
Datum des wissenschaftlichen Vortrags:  9. Februar 2004\\[2ex]

\vspace*{1cm}
\begin{tabular}{ll}
Dekan: & Prof. Dr. Hans Hagen \\[2ex]
Habilitationskommission: \\
Vorsitzender: & Prof. Dr. Otto Mayer \\
Berichterstatter: & Prof. Dr. Klaus E. Madlener \\
& Prof. Dr. Teo Mora \\
& Prof. Dr. Volker Weispfenning 
\end{tabular}
\end{titlepage}

%%% Local Variables: 
%%% mode: latex
%%% TeX-master: "testlauf"
%%% End: 

%% file: vorwort.tex
\chapter*{Vorwort}

Die vorliegende Arbeit ist die Quintessenz meiner Ideen und Erfahrungen, die ich in den letzten Jahren bei meiner Forschung auf dem Gebiet der Gr\"obnerbasen gemacht habe. Meine geistige Heimat war dabei die Arbeitsgruppe von Professor Klaus Madlener an der Technischen Universit\"at Kaiserslautern. Hier habe ich bereits im Studium Bekanntschaft mit der Theorie der Gr\"obnerbasen gemacht und mich w\"ahrend meiner Promotion mit dem Spezialfall dieser Theorie f\"ur Monoid- und Gruppenringe besch\"aftigt. 
Nach der Promotion konnte ich im Rahmen eines DFG-Forschungsstipendiums zus\"atzlich Problemstellungen und Denkweisen anderer Arbeitsgruppen kennenlernen - die Arbeitsgruppe von Professor Joachim Neub\"user in Aachen und die Arbeitsgruppe von Professor Theo Mora in Genua. 
Meine Aufenthalte in diesen Arbeitsgruppen haben meinen Blickwinkel f\"ur weitergehende Fragestellungen erweitert. 
An dieser Stelle m\"ochte ich mich bei allen jenen bedanken, die mich in dieser Zeit begleitet haben und so zum Entstehen und Gelingen dieser Arbeit beigetragen haben.

Mein besonderer Dank gilt meinem akademischen Lehrer Professor Klaus Madlener, der meine akademische Ausbildung schon seit dem Grundstudium begleitet und meine Denk- und Arbeitsweise wesentlich gepr\"agt hat. 
Durch ihn habe ich gelernt, mich intensiv mit diesem Thema zu besch\"aftigen und mich dabei nie auf nur einen Blickwinkel zu beschr\"anken. 
Insbesondere sein weitreichenden Literaturkenntnisse und die dadurch immer neu ausgel\"osten Fragen aus verschiedenen Themengebieten bewahrten meine Untersuchungen vor einer gewissen Einseitigkeit. Er hat mich gelehrt, selbst\"andig zu arbeiten, Ideen und Papiere zu hinterfragen, mir meine eigene Meinung zu bilden, diese zu verifizieren und auch zu vertreten.

Professor Teo Mora und Professor Volker Weispfenning danke ich f\"ur die \"Ubernahme der weiteren Begutachtungen dieser Arbeit. Professor Teo Mora danke ich insbesondere auch f\"ur die fruchtbare Zeit in seiner Arbeitsgruppe in Genua. Seine Arbeiten und seine Fragen haben meine Untersuchungen zum Zusammenhang zwischen Gr\"obnerbasen in Gruppenringen und dem Todd-Coxeter Ansatz f\"ur Gruppen und die Fragestellungen dieser Arbeit wesentlich gepr\"agt. 

Meinen Kollegen aus unserer Arbeitsgruppe danke ich f\"ur ihre Diskussionsbe\-reitschaft und ihre geduldige Anteilnahme an meinen Gedanken. Insbesondere Andrea Sattler-Klein, Thomas Dei{\ss}, Claus-Peter Wirth und Bernd L\"ochner haben immer an mich geglaubt und mich ermutigt, meinen Weg weiter zu gehen. Mit ihnen durfte ich nicht nur Fachliches sondern das Leben teilen.

Um eine solche Arbeit fertigzustellen braucht man jedoch nicht nur eine fachliche Heimat. Mein Mann Joachim hat nie an mir gezweifelt und mich immer unterst\"utzt. Auch nach unserem Umzug nach Rechberghausen hat er mein Pendeln nach Kaiserslautern und das Brachliegen unseres Haushalts mitgetragen. Ohne meine Eltern Irma und Helmut Weber und meine Patentante Anita Sch\"afer w\"are insbesondere nach der Geburt unserer Tochter Hannah diese Arbeit nie fertiggestellt worden. Sie haben Hannah liebevoll beh\"utet, so dass ich lesen, schreiben und arbeiten konnte, ohne ein schlechtes Gewissen zu haben. Hannah hat von Anfang an gelernt, dass Forschung ein Leben bereichern kann und die Zeit mit ihrer erweiterten Familie und die vielen Reisen nach Kaiserslautern genossen.

Gewidmet ist diese Arbeit meiner Mutter, die leider die Fertigstellung nicht mehr erleben durfte, und Hannah, die auch heute noch Geduld aufbringt, wenn ihre Mutter am Computer f\"ur ihre ``Schule'' arbeitet.

Rechberghausen, im August 2004 \hfill Birgit Reinert

%\pagebreak
\chapter*{}
\thispagestyle{empty}
\vspace*{15cm}
\hspace*{\fill}\parbox{6cm}{%
\begin{center}
{\em F\"ur meine Mutter} \\
{\em und Hannah}
\end{center}}
%%% Local Variables: 
%%% mode: latex
%%% TeX-master: t
%%% TeX-master: t
%%% End: 

%% file: introduction.tex
\chapter{Introduction}\label{introduction}
%%%%%%%%%%%%%%%%%%%%%%%%%%%%%%%%%%%%%%%%%%%%%%%%%%%%%%%%%%%%%%%%%%%%%%%%%%%%%%%
One of the amazing features of computers is the ability to do
 extensive computations impossible to be done by hand.
This enables to overcome the boundaries of constructive algebra as
 proposed by mathematicians as Kronecker.
He demanded that
 definitions of mathematical objects should be given in such a way that
 it is possible to decide in a finite number of steps whether a
 definition applies to an object.
While in the beginning computers
 were used to do incredible numerical calculations, a new dimension
 was added when they were used to do symbolical mathematical
 manipulations substantial to many fields in mathematics and
 physics.
These new possibilities led to
 open up whole new areas of mathematics and computer science.
In the wake of these developments has come a new access to abstract
 algebra in a computational fashion -- computer algebra.
One important contribution to this field which is the subject of this work
 is the theory of Gr\"obner bases -- the result of Buchberger's
 algorithm for manipulating systems of polynomials.
%%%%%%%%%%%%%%%%%%%%%%%%%%%%%%%%%%%%%%%%%%%%%%%%%%%%%%%%%%%%%%%%%%%%%%%%%%%%%%%
\section{The History of Gr\"obner Bases}
%%%%%%%%%%%%%%%%%%%%%%%%%%%%%%%%%%%%%%%%%%%%%%%%%%%%%%%%%%%%%%%%%%%%%%%%%%%%%%%
In 1965 Buchberger introduced the theory of Gr\"obner bases\footnote{Note
 that similar concepts appear in a paper of Hironaka where the notion of a 
 complete set of polynomials is called a standard basis \cite{Hi64}.} for
 polynomial ideals in commutative
 polynomial rings over fields  \cite{Bu65,Bu70}.
Let $\myk[X_1, \ldots, X_n]$ be a polynomial ring over a computable field $\myk$
 and ${\mathfrak i}$ an ideal in $\myk[X_1, \ldots, X_n]$.
Then the quotient $\myk[X_1, \ldots, X_n] / {\mathfrak i}$ is a $\myk$-algebra.
If this quotient is zero-dimensional the algebra has a finite basis consisting of
 power products $X_1^{i_1}\ldots X_n^{i_n}$.
This was the starting point for Buchberger's PhD thesis.
His advisor Wolfgang Gr\"obner wanted to compute the multiplication
 table and had suggested a procedure for zero-dimensional ideals,
 for which termination conditions were lacking.
The result of Buchberger's studies then was a terminating algorithm which 
 turned a basis of an ideal into a special basis which allowed to
 solve Gr\"obner's question of writing down an explicit multiplication for
 the multiplication table of the quotient in the zero-dimensional case and
 was even applicable to arbitrary polynomial ideals. 
Buchberger called these special bases of ideals Gr\"obner bases.
%%%%%%%%%%%%%%%%%%%%%%%%%%%%%%%%%%%%%%%%%%%%%%%%%%%%%%%%%%%%%%%%%%%%%%%%%%%%%%%
\section{Two Definitions of Gr\"obner Bases}\label{section.intro.definitions}
%%%%%%%%%%%%%%%%%%%%%%%%%%%%%%%%%%%%%%%%%%%%%%%%%%%%%%%%%%%%%%%%%%%%%%%%%%%%%%%
In literature there are two main ways to define Gr\"obner bases
 in polynomial rings over fields.
They both require an admissible\footnote{An
 ordering $\succeq$ on the set of terms is called an admissible 
 term ordering if for every term $s,t,u$,
 $s \succeq 1$ holds, and $s \succeq t$ implies $s \mm u \succeq t \mm u$.
 An ordering fulfilling the latter condition is also said to be
 compatible with the respective multiplication $\mm$.} ordering
 on the set of terms.
With respect to such an ordering, 
 given a polynomial $f$ the maximal term
 occurring in $f$ is called the head term denoted by $\hterm(f)$.

One way to characterize Gr\"obner bases in an algebraic fashion
 is to use the concept of term division:
A term $X_1^{i_1}\ldots X_n^{i_n}$ is said
 to divide another term $X_1^{j_1}\ldots X_n^{j_n}$ if and
 only if $i_l \leq j_l$ for all $1 \leq l \leq n$.
Then a set $G$ of polynomials is called a Gr\"obner basis
 of the ideal $\mathfrak{i}$ it generates if and only if for every $f$ in $\mathfrak{i}$
 there exists a polynomial $g \in G$ such that $\hterm(g)$ divides $\hterm(f)$.

Another way to define Gr\"obner bases in polynomial rings is
 to  establish a rewriting approach to the theory of polynomial
 ideals.
Polynomials can be used as rules by using the largest monomial according to the
 admissible ordering as a left hand side of a rule.
Then a term is reducible by a polynomial as a rule if the head term
 of the polynomial divides the term.
A Gr\"obner basis $G$ then is a  set of  polynomials
 such that every polynomial in the polynomial ring has a unique
 normal form with respect to this reduction relation using the polynomials
 in $G$ as rules (especially the
 polynomials in the ideal generated by $G$ reduce to zero using $G$).

Of course both definitions coincide for polynomial
 rings since the reduction relation defined by Buchberger can be compared to
 division of one polynomial by a
 set of finitely many polynomials.
%%%%%%%%%%%%%%%%%%%%%%%%%%%%%%%%%%%%%%%%%%%%%%%%%%%%%%%%%%%%%%%%%%%%%%%%%%%%%%%
\section{Applications of Gr\"obner Bases}\label{section.intro.applications}
%%%%%%%%%%%%%%%%%%%%%%%%%%%%%%%%%%%%%%%%%%%%%%%%%%%%%%%%%%%%%%%%%%%%%%%%%%%%%%%
The method of Gr\"obner bases allows to solve many problems related to
 polynomial ideals in a computational fashion.
It was shown by Hilbert (compare Hilbert's basis theorem) that every
 ideal in a polynomial ring has a finite generating set.
However, an arbitrary finite generating set need not provide much
insight into the nature of the ideal.
Let $f_1 = X_1^{2} + X_2$ and $f_2 = X_1^{2} + X_3$ be two polynomials
 in the polynomial ring\footnote{$\q$ denotes the rational numbers.} $\q[X_1,X_2,X_3]$.
Then  ${\mathfrak i} = \{ f_1 \mrm g_1 + f_2
\mrm g_2 \mid g_1, g_2 \in \q[X_1,X_2,X_3] \}$ is the ideal they generate 
and it is not hard to see that
the polynomial $X_2 - X_3$ belongs to ${\mathfrak i}$ since $X_2 - X_3 =
f_1 - f_2$.
But what can be said about the polynomial $f = X_3^{3} + X_1 + X_3$?
Does it belong to ${\mathfrak i}$ or not?

The problem to decide whether a given polynomial lies in a given
 ideal is called the membership problem for ideals.
In case the generating set is a Gr\"obner basis this problem becomes
immediately decidable, as the membership problem then reduces to checking whether
 the polynomial reduces to zero using the elements of the Gr\"obner basis for reduction.

In our example the set $\{ X_1^{2} + X_3, X_2 - X_3 \}$ is a
generating set of ${\mathfrak i}$ which is in fact a Gr\"obner basis.
Now returning to the polynomial $f = X_3^{3} + X_1 + X_3$ we find that
it cannot
belong to ${\mathfrak i}$ since neither $X_1^{2}$ nor $X_2$ is a divisor
of a term in  $f$ and hence $f$ cannot be reduced to zero using the
polynomials in the Gr\"obner basis as rules.

The terms $X_1^{i_1}X_2^{i_2}X_3^{i_3}$ which are not reducible
 by  the set $\{ X_1^{2} + X_3, X_2 - X_3 \}$ form a basis of
 the $\q$-algebra $\q[X_1,X_2,X_3]/{\mathfrak i}$.
By inspecting the head terms $X_1^{2}$ and $X_2$ of the Gr\"obner basis
 we find that the (infinite) set $\{ X_3^i, X_1X_3^i \mid i \in \n \}$ is 
 such a basis.
Moreover, an ideal is zero-dimensional, i.e.~this set is finite, if and only
 if for each variable $X_i$ the  Gr\"obner basis contains a polynomial with
 head term $X_i^{k_i}$ for some $k_i \in \n^+$.
Similarly the form of the Gr\"obner basis reveals whether the ideal is trivial:
 ${\mathfrak i} = \myk[X_1, \ldots, X_n]$ if and only if every\footnote{Notice
 that if one Gr\"obner basis contains an element from $\myk$ so will all the others.}
 Gr\"obner basis contains an element from $\myk$.

Further applications of Gr\"obner bases come from areas as widespread
 as robotics, computer vision, computer-aided design, geometric theorem proving,
 Petrie nets and many more.
More details can be found  e.g.~in Buchberger \cite{Bu87}, or the books of Becker and
 Weispfenning  \cite{BeWe92}, Cox, Little and O'Shea \cite{CoLiOS92}, and
 Adams and Loustaunau \cite{AdLo94}.
%%%%%%%%%%%%%%%%%%%%%%%%%%%%%%%%%%%%%%%%%%%%%%%%%%%%%%%%%%%%%%%%%%%%%%%%%%%%%%%
\section{Generalizations of Gr\"obner Bases}
%%%%%%%%%%%%%%%%%%%%%%%%%%%%%%%%%%%%%%%%%%%%%%%%%%%%%%%%%%%%%%%%%%%%%%%%%%%%%%%
In the last years, the method of Gr\"obner bases and its applications have been extended 
 from commutative polynomial rings over fields to various types of  algebras 
 over fields and other rings.
In general for such rings arbitrary finitely generated ideals will not have finite Gr\"obner bases.
Nevertheless, there are interesting classes for which every finitely generated  
 (left, right or even two-sided) ideal has a finite Gr\"obner basis which can be computed by appropriate 
 variants of completion based algorithms.

First successful generalizations were extensions  to commutative polynomial rings over
 coefficient domains other than fields.
It was shown by several authors including  Buchberger, Kandri-Rody, Kapur,
Narendran, Lauer, Stifter, and Weispfenning that Buchberger's approach remains valid for polynomial
 rings over the integers, or even Euclidean rings, and over
 regular rings  
 (see e.g.~\cite{Bu83,Bu85,KaKa84,KaKa88,KaNa85,La76,St87,We87}).
For regular rings Weispfenning has to deal with the situation that zero-divisors in the coefficient domain
 have to be considered.
He uses a technique he calls Boolean closure to repair this problem and this
 technique can be regarded as a special saturating process\footnote{Saturation techniques
 are used in various fields to enrich a generating set of a structure in such a way,
 that the new set still describes the same structure but allows more insight.
 For example symmetrization in groups can be regarded as such a saturating process.}.
We will later on see how such saturating techniques become important ingredients of Gr\"obner
 basis methods in many algebraic structures.

Since the development of computer algebra systems for commutative
 algebras made it possible to perform tedious calculations using computers,
 attempts to generalize such systems and especially Buchberger's ideas
 to non-commutative algebras followed.
Originating from special problems in physics, Lassner in \cite{La85}
 suggested how to extend existing computer algebra systems in order to
 additionally handle special classes of non-commutative algebras, e.g.~Weyl algebras.
He studied structures where the elements could be represented using
 the usual representations of polynomials in commutative variables and
 the non-commutative multiplication could be performed by a so-called
 ``twisted product'' which required only procedures involving
 commutative algebra operations and differentiation.
Later on together with Apel he extended Buchberger's algorithm to
 enveloping fields of Lie algebras  \cite{ApLa88}.
Because these ideas use representations of the elements by commutative polynomials,
 Dickson's Lemma\footnote{Dickson's Lemma in the context of commutative terms is as follows:
 For every infinite sequence of terms $t_s$, $s \in \n$, there exists an
 index $k \in \n$ such that for every index $i > k$ there exists an index 
 $j \leq k$ and a term $w$ such that $t_i = t_j w$.} can be carried over.
By this the existence and construction of finite Gr\"obner bases for finitely
 generated left ideals can be ensured using the same arguments as in the original
 approach.
%In \cite{Ga85} Galligo also studied algorithmic questions on ideals of
% differential operators.

On the other hand, Mora gave a concept of Gr\"obner bases for a class
 of non-commutative algebras by saving an other property of the
 commutative polynomial ring -- admissible orderings --
 while losing the validity of Dickson's Lemma.
The usual polynomial ring can be viewed as a monoid ring where the
 monoid is a finitely generated free commutative monoid.
Mora studied the class where the free commutative monoid is substituted
 by a free monoid -- the class of finitely generated free monoid rings (compare e.g.~\cite{Mo85,Mo94}).
The ring operations are mainly performed in the coefficient
 domain while the terms are treated like words, i.e., the variables no
 longer commute with each other and multiplication is concatenation.
The definitions of (one- and two-sided) ideals, reduction and Gr\"obner
 bases are carried over from the commutative  case to establish
 a similar theory of Gr\"obner bases in ``free non-commutative
 polynomial rings over fields''.
But these rings are no longer Noetherian if they are generated by more
 than one variable.
Mora presented a terminating completion procedure for finitely generated one-sided ideals
 and an enumeration procedure for finitely generated two-sided ideals
 with respect to some term ordering in free monoid rings.
For the special instance of ideals generated by bases of the restricted form
 $\{ \ell_i - r_i \mid \ell_i, r_i \mbox{ words}, 1 \leq i \leq n \}$, Mora's
 procedure coincides with Knuth-Bendix completion for string rewriting systems and the
 one-sided cases can be related to prefix respectively suffix rewriting \cite{MaRe95,MaRe97c}.
Hence many results known for finite string rewriting systems and their completion
 carry over to finitely generated ideals and the computation of their Gr\"obner bases.
Especially the undecidability of the word problem yields non-termination for Mora's
 general procedure (see also \cite{Mo87}).

Gr\"obner bases and Mora's procedure have been generalized to path algebras
 (see \cite{FaFeGr93,Ke97}); free non-commutative polynomial rings are in fact a particular
 instance of path algebras.

Another class of non-commutative rings  where the elements can
 be represented by the usual polynomials and which allow the construction of finite
 Gr\"obner bases for arbitrary ideals are the solvable polynomial rings, a class
 intermediate between commutative and general non-commutative
 polynomial rings.
They were studied by Kandri-Rody, Weispfenning and Kredel \cite{KaWe90,Kr93}.
Solvable polynomial rings can be described by ordinary polynomial
 rings $\myk[X_1, \ldots, X_n]$ provided with a ``new'' definition of multiplication which coincides with the
 ordinary multiplication except for the case that a variable $X_j$ is
 multiplied with a variable $X_i$ with lower index, i.e., $i<j$.
In the latter case multiplication can be defined by equations of the form
 $X_j \star X_i = c_{ij} X_iX_j + p_{ij}$
 where $c_{ij}$ lies in $\myk^*= \myk \backslash \{ 0 \}$ and $p_{ij}$ is a polynomial
 ``smaller'' than $X_iX_j$ with respect to a fixed admissible term
 ordering on the polynomial ring.

The more special case of twisted semi-group rings, where $c_{ij} = 0$
 is possible, has been studied in \cite{Ap88,Mo88}.

In \cite{We92} Weispfenning showed the existence of finite Gr\"obner
 bases for arbitrary finitely
 generated ideals in non-Noetherian skew polynomial rings over two
 variables $X,Y$ where a ``new'' multiplication $\star$ is introduced
 such that $X \star Y = XY$ and $Y \star X = X^eY$ for some fixed $e$
 in $\n^+$.

Ore extensions have been successfully studied by Pesch in his PhD Thesis \cite{Pe97}
 and his results on two-sided Gr\"obner bases are also presented in
 \cite{Pe98}.

Most of the results cited so far assume admissible well-founded orderings on
 the set of terms
 so that in fact the reduction relations can be defined by
 considering the head monomials mainly (compare the algebraic definition of Gr\"obner bases in Section
 \ref{section.intro.definitions}).
This is essential to characterize Gr\"obner bases in the
respective ring with respect to the corresponding reduction 
relation\footnote{These reduction relations are based on divisibility of terms,
 namely the term to be reduced is divisible by the head term of the polynomial
 used as rule for the reduction step.} in a
finitary manner and  to enable to decide whether a finite set is a Gr\"obner
basis by checking whether the s-polynomials are reducible to
zero\footnote{Note that we always assume that the reduction relation in the ring is effective.}.

There are rings combined with reduction relations where admissible well-founded orderings cannot
 be accomplished and, therefore, other concepts to characterize
 Gr\"obner bases have been developed.
For example
 in case the ring contains zero-divisors a well-founded ordering on the ring
 is no
 longer compatible with the ring multiplication\footnote{When studying monoid rings over reduction rings it is possible
 that the ordering on the ring is not compatible with scalar multiplication
 as well as with multiplication with monomials or polynomials.}.
This phenomenon has been studied for the case of zero-divisors in
 the coefficient domain by Kapur and Madlener \cite{KaMa89} and by
 Weispfenning for the special case of regular rings \cite{We87}.
In his PhD thesis \cite{Kr93}, Kredel described problems occurring when
 dropping the 
 axioms guaranteeing the existence
 of admissible orderings in the theory of solvable polynomial rings by
 allowing $c_{ij} = 0$ in the defining equations above.
He sketched the idea of using saturation techniques to repair some of them.
Saturation enlarges the generating sets of ideals in order to ensure that
 enough head terms exist to do all necessary reduction steps and 
 this process can often be related to additional special critical pairs.
Similar ideas can be found in the PhD thesis of Apel \cite{Ap88}.
For special cases, e.g.~for the Grassmann
(exterior) algebras, positive results can be  achieved (compare the paper
 of Stokes \cite{St90}).

Another important class of rings where reduction relations can be introduced
 and completion techniques can be applied to enumerate and sometimes
 compute Gr\"obner bases are  monoid and group rings.
They have been studied in detail by various authors, e.g.~free group rings (\cite{Ro93}),
 monoid and group rings (\cite{MaRe93b,MaRe96b,Re95,Re96,MaRe97a})
 (including finite and free monoids and finite, free, plain and
 polycyclic groups), and polycyclic group rings (\cite{Lo98}).
In this setting we again need saturation techniques to repair a severe defect
 due to the fact that in general we cannot expect the ordering on the set of terms
 (here of course now the monoid or group elements)
 to be both, well-founded and admissible.
Let $\f$ be the free group generated by one element $a$.
Then for the polynomial $a + 1$ in $\q[\f]$ we have $(a+1) \mrm a^{-1} = 1 + a^{-1}$,
 i.e.,~after multiplication with the inverse element
 $a^{-1}$ the largest term of the new polynomial no
 longer results from the largest one of the original polynomial.
Moreover, assuming our ordering is well-founded, it cannot be compatible with the
 group multiplication\footnote{Assuming $a \succ 1$ compatibility with multiplication
 would imply $1 \succ a^{-1}$ giving rise to an infinite descending chain $a^{-1} \succ
 a^{-2} \succ \ldots$ contradicting the well-foundedness of the ordering.
 On the other hand for  $1 \succ a$ compatibility with multiplication immediately gives
 us an infinite descending chain $a \succ
 a^2 \succ \ldots$.}.

All approaches cited in this section can be basically divided into two main streams:
One  extension was to study structures which still allow to present
 their elements by ordinary ``commutative'' polynomials.
The advantage of this generalization is that Dickson's Lemma, which
 is essential in proving termination for Buchberger's algorithm,
 carries over.
The other idea of generalization 
 was to view the polynomial ring as a special monoid ring
 and to try to extend Buchberger's approach to other monoid and 
 group rings.
Since then in general Dickson's Lemma no longer holds, other ways to
 prove termination, if possible, have to be established.
Notice that solvable rings, skew-polynomial rings and arbitrary quotients of
 non-commutative polynomial rings cannot be interpreted as monoid rings.
Hence to find a generalization which will subsume {\em all} results cited here, a more
 general setting is needed.
In his habilitation thesis \cite{Ap98}, Apel provides one generalization which basically extends
 the first one of these two in such a way that Mora's approach can be incorporated.
He uses an abstraction of graded structures which needs admissible well-founded orderings.
Hence he cannot deal with group rings and many cases of monoid rings where such
 orderings cannot exists.
On the other hand  he is much more interested in algebraic characterizations of Gr\"obner
 bases and the division algorithms associated to them.

In order to characterize structures where the well-founded ordering is no longer admissible, we
 extend Gr\"obner basis techniques to an abstract setting called function rings.
%%%%%%%%%%%%%%%%%%%%%%%%%%%%%%%%%%%%%%%%%%%%%%%%%%%%%%%%%%%%%%%%%%%%%%%%%%%%%%%
\section{Gr\"obner Bases in Function Rings -- A Guide for Introducing Reduction Relations
 to Algebraic Structures}
%%%%%%%%%%%%%%%%%%%%%%%%%%%%%%%%%%%%%%%%%%%%%%%%%%%%%%%%%%%%%%%%%%%%%%%%%%%%%%%
The aim of this work is to give a general setting which comprises {\em all}
 generalizations mentioned above and which is a basis for studying further
 structures in the light of introducing reduction relations and Gr\"obner basis techniques.
All structures mentioned so far can be viewed as rings of functions with
 finite support.
For such rings we introduce the familiar concepts of polynomials,
 (right) ideals, standard representations, standard bases, reduction relations
 and Gr\"obner bases.
A general characterization of Gr\"obner bases in an ``algorithmic fashion''
 is provided.
It is shown that in fact polynomial rings, solvable polynomial rings, free respectively finite monoid rings, and free, finite,
 plain, respectively polycyclic group rings are examples of our
 generalization where finite Gr\"obner bases can be computed. 
While most of the examples cited above are presented in the literature as rings over fields
 we will here also present the more general concept of function rings over
 reduction rings (compare \cite{Ma86,Re95,MaRe98}) and the impotant special case
 of function rings over the integers.
%%%%%%%%%%%%%%%%%%%%%%%%%%%%%%%%%%%%%%%%%%%%%%%%%%%%%%%%%%%%%%%%%%%%%%%%%%%%%%%
\section{Applications of Gr\"obner Bases Generalized to Function Rings}
%%%%%%%%%%%%%%%%%%%%%%%%%%%%%%%%%%%%%%%%%%%%%%%%%%%%%%%%%%%%%%%%%%%%%%%%%%%%%%%
For polynomial rings over fields many algebraic questions related to ideals can be
 solved using Gr\"obner bases and their associated reduction relations.
Hence the question arises which of these applications can be extended to more
 general settings.
While some questions e.g.~concerning algebraic geometry are strongly connected
 to polynomial rings over fields, many other applications carry over.
They include natural ones such as the membership problem for ideals,
 as well as special techniques such as elimination theory or the treatment
 of systems of linear equations.
%%%%%%%%%%%%%%%%%%%%%%%%%%%%%%%%%%%%%%%%%%%%%%%%%%%%%%%%%%%%%%%%%%%%%%%%%%%%%%%
\section{Organization of the Contents}
%%%%%%%%%%%%%%%%%%%%%%%%%%%%%%%%%%%%%%%%%%%%%%%%%%%%%%%%%%%%%%%%%%%%%%%%%%%%%%%
Chapter \ref{chapter.definitions} introduces some of the basic themes of
 this book. 
We need some basic notions from the theory of algebra as well as from
 the theory of rewriting systems.
Furthermore, as the aim of this book is to provide
 a systematic study of Gr\"obner basis methods, a short introduction
 to the original case of Gr\"obner bases in polynomial rings over fields is
 presented.

Chapter \ref{chapter.reduction.rings} concentrates on rings with reduction relations,
which are studied with regard to the existence of Gr\"obner bases.
They are called reduction rings in case they allow finite Gr\"obner bases for
 finitely generated ideals.
Moreover, special ring constructions are presented, which in many cases preserve the existence of
 Gr\"obner bases.
These constructions include quotients and sums of reduction rings as well as modules and polynomial
 rings over reduction rings.
Many structures with reduction relations allowing Gr\"obner bases can already be
 found in this setting.
For example knowing that the integers $\z$ for certain reduction relations 
 allow finite Gr\"obner bases, using the results of this chapter,
 we can conclude that the module $\z^k$ as well as the polynomial
 rings $\z[X_1, \ldots, X_n]$ and $\z^k[X_1, \ldots, X_n]$ allow the computation of finite
 Gr\"obner bases.

Chapter \ref{chapter.function.rings} is the heart of this book.
It establishes a generalizing framework for structures enriched with reduction
 relations and studied with respect to the existence of Gr\"obner bases in the literature.
Reduction relations are defined for the setting of function rings over fields
 and later on generalized to reduction rings.
Definitions for terms such as variations of standard representations, standard bases and
 Gr\"obner bases are given and compared to the known terms from the theory of Gr\"obner
 bases over polynomial rings.
It turns out that while completion procedures will still involve equivalents to s-polynomials
 or the more general concept of g- and m-polynomials for the ring case, these situations are no longer
 sufficient to characterize Gr\"obner bases.
Saturation techniques, which enrich the bases by additional polynomials, are needed.
Moreover, for function rings over reduction rings the characterizations no longer
 describe Gr\"obner bases but only weak\footnote{Weak Gr\"obner bases are
 bases such that
 any polynomial in the ideal they generate can be reduced to zero. For fields
 this property already characterizes Gr\"obner bases as the
 Translation Lemma holds. In general this is not true
 and while weak Gr\"obner bases allows to solve the ideal membership
 problem they no longer guarantee the existence of unique
 normal forms for elements of the quotient.}
 Gr\"obner bases, since the Translation Lemma\footnote{The Translation Lemma 
 establishes that if for two polynomials $f,g$ we have that $f-g$ reduces
 to zero, both polynomials reduce to the same normal form.} no longer holds. 
Since the ring of integers viewed as a reduction ring is of special interest
 in the literature and allows more insight into the respective chosen
 reduction relations, this special case is studied.

Chapter \ref{chapter.applications} outlines how some applications known for Gr\"obner bases
 in the literature carry over to function rings.
These applications include natural ones such as the ideal membership problem,
 representation problems, the ideal inclusion problem, the ideal triviality problem, and
 many more.
Another focus is on doing computations in quotient rings using Gr\"obner bases.
The powerful elimination methods are also generalized.
One of their applications to study polynomial mappings is outlined.
Finally solutions for linear equations over function rings in terms of Gr\"obner 
 bases are provided.

%%% Local Variables: 
%%% mode: latex
%%% TeX-master: "testlauf"
%%% TeX-master: "testlauf"
%%% End: 

%% file: algebra.tex
After introducing the necessary definitions required from
 algebra we focus on the subject of this book --- Gr\"obner bases.
One way of characterizing Gr\"obner bases is in terms of algebraic
 simplification or reduction.
The aim of this chapter is to introduce an abstract concept for the notion
 of reduction which is the basis of many syntactical methods for studying
 structures in mathematics or theoretical computer science in Section \ref{section.rewriting}.
It is the foundation for e.g.~term rewriting and string rewriting and we
 introduce a reduction relation for polynomials in the commutative
 polynomial ring over a field in a similar fashion.
Gr\"obner bases then arise naturally when doing completion in this setting in Section \ref{section.buchberger}. 

%%%%%%%%%%%%%%%%%%%%%%%%%%%%%%%%%%%%%%%%%%%%%%%%%%%%%%%%%%%%%%%%%%%%%%%
%%%%%%%%%%%%%%%% Section: A little Algebra %%%%%%%%%%%%%%%%%%%%%%%%%%%%
%%%%%%%%%%%%%%%%%%%%%%%%%%%%%%%%%%%%%%%%%%%%%%%%%%%%%%%%%%%%%%%%%%%%%%%
\section{Algebra}\label{section.algebra}%\RechterRand{Algebra}
%%%%%%%%%%%%%%%%%%%%%%%%%%%%%%%%%%%%%%%%%
Mathematical theories are closely related with the study of two objects, namely
 sets and functions.
Algebra can be regarded as the study of algebraic operations on sets, i.e., 
 functions that take elements from a set to the set itself.
Certain algebraic operations on sets combined with certain axioms are again 
 the objects of independent theories.
This chapter is a short introduction to some of the algebraic systems
  used later on:
 monoids, groups, rings, fields, ideals and modules.
\begin{definition}~\\
{\rm
A non-empty set of elements $\m$ together with a binary operation $\mm_{\m}$ is
 said to form a \index{monoid}\betonen{monoid}, if for all $\alpha,\beta,\gamma$ in
 $\m$      
  \begin{enumerate}
\item $\m$ is closed under $\mm_{\m}$, i.e.,  $\alpha \mm_{\m} \beta \in \m$,
\item the associative law
  holds for $\mm_{\m}$, i.e., 
  $\alpha \mm_{\m} ( \beta \mm_{\m} \gamma) =_{\m} (\alpha \mm_{\m} \beta) \mm_{\m} \gamma$, and
\item there exists  $1_{\m} \in \m$ such
  that $\alpha \mm_{\m} 1_{\m} =_{\m} 1_{\m} \mm_{\m} \alpha =_{\m} \alpha$.
  The  element $1_{\m}$ is called \index{identity}\betonen{identity}.
\dend
\end{enumerate}
}
\end{definition}
For simplicity of notation we will henceforth drop the index $\m$ and 
 write $\mm$ respectively $=$ if no confusion is likely to arise.
Furthermore, we will often
 talk about a monoid without mentioning its binary operation
 explicitly.
The monoid operation will often be called multiplication or  addition.
Since the algebraic operation is associative we can omit brackets, hence
 the product $\alpha_1 \mm \ldots \mm \alpha_n$ is uniquely defined.

\begin{example}~\\
{\rm
Let $\Sigma =\{ a_1, \ldots, a_n \}$ be a set of letters.
Then $\Sigma^*$ denotes the set of words over this alphabet.
For two words $u, v \in \Sigma^*$ we define $u \mm v = uv$, i.e.,~the
 word which arises from concatenating the two words $u$ and $v$.
Then $\Sigma^*$ is a monoid with respect to this binary operation and
 its identity element is the empty word, i.e.,~the word containing no letters. 
This monoid is called the \betonen{free monoid} over the alphabet $\Sigma$.
\exaend
}
\end{example}

For some $n$ in $\n$\footnote{In the following $\n$ denotes the
  set of natural numbers including zero and $\n^+ = \n \backslash \{ 0 \}$.}
 the product of $n$ times the same element $\alpha$ is called the
 \index{n-th power of an element}\betonen{n-th power of $\alpha$} and
 will be denoted by $\alpha^n$,
 where $\alpha^0 = 1$.
\begin{definition}~\\
{\rm
An element $\alpha$ of a monoid $\m$ is said to have 
 \index{infinite order}\index{order!infinite}\betonen{infinite
  order}\/
 in case for all $n,m \in \n$, $\alpha^n = \alpha^m$ implies $n=m$.
We say that $\alpha$ has \index{finite
  order}\index{order!finite}\betonen{finite order}\/ 
 in case the set $\{ \alpha^n \mid n \in \n^+ \}$
 is finite and the cardinality of this set is then called the \betonen{order} of
 $\alpha$.
\dend
}
\end{definition}
A subset of a monoid $\m$ which is again a monoid is called a \betonen{submonoid}
 of $\m$.
Other special subsets of monoids are (one-sided) ideals.
\begin{definition}~\\
{\rm
For a subset $S$ of a monoid $\m$ we call 
\begin{enumerate}
\item $\ideal{r}{\m}(S) = \{ \sigma \mm \alpha \mid \sigma \in S, \alpha \in \m \}$
       the right ideal,
\item $\ideal{l}{\m}(S) =\{ \alpha \mm \sigma \mid \sigma \in S, \alpha \in \m \}$
       the left ideal, and
\item $\ideal{}{\m}(S) =\{ \alpha \mm \sigma \mm \alpha' \mid \sigma \in S, \alpha, \alpha' \in \m \}$
       the ideal
\end{enumerate}
generated by $S$ in $\m$.
\dend
}
\end{definition}
A monoid $\m$ is called 
 \index{monoid!Abelian}\index{Abelian}\index{Abelian!monoid}\index{monoid!commutative}\index{commutative}\index{commutative!monoid}\betonen{commutative (Abelian)}\/
 if  we have $\alpha \mm \beta = \beta \mm \alpha$ for all elements $\alpha,\beta$ in $\m$.
A natural example for a commutative monoid are the integers together
 with multiplication or addition. 
Another example which will be of interest later on is the set
 of terms.
\begin{example}\label{exa.terms}~\\
{\rm
Let $X_1, \ldots , X_n$ be a set of (ordered) variables.
Then ${\cal T} = \{ X_1^{i_1} \ldots X_{n\phantom{1}}^{i_n} \mid i_1, \ldots  i_n \in \n\}$ 
 is called the set of 
 \betonen{terms} over these variables.
The multiplication $\mm$ is defined as $X_1^{i_1} \ldots X_{n\phantom{1}}^{i_n} \mm X_1^{j_1} \ldots X_{n\phantom{1}}^{j_n} =
  X_1^{i_1+j_1} \ldots X_{n\phantom{1}}^{i_n+j_n}$.
The identity is the empty term $1_{\cal T} = X_1^{0} \ldots X_{n\phantom{1}}^{0}$.
\exaend
}
\end{example}
A mapping $\phi$ from one monoid $\m_1$ to another monoid $\m_2$
 is called a \index{homomorphism}\betonen{homomorphism},
 if $\phi(1_{\m_1})=
 1_{\m_2}$ and for all $\alpha,\beta$ in $\m_1$, $\phi(\alpha \mm_{\m_1} \beta) =
 \phi(\alpha) \mm_{\m_2} \phi(\beta)$.
In case $\phi$ is surjective we call it an
 \index{epimorphism}\betonen{epimorphism}, in case $\phi$
 is injective a \index{monomorphism}\betonen{monomorphism}\/
 and in case it is both an \index{isomorphism}\betonen{isomorphism}.
The fact that two structures $S_1$, $S_2$ are isomorphic will be denoted by
 $S_1 \cong S_2$.

A monoid is called
 \index{left-cancellative}\betonen{left-cancellative}\/ 
 (respectively \index{right-cancellative}\betonen{right-cancellative}\/) if
for all $\alpha,\beta,\gamma$ in $\m$, $\gamma \mm \alpha = \gamma \mm \beta$
 (respectively $\alpha \mm \gamma = \beta \mm \gamma$) implies $\alpha = \beta$.
In case a monoid is both, left- and right-cancellative,
 it is called \betonen{cancellative}.
In case $\alpha \mm \gamma = \beta$ we say that $\alpha$ is a \betonen{left divisor}\/
of $\beta$  and $\gamma$ is called a \betonen{right
  divisor}\/ of $b$.
If $\gamma \mm \alpha \mm \delta =\beta$ then $\alpha$ is called a \betonen{divisor} of $\beta$.
A special class of monoids fulfill that for all $\alpha,\beta$ in $\m$
 there exist $\gamma,\delta$ in $\m$ such
 that $\alpha \mm \gamma = \beta$ and $\delta \mm \alpha = \beta$, i.e., right and left divisors
 always exist.
These structures are called groups and they can be specified by
 extending the definition of monoids and we do so by
 adding one further axiom.
\begin{definition}~\\
{\rm
A monoid $\m$ together with its binary operation $\mm$ is said to
 form a \index{group}\betonen{group}\/ if additionally 
\begin{enumerate}
\item[4.] for every $\alpha \in \m$ there exists an element $\inv{\alpha} \in \m$ (called
  \index{inverse}\betonen{inverse}\/ of $\alpha$) such
  that $\alpha \mm \inv{\alpha} = \inv{\alpha} \mm
  \alpha = 1$. 
\dend
\end{enumerate}
}
\end{definition}
Obviously, the integers form a group with respect to addition, 
but this is no longer
 true for multiplication.

A subset of a group $\g$ which is again a group is called a 
 \index{subgroup}\betonen{subgroup} of $\m$.
A subgroup $\h$ of a group $\g$ is called \index{normal}\betonen{normal} if
 for each $\alpha$ in $\g$ we have $\alpha\h = \h \alpha$ where 
 $\alpha\h = \{ \alpha \mm \beta \mid \beta \in \h \}$ and
 $\h \alpha = \{ \beta \mm \alpha \mid \beta \in \h\}$.

We end this section by briefly introducing some more algebraic structures that will be used throughout.
\begin{definition}~\\
{\rm
A nonempty set $\rr$ is called an 
 \betonen{(associative) ring (with
 unit element)}\/ if there are two binary operations $+$ (addition) and
 $\rmult$ (multiplication) such
 that for all $\alpha,\beta,\gamma$ in $\rr$
\begin{enumerate}
\item $\rr$ together with $+$ is an Abelian group with zero
  element $0$ and inverse $-\alpha$,
\item $\rr$ is closed under $\rmult$, i.e., $\alpha \rmult \beta \in \rr$, 
\item $\rmult$ is associative, i.e.,  $\alpha \rmult (\beta \rmult \gamma) =
 (\alpha \rmult \beta ) \rmult \gamma$,
\item the distributive laws hold, i.e., $\alpha \rmult (\beta + \gamma) =
     \alpha \rmult \beta + \alpha \rmult \gamma$
       and $(\beta + \gamma) \rmult \alpha = \beta \rmult \alpha + \gamma \rmult \alpha$,
\item there is an element $1 \in \rr$ (called
       \index{unit}\betonen{unit}\/) such that $1 \rmult \alpha = \alpha
       \rmult 1 = \alpha$.
\dend
\end{enumerate}
}
\end{definition}
A ring is called 
 \betonen{commutative (Abelian)}\/ if $\alpha \rmult \beta = \beta \rmult \alpha$ for
 all $\alpha,\beta$ in $\rr$.
The integers together with addition and multiplication are a well-known  example
 of a ring.
Other rings which will be of interest later on are monoid rings.
\begin{example}~\\
{\rm
Let $\z$ be the ring of integers and $\m$ a monoid.
Further let $\z[\m]$ denote the set of all
 mappings $f : \m \myr \z$ where the sets $\supp(f)=\{ \alpha \in \m \mid f(\alpha) \neq 0 \}$
 are finite.
We call $\z[\m]$ the \betonen{monoid ring} of $\m$ over $\z$.
The \betonen{sum} of two elements $f$ and $g$ is denoted by
 $f + g$ where $(f + g)(\alpha) = f(\alpha) + g(\alpha)$.
The \betonen{product} is denoted by $f \rmult g$ where
 $(f \rmult g)(\alpha) = \sum_{\beta \mm \gamma = \alpha} f(\beta) \rmult g(\gamma)$.
}
\end{example}
Polynomial rings are a special case of monoid rings namely over
 the set of terms as defined in Example \ref{exa.terms}.

A ring $\rr$ is said to contain 
 \betonen{zero-divisors}, if
 there exist not necessarily different elements $\alpha, \beta$ in $\rr$ such that $\alpha \neq 0$ and $\beta \neq 0$,
 but $ \alpha \rmult \beta =0$.
Then $\alpha$ is called a left zero-divisor and $\beta$ is called a right
 zero-divisor.
\begin{definition}~\\
{\rm
A commutative ring is called a \betonen{field}\/ if
 its non-zero elements form a group under multiplication.
\dend
}
\end{definition}
Similar to our proceeding in group theory we will now look at
 subsets of a ring $\rr$.
For a subset $U \subseteq \rr$ to be a \betonen{subring} of $R$ with the operations
 $+$ and $\rmult$ it is necessary and sufficient that
\begin{enumerate}
\item $U$ is a subgroup of $(\rr,+)$, i.e., for $a,b \in U$ we have $a-b \in U$, and
\item for all $\alpha,\beta \in U$ we have $\alpha \rmult \beta \in U$.
\end{enumerate}

We will now take a closer look at special subrings that play a
 role similar to normal subgroups in group theory.
\begin{definition}~\\
{\rm
A nonempty subset $\mswab{i}$ of a ring $\rr$ is called a
  \betonen{right (left) ideal}\/ of $\rr$, if
\begin{enumerate}
\item for all $\alpha, \beta \in \mswab{i}$ we have $\alpha - \beta \in \mswab{i}$, and
\item  for every $\alpha \in \mswab{i}$ and $\rho \in \rr$,
  the element $\alpha \rmult
  \rho$ (respectively $\rho \rmult \alpha$)
  lies in $\mswab{i}$.
\end{enumerate}
A subset that is both, a right and a left ideal, is called a
 \betonen{(two-sided) ideal}\/ of $\rr$.
\dend
}
\end{definition}
For each ring  the sets $\{ 0 \}$ and $\rr$ are trivial ideals.
Similar to subgroups, ideals can be described in terms of generating sets. 
\index{ideal!trivial}\index{trivial ideal}
\begin{lemma}\label{lem.idealconstruction}~\\
{\sl
Let $F$ be a non-empty subset of $\rr$. Then
\begin{enumerate}
\item  $\ideal{}{\rr}(F) = \{\sum_{i=1}^{n} \rho_i \rmult \alpha_i \rmult \sigma_i
         \mid \alpha_i \in F, \rho_i , \sigma_i \in \rr, n \in \n \}$ is an ideal of $\rr$,
\item  $\ideal{r}{\rr}(F) = \{\sum_{i=1}^{n} \alpha_i \rmult \rho_i
        \mid \alpha_i \in F, \rho_i \in \rr, n \in \n \}$ is a right ideal
        of $\rr$, and 
\item  $\ideal{l}{\rr}(F) = \{\sum_{i=1}^{n} \rho_i \rmult \alpha_i
        \mid \alpha_i \in F, \rho_i \in \rr, n \in \n \}$ is a left ideal
        of $\rr$.
\ohnebeweis
\end{enumerate}
}
\end{lemma}
Notice that the empty sum $\sum_{i=1}^0 \alpha_i$ is zero.

We will simply write $\ideal{}{}(F)$, $\ideal{r}{}(F)$ and
$\ideal{l}{}(F)$ if the context is clear.
Many algebraic problems for rings are related to ideals and we will
close this section by stating two of them\footnote{For more information on
such problems in the special case of commutative polynomial rings see e.g.~\cite{Bu87}.}.

\begin{tabbing}
XXXXXXXX\=XXXXX \kill
{\bf The Ideal Membership Problem} \\
\\
{\bf Given:} \> An element $\alpha \in \rr$ and a set of elements 
                $F \subseteq \rr$. \\
{\bf Question:} \> Is $\alpha$ in the ideal generated by $F$?
\end{tabbing}

\begin{definition}~\\
{\rm
Two elements $\alpha, \beta \in \rr$ are said to be
 \betonen{congruent modulo} $\ideal{}{}(F)$,
 denoted by $\alpha \equiv_{\ideal{}{}(F)} \beta$, if $\alpha = \beta + \rho$ for some $\rho \in
 \ideal{}{}(F)$, i.e., $\alpha - \beta \in \ideal{}{}(F)$.
\dend
}
\end{definition}
\begin{tabbing}
XXXXXXXX\=XXXXX \kill
{\bf The Congruence Problem} \\
\\
{\bf Given:} \> Two elements $\alpha, \beta \in \rr$ and a set of elements 
                $F \subseteq \rr$. \\
{\bf Question:} \> Are $\alpha$ and $\beta$ congruent modulo the ideal generated by $F$?
\end{tabbing}

Note that both problems can similarly be specified for left and right ideals.

We have seen that a non-empty subset of $\rr$ is an ideal if it is closed under
 addition and closed under multiplication with arbitrary elements
 of $\rr$.
Modules now can be viewed as a natural generalization of the concept of ideals
 to arbitrary commutative groups.
\begin{definition}~\\
{\rm
Let $\rr$ be a ring.
A \betonen{left $\rr$-module} $M$ is an additive commutative group with an additional
 operation $\skm : \rr \times M \myr M$, called scalar multiplication,
 such that for all $\alpha, \beta \in \rr$ and $a,b \in M$, the following hold:
 \begin{enumerate}
 \item $\alpha \skm (a+b) = \alpha \skm a + \alpha \skm b$,
 \item $(\alpha + \beta) \skm a = \alpha \skm a + \beta \skm a$,
 \item $(\alpha \rmult \beta) \skm a = \alpha \skm (\beta \skm a)$, and
 \item $1 \skm a = a$.
 \dend
 \end{enumerate}
}
\end{definition}
We can define right $\rr$-modules and (two-sided) $\rr$-modules (also called $\rr$-bimodules)
 in a similar fashion.

Notice that a (left, right) ideal $\mathfrak{i} \subseteq \rr$ forms a (left, right) $\rr$-module with
 respect to the addition and multiplication in $\rr$.
This obviously holds for the trivial (left, right) ideals $\{ 0 \}$ and $\rr$ of $\rr$.

Another example of (left, right) $\rr$-modules we will study are the finite direct
 products of the ring called free (left, right) $\rr$-modules $\rr^k$, $k \in \rr$.

An additive subset of a (left, right) $\rr$-module is called a \betonen{(left, right) submodule}
 if it is closed under scalar multiplication with elements of $\rr$.
For a subset $F \subseteq M$ let $\langle F \rangle$ denote the
 submodule generated by $F$ in $M$.

\begin{tabbing}
XXXXXXXX\=XXXXX \kill
{\bf The Submodule Membership Problem} \\
\\
{\bf Given:} \> An element $a \in M$ and a set of elements 
                $F \subseteq M$. \\
{\bf Question:} \> $a \in \langle F \rangle$?
\end{tabbing}

Similar to the congruence problem for ideals we can specify the congruence
 problem for submodules as follws:

\begin{definition}~\\
{\rm
Two elements $a, b \in \rr$ are said to be
 \betonen{congruent modulo} the submodule $\langle F \rangle$ for some
 $F \subseteq M$,
 denoted by $a \equiv_{\langle F \rangle} b$, if $a - b \in \langle F \rangle$.
\dend
}
\end{definition}

\begin{tabbing}
XXXXXXXX\=XXXXX \kill
{\bf The Congruence Problem for submodules} \\
\\
{\bf Given:} \> Two elements $a, b\in \rr$ and a set of elements 
                $F \subseteq M$. \\
{\bf Question:} \> $a \equiv_{\langle F \rangle} b$?
\end{tabbing}

%%% Local Variables: 
%%% mode: latex
%%% TeX-master: "testlauf"
%%% End: 

%% file: rewriting.tex
%%%%%%%%%%%%%%%%%%%%%%%%%%%%%%%%%%%%%%%%%%%%%%%%%%%%%%%%%%%%%%%%%%%%%%
%%%%%%%%%%%%%%%%% Section: Reduction %%%%%%%%%%%%%%%%%%%%%%%%%%%%%%%%%
%%%%%%%%%%%%%%%%%%%%%%%%%%%%%%%%%%%%%%%%%%%%%%%%%%%%%%%%%%%%%%%%%%%%%%
\section{The Notion of Reduction}\label{section.rewriting}
%\rechterRand{The Notion of Reduction}
%
This section summarizes some important notations and definitions of
reduction relations and basic properties related to them, as can be 
found more explicitly for example in the work of Huet or Book and Otto
(\cite{Hu80,Hu81,BoOt93}).

Let ${\cal E}$ be a set of elements and $\myr$ a binary relation  on
${\cal E}$ called 
 \index{reduction!relation}\index{reduction}\betonen{reduction}.
For $a,b \in {\cal E}$ we will write $a \red{}{\myr}{}{} b$ in case
$(a,b) \in\;\; \myr$.
A pair $({\cal E},\myr)$ will be called a \index{reduction system}\betonen{reduction system}.
Then we can expand the binary relation as follows:
\begin{tabbing}
%\hspace{8mm}\=\hspace{4cm}\=\hspace{6cm}\= \kill
\hspace{3mm}\=\hspace{1cm}\=\hspace{5mm}\=\hspace{2.3cm}\=\hspace{6cm}\= \kill
\> $\red{0}{\myr}{}{}$\> \>  \> denotes the identity on ${\cal E}$, \\
\> $\red{}{\longleftarrow}{}{}$\phantom{$\red{+}{\lr}{}{}$} \> \> \> denotes the inverse relation for $\myr$, \\
\> $\red{n+1}{\myr}{}{}$\> $:=$ \> $\red{n}{\myr}{}{} \circ \myr$ \> where
                $\circ$ denotes composition of relations and $n \in
                \n$, \\
\> $\red{\leq n}{\myr}{}{}$\> $:=$ \> $\;\!\bigcup_{0 \leq i \leq n}\!\!
\red{i}{\myr}{}{}$, \\
\> $\red{+}{\myr}{}{}$ \> $:=$ \> $\;\!\bigcup_{n>0}\! \red{n}{\myr}{}{}$ \>  denotes the transitive closure
                                               of $\myr$, \\
\> $\red{*}{\myr}{}{}$ \> $:=$ \> $\red{+}{\myr}{}{} \cup \red{0}{\myr}{}{}$ \>  denotes 
   the reflexive transitive closure of $\myr$, \\
\> $\red{}{\lr}{}{}$ \> $:=$ \> $\;\!\longleftarrow \cup \myr$\phantom{$\red{+}{\lr}{}{}$} \>  denotes the
   symmetric closure of $\myr$, \\
\> $\red{+}{\lr}{}{}$ \> \>\>  denotes the symmetric transitive closure
   of $\myr$, \\
\> $\red{*}{\lr}{}{}$\phantom{$\red{+}{\lr}{}{}$} \>\>\>  denotes the reflexive symmetric transitive closure
   of $\myr$.   
\end{tabbing}
%Obviously $\red{*}{\lr}{}{}$ is an equivalence
%relation on ${\cal E}$ and $$\red{*}{\myr}{}{}$ can be viewed as a
%reduction relation on ${\cal E}$.
%
A well-known decision problem related to a reduction system is the
word problem.
\begin{definition}~\\
{\rm
The \index{word problem}
 \betonen{word problem} for a reduction system $({\cal E}, \myr)$ is to decide for 
 $a,b$ in ${\cal E}$, whether $a
 \red{*}{\lr}{}{} b$ holds.
\dend
}
\end{definition}
Instances of this problem are well-known in the literature and undecidable in
general.
In the following we will outline sufficient conditions such that a reduction system $({\cal E}, \myr)$ has
solvable word problem.

An element $a \in {\cal E}$ is said to be
 \index{reducible}\betonen{reducible}\/
 (with respect to $\myr$) if there
 exists an element $b \in {\cal E}$ such that $a \myr b$.
All elements $b \in {\cal E}$ such that $a \red{*}{\myr}{}{} b$ are
 called \index{successor}\betonen{successors}\/
 of $a$ and in case  $a \red{+}{\myr}{}{} b$
 they are called 
 \index{proper!successor}\index{successor!proper}\betonen{proper successors}.
An element which has no proper successors  is called
 \index{irreducible}\betonen{irreducible}.
In case $a \red{*}{\myr}{}{} b$ and $b$ is irreducible, $b$ is called a
 \index{normal form}\betonen{normal form}\/ of $a$.
Notice that for an element $a$ in ${\cal E}$ there can be no, one or many normal
forms.
\begin{definition}~\\
{\rm
A reduction system $({\cal E}, \myr)$ is said to be 
 \betonen{Noetherian}
 (or \betonen{terminating}\/) in case
 there are no infinitely
 descending reduction chains $a_0 \myr a_1
 \myr \ldots\; $, with  $a_i \in {\cal E}$, $i \in \n$.
\dend
}
\end{definition}
In case a reduction system $({\cal E}, \myr)$ is Noetherian every element in ${\cal E}$ has
at least one normal form.
\begin{definition}~\\
{\rm
A reduction system $({\cal E}, \myr)$ is called 
 \index{reduction system!confluent}\index{confluence}\betonen{confluent},
 if for all $a, a_1, a_2
 \in {\cal E}$, $a \red{*}{\myr}{}{} a_1$ and $a \red{*}{\myr}{}{} a_2$
implies the existence of $a_3 \in {\cal E}$ such that $a_1
\red{*}{\myr}{}{} a_3$ and $a_2 \red{*}{\myr}{}{} a_3$, and $a_1$, $a_2$ are called 
\index{joinable}\betonen{joinable}.
\dend
}
\end{definition}
In case a reduction system $({\cal E}, \myr)$ is confluent every element has at most one
normal form.
We can combine these two properties to give sufficient conditions for
the solvability of the word problem.
\begin{definition}~\\
{\rm
A reduction system $({\cal E}, \myr)$ is said to be 
 \index{reduction system!complete}\index{complete}\betonen{complete}\/
 (or 
  \index{reduction system!convergence}\index{convergent}\betonen{convergent}\/) in case it
 is both, Noetherian and confluent.
\dend
}
\end{definition}
Complete reduction systems with effective or computable\footnote{By effective or computable we mean that 
 given an element we can always construct a successor in case one exists.} 
 reduction relations have solvable word problem, as every
 element has a unique normal form and two elements are equal if and
 only if their normal forms are equal.
Of course we cannot always expect $({\cal E}, \myr)$ to be complete.
Even worse, both properties -- termination and confluence -- are undecidable in general.
Nevertheless, there are weaker conditions which guarantee
 completeness.
\begin{definition}~\\
{\rm
A reduction system $({\cal E}, \myr)$ is said to be 
 \index{local confluence}\index{confluence!local}\betonen{locally confluent}, 
 if for all $a, a_1, a_2
 \in {\cal E}$, $a \red{}{\myr}{}{} a_1$ and $a \red{}{\myr}{}{} a_2$
 implies the existence of an element $a_3 \in {\cal E}$ such that $a_1
 \red{*}{\myr}{}{} a_3$ and $a_2 \red{*}{\myr}{}{} a_3$.
\dend
} 
\end{definition}
I.e.~local confluence is a special instance of confluence, namely a
 {\em localization} of confluence to one-reduction-step successors of elements only.
The next lemma gives an important connection between local confluence and
 confluence.
\index{Newman's Lemma}
\begin{lemma}[Newman]~\\
{\sl
Let $({\cal E}, \myr)$ be a Noetherian reduction system.
%\\
Then $({\cal E}, \myr)$ is confluent if and only if $({\cal E}, \myr)$ is
locally confluent.
\lemend
}
\end{lemma}
To prove Newman's lemma we need  the concept of Noetherian
 induction which is based on the following definition.
\begin{definition}~\\
{\rm
Let $({\cal E}, \myr)$ be a reduction system.
A predicate ${\cal P}$ on ${\cal E}$ is called \betonen{$\myr$-complete}, in case for
 every $a \in {\cal E}$ the following implication holds:  if ${\cal
   P}(b)$ is true for all proper successors of $a$,
 then ${\cal P}(a)$ is true.
\dend
}
\end{definition}
{\sl
{\bf The Principle of Noetherian Induction}:\index{Noetherian!induction} \\
In case  $({\cal E}, \myr)$ is a Noetherian reduction system and ${\cal P}$ is
a predicate that is $\myr$-complete, then  for all $a \in {\cal E}$,
${\cal P}(a)$ is true.
}

\Ba{of Newman's lemma}~\\
Suppose, first, that  the reduction system $({\cal E}, \myr)$ is confluent.
This immediately implies the local
 confluence of  $({\cal E}, \myr)$ as a special case.
To show the converse, since  $({\cal E}, \myr)$ is Noetherian we can apply the principle of
Noetherian induction to the following predicate:
\begin{center}
${\cal P}(a)$ \\
if and only if \\
for all $a_1,a_2 \in {\cal E}$, $a
\red{*}{\myr}{}{} a_1$ and $a \red{*}{\myr}{}{} a_2$ implies that $a_1$
 and $a_2$ are joinable.
\end{center}
All we have to do now is to show that ${\cal P}$ is $\myr$-complete.
Let $a \in {\cal E}$ and let ${\cal P}(b)$ be true for all proper successors $b$
of $a$.
We have to prove that ${\cal P}(a)$ is true.
Suppose $a \red{*}{\myr}{}{} a_1$ and $a \red{*}{\myr}{}{} a_2$.
In case $a = a_1$ or $a = a_2$ there is nothing to show.
Therefore, let us assume $a \neq a_1$ and $a \neq a_2$, i.e.,
 $a \myr \tilde{a}_1 \red{*}{\myr}{}{} a_1$ and $a \myr \tilde{a}_2
 \red{*}{\myr}{}{} a_2$. 
Then we can deduce the following  figure
%\vspace{5mm}%gepfuscht
\begin{diagram}[size=1.6em]
  & & & & a  & & & &  \\
  & & &\ldTo &    &\rdTo & & &  \\
  & &\tilde{a}_1 & &    & &\tilde{a}_2 & &  \\
  &\ldTo^{*} & &\rdTo^{*} &     &\ldTo^{*}& &\rdTo^{*} &  \\
a_1  & & & & b_0 & & & & a_2  \\
  &\rdTo^{*} & &\ldTo^{*} & & & &\ldTo(4,4)^{*}  &  \\
  & & b_1 & & & && & \\
  & & &\rdTo^{*} & & & & &  \\
  & & & & b  \\
\end{diagram}
\auskommentieren{
\begin{center}
\setlength{\unitlength}{0.0125in}%
\begin{picture}(160,171)(240,540)
\thicklines
\put(315,695){\vector(-1,-1){ 28}}
\put(275,655){\vector(-1,-1){ 28}}
\put(355,655){\vector(-1,-1){ 28}}
\put(315,615){\vector(-1,-1){ 28}}
\put(325,695){\vector( 1,-1){ 28}}
\put(285,655){\vector( 1,-1){ 28}}
\put(245,615){\vector( 1,-1){ 28}}
\put(365,655){\vector( 1,-1){ 28}}
\put(285,575){\vector( 1,-1){ 28}}
\put(395,615){\vector(-1,-1){ 68}}
\put(316,700){\makebox(0,0)[lb]{\raisebox{0pt}[0pt][0pt]{\twlrm $a$}}}
\put(276,660){\makebox(0,0)[lb]{\raisebox{0pt}[0pt][0pt]{\twlrm $\tilde{a}_1$}}}
\put(356,660){\makebox(0,0)[lb]{\raisebox{0pt}[0pt][0pt]{\twlrm $\tilde{a}_2$}}}
\put(236,620){\makebox(0,0)[lb]{\raisebox{0pt}[0pt][0pt]{\twlrm $a_1$}}}
\put(316,620){\makebox(0,0)[lb]{\raisebox{0pt}[0pt][0pt]{\twlrm $b_0$}}}
\put(396,620){\makebox(0,0)[lb]{\raisebox{0pt}[0pt][0pt]{\twlrm $a_2$}}}
\put(276,580){\makebox(0,0)[lb]{\raisebox{0pt}[0pt][0pt]{\twlrm $b_1$}}}
\put(316,540){\makebox(0,0)[lb]{\raisebox{0pt}[0pt][0pt]{\twlrm $b$}}}
\put(255,645){\makebox(0,0)[lb]{\raisebox{0pt}[0pt][0pt]{\twlrm $*$}}}
\put(299,645){\makebox(0,0)[lb]{\raisebox{0pt}[0pt][0pt]{\twlrm $*$}}}
\put(335,645){\makebox(0,0)[lb]{\raisebox{0pt}[0pt][0pt]{\twlrm $*$}}}
\put(379,645){\makebox(0,0)[lb]{\raisebox{0pt}[0pt][0pt]{\twlrm $*$}}}
\put(259,605){\makebox(0,0)[lb]{\raisebox{0pt}[0pt][0pt]{\twlrm $*$}}}
\put(295,605){\makebox(0,0)[lb]{\raisebox{0pt}[0pt][0pt]{\twlrm $*$}}}
\put(299,565){\makebox(0,0)[lb]{\raisebox{0pt}[0pt][0pt]{\twlrm $*$}}}
\put(355,585){\makebox(0,0)[lb]{\raisebox{0pt}[0pt][0pt]{\twlrm $*$}}}
\end{picture}
\end{center}}
where $b_0$ exists, as $({\cal E}, \myr)$ is locally confluent and
$b_1$ and $b$  exist by our induction hypothesis since $a_1$, $b_0$ as well
as $a_2$, $b_1$ are proper successors of $a$.
Hence $a_1$ and $a_2$ must be joinable, i.e., the reduction system $({\cal E}, \myr)$ is confluent.
\\
\qed
Therefore, if the reduction system is terminating,
 a check for confluence can be reduced to a check for local
 confluence.
The concept of completion then is based
 on two steps:
 \begin{enumerate}
 \item Check the system for local confluence. \\
       If it is locally confluent, then it is also complete.
 \item Add new relations arising from situations where the system is not locally
       confluent.
 \end{enumerate}
For many reduction systems, e.g.~string rewriting systems or term rewriting systems,
 the check for local confluence again can be localized, often to finite test sets 
 of so-called critical pairs.
The relations arising from such critical situations are either confluent or give rise to new
 relations which stay within the congruence described by the reduction system.
Hence adding them in order to increase the descriptive power of the reduction system
 is correct.
This can be done until a complete set is reached. 
If fair strategies are used in the test for local confluence, the limit system will be 
 complete.

We close this section by providing sufficient conditions to ensure a reduction system $({\cal E}, \myr)$ to be
Noetherian.
\begin{definition}~\\
{\rm
A binary relation $\succeq$ on  a set $M$
 is said to be a 
 \index{partial ordering}\index{ordering!partial}\betonen{partial ordering},
 if for all $a,b,c$ in $M$:
\begin{enumerate}
\item $\succeq$ is reflexive, i.e., $a \succeq a$,
\item $\succeq$ is transitive, i.e., $a \succeq b$ and $b \succeq c$
  imply $a \succeq c$, and
\item $\succeq$ is anti-symmetrical, i.e., $a \succeq b$ and $b \succeq
  a$ imply $a=b$.
\dend
\end{enumerate}
}
\end{definition}
A partial ordering is called 
 \index{total ordering}\index{ordering!total}\betonen{total},
 if for all $a, b \in M$
 either $a \succeq b$ or $b \succeq a$ holds.
Further a partial ordering $\succeq$ defines a transitive irreflexive
 ordering $\succ$, where $a \succ b$ if and only if $a \succeq b$ and $a
\neq b$, which is often called a \betonen{proper} or \betonen{strict} ordering.
We call a partial ordering $\succeq$ 
 \index{well-founded
   ordering}\index{ordering!well-founded}\betonen{well-founded},
 if the corresponding strict
 ordering $\succ$ allows  no infinite descending chains $a_0 \succ a_1 \succ
 \ldots\;$, with $a_i \in M$, $i \in \n$.
Now we can give a sufficient condition for a reduction system to be terminating. 
\begin{lemma}~\\
{\sl
Let $({\cal E}, \myr)$ be a reduction system and
suppose there exists a  partial ordering  $\succeq$ on ${\cal E}$ which is
well-founded such that $\myr\;\;\subseteq\;\;\succ$.
Then $({\cal E}, \myr)$ is Noetherian.
}
\end{lemma}
\Ba{}~\\
Suppose the reduction system  $({\cal E}, \myr)$ is not Noetherian.
Then there is an infinite sequence $a_0 \myr a_1 \myr \ldots\;$,  $a_i
\in {\cal E}$, $i \in \n$.
As  $\myr \;\subseteq\; \succ$ this sequence gives us an  infinite sequence
$a_0 \succ a_1 \succ \ldots\;$, with $a_i \in {\cal E}$, $i \in \n$ 
 contradicting our assumption that
$\succeq$ is well-founded on ${\cal E}$.
\\
\qed
%

%%% Local Variables: 
%%% mode: latex
%%% TeX-master: "testlauf"
%%% TeX-master: "testlauf"
%%% End: 

%% file: groebner.tex
%%%%%%%%%%%%%%%%%%%%%%%%%%%%%%%%%%%%%%%%%%%%%%%%%%%%%%%%%%%%%%%%%%%%%%%%
%%%%%%%%%%%%%%%%% Section: Groebner Bases %%%%%%%%%%%%%%%%%%%%%%%%%%%%%%
%%%%%%%%%%%%%%%%%%%%%%%%%%%%%%%%%%%%%%%%%%%%%%%%%%%%%%%%%%%%%%%%%%%%%%%%
\section{Gr\"obner Bases in Polynomial
  Rings}\label{section.buchberger}
%\RechterRand{Gr\"obner Bases in Polynomial Rings}
%
The main interest in this section is the study of ideals in polynomial
 rings over fields. 
Let $\myk [X_1, \ldots , X_n]$ denote a polynomial ring over the
 (ordered) variables $X_1, \ldots , X_n$ and the computable field $\myk$.
By ${\cal T} = \{ X_1^{i_1} \ldots X_{n\phantom{1}}^{i_n} \mid i_1, \ldots  i_n \in \n\}$ 
 we define the set of 
 \betonen{terms} in this structure.
A polynomial then is a formal sum $\sum_{i=1}^{n} \alpha_i \skm t_i$
 with non-zero coefficients $\alpha_i \in \myk \backslash \{ 0 \}$ and terms $t_i \in \myt$.
The products $\alpha \skm t$ for $\alpha \in \myk$, $t \in \myt$ are called
 monomials and will often be denoted as $m = \alpha \skm t$.
We recall that a subset
 $F$ of $\myk[X_1, \ldots , X_n]$
 generates an ideal
 $\ideal{}{}(F) = \{ \sum_{i=1}^{k} f_i \mrm g_i \mid k
 \in \n, f_i \in F, g_i \in \myk[X_1, \ldots, X_n] \}$
 and $F$ is called a basis of this ideal. 
It was shown by Hilbert using non-constructive arguments
 that every ideal in $\myk [X_1, \ldots ,X_n]$ in fact
 has a finite basis, but such a generating set need not allow
 algorithmic solutions for the membership or congruence problem
 related to the ideal as we have seen in the introduction.
It was Buchberger who developed a special type of basis, namely the
 Gr\"obner basis, which allows algorithmic solutions for several
 algebraic problems concerning ideals.
He introduced a reduction relation to $\myk [X_1, \ldots , X_n]$ by transforming
 polynomials into ``rules'' and gave a terminating procedure to
 ``complete'' an ideal basis interpreted as a reduction system.
This procedure is called Buchberger's algorithm in the literature. 
We will give a sketch of his approach below.

Let $\succeq$ be a total well-founded ordering on the set of terms
 ${\cal T}$, which is
 admissible, i.e., $t \succeq 1$, and $s \succ t$ implies $s \mm u \succ
 t \mm u$ for all $s,t,u$ in ${\cal T}$.
The latter property is called compatibility with the multiplication $\mm$.
In this context $\mm$ denotes the multiplication in ${\cal T}$, i.e.,
 $X_1^{i_1} \ldots X_{n\phantom{1}}^{i_n} \mm X_1^{j_1} \ldots X_{n\phantom{1}}^{j_n} =
  X_1^{i_1+j_1} \ldots X_{n\phantom{1}}^{i_n+j_n}$.
With respect to this multiplication we say that a term $s=X_1^{i_1}
\ldots X_{n\phantom{1}}^{i_n}$ divides a term $t=X_1^{j_1} \ldots X_{n\phantom{1}}^{j_n}$, if
for all $1 \leq l \leq n$ we have $i_l \leq j_l$.
The \betonen{least common
  multiple} $\lcm(s,t)$ of the terms $s$ and $t$  is the term
 $X_1^{\max \{i_1,j_1\}} \ldots X_{n\phantom{1}}^{\max \{i_n,j_n\}}$.
Note that ${\cal T}$ can be interpreted
 as the free commutative monoid
 generated by $X_1, \ldots, X_n$ with the same multiplication $\mm$ as
 defined above and identity $1 = X_1^{0} \ldots
 X_{n\phantom{1}}^{0}$ (recall Example \ref{exa.terms}).
We proceed to give an example for a total well-founded admissible
ordering on the set of terms ${\cal T}$.
\begin{example}~\\
{\rm
A \betonen{total degree ordering} $\succ$
 on ${\cal T}$ is specified as follows: 
$X_1^{i_1} \ldots X_{n\phantom{1}}^{i_n} \succ  X_1^{j_1} \ldots X_{n\phantom{1}}^{j_n}$
 if and only if
$ \sum_{s=1}^{n} i_s > \sum_{s=1}^{n} j_s$
 or
$\sum_{s=1}^{n} i_s = \sum_{s=1}^{n} j_s$ and there exists $k$
such that  $i_k > j_k$ and  $i_s=j_s, 1 \leq s < k$.
\dend
}
\end{example}
Henceforth, let $\succeq$ denote a total admissible
ordering on ${\cal T}$ which is of course well-founded.
%Man kann keine beliebige absteigende Kette bilden, 
%da sonst durch Dickson Widerspruch zu admissible entstehen wuerde!
%
\begin{definition}~\\
{\rm
Let $p = \sum_{i=1}^{k} \alpha_i \skm t_i$ be a non-zero polynomial in $\myk[X_1, \ldots , X_n]$ such
that $\alpha_i \in \myk^* = \myk \backslash \{ 0 \}$, $t_i \in {\cal T}$ and $t_1 \succ \ldots \succ t_n$.
Then we let $\hm(p) = \alpha_{1} \cd t_{1}$ denote the
 \index{head monomial}\betonen{head monomial},
 $\hterm(p) = t_{1}$ the \index{head term}\betonen{head term}\/
 and $\hc(p) = \alpha_{1}$ the 
 \index{head coefficient}\betonen{head coefficient}\/ of $p$.
 $\reductum(p) = p - \hm(p)$ stands for the
 \index{reduct}\betonen{reductum}\/ of $p$.
 We call $p$  \betonen{monic}\index{monic} in case $\hc(p) = 1$.
These definitions can be extended to sets $F$ of polynomials by setting
 $\hterm(F) = \{ \hterm(f) \mid f \in F \}$,
 $\hc(F) = \{ \hc(f) \mid f \in F \}$, respectively
 $\hm(F) = \{ \hm(f) \mid f \in F \}$. 
\mbox{ }\dend
}
\end{definition}
Using the notions of this definition we can recursively extend 
 $\succeq$ from ${\cal T}$ to a partial well-founded admissible ordering
 $\geq$ on $\myk[X_1,\ldots, X_n]$.
\begin{definition}\label{def.ordering.polyring}~\\
{\rm
Let $p,q$ be two polynomials in $\myk[X_1,\ldots, X_n]$.
Then we say $p$ is \betonen{greater} than $q$ with respect to a total well-founded admissible
ordering $\succeq$ on ${\cal T}$, i.e., $p > q$, if
\begin{enumerate}
\item $\hterm(p) \succ \hterm(q)$ or
\item $\hm(p) = \hm(q)$ and $\reductum(p) > \reductum(q)$.
\dend
\end{enumerate}
}
\end{definition}
Now one first specialization of right ideal bases in terms of the representations
 they allow can be given according to standard representations as introduced
 e.g.~in \cite{BeWe92} for polynomial rings over fields.
\begin{definition}\label{def.standard.rep.polyring}~\\
{\rm
Let $F$ be a set of polynomials in $\myk[X_1,\ldots, X_n]$
 and $g$ a non-zero polynomial in $\ideal{}{}(F)\subseteq \myk[X_1,\ldots, X_n]$.
A representations of the form
 \begin{eqnarray}
  g & = & \sum_{i=1}^n f_i \rmult m_i, 
   f_i \in F, m_i = \alpha_i \skm t_i, \alpha_i \in \myk, t_i \in \myt, n \in \n\label{eqn.standard.poly}
 \end{eqnarray}
 where additionally $\hterm(g) \succeq \hterm(f_i \rmult m_i)$
 holds for
 $1 \leq i \leq n$ is called a \betonen{standard representation}
 of $g$ in terms of $F$.
If every $g \in \ideal{}{}(F) \backslash \{ 0 \}$ has such a representation in terms of $F$,
 then $F$ is called a \betonen{standard basis} of $\ideal{}{}(F)$.
\dend
}
\end{definition}
What distinguishes an arbitrary representation from a
 standard representation is the fact that the former may
 contain polynomial multiples
 with head terms larger than the head term of the represented polynomial.
For example\label{exa.qxyz} let $f_1 = X_1 + X_2$, $f_2 = X_1 + X_3$ and
 $F = \{ f_1, f_2 \}$ in $\q[X_1, X_2]$ with $X_1 \succ X_2 \succ X_3$.
Then for the polynomial $g = X_2 - X_3$ we have the representation
 $g = f_1 + (-1) \skm f_2$ which is no standard one as 
 $\hterm(g) = X_2 \pred \hterm (f_1) = \hterm(f_2) = X_1$.
Obviously the larger head terms have to vanish in the sum.
Therefore, in order to change an arbitrary representation into one fulfilling our
 additional condition (\ref{eqn.standard.poly})
 we have to deal with special sums of polynomials related to such situations.
\begin{definition}\label{def.critical.situations.polynomialring}~\\
{\rm
Let $F$ be a set of polynomials in $\myk[X_1,\ldots, X_n]$
 and $t$ an element in $\myt$.
Then we define the set of \betonen{critical situations}
 ${\cal C}(t,F)$ related to $t$ and $F$ to contain all tuples of the form
 $(t, f_1, \ldots, f_k, m_1, \ldots, m_k)$, $k \in \n$, $f_1, \ldots, f_k \in F$\footnote{Notice that $f_1, \ldots, f_k$ are not
 necessarily different polynomials from $F$.}, $m_i = \alpha_i \skm t_i$,
  such that
 \begin{enumerate}
 \item $\hterm(f_i \rmult m_i) = t$, $1 \leq i \leq k$, and
 \item $\sum_{i=1}^k \hm(f_i \rmult m_i) = 0$.
 \end{enumerate}
We set ${\cal C}(F) = \bigcup_{t \in \myt} {\cal C}(t, F)$.
\dend
}
\end{definition}
In our example the tuple $(X_1, f_1, f_2, 1, -1)$ is an elements of
 the critical set ${\cal C}(X_1, F)$.
We can characterize standard bases using these special sets.
\begin{theorem}\label{theo.standard.basis.polyring}~\\
{\sl
Let $F$ be a set of polynomials in $\myk[X_1,\ldots, X_n] \backslash \{ 0\}$.
Then $F$ is a standard basis of $\ideal{}{}(F)$ if and only if
 for every tuple 
 $(t, f_1, \ldots, f_k, m_1, \ldots, m_k)$
 in ${\cal C}(F)$ as specified in Definition
 \ref{def.critical.situations.polynomialring}
 the polynomial $\sum_{i=1}^k f_i \rmult m_i$
 has a standard representation
 with respect to $F$.
\theoend
}
\end{theorem}
\Ba{}~\\
In case $F$ is a standard basis
 since these polynomials are all elements of $\ideal{}{}(F)$ they must
 have standard representations with respect to $F$.
\\
To prove the converse, it remains to show that every element in
 $\ideal{}{}(F)$ has a standard representation
 with respect to $F$.
Hence, let $g = \sum_{j=1}^m f_{j} \rmult m_j$ be an arbitrary
 representation of a non-zero  polynomial $g\in \ideal{}{}(F)$ such that
 $f_j \in F$, and $m_j = \alpha_j \skm t_j$ with
 $\alpha_j \in \myk$, $t_j \in \myt$.
Depending on this  representation of $g$ and the
 well-founded total ordering $\succeq$ on $\myt$ we define
 $t = \max_{\succeq} \{ \hterm(f_{j} \rmult t_{j}) \mid 1\leq j \leq m \}$ and
 $K$ as the number of polynomials $f_j \rmult t_j$ with head term $t$.
%\\
Then $t \succeq \hterm(g)$ and 
 in case $\hterm(g) = t$ this immediately implies that this representation is
 already a standard representation. 
%\\
Else we proceed by induction
 on the term $t$.
%\\
Without loss of generality let $f_1, \ldots, f_K$ be the polynomials
 in the corresponding representation
 such that  $t=\hterm(f_i \rmult t_i)$, $1 \leq i \leq K$.
Then the tuple $(t, f_1, \ldots, f_K, m_1, \ldots, m_K)$
 is in ${\cal C}(F)$ and let $h = \sum_{i=1}^K f_i \rmult m_i$.
%\\
We will now change our representation of $g$ in such a way that for the new
 representation of $g$ we have a smaller maximal term.
%\\
Let us assume $h$ is not $0$\footnote{In case  $h =0$,
 just substitute the empty sum for the representation of $h$
 in the equations below.}. 
%\\
By our assumption, $h$ has a standard representation
 with respect to $F$, say $\sum_{j=1}^n h_j \rmult n_j$, 
 where $h_j \in F$, and $n_j = \beta_j \skm s_j$ with
 $\beta_j \in \myk$, $s_j \in \myt$ 
 and  all terms occurring in the sum are bounded by
 $t \succ \hterm(h)$ as $\sum_{i=1}^K \hm(f_i \rmult m_i) = 0$.
%\\
This gives us: 
\begin{eqnarray}
  g   & = & \sum_{i=1}^K f_i \rmult m_i  +
             \sum_{i=K+1}^m f_i \rmult m_i
             \nonumber\\                                                           
  & = & \sum_{j=1}^n h_j \rmult n_j + \sum_{i=K+1}^m f_i \rmult m_i
              \nonumber
\end{eqnarray}
which is a representation of $g$ where the maximal term is smaller than $t$.
\\
\qed
In fact for the case of polynomial rings over fields one
 can show that it is sufficient to consider critical sets for
 subsets of $F$ of size 2 and we can restrict the terms to the
 least common multiples of the head terms of the respective two polynomials. 
These sets then correspond to the concept of s-polynomials used to characterize
 Gr\"obner bases which
 will be introduced later on.

Reviewing our example on page \pageref{exa.qxyz} we find that the set 
 $F = \{ X_1 + X_2, X_1 + X_3 \}$ is no standard basis as
 the polynomial $g = X_2 - X_3$ has no standard representation
 although it is an elements of $\ideal{}{}(F)$.
However the set $F \cup \{ g \}$ then is a standard basis of 
 $\ideal{}{}(F)$.

In the literature standard representations in $\myk[X_1,\ldots, X_n]$
 are closely related to reduction relations based on the divisibility of terms
 and standard bases are in fact Gr\"obner bases.
Here we want to introduce Gr\"obner bases in terms of rewriting.
Hence we continue by introducing the concept of reduction to
 $\myk[X_1,\ldots, X_n]$.

We can split a non-zero polynomial $p$ into a \index{rule}\betonen{rule}
 $\hm(p) \myr - \reductum(p)$
       and we have $\hm(p) > - \reductum(p)$.
Therefore, a set of polynomials gives us a binary relation $\myr$ on
 $\myk[X_1, \ldots , X_n]$ which induces a one-step reduction relation as follows.
\begin{definition}\label{def.buchberger.red}~\\
{\rm 
Let $p, f$ be two polynomials  in $\myk[X_1, \ldots , X_n]$. 
We say $f$ 
 \index{Buchberger's
   reduction}\index{reduction!Buchberger}\betonen{reduces}
 $p$ to $q$ at a monomial $m = \alpha \skm t$ of $p$
 in one step, denoted by $p \red{}{\myr}{b}{f} q$, if
\begin{enumerate}
\item[(a)] $\hterm(f) \mm u = t$ for some $u \in {\cal
    T}$, i.e., $\hterm(f)$ divides $t$, and
\item[(b)] $q = p - \alpha \skm \hc(f)^{-1} \skm f \mrm u$.
\end{enumerate}
We write $p \red{}{\myr}{b}{f}$ if there is a polynomial $q$ as defined
above and $p$ is then called reducible by $f$. 
Further, we can define $\red{*}{\myr}{b}{}, \red{+}{\myr}{b}{}$, and
 $\red{n}{\myr}{b}{}$ as usual.
Reduction by a set $F \subseteq \myk[X_1, \ldots , X_n]$ is denoted by
 $p \red{}{\myr}{b}{F} q$ and abbreviates $p \red{}{\myr}{b}{f} q$
 for some $f \in F$,
 which is also written as  $p \red{}{\myr}{b}{f \in F} q$.
\dend
}
\end{definition} 
Note that if $f$ \betonen{reduces}
 $p$ to $q$ at a monomial $m = \alpha \skm t$ then $t$ is no longer among
 the terms of $q$.
We call a set of polynomials $F \subseteq \myk[X_1, \ldots , X_n]$
 \betonen{interreduced}, if no $f \in F$ is reducible by a polynomial
 in $F \backslash \{ f \}$.

In the classical case of polynomial rings over fields the existence
 of a standard representation
 for a polynomial immediately implies reducibility of the head monomial
 of the polynomial by any
 reduction relation based on divisibility of terms, hence by the reduction relation
 defined here.
This is due to the fact that if a polynomial $g$ has a standard representation
 in terms of a set of polynomials $F$ for at least one polynomial
 $f$ in $F$ and some term $t$ in ${\cal T}$ we have
 $\hterm(g) = \hterm(f \rmult t) = \hterm(f) \mm t$
 and hence $g$ is reducible at the monomial $\hm(g)$ by $f$.
Notice that this is no longer true for polynomial rings over the integers.
Let $F = \{ 3 \skm X^2 + X, 2 \skm X^2 + X \}$ be a subset
 of $\z[X]$.
Then the polynomial $g = (3 \skm X^2 + X) - (2 \skm X^2 + X)
 = X^2$
 has a standard representation in terms of $F$ but neither $3 \skm X^2$
 nor $2 \skm X^2$ are divisors of the monomial $X^2$ as neither $3$ nor $2$
 devide $1$ in $\z$. 

Notice that we have $\myr \;\subseteq\;\; >$ and indeed one can show that
 our reduction relation  on $\myk[X_1, \ldots, X_n]$ is Noetherian.
Therefore, we can restrict ourselves to ensuring local confluence when
 describing a completion procedure to compute Gr\"obner bases later on.
But first we have to provide a definition of Gr\"obner bases in the 
 context of rewriting.
\begin{definition}\label{def.buchberger.gb}~\\
{\rm
A set $G \subseteq \myk[X_1,\ldots,X_n]$ is said to be a 
 \index{Gr\"obner basis!Buchberger}\index{Gr\"obner basis}\betonen{Gr\"obner basis}
 of the ideal it generates, if
\begin{enumerate}
\item $\red{*}{\lr}{b}{G} = \;\;\equiv_{\ideal{}{}(G)}$, and
\item $\red{}{\myr}{b}{G}$ is confluent.
\dend
\end{enumerate}
}
\end{definition}
The first statement expresses that the reduction relation describes the ideal
congruence.
It holds for any basis of an ideal in $\myk[X_1,\ldots,X_n]$ and is hence normally
 omitted in the definitions provided in the literature.
However, when generalizing the concept of Gr\"obner bases to other structures it
 is no longer guaranteed and hence we have included it in our definition.
The second statement ensures the existence of unique normal forms.
If we additionally require a Gr\"obner basis to be interreduced, such
a basis is unique in case we assume that the polynomials  are monic, i.e., their
 head coefficients are $1$.
The following lemma gives some properties of the reduction relation, which are
 essential in giving a constructive description of a Gr\"obner basis not only
 in the setting of commutative polynomial rings over fields.
\begin{lemma}\label{lem.buchberger.confluence}~\\
{\sl
Let $F$ be a set of polynomials  and $p,q,h$ some polynomials in $\myk[X_1, \ldots, X_n]$.
Then the following statements hold:
\begin{enumerate}
\item
Let $p-q \red{}{\myr}{b}{F} h$.
%\\
Then there are polynomials  $p',q' \in \myk[X_1, \ldots, X_n]$ such that 
 $p  \red{*}{\myr}{b}{F} p'$, $q  \red{*}{\myr}{b}{F} q'$ and $h=p'-q'$.
\item
Let $0$ be a  normal form of $p-q$ with respect to $F$.
%\\
Then there exists a polynomial  $g \in \myk[X_1, \ldots, X_n]$ such that
 $p  \red{*}{\myr}{b}{F} g$ and $q  \red{*}{\myr}{b}{F} g$.
\item
 $p \red{*}{\lr}{b}{F} q \mbox{ if and only if } p - q \in \ideal{}{}(F)$.
\item $p \red{*}{\myr}{b}{F} 0$ implies $\alpha \skm p \mrm u \red{*}{\myr}{b}{F}
  0$ for all $\alpha \in \myk$ and $u \in {\cal T}$.
\item $\alpha \skm p \mrm u \red{}{\myr}{b}{p} 0$ for all $\alpha \in \myk^*$ and $u \in {\cal T}$.
\lemend
\end{enumerate}
}
\end{lemma}
The second statement of this lemma is often called the \betonen{Translation Lemma} in the
 literature.
Statement 3 shows that Buchberger's reduction relation always captures the ideal congruence.
Statement 4 is connected to the important fact that
 reduction steps are preserved under multiplication with 
 monomials.

The set $F = \{ X_1 + X_2, X_1 + X_3 \}$ of polynomials in $\q [X_1,X_2,X_3]$ from
 page \pageref{exa.qxyz} is an example of an ideal basis which is not complete,
 i.e.~the reduction relation is not complete\footnote{Note
 that we call a set of polynomials complete (confluent,
 etc.) if the reduction relation induced by these polynomials used as rules is
 complete (confluent, etc.).}.
This follows as the polynomial $X_1$ can be reduced by $\red{}{\myr}{b}{F}$
 to $-X_2$ as well as to $-X_3$ and the latter two polynomials cannot be joined
 using $\red{}{\myr}{b}{F}$.

Of course we cannot expect an arbitrary ideal basis to be complete.
But Buchberger was able to show that in order to ``complete'' a given
 basis one only has to add finitely many special polynomials which arise from
 critical situations as described in the context of reduction systems in the previous
 section and Definition \ref{def.critical.situations.polynomialring}.

The term $X_1$ in our example describes such a critical situation which is in fact
 the only one relevant for completing the set $F$.
\begin{definition}\label{def.buchberger.spol}~\\
{\rm
The \index{s-polynomial!Buchberger}\betonen{s-polynomial}\/ for two
 non-zero polynomials
 $p,q \in \myk[X_1, \ldots,X_n]$ is defined as
$$ \spol{}(p,q) = \hc(p)^{-1} \skm p \mrm u - \hc(q)^{-1} \skm q \mrm v,$$
 where $\lcm(\hterm(p),\hterm(q))=\hterm(p) \mm u = \hterm(q)\mm v$
 for some $u,v \in {\cal T}$.
\dend
}
\end{definition}
An s-polynomial will be called \betonen{non-trivial} in case it is not zero and
notice that for non-trivial s-polynomials we always have $\hterm(\spol{}(p,q))
\pred \lcm(\hterm(p),\hterm(q))$.
The s-polynomial for $p$ and $q$ belongs to the set of critical situations
 ${\cal C}(\lcm(\hterm(p),\hterm(q)),\{p,q\})$.

In our example we find $\spol{}(X_1 + X_2, X_1 + X_3) = X_1 + X_2 -( X_1 + X_3) =
 X_2-X_3$.
 
Why are s-polynomials related to testing for local confluence? 
To answer this question we have to look at critical situations related to the
 reduction relation as defined in Definition \ref{def.buchberger.red}.
Given two polynomials $p,q \in \myk[X_1, \ldots,X_n]$ the smallest situation
 where both of them can be applied as rules is the least common multiple of their head
 terms.
Let $\lcm(\hterm(p),\hterm(q)) = \hterm(p) \mm u = \hterm(q)\mm v = t$
 for some $u,v \in {\cal T}$.
This gives us the following situation:
%\vspace{-2mm}
\begin{diagram}[size=1.6em]
   & & \lcm(\hterm(p),\hterm(q)) = t  & &   \\
   &\ldTo\phantom{XX} q &    &p\phantom{XX}\rdTo&  \\
    t - \hc(q)^{-1} \skm q \mrm v & &    & & t - \hc(p)^{-1} \skm p \mrm u   \\
   =p'  & &    & & =q'
\end{diagram}

Then we get $p' - q' = t - \hc(q)^{-1} \skm q \mrm v - (t - \hc(p)^{-1} \skm p \mrm u) = 
 \hc(p)^{-1} \skm p \mrm u - \hc(q)^{-1} \skm q \mrm v = \spol{}(p,q)$,
 i.e., the s-polynomial is derived from the two one-step successors by subtraction.
Now by Lemma \ref{lem.buchberger.confluence} we know that $\spol{}(p,q) \red{*}{\myr}{b}{F} 0$
 implies the existence of a common normal form for the polynomials $p'$ and $q'$.
Since the reduction relation based on Definition \ref{def.buchberger.red} is terminating,
 the confluence test can hence be reduced to checking whether all s-polynomials  reduce to zero.
The following theorem now gives a constructive characterization of
Gr\"obner bases based on these ideas.
\begin{theorem}\label{theo.buchberger.completion}~\\
{\sl
For a set of polynomials $F$ in $\myk[X_1, \ldots,X_n]$,
 the following statements are equivalent:
\begin{enumerate}
\item $F$ is a  Gr\"obner basis.
\item For all polynomials $g \in \ideal{}{}(F)$ we have $g \red{*}{\myr}{b}{F} 0$.
\item For all polynomials $f_{k}, f_{l} \in F$  we have 
  $ \spol{}(f_{k}, f_{l}) \red{*}{\myr}{b}{F} 0$.
\theoend
\end{enumerate}
}
\end{theorem}
\Ba{}~\\
\mbox{$1 \R 2:$ }%\\
Let $F$ be a Gr\"obner basis and $g \in \ideal{}{}(F)$.
Then $g$ is congruent to $0$ modulo the ideal generated by
 $F$, i.e., $g \red{*}{\lr}{b}{F} 0$.
Thus, as $0$  is irreducible and $G$ is confluent,
 we get $g  \red{*}{\myr}{b}{F} 0$.
\\
\mbox{$2 \R 1:$ }%\\
By Lemma \ref{lem.buchberger.confluence} 3 we know 
 $\red{*}{\lr}{b}{G} = \;\;\equiv_{\ideal{}{}(G)}$.
Hence it remains to show that reduction with respect to $F$ is confluent.
Since our reduction is terminating it is sufficient to show local
confluence.
Thus, suppose there are three different polynomials $g,h_1,h_2$
 such that $g \red{}{\myr}{b}{F} h_1$ and  $g \red{}{\myr}{b}{F} h_2$.
Then we know $h_1 \equiv_{\ideal{}{}(F)} g \equiv_{\ideal{}{}(F)} h_2$ and hence 
 $h_1 - h_2 \in \ideal{}{}(F)$.
Now by lemma \ref{lem.buchberger.confluence} (the translation lemma), 
 $h_1 - h_2  \red{*}{\myr}{b}{F} 0$ implies the
 existence of  a polynomial  $h \in \myk[X_1, \ldots, X_n]$ such that
 $h_1  \red{*}{\myr}{b}{F} h$ and $h_2  \red{*}{\myr}{b}{F} h$.
Hence, $h_1$ and $h_2$ are joinable.
\\
\mbox{$2 \R 3:$ }%\\ 
By definition \ref{def.buchberger.spol} the s-polynomial for  two
 non-zero polynomials
 $f_k,f_l \in \myk[X_1, \ldots,X_n]$ is defined as
 $$\spol{}(f_k,f_l) = \hc(f_k)^{-1} \skm f_k \mrm u -
                   \hc(f_l)^{-1} \skm f_l \mrm v,$$
 where $\lcm(\hterm(p),\hterm(q))=\hterm(p) \mm u = \hterm(q)\mm v$  
      and, hence, $\spol{}(f_{k}, f_{l}) \in \ideal{}{}(F)$.
Therefore,  $\spol{}(f_{k}, f_{l})  \red{*}{\myr}{b}{F} 0$ follows immediately.
\\
\mbox{$3 \R 2:$ }%\\  
We have to show that every 
 $g \in \ideal{}{}(F) \backslash \{ 0 \}$ is $\red{}{\myr}{b}{F}$-reducible
 to zero.
%\\
Remember that for
 $h \in \ideal{}{}(F)$, $ h \red{}{\myr}{b}{F} h'$ implies $h' \in \ideal{}{}(F)$.
As  $\red{}{\myr}{b}{F}$ is Noetherian, thus
 it suffices to show that every  $g \in \ideal{}{}(F) \backslash \{ 0 \}$ 
 is $\red{}{\myr}{b}{F}$-reducible.
%\\
Let $g = \sum_{j=1}^m \alpha_{j} \skm f_{j} \mrm w_{j}$ be an arbitrary
 representation of $g$ with $\alpha_{j} \in \myk^*$, $f_j \in F$, and  $w_{j} \in {\cal T}$.
%\\
Depending on this representation of $g$ and a total well-founded
admissible  ordering 
 $\succeq$ on ${\cal T}$ we define
 $t = \max \{ \hterm(f_{j}) \mm w_{j} \mid j \in \{ 1, \ldots, m \}  \}$ and
 $K$ is the number of polynomials $f_j \mrm w_j$ containing $t$ as a term.
%\\
Then $t \succeq \hterm(g)$ and
 in case $\hterm(g) = t$ this immediately implies that $g$ is
 $\red{}{\myr}{b}{F}$-reducible. 
%\\
Thus we will prove that $g$ has a representation where every
occurring term is less or equal to $\hterm(g)$, i.e., there exists a
representation such that $t = \hterm(g)$\footnote{Such
  representations are often called standard representations in the
  literature (compare \cite{BeWe92}).}.
This will be done by induction
 on $(t,K)$, where
 $(t',K')<(t,K)$ if and only if $t' \prec t$ or $(t'=t$ and
 $K'<K)$\footnote{Note that this ordering is well-founded since $\succ$
                  is well-founded on ${\cal T}$ and $K \in\n$.}.
%\\
In case $t \succ \hterm(g)$ there are  two polynomials $f_k,f_l$ in the corresponding 
 representation\footnote{Not necessarily $f_l \neq f_k$.}
 such that  $\hterm(f_k) \mm w_k = \hterm(f_l) \mm w_l=t$.
By definition \ref{def.buchberger.spol} we have an s-polynomial
 \mbox{$\spol{}(f_k,f_l) = \hc(f_k)^{-1} \skm  f_k
 \mrm z_k -\hc(f_l)^{-1}\skm f_l \mrm z_l$} such that
 $\hterm(f_k) \mm z_k = \hterm(f_l) \mm z_l = \lcm(\hterm(f_k),\hterm(f_l))$.
Since $\hterm(f_k) \mm w_k = \hterm(f_l) \mm w_l$ there exists an element 
 $z \in {\cal T}$ such that
 $w_k = z_k \mm z$ and $w_l = z_l \mm z$.
%\\
We will now change our representation of $g$ by using the additional
information on this s-polynomial in such a way that for the new
representation of $g$ we either have a smaller maximal term or the 
occurrences of the term $t$
are decreased by at least 1.
%\\
Let us assume that $\spol{}(f_k,f_l)$ is not trivial\footnote{In case 
               $\spol{}(f_k,f_l) = 0$, just substitute $0$
               for the sum $\sum_{i=1}^n \delta_i \skm h_i \mrm
               v_i$ in the equations below.}.
%\\
Then the reduction sequence $\spol{}(f_k,f_l) \red{*}{\myr}{b}{F} 0 $
results in a representation of the form
 $\spol{}(f_k,f_l) =\sum_{i=1}^n \delta_i \skm h_i \mrm v_i$,
 where $\delta_i  \in \myk^*,h_i \in F,v_i \in {\cal T}$.
As the $h_i$ are due to the reduction of the s-polynomial,
 all terms occurring in the sum are bounded by the term $\hterm(\spol{}(f_k,f_l))$.
Moreover, since $\succeq$ is admissible on ${\cal T}$ this 
 implies that all terms of the sum 
 $\sum_{i=1}^n \delta_i \skm h_i \mrm v_i \mrm z$ are bounded by 
 $\hterm(\spol{}(f_k,f_l)) \mm z \pred t$, i.e., they are strictly bounded
 by $t$\footnote{This can also be concluded by statement four of lemma
   \ref{lem.buchberger.confluence} since $\spol{}(f_k,f_l)
   \red{*}{\myr}{b}{F} 0 $ implies $\spol{}(f_k,f_l) \mrm z
   \red{*}{\myr}{b}{F} 0 $ and $\hterm(\spol{}(f_k,f_l) \mrm z) \prec t$.}.
%\\
We can now do the following transformations: 
\begin{eqnarray}
 &  & \alpha_{k} \skm f_{k} \mrm w_{k} + \alpha_{l} \skm f_{l} \mrm w_{l}  \nonumber\\  
 &  &  \nonumber\\                                                             
 & = &  \alpha_{k} \skm f_{k} \mrm w_{k} +
         \underbrace{ \alpha'_{l} \skm \beta_k \skm f_{k} \mrm w_{k}
                     - \alpha'_{l} \skm \beta_k \skm f_{k} \mrm w_{k}}_{=\, 0} 
        + \alpha'_{l}\skm \beta_l  \skm f_{l} \mrm w_{l} \nonumber\\ 
 & = & (\alpha_{k} + \alpha'_{l} \skm \beta_k) \skm f_{k} \mrm w_{k} - \alpha'_{l} \skm 
        \underbrace{(\beta_k \skm f_{k} \mrm w_{k}
        -  \beta_l \skm f_{l} \mrm w_{l})}_{=\, \spol{}(f_k,f_l) \mrm z} \nonumber\\
 & = & (\alpha_{k} + \alpha'_{l} \skm \beta_k) \skm f_{k} \mrm w_{k} -
        \alpha'_{l} \skm (\sum_{i=1}^n \delta_{i} \skm h_{i} \mrm (v_{i} \mm z)
        ) \label{s.buchberger}
\end{eqnarray}
where, $\beta_k=\hc(f_k)^{-1}$, $\beta_l=\hc(f_l)^{-1}$, and  $\alpha'_l \skm \beta_l = \alpha_l$.
By substituting (\ref{s.buchberger}) in our representation of $g$
 either $t$ disappears or $K$ is decreased.
\\
\qed
The second item of this theorem immediately implies the correctness of the algebraic
 definition of Gr\"obner bases, which is equivalent to Definition \ref{def.buchberger.gb}.
\begin{definition}\label{def.buchberger.gb.algebraisch}~\\
{\rm
A set $G$ of polynomials in $\myk[X_1,\ldots,X_n] \backslash \{ 0 \}$ is said to be a 
 \betonen{Gr\"obner basis}, if
 $\hterm(\ideal{}{}(G)) = \{ \hterm(g) \mrm t \mid g \in G, t \in \myt \}$.
\dend
}
\end{definition}
\begin{remark}~\\
{\rm
A closer inspection of the proof of $3 \R 2$ given above reveals a
 concept which is essential in the proofs of similar theorems for
 specific function rings in the following chapters.
The heart of this proof consists in transforming an arbitrary
 representation of an element $g$ belonging to the ideal generated by
 the set $F$ in such a way that we can deduce a top reduction
 sequence for $g$ to zero, i.e., a reduction sequence where the
 reductions only take place at the respective head term.
Such a representation of $g$ then is a standard representation
 and Gr\"obner bases are standard bases.
\remend
}
\end{remark}
As a consequence of Theorem \ref{theo.buchberger.completion} it is
decidable whether a finite set of polynomials is a Gr\"obner basis.
Moreover, this theorem gives rise to the following completion procedure
 for sets of polynomials.

\procedure{Buchberger's Algorithm}%
{\vspace{-4mm}\begin{tabbing}
XXXXX\=XXXX \kill
\removelastskip
{\bf Given:} \> A finite set of polynomials $F  \subseteq \myk[X_1, \ldots, X_n]$. \\
{\bf Find:} \> $\gb(F)$, a  Gr\"obner basis of $F$.
\end{tabbing}
\vspace{-7mm}
\begin{tabbing}
XX\=XX\= XXX \= XXX \=\kill
$G$ := $F$; \\
$B$ := $\{ (q_{1}, q_{2}) \mid q_{1}, q_{2} \in G, q_{1} \neq q_{2}
\}$; \\
{\bf while} $B \neq \emptyset$ {\bf do} \\
\>      $(q_{1}, q_{2}) := {\rm remove}(B)$; \\
\>     {\rm\kommentar \% Remove an element from
                     the set $B$}\\
\>      $h := {\rm normalform}(\spol{}(q_{1}, q_{2}),\red{}{\myr}{b}{G})$ \\
\>     {\rm\kommentar \% Compute a normal form of
                     $\spol{}(q_{1}, q_{2})$ with respect to
                     $\red{}{\myr}{b}{G}$} \\
\>      {\bf if} \>$h \neq 0$\\
\> \>     {\bf then}  \>      $B := B \cup \{ (f,h) \mid f \in G \}$; \\
\> \>                 \>      $G := G \cup \{ h \}$; \\
\>      {\bf endif}\\
{\bf endwhile} \\
$\gb (F):= G$
\end{tabbing}}

Applying this procedure to our example $F = \{ X_1 + X_2, X_1 + X_3\}$
 from page \pageref{exa.qxyz} gives us $h = X_2-X_3$ and $G= F \cup \{ h \}$
 is a Gr\"obner basis as all other critical situations are resolvable.

Termination of the procedure can be shown by using a slightly different
 characterization of Gr\"obner bases (see Section \ref{section.intro.definitions}): 
A subset $G$ of $\ideal{}{\myk[X_1, \ldots, X_n]}(F)$ is a Gr\"obner basis of
 $\ideal{}{\myk[X_1, \ldots, X_n]}(F)$ if and only if
 $\hterm(\ideal{}{\myk[X_1, \ldots, X_n]}(F) \backslash \{ 0 \}) = \ideal{}{\cal T}(\hterm(G))$, i.e., the set of the head
 terms of the polynomials in the ideal generated by $F$ in $\myk[X_1, \ldots, X_n]$
 coincides with
 the ideal (in ${\cal T}$) generated by the head terms of the
 polynomials in $G$.
Reviewing the procedure, we find that every polynomial added 
 in the {\bf while} loop has
 the property that its head term cannot be divided by the head terms of
 the polynomials already in $G$.
By Dickson's Lemma or Hilbert's Basis Theorem,
 the head terms of the polynomials in $G$ will at some step form a
 basis for the set of head terms of the polynomials of the ideal generated by
 $F$ which itself is  the ideal in ${\cal T}$
 generated by the head terms of the polynomials in $G$.
From this time on for every new polynomial $h$ computed by the
 algorithm the head term $\hterm(h)$ must lie in this ideal.
Therefore, its head term must be divisible by at least
 one of the head terms of the polynomials in $G$, i.e., $\hterm(h)$ and hence $h$ cannot be
 in normal form with respect to $G$ unless it is zero.

%%% Local Variables: 
%%% mode: latex
%%% TeX-master: "testlauf"
%%% End: 

%% file: reductionrings.tex
In this chapter we proceed to distinguish sufficient conditions, which allow to define 
 a reduction relation for a ring in such a way that every finitely generated ideal in the ring
 has a finite Gr\"obner basis with respect to that reduction relation.
Such rings will be called reduction rings.
Often additional conditions can be given
 to ensure effectivity for the ring operations, the reduction relation and the computation of the Gr\"obner bases
 -- the ring is then called an effective reduction ring.
Naturally the question arises, when and how the property of being a
 reduction ring is preserved under various ring constructions.
This can be studied from an existential as well as from a constructive point of view.
One main goal of studying abstract reduction rings is to provide universal methods for 
 constructing new reduction rings
 without having to generalize the whole setting individually for each new
 structure:
e.g.~knowing that the integers $\z$ are a  reduction ring and that the property
 lifts to polynomials in one variable, we find that $\z[X]$ is again a 
 reduction ring and we can immediately conclude that also $\z[X_1, \ldots, X_n]$
 is a  reduction ring.
Similarly, as sums of reduction rings are again reduction rings, we can directly
 conclude that $\z^k[X_1, \ldots, X_n]$ or even $(\z[Y_1, \ldots, Y_m])^k[X_1, \ldots, X_n]$
 are reduction rings.
Moreover, since $\z$ is an effective reduction ring it can be shown that these new reduction rings
 again are effective.
Commutative effective reduction rings have been studied by Buchberger, Madlener, and Stifter in \cite{Bu83,Ma86,St87}.

On the other hand, many rings of interest are non-commutative, e.g.~rings of
 matrices, the ring of quaternions, Bezout rings and various monoid rings, and since in many cases 
 they can be regarded as reduction rings, they are again candidates
 for applying ring constructions. More interesting examples of
 non-commutative reduction rings have been studied by Pesch in \cite{Pe97}. 

A general framework for reduction rings and ring constructions
 including the non-commutative case was presented at the Linz conference ``33 years of
 Gr\"obner Bases'' in \cite{MaRe98}.
Here we extend this framework by giving more details and insight.
Additionally, we add a section on modules over reduction rings, as this
 concept arises naturally as a generalization of ideals in rings.

Of course there are also rings of interest, which can be enriched by a reduction relation,
 but will not allow finite Gr\"obner bases for all ideals.
Monoid and group rings provide such a setting.
For such structures still many of the properties studied here are of interest and can
 be shown in weaker forms, e.g.~provided a monoid ring with a reduction relation we
 can define a reduction relation for the polynomial ring with one variable over the monoid ring. 

The chapter is organized as follows:
In Section \ref{section.reductionrings} we introduce axioms for specifying reduction
 relations in rings and give two concepts involving special forms of ideal bases
 -- weak reduction rings and reduction rings.
In Section \ref{section.quotients} -- \ref{section.polyrings} we study quotients, sums, modules, and 
 polynomial rings  of these structures.

%%%%%%%%%%%%%%%%%%%%%%%%%%%%%%%%%%%%%%%%%%%%%%%%%%%%%%%%%%%%%%%%%%%%%%%
\section{Reduction Rings}\label{section.reductionrings}
%%%%%%%%%%%%%%%%%%%%%%%%%%%%%%%%%%%%%%%%%%%%%%%%%%%%%%%%%%%%%%%%%%%%%%%
Let $\rr$ be a  ring with unit $1$ and  a (not necessarily effective) reduction relation $\Longrightarrow_B \subseteq \rr \times \rr$
associated with  subsets $B \subseteq\rr$ satisfying the following
axioms:
\begin{enumerate}
\item[(A1)] $ \Longrightarrow_B\; =\; \bigcup_{\beta \in B} \Longrightarrow_{\beta}$,\\
  $\Longrightarrow_B$ is terminating for all {\em finite}
  subsets $B \subseteq\rr$.
\item[(A2)] $\alpha \Longrightarrow_{\beta} \gamma$ implies $\alpha - \gamma
  \in \ideal{}{\rr}(\beta)$.
\item[(A3)] $\alpha \Longrightarrow_{\alpha} 0$ for all $\alpha \in \rr
  \backslash \{ 0 \}$.
\end{enumerate}
Part one of Axiom (A1) states how a reduction relation using sets is defined 
 in terms of a reduction relation using elements of $\rr$ and
 is hence applicable to arbitrary sets $B \subseteq \rr$.
However, Axiom (A1) does {\em not} imply termination of reduction with
 respect to arbitrary sets:
Just assume for example the ring $\rr = \q[\{ X_i \mid i \in \n \}]$,
 i.e., the polynomial ring with infinitely many indeterminates, and the 
 reduction relation based on divisibility of head terms
 with respect to the length-lexicographical
 ordering induced by $X_1 \succ X_2 \succ \ldots$.
Then although reduction when using a finite set of polynomials is terminating,
 this is no longer true for infinite sets.
For example the infinite set
 $\{ X_i - X_{i+1} \mid i \in \n \}$ gives rise to an infinite
 reduction sequence 
 $X_1 \red{}{\R}{}{X_1 - X_2} X_2 \red{}{\R}{}{X_2 - X_3} X_3 \ldots$.
This phenomenon of course has many consequences.
Readers familiar with Gr\"obner bases in polynomial rings know that when proving
 that a set of polynomials is a Gr\"obner basis if and only if all ideal 
 elements reduce to zero using the set, this is shown by proving that every
 ideal element is reducible by some element in the set (compare Theorem \ref{theo.buchberger.completion}).
Unfortunately, this only implies reducibility to zero in case
 the reduction relation is
 terminating.
Without this property other methods have to be applied.

In order to ensure termination for arbitrary subsets of $\rr$
 it is possible to give a more restricted form of Axiom (A1):
\begin{enumerate}
\item[(A1')] $ \Longrightarrow_B\; =\; \bigcup_{\beta \in B} \Longrightarrow_{\beta}$,\\
  $\Longrightarrow_B$ is terminating for {\em all} subsets $B \subseteq\rr$.
\end{enumerate}
Then of course reduction sequences are always terminating and many additional restrictions, which we
 have to add later, are no longer necessary. 
Still we prefer the more general formulation of the axiom since it allows to
 state more clearly why and where termination is needed and how it can be achieved.

Axiom (A2) states how reduction steps are related to the ideal congruence, namely that
 one reduction step using an element $\beta \in \rr$ is captured by the congruence generated by
 $\ideal{}{\rr}(\beta)$.
We will later on see that this extends to 
 the reflexive transitive symmetric closure $\red{*}{\Longleftrightarrow}{}{B}$ 
 of any reduction relation $\Longrightarrow_B$ for arbitrary sets $B \subseteq \rr$.

Notice that in case $\rr$ is commutative (A2) implies
  $\gamma = \alpha - \beta \skm \rho$ for some $\rho$ in $\rr$.
In the non-commutative case using a single element $\beta$ for reduction 
 $\alpha - \gamma  \in \ideal{}{\rr}(\beta)$ only implies
 \mbox{$\gamma = \alpha - \sum_{i=1}^k \rho_{i1} \skm \beta \skm \rho_{i2}$}
 for some $\rho_{i1}, \rho_{i2} \in \rr$, $1 \leq i \leq k$, 
 hence possibly involving $\beta$ more than once with different multipliers.
This provides a large range of possibilities for defining reduction steps, e.g.~by
 subtracting one or more appropriate multiples of $\beta$ from $\alpha$.
Notice further that on the converse Axiom (A2) does {\em not} provide any information 
 on how $\alpha$, $\gamma \in \rr$ with $\alpha - \gamma \in \ideal{}{\rr}(\beta)$ are
 related with respect to the reduction relation $\Longrightarrow_{\{ \beta \}}$.
As a consequence many properties of specialized reduction relations as 
 known from the literature, e.g.~the useful Translation Lemma, cannot be shown to
 hold in this general setting.

We can define \betonen{one-sided} (right or left) reduction relations in rings by refining Axiom
 (A2) as follows: 
\begin{enumerate}
\item[(A2r)]$\alpha \Longrightarrow_{\beta} \gamma$ implies $\alpha - \gamma
  \in \ideal{r}{\rr}(\beta)$, respectively
\item[(A2l)]$\alpha \Longrightarrow_{\beta} \gamma$ implies $\alpha - \gamma
  \in \ideal{l}{\rr}(\beta)$.
\end{enumerate}
In these special cases again we always get $\gamma = \alpha - \beta \skm \rho$ respectively 
 $\gamma = \alpha - \rho \skm \beta$ for some $\rho \in \rr$.

Remember that Axiom (A2) while not specific on the exact form of the reduction step ensures
 that reduction steps ``stay'' within the ideal congruence.
Let us now study the situation for a set $B \subseteq \rr$ and
 let $\equiv_{{\mathfrak i}}$ denote the congruence
 generated by the ideal ${\mathfrak i}=\ideal{}{}(B)$,
 i.e., $\alpha \equiv_{{\mathfrak i}} \beta$ if and only if
 $\alpha - \beta \in {\mathfrak i}$.
Then (A1)\footnote{We only need the first part of Axiom (A1), namely how $\Longrightarrow_B$ is
 defined, and hence we do not have to restrict ourselves to finite sets.} and (A2) immediately imply
 $\red{*}{\Longleftrightarrow}{}{B} \subseteq \;\;\equiv_{{\mathfrak i}}$.
Hence, in case the reduction relation is effective one method for deciding
 the membership problem for a finitely generated ideal ${{\mathfrak i}}$
 is to transform a finite generating set $B$ into a finite set $B'$ such that 
 $B'$ still generates ${\mathfrak i}$ and
 $\Longrightarrow_{B'}$ is confluent on ${\mathfrak i}$.
Notice that $0$ has to be irreducible\footnote{$0$ cannot be reducible by itself since this would contradict
 the termination property in (A1).
Similarly, $0 \R_{\beta} 0$ and $0 \R_{\beta} \gamma$, both $\beta$ and
 $\gamma$ not equal $0$, give rise to infinite reduction sequences
 again contradicting (A1).} for all $\R_{\alpha}$, 
 $\alpha \in \rr$.
Therefore, $0$ has to be {\em the} normal form of the ideal elements.
Hence the goal is to achieve $\alpha \in {{\mathfrak i}}$ if and only if
 $\alpha \red{*}{\Longrightarrow}{}{B'} 0$.
In particular ${\mathfrak i}$ is one
 equivalence class of $\red{*}{\Longleftrightarrow}{}{B'}$.
The different definitions of reduction  relations for rings existing in literature show
 that for deciding the membership problem of an ideal ${\mathfrak i}$ it is not necessary to enforce
 $\red{*}{\Longleftrightarrow}{}{B'} = \;\equiv_{{\mathfrak i}}$.
For example the D-reduction notion given by Pan in \cite{Pa85} does not
 have this property but is still sufficient 
 to decide $\equiv_{{\mathfrak i}}$-equivalence of two elements because 
 $\alpha \equiv_{{\mathfrak i}} \beta$ if and only if
 $\alpha - \beta \in {{\mathfrak i}}$.
It may even happen that D-reduction is not only confluent on
 ${{\mathfrak i}}$ but
 confluent everywhere and still $\alpha \equiv_{{\mathfrak i}} \beta$ does not
 imply that the normal forms with respect to D-reduction are the same.
This phenomenon is illustrated in the next example.
\begin{example}\label{exa.Z}~\\
{\rm
Let us look at different ways of introducing reduction relations for the ring of integers $\z$.
For $\alpha, \beta, \gamma \in \z$ we define:
\begin{itemize}
\item $\alpha \red{}{\R}{}{\beta} \gamma$ if and only if 
       $\alpha = \kappa \skm | \beta | + \gamma$ where
       $0 \leq \gamma < |\beta |$ and $\kappa \in \z$ (division with remainder), 
\item $\alpha \red{}{\R}{D}{\beta} 0$ if and only if
       $\alpha = \kappa \skm \beta$, i.e.~$\beta$ is a proper divisor of
       $\alpha$ (D-reduction).
\end{itemize}
Then for example we have $5 \red{}{\R}{}{4} 1$ but $5 \nred{}{\R}{D}{4}$. \\
It is easy to show that both reduction relations satisfy (A1) -- (A3).
Moreover, all elements in $\z$ have unique normal forms.
An element belongs to $\ideal{}{}(4)$
 if and only if it is reducible to zero using $4$.
For $\R$-reduction the normal forms are unique representatives of the quotient
 $\z/\ideal{}{}(4)$.
This is no longer true for $\R^D$-reduction, since e.g.~$3 \equiv_{\ideal{}{}(4)} 7$ since
 $7 = 3 + 4$,
 but both are $\R^D$-irreducible.
On the other hand, as $\R_{\alpha}^D$ is only applicable to multiples
 $\kappa \skm \alpha$ and then reduces them to zero, $\R_{4}^D$ is confluent
 everywhere on $\z$.
\exaend
}
\end{example}
Since confluence of a reduction relation on the ideal is already sufficient to solve its
 membership problem, bases with this property called weak Gr\"obner bases have been studied in the literature.
We proceed here by defining such weak Gr\"obner bases in our context.
\begin{definition}\label{def.weak.gb.reduction.ring}~\\
{\rm
A subset $B$ of  $\rr$ is called a 
 \betonen{weak Gr\"obner basis} of the ideal ${\mathfrak i}= \ideal{}{}(B)$ it
 generates, if 
 $\red{}{\R}{}{B}$ is terminating and $\alpha \red{*}{\R}{}{B} 0$
 for all $\alpha \in {\mathfrak i}$.
\dend
}
\end{definition}
Notice that in Theorem \ref{theo.buchberger.completion} this property was one way of
 characterizing Gr\"obner bases in
 $\myk[X_1, \ldots, X_n]$.
We will later on see why in polynomial rings the terms weak Gr\"obner basis and Gr\"obner basis 
 coincide.
\begin{definition}\label{def.weakredring}~\\
{\rm
A ring $(\rr, \R)$ satisfying (A1) -- (A3) 
 is called a \betonen{weak reduction ring}
 if every finitely generated ideal in $\rr$ has a finite weak Gr\"obner basis.
\dend
}
\end{definition}
As stated before such a weak Gr\"obner basis is sufficient to decide the ideal membership
 problem in case the reduction relation is effective.
However, if we want unique normal forms for {\em all} elements in $\rr$ such that each congruence
 has one unique representative we need a stronger kind of ideal basis.
\begin{definition}\label{def.gb}~\\
{\rm
A subset $B$ of  $\rr$ is called a 
 \betonen{Gr\"obner basis} of the ideal ${\mathfrak i}= \ideal{}{}(B)$
 it generates, if
 $\red{*}{\Longleftrightarrow}{}{B} = \;\;\equiv_{{\mathfrak i}}$
 and $\red{}{\R}{}{B}$ is complete\footnote{Notice that in the literature definitions of Gr\"obner
 bases normally only require that $\red{}{\R}{}{B}$ is ``confluent''. This is due to the fact that in these
 cases $\red{}{\R}{}{B}$ is terminating. In our context, however for arbitrary
 sets $B \subseteq \rr$ we have seen that $\red{}{\R}{}{B}$ need not be Noetherian. Hence we have
 to incorporate this additional requirement into our definition, which is done by demanding completeness.
 Hence here we have a point where the weaker form (A1) demands more care in defining the
 term ``Gr\"obner basis''.
 In rings where the reduction relation using an arbitrary set of elements is always Noetherian, the weaker demand
 for (local) confluence is of course sufficient.}.
\dend
}
\end{definition}
Of course Gr\"obner bases are also weak Gr\"obner bases.
This can be shown by induction on $k$, where for 
 $\alpha \in \ideal{}{}(B)$
 we have $\alpha \red{k}{\Longleftrightarrow}{}{B} 0$.
In case $k = 1$ we immediately get that $\alpha \red{}{\R}{}{B} 0$ must hold
 as  $0$ is irreducible.
In case $k>1$ we find $\alpha \red{}{\Longleftrightarrow}{}{B} \beta \red{k-1}{\Longleftrightarrow}{}{B} 0$
 and by our induction hypothesis $\beta \red{*}{\R}{}{B} 0$ must hold.
Now either $\alpha \red{}{\R}{}{B} \beta$ and we are done or
 $\beta \red{}{\R}{}{B} \alpha$.
In the latter case the completeness of our reduction relation combined with the irreducibility  of zero then must
 yield $\alpha \red{*}{\R}{}{B} 0$ and we are done.

The converse is not true.
To see this let us review the definition of $\R^D$-reduction for $\z$ as presented in Example \ref{exa.Z}.
Then the set $\{ 2 \}$ is a weak Gr\"obner basis of the ideal $2\skm \z =\{ 2 \skm 
 \alpha \mid \alpha \in \z \} $ as for every
 $\alpha \in ( 2\skm \z ) \backslash \{ 0 \}$ we  have $\alpha \R_{\{2\}}^{D} 0$.
On the other hand elements in $\z \backslash (2\skm\z)$ are irreducible and hence $3$ and $5$ are in normal form
 with respect to $\R_{\{2\}}^{D}$. 
Therefore, $3 \nred{*}{\Longleftrightarrow}{D}{\{2\}} 5$ although $5 \equiv_{2\skm\z} 3$ as $5 = 3 + 1 \skm 2$.

However, for many rings as e.g.~polynomial rings over fields, weak Gr\"obner
 bases are also Gr\"obner bases.
This is due to the fact that many
 rings with reduction relations studied in the literature fulfill
 a certain property for the reduction relation called the
 Translation Lemma (compare Lemma \ref{lem.buchberger.confluence} (2)).
Rephrased in our context
 the Translation Lemma states that for a set $F \subseteq \rr$
 and for all $\alpha, \beta \in \rr$, $\alpha - \beta \red{*}{\R}{}{F} 0$
 implies the existence of $\gamma \in \rr$ such that $\alpha \red{*}{\R}{}{F} \gamma$ and 
 $\beta \red{*}{\R}{}{F} \gamma$.
As mentioned before, the validity of this lemma for a reduction relation in a
 ring has consequences on the relation between weak Gr\"obner
 bases and Gr\"obner bases.

\begin{theorem}\label{theo.weak+translation}~\\
{\sl
Let $\rr$ be a ring with a reduction relation $\R$ fulfilling (A1) -- (A3).
If additionally the Translation Lemma holds for the reduction relation $\R$ in $\rr$, 
 then weak Gr\"obner bases are also Gr\"obner bases.
\theoend
}
\end{theorem}
\Ba{}~\\
Let $\rr$ be a  ring where the Translation Lemma holds for the reduction relation $\R$.
Further let $B$ be a weak Gr\"obner basis of the ideal ${\mathfrak i}= \ideal{}{}(B)$.
In order to prove that $B$ is in fact a Gr\"obner basis we 
 have to show two properties:
\begin{enumerate}
\item $\red{*}{\Longleftrightarrow}{}{B} = \;\;\equiv_{{\mathfrak i}}$:\\
      The inclusion $\red{*}{\Longleftrightarrow}{}{B} \subseteq \;\;\equiv_{{\mathfrak i}}$
       follows by (A1) and (A2).
      To see the converse let $\alpha \equiv_{{\mathfrak i}} \beta$.
      Then $\alpha - \beta \in {\mathfrak i}$, and $\alpha - \beta \red{*}{\R}{}{B} 0$, as
       $B$ is a weak Gr\"obner basis.
      But then the Translation Lemma yields that $\alpha$ and $\beta$ are joinable by $\red{}{\R}{}{B}$
      and hence $\alpha \red{*}{\Longleftrightarrow}{}{B} \beta$.
\item $\red{}{\R}{}{B}$ is complete:\\
      Since $\red{}{\R}{}{B}$ is terminating it suffices to show local confluence.
      Let $\alpha, \beta_1, \beta_2 \in \rr$ such that $\alpha \red{}{\R}{}{B} \beta_1$ and
       $\alpha \red{}{\R}{}{B} \beta_2$.
      Then again $\beta_1 - \beta_2 \in {\mathfrak i}$, and $\beta_1 - \beta_2 \red{*}{\R}{}{B} 0$, since
       $B$ is a weak Gr\"obner basis.
      As before the Translation Lemma yields that $\beta_1$ and $\beta_2$ are joinable 
       by $\red{}{\R}{}{B}$ and we are done.
\end{enumerate}
\qed
On the other hand, looking at proofs of variations of the Translation Lemma in the literature
 we find that in order to show this property 
 for a ring with a reduction relation we need more information on the reduction step as is provided 
 by the very general form of Axiom (A2).
Hence in this general setting weak Gr\"obner bases and Gr\"obner bases have to be distinguished.

Rings where finitely generated ideals have finite Gr\"obner bases are of particular interest.
\begin{definition}\label{def.redring}~\\
{\rm
A ring $(\rr, \R)$ satisfying (A1) -- (A3) 
 is called a \betonen{reduction ring}
 if every finitely generated ideal in $\rr$ has a finite Gr\"obner basis.
\dend
}
\end{definition}
The connection between weak reduction rings and reduction rings follows from Theorem
 \ref{theo.weak+translation}.
\begin{corollary}~\\
{\sl
Let $(\rr, \R)$ be a weak reduction ring.
If additionally the Translation Lemma holds, 
 then $(\rr, \R)$ is a reduction ring.
\corend
}
\end{corollary}
To simplify notations sometimes  we will identify $(\rr, \R)$ with $\rr$ in case $\R$
 is known or irrelevant.
The notion of \betonen{one-sided weak reduction rings} and \betonen{one-sided reduction rings}
  is straightforward\footnote{An example for a one-sided weak reduction ring 
  which is not a one-sided reduction ring can be given using the two 
  different reduction relations $\red{}{\R}{}{}$ and $\red{}{\R}{D}{}$
  for the integers provided in Example \ref{exa.Z}.
  Then the free monoid ring $\z[\{ a,b \}]$ with prefix reduction induced
   by $\red{}{\R}{}{}$ is a one-sided reduction ring while for prefix
   reduction induced by $\red{}{\R}{D}{}$ we get a one-sided weak reduction ring.}.

\betonen{Effective} or \betonen{computable weak reduction rings} 
 and \betonen{effective} or \betonen{computable reduction rings}
 can be defined similar to 
 Buchberger's commutative reduction rings (see \cite{Bu83,St87}), 
 in our case by demanding that
 the ring operations are computable, the reduction relation is effective, and, 
 additionally, Gr\"obner bases can be computed.
Procedures which compute Gr\"obner bases are normally completion procedures based
 on effective tests for local confluence
 to decide whether a finite set is a Gr\"obner basis and to enrich that set if not.
But of course other procedures are also possible, e.g.~when using division with
 remainders as reduction relation in $\z$ the Euclidean algorithm can be used for
 computing Gr\"obner bases of ideals.

Notice that Definition \ref{def.redring} does not imply that Noetherian rings
 satisfying the Axioms (A1), (A2) and (A3) are indeed reduction rings.
This is due to the fact that while of course all ideals then have
 finite bases, the property of being a Gr\"obner basis
 strongly depends on the reduction ring which is of course itself
 strongly dependent on the reduction relation chosen for the ring.
Hence the existence of finite ideal bases does not imply the existence
 of finite Gr\"obner bases as the following example shows:
Given an arbitrary Noetherian ring $\rr$ we can associate a (very simple) reduction relation
 to elements of $\rr$ by defining for any $\alpha \in \rr \backslash \{ 0 \}$,
 $\alpha \Longrightarrow_{\beta}$ if and only if $\alpha = \beta$.
Additionally we define $\alpha \Longrightarrow_{\alpha} 0$.
Then the Axioms (A1), (A2)  and (A3) are fulfilled but although every ideal in 
 the Noetherian ring $\rr$
 has a finite basis (in the sense of a generating set), 
 infinite ideals will not have finite Gr\"obner bases, as for any ideal
 ${\mathfrak i} \subseteq \rr$
 in this setting the set ${\mathfrak i} \backslash \{ 0 \}$ is the only
 possible Gr\"obner basis.

Another interesting question concerns which changes to ideal bases preserve the
 property of being a Gr\"obner basis.
Extensions of (weak) Gr\"obner bases by ideal elements are not critical\footnote{Extensions
 of (weak) Gr\"obner bases by elements not belonging to the ideal make no sense in our context as
 then the reduction relation no longer is a proper means for describing the original ideal congruence.}.
\begin{remark}\label{rem.add.elements.to.basis}~\\
{\rm
If $B$ is a finite (weak) Gr\"obner basis of ${\mathfrak i}$ and $\alpha \in {\mathfrak i}$,
 then $B' = B \cup \{ \alpha \}$ is again a (weak) Gr\"obner basis of ${\mathfrak i}$:
First of all we find 
$ \red{*}{\Longleftrightarrow}{}{B} \subseteq \red{*}{\Longleftrightarrow}{}{B'}
 \subseteq\;\;\equiv_{{\mathfrak i}}\;\; = \red{*}{\Longleftrightarrow}{}{B}$.
Moreover, since $B'$ is again a finite set, $\R_{B'}$ is terminating.
Finally $\R_{B'}$ inherits its confluence from $\R_{B}$ since
 $\beta \R_{\alpha} \gamma$ implies $\beta \equiv_{{\mathfrak i}} \gamma$, and hence
 $\beta$ and $\gamma$ have the same normal form with
 respect to $\R_{B}$.
\remend
}
\end{remark}
Hence, if $B$ is a finite Gr\"obner basis of an ideal ${\mathfrak i}$ and $\beta \in B$
 is reducible by $B \backslash \{ \beta \}$ to $\alpha$, then
 $B \cup \{ \alpha \}$ is again a Gr\"obner basis of ${\mathfrak i}$.
The same is true for weak Gr\"obner bases.

Removing elements from a set is critical as we might decrease
 the set of elements which are reducible with respect to the set.
Hence if the set is a Gr\"obner basis, after removing elements
 the ideal elements might no longer reduce to zero using the
 remaining set.
Reviewing the example presented in Section \ref{section.intro.applications} we find that
 while the set $\{ X_,^2 + X_2, X_1^2 + X_3, X_2 - X_3 \}$ is a Gr\"obner basis in $\q[X_1, X_2, X_3]$
 the subset $\{ X_,^2 + X_2, X_1^2 + X_3 \}$, although it generates the same ideal, is none.
In order to remove $\beta$ from a Gr\"obner basis 
 $B$ without losing
 the Gr\"obner basis property it is important for the reduction relation $\R$  to
 satisfy an  additional axiom:
\begin{enumerate}
\item[(A4)] $\alpha \Longrightarrow_{\beta}$ and
             $\beta \Longrightarrow_{\gamma} \delta$ 
             imply $\alpha \Longrightarrow_{\gamma}$ or 
             $\alpha \Longrightarrow_{\delta}$.
\end{enumerate}
It is not easy to give a simple example for a ring with a reduction relation fulfilling
 (A1) -- (A3) but not (A4) as the reduction rings we have introduced so far all
 satisfy (A4)\footnote{An example using a right reduction relation in a 
 monoid ring can be found in Example 3.6
 in \cite{MaRe95}:
 Let  $\Sigma = \{ a, b, c \}$  and
 $T = \{ a^2 \myr 1, b^2 \myr 1, c^2 \myr 1 \}$  be a
 monoid presentation of $\m$ with a length-lexicographical ordering induced
 by $a \succ b \succ c$.
%\\
For $p, f \in \myk[\m]$ a (right) reduction relation is defined by 
 $p \red{}{\myr}{s}{f} q$ at a monomial $\alpha \skm t$, if
\begin{enumerate}
\item[(a)] $\hterm(f \mrm w) = t$ for some $w \in \m$, and
\item[(b)] $q = p - \alpha \skm \hc(f \mrm w)^{-1} \skm f \mrm w$.
\end{enumerate}
Looking at $p = ba + b, q = bc + 1$ and $r = ac + b \in \q[\g]$ we get
 $p \red{}{\myr}{s}{q} p - q \mrm ca = -ca + b$ and $q \red{}{\myr}{s}{r} q - r \mrm c
 = -a + 1 = q_{1}$, but $p \nred{}{\myr}{s}{\{r,q_1\}}$.
%\\
Trying to reduce  $ba$ by $r$ or $q_{1}$ we get
 $r \mrm a = \underline{aca} + ba, r \mrm caba = ba + \underline{bcaba}$ and
 $q_{1} \mrm aba = -ba + \underline{aba},q_{1} \mrm ba = -\underline{aba} + ba$
 all violating condition (a).
Trying to reduce $b$ we get the same problem as $r \mrm cab = b + \underline{bcab},
 q_1 \mrm ab = -b+\underline{a}$
 and $q_1 \mrm b = -\underline{ab}+b$.}.
\begin{lemma}\label{lem.BB'}~\\
{\sl
Let $(\rr, \R)$ be a reduction ring satisfying (A4).
Further let $B \subseteq \rr$ be a  (finite) Gr\"obner basis of a finitely generated ideal in $\rr$
 and $B' \subseteq B$ such that for all $\beta \in B$, $\beta \red{*}{\R}{}{B'} 0$ holds.
Then $B'$ is a Gr\"obner basis of $\ideal{}{\rr}(B)$.
In particular, for all $\alpha \in \rr$, $\alpha \red{*}{\R}{}{B} 0$ implies
 $\alpha \red{*}{\R}{}{B'} 0$.
\lemend
}
\end{lemma}
\Ba{}~\\
In this proof let $\alpha\rnf{B}$ denote a normal form of $\alpha$ with respect to $\R_B$
 and let $\irr(\R_B)$ denote the $\R_B$-irreducible elements in $\rr$.
Notice that by the Axioms (A1) and  (A4) and our assumptions on $B'$, all elements reducible
 by $B$ are also reducible by $B'$:
We show a more general claim by induction on $n$:
If $\alpha, \beta \in \rr$ such that $\alpha \R_{\beta}$ and $\beta \red{n}{\R}{}{B'} 0$,
 then $\alpha \R_{B'}$.
The base case $n=1$ is a direct consequence of (A4), as $\alpha \R_{\beta}$ and $\beta \red{}{\R}{}{\beta' \in B'} 0$
 immediately imply $\alpha \R_{\beta' \in B'}$.
In the induction step we find $\beta \R_{\beta' \in B'} \delta \red{n-1}{\R}{}{B'} 0$ and either
 $\alpha \R_{\beta' \in B'}$ or $\alpha \R_{\delta}$ and our induction hypothesis yields $\alpha \R_{B'}$.
\\
Hence we can conclude
 $\irr(\R_{B'}) \subseteq \irr(\R_B)$.
We want to show that $B'$ is a Gr\"obner basis of $\ideal{}{\rr}(B)$:
Assuming $\alpha \red{*}{\R}{}{B} \alpha\rnf{B}$ but $\alpha \red{*}{\R}{}{B'} \alpha\rnf{B'} \neq \alpha\rnf{B}$,
  we find $\alpha\rnf{B'} \in \ideal{}{\rr}(B)$ and
  $\alpha\rnf{B'} \in \irr(\R_{B'}) \subseteq \irr(\R_B)$, contradicting the confluence of $\R_B$.
Hence, $\alpha\rnf{B'} = \alpha\rnf{B}$, 
 implying that $\R_{B'}$ is also confluent, as $\alpha\rnf{B}$ is unique.
Now it remains to show that
 $\red{*}{\Longleftrightarrow}{}{B} \subseteq 
 \red{*}{\Longleftrightarrow}{}{B'}$ holds.
This follows immediately, as for 
 $\alpha \red{*}{\Longleftrightarrow}{}{B} \beta$ we have
 $\alpha\rnf{B'} = \alpha\rnf{B} = \beta\rnf{B}  = \beta\rnf{B'}$ 
 which implies $\alpha \red{*}{\Longleftrightarrow}{}{B'} \beta$.
\\ \qed
This result carries over for weak Gr\"obner bases.
\begin{corollary}~\\
{\sl
Let $(\rr, \R)$ be a weak reduction ring satisfying (A4).
Further let $B \subseteq \rr$ be a (finite) weak Gr\"obner basis of a finitely generated ideal in $\rr$
 and $B' \subseteq B$ such that for all $\beta \in B$, $\beta \red{*}{\R}{}{B'} 0$ holds.
Then $B'$ is a weak Gr\"obner basis of $\ideal{}{\rr}(B)$.
In particular, for all $\alpha \in \rr$, $\alpha \red{*}{\R}{}{B} 0$ implies 
 $\alpha \red{*}{\R}{}{B'} 0$.
\corend
}
\end{corollary}
\Ba{}~\\
As in the proof of Lemma \ref{lem.BB'} we can conclude $\irr(\R_{B'}) \subseteq \irr(\R_B)$.
Hence assuming that $\alpha \red{*}{\R}{}{B} 0$ while $\alpha \red{*}{\R}{}{B'} \alpha\rnf{B'} \neq 0$
 would imply $\alpha\rnf{B'} \in  \irr(\R_B)$.
As $B' \subseteq B$ this would give us a contradiction since then $\alpha \in \ideal{}{\rr}(B)$
 would have two different normal forms at least one of them not equal
 to zero with respect to $B$ contradicting the fact that $B$ is supposed
 to be a weak Gr\"obner basis.
\\ \qed
Remark \ref{rem.add.elements.to.basis} and Lemma \ref{lem.BB'} are closely related to
 interreduction and reduced (weak) Gr\"obner bases.
We call a (weak) Gr\"obner basis $B \subseteq \rr$ \betonen{reduced} if no element $\beta \in B$
 is reducible by $\R_{B \backslash \{ \beta \} }$.

The results of this section carry over to rings with appropriate one-sided reduction relations.

In the remaining sections of this chapter we study the question
  which ring constructions preserve the property of being a (weak) reduction ring.
%
%%%%%%%%%%%%%%%%%%%%%%%%%%%%%%%%%%
\section{Quotients of Reduction Rings}\label{section.quotients}
Let $\rr$ be a ring with a reduction relation $\R$ fulfilling (A1) -- (A3) 
 and ${\mathfrak i}$ a finitely generated
 ideal in $\rr$ with a finite Gr\"obner basis $B$.
Then every element
 $\alpha \in \rr$ has a unique normal form $\alpha\rnf{B}$ with respect to $\R_{B}$.
We choose 
 the set of $\R_{B}$-irreducible elements of $\rr$ as representatives 
 for the elements in the \betonen{quotient} $\rr/{\mathfrak i}$.
Addition is defined by $\alpha + \beta :=  (\alpha + \beta)\rnf{B}$
 and multiplication by $\alpha \skm \beta :=  (\alpha \skm \beta)\rnf{B}$.
Then a natural reduction relation can be defined on the quotient
 $\rr/{\mathfrak i}$ as follows:
\begin{definition}\label{def.quotientred}~\\
{\rm
Let $\alpha, \beta, \gamma \in \rr/{\mathfrak i}$.
We say $\beta$ \betonen{reduces} $\alpha$ to $\gamma$ in one step,
 denoted by $\alpha \myr_{\beta} \gamma$, if there exists $\gamma' \in \rr$
 such that $\alpha \R_{\beta} \gamma'$ and $(\gamma')\rnf{B} = \gamma$.
\dend
}
\end{definition}
First we ensure that the Axioms (A1) -- (A3)  hold for the reduction relation in
 $\rr/{\mathfrak i}$ based on Definition \ref{def.quotientred}:
$\myr_{S}\; = \bigcup_{s \in S} \myr_s$ is terminating for
 all finite $S \subseteq \rr/{\mathfrak i}$ since otherwise
 $\R_{B\cup S}$ would not be terminating in $\rr$ although $B \cup S$ is finite.
Hence (A1) is satisfied.
If $\alpha \myr_{\beta} \gamma$ for some
  $\alpha, \beta, \gamma \in \rr/{\mathfrak i}$ we know
 $\alpha \R_{\beta} \gamma' \red{*}{\R}{}{B} \gamma$, i.e., $\alpha - \gamma \in
 \ideal{}{\rr}(\{ \beta \} \cup B)$, and hence $\alpha - \gamma \in
 \ideal{}{\rr/{\mathfrak i}}(\beta)$.
Therefore, (A2) is also fulfilled.
Finally Axiom (A3) holds since $\alpha \R_{\alpha} 0$ for all $\alpha \in \rr \backslash \{ 0 \}$
 implies $\alpha \myr_{\alpha} 0$.

Moreover, in case (A4) holds  in $\rr$ this is also true for $ \rr/{\mathfrak i}$:
For $\alpha, \beta, \gamma, \delta \in \rr/{\mathfrak i}$ we have that
 $\alpha \myr_{\beta}$ and $\beta \myr_{\gamma} \delta$ imply
 $\alpha \R_{\beta}$ and $\beta \R_{\gamma} \delta' \red{*}{\R}{}{B} \delta$
 and since $\alpha$ is $\R_{B}$-irreducible\footnote{Remember that in the proof
 of Lemma \ref{lem.BB'} we have shown that $\alpha \R_{\beta}$ and
 $\beta \red{*}{\R}{}{B'} 0$
 imply $\alpha \R_{B'}$. This carries over to our situation in the form that
 $\alpha \R_{\beta}$ and $\beta \R_{\gamma} \delta' \red{*}{\R}{}{B} \delta$
 implies $\alpha \R_{\{\gamma, \delta',\delta\}\cup B}$ and using induction
 to $\alpha \R_{\{\gamma,\delta\}\cup B}$.} this implies
 $\alpha \R_{\{\gamma,\delta\}}$ and hence $\alpha \myr_{\{\gamma,\delta\}}$.
\begin{theorem}\label{theo.quotients}~\\
{\sl
If $(\rr,\R)$ is a reduction ring with (A4), then for
 every finitely generated ideal ${\mathfrak i}$ the quotient
 $(\rr/{\mathfrak i},\myr)$ again is a reduction ring with (A4).
\theoend
}
\end{theorem}
\Ba{}~\\
Since reduction in $\rr/{\mathfrak i}$ as defined above inherits
 (A1) -- (A4) from $\rr$, 
 it remains to show that every finitely generated ideal  
 ${\mathfrak j} \subseteq \rr/{\mathfrak i}$ has a finite Gr\"obner basis.
Let ${\mathfrak j}_{\rr} = \{ \alpha \in \rr \mid \alpha\rnf{B}
 \in {\mathfrak j} \}$ be an ideal\footnote{${\mathfrak j}_{\rr}$ is an ideal
 in $\rr$ since
 \begin{enumerate}
 \item $0 \in {\mathfrak j}_{\rr}$ as $0 \in {\mathfrak j}$.
 \item $\alpha,\beta \in {\mathfrak j}_{\rr}$ implies 
       $\alpha\rnf{B},\beta\rnf{B} \in {\mathfrak j}$,
        hence $\alpha\rnf{B}+\beta\rnf{B}=(\alpha +\beta)\rnf{B} \in
        {\mathfrak j}$ and $\alpha +\beta \in {\mathfrak j}_{\rr}$.
 \item $\alpha \in {\mathfrak j}_{\rr}$ and $\gamma \in \rr$ implies
       $\alpha\rnf{B} \in {\mathfrak j}$ and $\gamma \skm \alpha\rnf{B}
       = (\gamma \skm \alpha)\rnf{B} \in {\mathfrak j}$,
       $ \alpha\rnf{B} \skm \gamma = (\alpha\skm \gamma)\rnf{B} \in {\mathfrak j}$,
        hence $\gamma \skm \alpha,\alpha\skm \gamma \in {\mathfrak j}_{\rr}$. 
 \end{enumerate}} in $\rr$ corresponding to ${\mathfrak j}$.
Then ${\mathfrak j}_{\rr}$ is {\em finitely} generated as an ideal in $\rr$ by
 its finite basis in $\rr/{\mathfrak i}$ viewed as elements of $\rr$ and the finite basis of ${\mathfrak i}$.
Hence ${\mathfrak j}_{\rr}$ has a finite Gr\"obner basis in $\rr$, say $G_{\rr}$.
Then $G = \{  \alpha\rnf{B} \mid  \alpha \in G_{\rr} \} \backslash \{ 0 \}$
 is a finite Gr\"obner basis of ${\mathfrak j}$:
If $\alpha \in {\mathfrak j}$ we have $\alpha \red{*}{\myr}{}{G} 0$ and
 $\ideal{}{\rr/{\mathfrak i}}(G) = {\mathfrak j}$, as every element which is
 reducible with an element $\beta \in G_{\rr}$ is also reducible with an
 element of $G \cup B$ because (A4) holds.
Since $G \cup B$ is also a Gr\"obner basis of ${\mathfrak j}_{\rr}$ and
 $\red{}{\myr}{}{G} \subseteq \red{*}{\R}{}{G \cup B}$, when restricted to
 elements in $\rr/{\mathfrak i}$ we have
 $\irr(\myr_{G}) = \irr(\R_{G \cup B})$ and $\myr_{G}$ is
 confluent.
Furthermore, since
 $\equiv_{\mathfrak j} \;=\; \equiv_{{\mathfrak j}_{\rr}}$
 when restricted to $\rr/{\mathfrak i}$ we get 
 $\red{*}{\lr}{}{G} =\;\equiv_{\mathfrak j}$ on $\rr/{\mathfrak i}$
 implying that $\rr/{\mathfrak i}$ is a reduction ring.
\\ \qed
In Example \ref{exa.Z} we have seen how to associate the integers 
 with a reduction relation $\R$ and in fact $(\z, \R)$ is
 a reduction ring.
Theorem \ref{theo.quotients} then states that for every $m \in \z$ the quotient
 $\z / \ideal{}{}(m)$ again is a reduction ring with respect to the reduction
 relation defined analogue to Definition \ref{def.quotientred}. 
In particular reduction rings with
 zero divisors can be constructed in this way.

Of course if we only assume that $\rr$ is a weak reduction ring
 we no longer have unique normal forms for the elements in the quotient.
Still comparing elements is possible as $\alpha = \beta$ in  $\rr/{\mathfrak i}$
 if and only if $\alpha - \beta \in {\mathfrak i}$ if and only if 
 $\alpha - \beta \red{*}{\R}{}{B} 0$ for a weak Gr\"obner basis $B$ of
 ${\mathfrak i}$.
Hence the elements in the quotient are no longer given by unique elements
 but by the respective sets of all representatives with respect to the weak Gr\"obner basis
 chosen for the ideal\footnote{Such an element $\alpha$ in the quotient can be represented by any element which is equivalent to it. When doing computations then of course to decide whether $\alpha = \beta$ in  $\rr/{\mathfrak i}$ one has to check if $\alpha - \beta \red{*}{\R}{}{B} 0$ for a weak Gr\"obner basis $B$ of
 ${\mathfrak i}$.}.
\begin{corollary}~\\
{\sl
If $(\rr,\R)$ is a weak reduction ring with (A4), then for
 every finitely generated ideal ${\mathfrak i}$ the quotient
 $(\rr/{\mathfrak i},\myr)$ again is a weak reduction ring with (A4).
\corend
}
\end{corollary}
\Ba{}~\\
It remains to show that every finitely generated ideal
 ${\mathfrak j} \subseteq \rr/{\mathfrak i}$ has a
 finite weak Gr\"obner basis.
Let $B$ be a finite weak Gr\"obner basis of ${\mathfrak i}$
 in $\rr$ and $B_{\mathfrak j}$ a finite generating set
 for the ideal ${\mathfrak j}$ in $\rr/{\mathfrak i}$.
\\
Let ${\mathfrak j}_{\rr} = \bigcup_{\alpha \in  {\mathfrak j}} 
 \{ \beta \in \rr \mid \beta \Longleftrightarrow^{*}_{B} \alpha \}$,
 be an ideal in $\rr$ corresponding to ${\mathfrak j}$.
Then ${\mathfrak j}_{\rr}$ is  {\em finitely} generated by
 the set $B \cup \tilde{B}_{\mathfrak j}$ where for each element $\alpha \in B_{\mathfrak j}$
 the set $\tilde{B}_{\mathfrak j}$ contains some 
 $\tilde{\alpha} \in \{ \beta \in \rr \mid \beta \Longleftrightarrow^{*}_{B} \alpha \}$. 
Moreover, ${\mathfrak j}_{\rr}$ has a finite weak Gr\"obner basis, say $G_{\rr}$.
Then the set $G = \{  \alpha\rnf{B} \mid  \alpha \in G_{\rr} \} \backslash \{ 0 \}$
 containing for each $\alpha \in G_{\rr}$ one not necessarily unique
 normal form $\alpha\rnf{B}$ is a finite weak Gr\"obner basis of ${\mathfrak j}$:
If $\alpha \in {\mathfrak j}$ we have $\alpha \red{*}{\myr}{}{G} 0$ and
 $\ideal{}{\rr/{\mathfrak i}}(G) = {\mathfrak j}$, as every element in
 ${\mathfrak j}$ (i.e.~in particular irreducible with respect to $B$)
 which is
 reducible with an element $\beta \in G_{\rr}$ is also reducible with an
 element of $G$ because (A4) holds\footnote{Since $\alpha \in {\mathfrak j}$ is irreducible by
 $B$, we have $\alpha \R_{\beta} \delta' \red{*}{\R}{}{G_{\rr}} \delta$
 and $\beta \not\in B$. Then looking at the situation $\alpha \R_{\beta}$
 and $\beta \red{*}{\R}{}{G_{\rr}} \beta\rnf{B}$, (A4)
 yields $\alpha \R_{\beta\rnf{B}}$.}.
\\ \qed

Now if $(\rr, \R)$ is an effective reduction ring, then $B$ can be computed and
 addition and  multiplication in $\rr/{\mathfrak i}$ as well as the reduction relation
 based on  Definition \ref{def.quotientred}
 are computable operations.
Moreover, Theorem \ref{theo.quotients} can be generalized:
\begin{corollary}\label{cor.quotients}~\\
{\sl
If $(\rr,\R)$ is an effective reduction ring with (A4), then for
 every finitely generated ideal ${\mathfrak i}$ the quotient
 $(\rr/{\mathfrak i},\myr)$ again is an effective reduction ring with (A4).
\corend
}
\end{corollary}
\Ba{}~\\
Given $\rr$, $B$ and a finite generating set $F$ for an ideal ${\mathfrak j}$ in $\rr/{\mathfrak i}$
 we can compute a finite Gr\"obner basis for ${\mathfrak j}$ using the method for computing Gr\"obner bases
 in $\rr$:
Compute a Gr\"obner basis $G_{\rr}$ of the ideal generated by $B \cup F$ in $\rr$.
Then the set $G = \{ \alpha\rnf{B} \mid \alpha \in G_{\rr} \}$, where 
 $\alpha\rnf{B}$ is the normal form of $g$ with respect to $\R_{B}$ in $\rr$ and hence an
 element of $\rr/{\mathfrak i}$, is a Gr\"obner basis of ${\mathfrak j}$ in $\rr/{\mathfrak i}$.
\\ \qed

The same is true for effective weak reduction rings.

Finally the results carry over to the case of one-sided reduction
 rings with (A4) provided that the two-sided ideal has a finite
  right respectively left Gr\"obner basis.
 
%%%%%%%%%%%%%%%%%%%%%%%%%%%%%%%%%%
\section{Sums of Reduction Rings}\label{section.sums}
Let $\rr_1,\rr_2$ be rings with reduction relations
 $\red{}{\R}{1}{}$ respectively $\red{}{\R}{2}{}$ fulfilling (A1) -- (A3).
Then $\rr = \rr_1 \times \rr_2 = \{ (\alpha_1, \alpha_2) \mid
 \alpha_1 \in \rr_1, \alpha_2 \in \rr_2 \}$ is called the \betonen{direct sum} 
 of $\rr_1$ and $\rr_2$.
Addition and multiplication are defined  component wise, the
 unit is $(1_1, 1_2)$ where $1_i$ is the respective unit in $\rr_i$.
A natural reduction relation can be defined on $\rr$ as follows:
\begin{definition}\label{def.sumred}~\\
{\rm
Let $\alpha = (\alpha_1, \alpha_2)$, $\beta = (\beta_1, \beta_2)$,
 $\gamma = (\gamma_1, \gamma_2) \in \rr$.
We say that
 $\beta$ \betonen{reduces} $\alpha$ to $\gamma$ in one step, 
 denoted by $\alpha \myr_{\beta} \gamma$, if either
 $(\alpha_1 \red{}{\R}{1}{\beta_1} \gamma_1$ and
 $\alpha_2 = \gamma_2)$ or
 $(\alpha_1 = \gamma_1$ and
 $\alpha_2  \red{}{\R}{2}{\beta_2} \gamma_2)$ or
 $(\alpha_1 \red{}{\R}{1}{\beta_1} \gamma_1$ and
 $\alpha_2  \red{}{\R}{2}{\beta_2} \gamma_2)$.
\dend
}
\end{definition}
Again we have to prove that the Axioms (A1) -- (A3) hold for the reduction relation in  $\rr$:
$\myr_{B} = \bigcup_{\beta \in B} \myr_{\beta}$ is terminating for
 finite sets $B \subseteq \rr$ since this property is inherited
 from the termination of the respective reduction relations in $\rr_i$.
Hence (A1) holds.
(A2) is satisfied since $\alpha \myr_{\beta} \gamma$ implies $\alpha - \gamma \in \ideal{}{\rr}(\beta)$.
(A3) is true as $\alpha \myr_{\alpha} (0_1,0_2)$ holds for all
 $\alpha \in \rr \backslash \{ (0_1,0_2) \}$.
Moreover, it is easy to see that if condition (A4) holds for $\red{}{\R}{1}{}$ and
 $\red{}{\R}{2}{}$ then this is inherited by $\myr$.
\begin{theorem}\label{theo.sum}~\\
{\sl
If $(\rr_1,\red{}{\R}{1}{})$, $(\rr_2,\red{}{\R}{2}{})$ are reduction rings, then $(\rr = \rr_1 \times \rr_2, \myr)$
  is again a reduction ring.
\theoend
}
\end{theorem}
\Ba{}~\\
Since the reduction relation in $\rr$ as defined above inherits
 (A1) -- (A3) respectively (A4) from the reduction relations in the $\rr_i$, 
it remains to show that every finitely generated ideal
 ${\mathfrak i} \subseteq \rr$ has a finite Gr\"obner basis.
To see this notice that the restrictions
 ${\mathfrak i}_1 = \{ \alpha_1 \mid (\alpha_1, \alpha_2) \in {\mathfrak i}
 \mbox{ for some } \alpha_2 \in \rr_2 \}$ and
 ${\mathfrak i}_2 = \{ \alpha_2 \mid (\alpha_1, \alpha_2) \in {\mathfrak i}
 \mbox{ for some } \alpha_1 \in \rr_1 \}$ are finitely
 generated ideals in $\rr_1$ respectively
 $\rr_2$ and hence have finite Gr\"obner bases $B_1$  respectively $B_2$.
We claim that $B =  \{ (\beta_1, 0_2), (0_1, \beta_2) \mid \beta_1 \in B_1,
 \beta_2 \in B_2 \}$ is a finite Gr\"obner basis of ${\mathfrak i}$.
Notice that ${\mathfrak i} = {\mathfrak i}_1 \times {\mathfrak i}_2$.
Then $\ideal{}{}(B) = {\mathfrak i}$ and $\alpha \in {\mathfrak i}$ implies
 $\alpha \red{*}{\myr}{}{B} (0_1,0_2)$ due to the fact that for 
 $\alpha = (\alpha_1, \alpha_2)$ we have $\alpha_1 \in {\mathfrak i}_1$
 and $\alpha_2 \in {\mathfrak i}_2$ implying $\alpha_1 \red{*}{\R}{1}{B_1} 0_1$
 and $\alpha_2 \red{*}{\R}{2}{B_2} 0_2$.
Similarly $\myr_B$ is confluent because $\red{}{\R}{1}{B_1}$ and
 $\red{}{\R}{2}{B_2}$ are confluent.
Finally $\red{*}{\lr}{}{B} = \;\equiv_{\mathfrak i}$ since
 $(\alpha_1,\alpha_2) \equiv_{\mathfrak i} (\beta_1, \beta_2)$
 implies $\alpha_1 \equiv_{{\mathfrak i}_1} \beta_1$ respectively
 $\alpha_2 \equiv_{{\mathfrak i}_2} \beta_2$ and hence
 $\alpha_1 \red{*}{\Longleftrightarrow}{1}{B_1} \beta_1$
 respectively
 $\alpha_2 \red{*}{\Longleftrightarrow}{2}{B_2} \beta_2$.
\\ \qed

Special regular rings as introduced by Weispfenning in \cite{We87}
 provide examples of such sums of reduction rings, e.g.~any direct sum
 of fields.

\begin{corollary}~\\
{\sl
If $(\rr_1,\red{}{\R}{1}{})$, $(\rr_2,\red{}{\R}{2}{})$ are weak reduction rings,
 then $(\rr = \rr_1 \times \rr_2, \myr)$
  is again a weak reduction ring.
\corend
}
\end{corollary}
\Ba{}~\\
Reviewing the proof of Theorem \ref{theo.sum} it remains to show that 
 every finitely generated ideal ${\mathfrak i} \subseteq \rr$ has a finite weak Gr\"obner basis.
Again we look at  the restrictions
 ${\mathfrak i}_1 = \{ \alpha_1 \mid (\alpha_1, \alpha_2) \in {\mathfrak i}
 \mbox{ for some } \alpha_2 \in \rr_2 \}$ and
 ${\mathfrak i}_2 = \{ \alpha_2 \mid (\alpha_1, \alpha_2) \in {\mathfrak i}
 \mbox{ for some } \alpha_1 \in \rr_1 \}$ which are finitely
 generated ideals in $\rr_1$ respectively
 $\rr_2$ and hence have finite weak Gr\"obner bases $B_1$  respectively $B_2$.
We claim that $B =  \{ (\beta_1, 0_2), (0_1, \beta_2) \mid \beta_1 \in B_1,
 \beta_2 \in B_2 \}$ is a finite weak Gr\"obner basis of ${\mathfrak i}$.
As before ${\mathfrak i} = {\mathfrak i}_1 \times {\mathfrak i}_2$ and
 $\ideal{}{}(B) = {\mathfrak i}$.
Then $\alpha \in {\mathfrak i}$ implies
 $\alpha \red{*}{\myr}{}{B} (0_1,0_2)$ due to the fact that for 
 $\alpha = (\alpha_1, \alpha_2)$ we have $\alpha_1 \in {\mathfrak i}_1$
 and $\alpha_2 \in {\mathfrak i}_2$ implying $\alpha_1 \red{*}{\R}{1}{B_1} 0_1$
 and $\alpha_2 \red{*}{\R}{2}{B_2} 0_2$ as $B_1$ and $B_2$ are respective
 weak Gr\"obner bases, and we are done.
\\ \qed

Now if $(\rr_1,\red{}{\R}{1}{})$, $(\rr_2,\red{}{\R}{2}{})$ are effective reduction rings, 
 then  addition and  multiplication in $\rr$ as well as the reduction relation based on 
 Definition \ref{def.sumred} are computable operations.
Moreover, Theorem \ref{theo.sum} can be generalized:
\begin{corollary}\label{cor.sum}~\\
{\sl
If $(\rr_1,\red{}{\R}{1}{})$, $(\rr_2,\red{}{\R}{2}{})$ are effective reduction rings, 
 then $(\rr = \rr_1 \times \rr_2, \myr)$
  is again an effective reduction ring.
\corend
}
\end{corollary}
\Ba{}~\\                                                                            
Given a finite generating set 
 $F = \{ (\alpha_i, \beta_i) \mid 1 \leq i \leq k, \alpha_i \in \rr_1, \beta_i \in \rr_2 \}$
 a Gr\"obner basis of the ideal generated by $F$ can be computed using the respective
 methods for Gr\"obner basis computation in $\rr_1$ and $\rr_2$:
Compute $B_1$ a Gr\"obner basis of the ideal generated by $\{ \alpha_1, \ldots, \alpha_k \}$ in $\rr_1$
 and $B_2$ a Gr\"obner basis of the ideal generated by $\{ \beta_1, \ldots, \beta_k \}$ in $\rr_2$.
Then $B =  \{ (\gamma_1, 0_2), (0_1, \gamma_2) \mid \gamma_1 \in B_1,
 \gamma_2 \in B_2 \}$ is a finite Gr\"obner basis of the ideal generated by $F$ in $\rr$.
\\ \qed

A similar result holds for effective weak reduction rings.

Due to the ``simple'' multiplication used when defining direct
 sums, Theorem \ref{theo.sum} and Corollary \ref{cor.sum} extend directly to one-sided reduction rings.
More complicated multiplications are possible and have to be treated individually. 
%%%%%%%%%%%%%%%%%%%%%%%%%%%%%%%%%%
\section{Modules over Reduction Rings}\label{section.modules}
Another structure which can be studied by reduction techniques are modules and
 their submodules.
Given a  ring $\rr$ with unit $1$ and a natural number $k$, 
 let $\rr^k = \{ {\bf a}=(\alpha_1, \ldots, \alpha_k) \mid \alpha_i \in \rr \}$
 be the set of all vectors of length $k$ with coordinates in $\rr$.
Obviously $\rr^k$ is an additive commutative group with respect to ordinary vector addition and we
 denote the zero by ${\bf 0}$.
Moreover,
$\rr^k$ is an \betonen{$\rr$-module} for scalar
 multiplication defined as $\alpha \mmult (\alpha_1, \ldots, \alpha_k) = 
 (\alpha \skm \alpha_1, \ldots, \alpha \skm \alpha_k)$ and $ (\alpha_1, \ldots, \alpha_k) \mmult \alpha = 
(\alpha_1 \skm \alpha , \ldots, \alpha_k \skm \alpha)$.
Additionally $\rr^k$ is called \betonen{free} as it has a basis\footnote{Here the term {\em basis} is
 used in the meaning of being a linearly independent set of generating vectors.}.
One such basis is the set of unit vectors ${\bf e}_1 = (1, 0, \ldots, 0), {\bf e}_2 = (0, 1, 0, \ldots, 0),
 \ldots, {\bf e}_k = (0, \ldots, 0, 1)$.
Using this basis the elements of $\rr^k$ can be written uniquely as ${\bf a} = \sum_{i=1}^k \alpha_i\mmult {\bf e}_i$
 where ${\bf a} =  (\alpha_1, \ldots, \alpha_k)$.
\begin{definition}\label{def.redmodule}~\\
{\rm
A subset of $\rr^k$ which is again an $\rr$-module is called a \betonen{submodule} of $\rr^k$.
\dend
}
\end{definition}
For example any ideal of $\rr$ is an $\rr$-module and even a submodule of the $\rr$-module $\rr^1$.
Provided a set of vectors $S = \{ {\bf a}_1, \ldots, {\bf a}_n \}$ the set
 $\{ \sum_{i=1}^n  \sum_{j=1}^{m_i} \beta_{ij} \mmult {\bf a}_i \mmult {\beta_{ij}}' \mid \beta_{ij}, {\beta_{ij}}' \in \rr \}$ is a submodule of $\rr^k$.
This set is denoted as $\langle S \rangle$ and $S$ is called its generating set.

Now similar to the case of modules over commutative polynomial rings,
 being Noetherian is inherited by $\rr^k$ from $\rr$.
\begin{theorem}\label{theo.module.noetherian}~\\  
{\sl
Let $\rr$ be a Noetherian ring.
Then every submodule in $\rr^k$ is also finitely generated. 
\theoend
} 
\end{theorem}
\Ba{}~\\
Let ${\cal S}$ 
 be a submodule of $\rr^k$.
We show our claim by induction on $k$.
For $k=1$ we find that ${\cal S}$ is in fact an ideal in $\rr$ and hence by our 
 hypothesis must be finitely generated.
For $k>1$ let us look at the set ${\mathfrak i} = \{ \beta_1 \mid (\beta_1, \ldots, \beta_k) \in {\cal S} \}$
 which is again
 an ideal in $\rr$ and hence finitely generated by some set
 $\{\gamma_1, \ldots, \gamma_s \mid \gamma_i \in \rr \}$.
Choose\footnote{In this step we need the Axiom of Choice and hence the construction is
 not constructive.} $H = \{ {\bf c}_1, \ldots, {\bf c}_s \} \subseteq {\cal S}$
 such that the first coordinate of 
 ${\bf c}_i$ is $\gamma_i$.
Similarly the set ${\cal M} = \{ (\beta_2, \ldots, \beta_{k})
 \mid (0, \beta_2, \ldots, \beta_k) \in {\cal S} \}$
 is a submodule in
 $\rr^{k-1}$ and therefore  finitely generated by
 our induction hypothesis.
Let $\{ (\delta_2^i, \ldots, \delta^i_{k}) \mid 1 \leq i \leq w\}$ be such a 
 finite generating set.
Then ${\bf d}_i = (0,\delta_2^i, \ldots, \delta^i_{k} ) \in {\cal S}$, $1 \leq i \leq w$
 and the set $G = \{{\bf c}_1, \ldots, {\bf c}_s \} \cup 
 \{ {\bf d}_i \mid  1 \leq i \leq w \}$ is a
 finite generating set for ${\cal S}$.
To see this assume ${\bf t} = (\tau_1, \ldots, \tau_k)\in {\cal S}$.   
Then $\tau_1 = \sum_{i=1}^s \sum_{j=1}^{n_i} \zeta_{ij} \skm \gamma_i \skm {\zeta_{ij}}'$
 for some $\zeta_{ij}, {\zeta_{ij}}' \in \rr$ and
 ${\bf t}' = {\bf t} - \sum_{i=1}^s \sum_{j=1}^{n_i} \zeta_{ij} \mmult {\bf c}_i \mmult {\zeta_{ij}}' \in {\cal S}$ 
 with first coordinate $0$.
Hence ${\bf t}' =  \sum_{i=1}^w \sum_{j=1}^{m_i} \eta_{ij} \mmult {\bf d}_i \mmult {\eta_{ij}}'$ 
 for some $\eta_{ij}, {\eta_{ij}}' \in \rr$ giving rise to
$${\bf t} = {\bf t}' + \sum_{i=1}^s \sum_{j=1}^{n_i} \zeta_{ij} \mmult {\bf c}_i \mmult {\zeta_{ij}}' = 
 \sum_{i=1}^w \sum_{j=1}^{m_i}  \eta_{ij} \mmult {\bf d}_i \mmult {\eta_{ij}}' + 
 \sum_{i=1}^s  \sum_{j=1}^{n_i} \zeta_{ij} \mmult {\bf c}_i \mmult {\zeta_{ij}}'.$$
\\ \qed
We will now study submodules of
 modules using reduction relations.
Let $\R$ be a reduction relation on $\rr$ fulfilling (A1) -- (A3).
A natural reduction relation on $\rr^k$ can be defined using the representations as polynomials with respect to
 the basis of unit vectors as follows:
\begin{definition}\label{def.module.red}~\\ 
{\rm  
Let ${\bf a} = \sum_{i=1}^k \alpha_i \mmult {\bf e}_i$, ${\bf b} = \sum_{i=1}^k \beta_i \mmult {\bf e}_i \in \rr^k$.
We say that
 ${\bf b}$ \betonen{reduces} ${\bf a}$ to ${\bf c}$ at $\alpha_s \mmult {\bf e}_s$ in one step, 
 denoted by ${\bf a} \myr_{\bf b} {\bf c}$, if 
\begin{itemize}
\item[(a)] $\beta_j = 0$ for $1 \leq j < s$,
\item[(b)] $\alpha_s \R_{\beta_s} \gamma_s$ with  
           $\alpha_s = \gamma_s + \sum_{i=1}^n \delta_i \skm \beta_s \skm {\delta_i}'$, 
           $\delta_i,{\delta_i}' \in \rr$, and
\item[(c)] ${\bf c} = {\bf a} -  \sum_{i=1}^n \delta_i \mmult {\bf b} \mmult {\delta_i}' =
            (\alpha_1, \ldots , \alpha_{s-1}, \gamma_s, 
             \alpha_{s+1} - \sum_{i=1}^n \delta_i \skm \beta_{s+1} \skm {\delta_i}', \ldots , 
             \alpha_k - \sum_{i=1}^n \delta_i \skm \beta_k \skm {\delta_i}')$.
\dend
\end{itemize}
}
\end{definition}
The Axioms (A1) -- (A3) hold for this reduction relation on  $\rr^k$:
$\myr_{B} = \bigcup_{{\bf b} \in B} \myr_{\bf b}$ is terminating for
 finite $B \subseteq \rr^k$ since this property is inherited
 from the termination of the respective reduction relation $\R$ in $\rr$.
Hence (A1) holds.
(A2) is satisfied now of course in the context of submodules since ${\bf a} \myr_{\bf b} {\bf c}$ implies 
 ${\bf a} - {\bf c} \in  \langle \{ {\bf b } \} \rangle$.
(A3) is true as ${\bf a} \myr_{\bf a} {\bf 0}$ holds for all
 ${\bf a} \in \rr^k \backslash \{ {\bf 0}  \}$.
Moreover, it is easy to see that if condition (A4) holds for $\red{}{\R}{}{}$ then this is inherited by 
 $\myr$ as defined in Definition \ref{def.module.red} for $\rr^k$.
First we show how the existence of weak Gr\"obner bases carries over for Noetherian $\rr$.
\begin{definition}~\\
{\rm
A subset $B$ of  $\rr^k$ is called a 
 \betonen{weak Gr\"obner basis} of the submodule ${\cal S} = \langle B \rangle$, if
  $\red{}{\myr}{}{B}$ is terminating and ${\bf a} \red{*}{\myr}{}{B} {\bf 0}$ for all
 ${\bf a} \in {\cal S}$. 
\dend
}
\end{definition}
\begin{theorem}\label{theo.module.weak}~\\  
{\sl
Let $\rr$ be a Noetherian ring with reduction relation
 $\red{}{\R}{}{}$ fulfilling (A1) -- (A3).
If in $\rr$ every ideal has a finite weak Gr\"obner basis, then 
 the same holds for submodules in $(\rr^k, \myr)$. 
\theoend
} 
\end{theorem}
\Ba{}~\\
Let ${\cal S}$ 
 be a submodule of $\rr^k$.
We show our claim by induction on $k$.
For $k=1$ we find that ${\cal S}$ is in fact an ideal\footnote{At this point we 
 could also proceed with a much weaker hypothesis, namely instead of requiring
 $\rr$ to be Noetherian assuming that ${\cal S}$ is finitely generated. 
 Then still the fact that $\rr$ is supposed to be a weak reduction ring would imply
 the existence of a finite weak Gr\"obner basis for ${\cal S}$.}
 in $\rr$ and hence by our 
 hypothesis must have a finite weak Gr\"obner basis.
For $k>1$ let us look at the set
 ${\mathfrak i} = \{ \beta_1 \mid (\beta_1, \ldots, \beta_k) \in {\cal S} \}$
 which is again an ideal\footnote{Here it still would
 be sufficient to require that ${\cal S}$ is finitely generated as the first
 coordinates of a finite generating set for ${\cal S}$ then would generate 
 ${\mathfrak i}$ hence implying that the ideal is finitely generated as well.}.
Hence ${\mathfrak i}$ must have a finite weak Gr\"obner basis
 $\{\gamma_1, \ldots, \gamma_s \mid \gamma_i \in \rr \}$.
Choose $H = \{ {\bf c}_1, \ldots, {\bf c}_s \} \subseteq {\cal S}$
 such that the first coordinate of  ${\bf c}_i$ is $\gamma_i$.
Similarly the set ${\cal M} = \{ (\beta_2, \ldots, \beta_{k}) \mid
 (0, \beta_2, \ldots, \beta_k) \in {\cal S} \}$
 is a submodule\footnote{{\em Now} we really need that $\rr^{k-1}$ is
 Noetherian. Assuming that ${\cal S}$ is finitely generated would {\em not}
 help to deduce that ${\cal M}$ is finitely generated.} in
 $\rr^{k-1}$  
 which by our induction hypothesis must have a finite weak Gr\"obner basis 
 $\{ (\delta_2^i, \ldots, \delta^i_{k}) \mid 1 \leq i \leq w\}$.
Then the set $G = \{{\bf c}_1, \ldots, {\bf c}_s \} \cup 
 \{ {\bf d}_i = (0,\delta_2^i, \ldots, \delta^i_{k} ) \mid  1 \leq i \leq w \}$ is a
 weak Gr\"obner basis for ${\cal S}$.
\\
That $G$ is a generating set for ${\cal S}$ follows as in the proof of Theorem
 \ref{theo.module.noetherian}.
It remains to show that $G$ is in fact a weak Gr\"obner basis, i.e., for every
 ${\bf t} = (\tau_1, \ldots, \tau_k)\in {\cal S}$ we have
 ${\bf t} \red{*}{\myr}{}{G} {\bf 0}$.
Since $\tau_1 \red{*}{\R}{}{\{\gamma_1, \ldots, \gamma_s\}} 0$ with
 $\tau_1 = \sum_{i=1}^s \sum_{j=1}^{n_i} \zeta_{ij} \skm \gamma_i \skm {\zeta_{ij}}'$,
 by the definition of $G$ we get
 ${\bf t} \red{*}{\myr}{}{\{{\bf c}_1, \ldots, {\bf c}_s \} }
 {\bf t} - \sum_{i=1}^s \sum_{j=1}^{n_i} \zeta_{ij} \mmult {\bf c}_i \mmult {\zeta_{ij}}' = {\bf t}'$
 where ${\bf t}' = (0, {\tau_2}', \ldots, {\tau_k}') \in {\cal M}$.
Hence, as $({\tau_2}',\ldots, {\tau_k}')
 \red{*}{\myr}{}{\{(\delta_2^i, \ldots, \delta^i_{k} ) \mid  1 \leq i \leq w \}}
 {\bf 0}$, we get ${\bf t} \red{*}{\myr}{}{G} {\bf 0}$ and are done. 
\\ \qed
Now we turn our attention to Gr\"obner bases of submodules in $\rr^k$.
\begin{definition}~\\
{\rm
A subset $B$ of  $\rr^k$ is called a 
 \betonen{Gr\"obner basis} of the submodule ${\cal S} = \langle B \rangle$, if
 $\red{*}{\lr}{}{B} = \;\;\equiv_{{\cal S}}$
 and $\red{}{\myr}{}{B}$ is complete. 
\dend
}
\end{definition}
\begin{theorem}\label{theo.module.GB}~\\  
{\sl
Let $\rr$ be a Noetherian
 ring with reduction relation $\red{}{\R}{}{}$ fulfilling (A1) -- (A3).
If in $\rr$ every  ideal has a finite Gr\"obner basis, then 
 the same holds for submodules in $(\rr^k, \myr)$. 
\theoend
} 
\end{theorem}
\Ba{}~\\
The candidate for the Gr\"obner basis can be built similar to
 the set $G$ in the proof of Theorem
  \ref{theo.module.weak} now of course using Gr\"obner bases in the
 construction instead of weak Gr\"obner bases:
Let ${\cal S}$ be a submodule of $\rr^k$.
We show our claim by induction on $k$.
For $k=1$ we find that ${\cal S}$ is in fact an ideal in $\rr$ and hence by our 
 hypothesis must have a finite Gr\"obner basis.
For $k>1$ let us look at the set ${\mathfrak i} = \{ \beta_1 \mid (\beta_1, \ldots, \beta_k) \in {\cal S} \}$ which is again an ideal in $\rr$.
Hence ${\mathfrak i}$ must have a finite  Gr\"obner basis 
 $\{\gamma_1, \ldots, \gamma_s \mid \gamma_i \in \rr \}$ by our assumption.
Choose $H = \{ {\bf c}_1, \ldots, {\bf c}_s \} \subseteq {\cal S}$ such that the first coordinate of 
 ${\bf c}_i$ is $\gamma_i$.
Similarly the set ${\cal M} = \{ (\beta_2, \ldots, \beta_{k}) \mid (0, \beta_2, \ldots, \beta_k) \in {\cal S} \}$
 is a submodule in
 $\rr^{k-1}$  finitely generated as $\rr^{k-1}$ is Noetherian.
Hence by our induction hypothesis ${\cal M}$ then must have a finite Gr\"obner basis 
 $\{ (\delta_2^i, \ldots, \delta^i_{k})\mid 1 \leq i \leq w\}$.
Let $G = \{{\bf c}_1, \ldots, {\bf c}_s \} \cup 
 \{ {\bf d}_i = (0,\delta_2^i, \ldots, \delta^i_{k} )\mid  1 \leq i \leq w \}$.
Since $G$ generates ${\cal S}$ (see the proof of Theorem \ref{theo.module.weak})
 it remains to show that it is a Gr\"obner basis.
\\
By the definition of the reduction relation in $\rr^k$ we immediately find
 $\red{*}{\longleftrightarrow}{}{G} \subseteq \;\;\equiv_{{\cal S}}$.
To see the converse let ${\bf r} = (\rho_1, \ldots , \rho_k) \equiv_{{\cal S}}
 {\bf s} = (\sigma_1, \ldots, \sigma_k)$.
Then as $\rho_1 \equiv_{ \{ \beta_1 \mid {\bf b} = (\beta_1, \ldots, \beta_k) \in {\cal S} \}} \sigma_1$ 
 by the definition of $G$ we get 
 $\rho_1 \red{*}{\Longleftrightarrow}{}{\{\gamma_1, \ldots, \gamma_s\}} \sigma_1$.
But this gives us 
 ${\bf r} \red{*}{\longleftrightarrow}{}{H} 
  {\bf r} + \sum_{i=1}^s \sum_{j=1}^{m_i} \chi_{ij} \mmult {\bf c}_i \mmult {\chi_{ij}}' =
  {\bf  r}' = (\sigma_1,{\rho_2}', \ldots, {\rho_k}')$
 and we get
 $(\sigma_1,{\rho_2}', \ldots, {\rho_k}')\equiv_{\cal S} (\sigma_1, \ldots, \sigma_k)$.
Hence $(\sigma_1,{\rho_2}', \ldots, {\rho_k}') - (\sigma_1, \ldots, \sigma_k)
  = (0, {\rho_2}' - \sigma_2, \ldots, {\rho_k}' - \sigma_k) \in {\cal S}$, implying
  $({\rho_2}' - \sigma_2, \ldots, {\rho_k}' - \sigma_k) \in {\cal M}$.
Now we have to be more careful since we cannot conclude that
  $({\rho_2}', \ldots, {\rho_k}'),(\sigma_2, \ldots, \sigma_k) \in {\cal M}$.
But we know $(\sigma_1,{\rho_2}', \ldots, {\rho_k}') = 
 (\sigma_1, \ldots, \sigma_k) + (0, {\rho_2}' - \sigma_2, \ldots, {\rho_k}' - \sigma_k) = 
 (\sigma_1, \ldots, \sigma_k) + \sum_{i=1}^w \sum_{j=1}^{n_i} \eta_{ij} \mmult {\bf d}_i \mmult {\eta_{ij}}'$ 
 where $(0, {\rho_2}' - \sigma_2, \ldots, {\rho_k}' - \sigma_k) = 
 \sum_{i=1}^w \sum_{j=1}^{n_i} \eta_{ij} \mmult {\bf d}_i \mmult {\eta_{ij}}'$
 for  $\eta_{ij},{\eta_{ij}}' \in \rr$, i.e., $(\sigma_1,{\rho_2}', \ldots, {\rho_k}') 
 \equiv_{\langle {\bf d}_1, \ldots, {\bf d}_w\rangle} (\sigma_1, \ldots, \sigma_k)$.
Hence, as $\{ (\delta_2^i, \ldots, \delta^i_{k})\mid 1 \leq i \leq w\}$
 is a Gr\"obner basis of ${\cal M}$ both vectors $(\sigma_1,{\rho_2}', \ldots, {\rho_k}')$ and
 $(\sigma_1, \ldots, \sigma_k)$ must have a common normal form using
 $\{ {\bf d}_i = (0,\delta_2^i, \ldots, \delta^i_{k})
 \mid  1 \leq i \leq w \}$ for reduction\footnote{The elements in this set cannot influence the
 first coordinate which is $\sigma_1$ for both vectors.} and we are done.
\\
The same argument applies to show local confluence. 
Let us assume there are ${\bf r}$, ${\bf s}_1$, ${\bf s}_2 \in \rr^k$
 such that ${\bf r} \myr_{G} {\bf s}_1$ and ${\bf r} \myr_{G} {\bf s}_2$.
Then by the definition of $G$, the first coordinates $\sigma^1_1$ and $\sigma_1^2$ 
 of ${\bf s}_1$ respectively ${\bf s}_2$ are joinable
 by $\{\gamma_1, \ldots, \gamma_s \}$ to some element, say $\sigma$, giving rise to the elements
 ${\bf r}_1 = {\bf s}_1 + \sum_{i=1}^s \sum_{j=1}^{n_i} \chi_{ij} \mmult {\bf c}_i \mmult {\chi_{ij}}'$ and
 ${\bf r}_2 = {\bf s}_2 + \sum_{i=1}^s \sum_{j=1}^{m_i} \psi_{ij} \mmult {\bf c}_i \mmult {\psi_{ij}}'$ with first
 coordinate $\sigma$.
Again  we know $(\sigma,{\rho^1_2}, \ldots, {\rho^1_k}) = 
 (\sigma,{\rho^2_2}, \ldots, {\rho^2_k}) + (0,{\rho^1_2} - {\rho^2_2}, \ldots, {\rho^1_k} - {\rho^2_k})$
 with $({\rho^1_2} - {\rho^2_2}, \ldots, {\rho^1_k} - {\rho^2_k}) \in {\cal M}$.
Hence $(\sigma,{\rho^1_2}, \ldots, {\rho^1_k}) = 
 (\sigma,{\rho^2_2}, \ldots, {\rho^2_k})
  + \sum_{i=1}^w \sum_{j=1}^{n_i} \eta_{ij} \mmult {\bf d}_i \mmult {\eta_{ij}}'$
 for $\eta_{ij},{\eta_{ij}}' \in \rr$, i.e., $(\sigma_1,{\rho_2}', \ldots, {\rho_k}') 
 \equiv_{\langle {\bf d}_1, \ldots, {\bf d}_w\rangle} (\sigma_1, \ldots, \sigma_k)$.
As again  $\{ (\delta_2^i, \ldots, \delta^i_{k})\mid 1 \leq i \leq w\}$
 is a Gr\"obner basis of ${\cal M}$
 both vectors must have a common normal with respect to reduction using
 $\{ {\bf d}_i = (0,\delta_2^i, \ldots, \delta^i_{k})
 \mid  1 \leq i \leq w \}$.
\\ \qed
Let us close this section with a remark on why the additional property of
 being Noetherian is so important.
In the proofs of Theorem \ref{theo.module.weak} and \ref{theo.module.GB}
 in the induction step the ``projection'' of ${\cal S}$ on $\rr^{k-1}$
 plays an essential role.
If this projection is defined as ${\cal M} = \{ (\beta_2, \ldots, \beta_{k})
 \mid (0, \beta_2, \ldots, \beta_k) \in {\cal S} \}$ we have to show
 that this module is again finitely generated.
In assuming Noetherian for $\rr$ this then follows as ${\cal M}$ is a
 submodule of $\rr^{k-1}$ which is again Noetherian.
Assuming that ${\cal S}$ is finitely generated by some set
 $\{{\bf a}_1, \ldots, {\bf a}_n \}$ does not improve the situation as
 in general we cannot extract a finite generating set for ${\cal M}$ from this set\footnote{
   Another idea might be to look at an other projection of ${\cal S}$: 
   ${\cal M'} = \{ (\beta_2, \ldots, \beta_{k}) \mid \mbox{ there exists } 
   \beta_1 \in \rr \mbox{ such that } 
   (\beta_1, \beta_2, \ldots, \beta_k) \in {\cal S} \}$.
   ${\cal M'}$ then is again a module now finitely generated by
   $(\alpha^1_2, \ldots, \alpha^1_{k}), \ldots, (\alpha^n_2, \ldots, \alpha^n_{k})$.
   Unfortunately in this case having a Gr\"obner basis for this module is
   of no use as we can no longer lift this special basis to $\rr^k$.
   The trick with adding $0$ as the first coordinate will no longer work
   as for some $(\gamma_2, \ldots, \gamma_k) \in {\cal M'}$ we only know
   that there exists some $\gamma \in \rr$ such that 
   $(\gamma,\gamma_2, \ldots, \gamma_k) \in {\cal S}$ and we cannot enforce
   that $\gamma = 0$.
   However, if we lift the set by adding appropriate elements $\gamma \in\rr$
   as first coordinates, then the resulting set does not lift the Gr\"obner basis
   properties for the reduction relation.
   Especially in the induction step the first coordinate of the vector being modified
   can no longer be expected to be left unchanged which is the case when using
   vectors with first coordinate $0$ for reduction.}.
The situation improves if we look at one-sided reduction rings $\rr$
 and demand that in $\rr$ all  (left respectively right) syzygy modules have finite
 bases.

$\rr^k$ is a right $\rr$-module 
 with scalar multiplication $(\alpha_1, \ldots, \alpha_k) \mmult \alpha = 
(\alpha_1 \skm \alpha , \ldots, \alpha_k \skm \alpha)$. 
Provided a finite subset $\{ \alpha_1, \ldots, \alpha_n \} \subseteq \rr$ the set  
  of solutions of the equation $\alpha_1 \skm X_1 + \ldots + \alpha_n \skm X_n = 0$ 
 is a submodule of the right $\rr$-module $\rr^n$. 
It is called the (first) module of syzygies of $\{ \alpha_1, \ldots, \alpha_n \}$ in the literature. 
We will see that these special modules can be used to characterize Gr\"obner bases of submodules  
 in $\rr^k$.   

A reduction relation can be defined similarly to Definition \ref{def.module.red}. 
\begin{definition}\label{def.module.red.right}~\\  
{\rm   
Let ${\bf a} = \sum_{i=1}^k {\bf e}_i \mmult \alpha_i$, 
    ${\bf b} = \sum_{i=1}^k {\bf e}_i \mmult \beta_i \in \rr^k$.
We say that
 ${\bf b}$ \betonen{right reduces} ${\bf a}$ to ${\bf c}$ at the monomial
 ${\bf e}_s \mmult \alpha_s$ in one step, 
 denoted by ${\bf a} \red{}{\myr}{r}{\bf b} {\bf c}$, if 
\begin{itemize}
\item[(a)] $\beta_j = 0$ for $1 \leq j < s$,
\item[(b)] $\alpha_s \R_{\beta_s} \gamma_s$ with  
           $\alpha_s = \gamma_s +  \beta_s \skm {\delta}$, 
           $\delta \in \rr$, and
\item[(c)] ${\bf c} = {\bf a} -  {\bf b} \mmult \delta =
            (\alpha_1, \ldots , \alpha_{s-1}, \gamma_s, 
             \alpha_{s+1} - \beta_{s+1} \skm \delta, \ldots , 
             \alpha_k - \beta_k \skm \delta)$.
\dend
\end{itemize}
}
\end{definition}
\begin{theorem}~\\
{\sl
Let $\rr$ be a ring with a right reduction relation $\red{}{\R}{}{}$ fulfilling (A1) -- (A3).
Additionally let every right module of syzygies in $\rr$ have a finite basis.
If every finitely generated right ideal in $\rr$ has a finite Gr\"obner basis, 
 then the same holds for every finitely generated right submodule in $(\rr^k, \myr)$.
\theoend
}
\end{theorem}
\Ba{}~\\
Again the candidate for the right Gr\"obner basis can be built similar to
 the set $G$ in the proofs of Theorem
  \ref{theo.module.weak} and \ref{theo.module.GB}:
Let ${\cal S}$ be a right submodule of $\rr^k$ which is finitely generated by a set
 $\{{\bf a}_1, \ldots, {\bf a}_n \}$.
We show our claim by induction on $k$.
For $k=1$ we find that ${\cal S}$ is in fact a finitely generated right  ideal in $\rr$ and hence by our 
 hypothesis must have a finite right Gr\"obner basis.
For $k>1$ let us look at the set ${\mathfrak i} = \{ \beta_1 \mid (\beta_1, \ldots, \beta_k) \in {\cal S} \}$ which is again a right  ideal in $\rr$ finitely generated by $\{ \alpha_1^1, \ldots, \alpha_1^n \}$ where
 ${\bf a}_i = (\alpha_1^i, \ldots, \alpha_k^i)$.
Hence ${\mathfrak i}$ must have a finite  right  Gr\"obner basis 
 $\{\gamma_1, \ldots, \gamma_s \mid \gamma_i \in \rr \}$ by our assumption.
Choose $H = \{ {\bf c}_1, \ldots, {\bf c}_s \} \subseteq {\cal S}$ such that the first coordinate of 
 ${\bf c}_i$ is $\gamma_i$.
On the other hand the right syzygy module 
 $\{ (\psi_1 , \ldots , \psi_n) \mid \sum_{i=1}^n \alpha_1^i \skm \psi_i = 0, \psi_i \in \rr \}$ has a finite basis 
 $B = \{ (\beta_1^j, \ldots, \beta_n^j) \mid 1 \leq j \leq m \} \subseteq \rr^n$.
Then the set
 $\{ \sum_{i=1}^n {\bf a}_i \mmult \beta_i^j \mid 1 \leq j \leq m \} \cup \{ {\bf a}_i \mid \alpha_1^i = 0, 1 \leq i \leq n\}$
 is a finite generating set for the submodule
 ${\cal M} = \{ (\beta_2, \ldots, \beta_{k}) \mid (0, \beta_2, \ldots, \beta_k) \in {\cal S} \}$
 of  $\rr^{k-1}$.
To see this let $(0, \beta_2, \ldots, \beta_k) \in {\cal S}$.
Then $(0, \beta_2, \ldots, \beta_k) = \sum_{i=1}^n {\bf a}_i \mmult \zeta_i$, $\zeta_i \in \rr$
 implies $\sum_{i=1}^n \alpha_1^i \skm \zeta_i = 0$ and hence $(\zeta_1 , \ldots , \zeta_n)$
 lies in the right syzygy module and we are done.
Hence by our induction hypothesis ${\cal M}$ then must have a finite right  Gr\"obner basis 
 $\{ (\delta_2^i, \ldots, \delta^i_{k})\mid 1 \leq i \leq w\}$.
Let $G = \{{\bf c}_1, \ldots, {\bf c}_s \} \cup 
 \{ {\bf d}_i = (0,\delta_2^i, \ldots, \delta^i_{k} )\mid  1 \leq i \leq w \}$.
Since $G$ generates ${\cal S}$ it remains to show that it is a right  Gr\"obner basis.
By the definition of the reduction relation in $\rr^k$ we immediately find
 $\red{*}{\longleftrightarrow}{}{G} \subseteq \;\;\equiv_{{\cal S}}$.
To see the converse let ${\bf r} = (\rho_1, \ldots , \rho_k) \equiv_{{\cal S}}
 {\bf s} = (\sigma_1, \ldots, \sigma_k)$.
Then as $\rho_1 \equiv_{ \{ \alpha_1 \mid {\bf a} = (\alpha_1, \ldots, \alpha_k) \in {\cal S} \}} \sigma_1$ 
 by the definition of $G$ we get 
 $\rho_1 \red{*}{\longleftrightarrow}{}{\{\gamma_1, \ldots, \gamma_s\}} \sigma_1$.
But this gives us 
 ${\bf r} \red{*}{\Longleftrightarrow}{}{H} 
  {\bf r} + \sum_{i=1}^s {\bf c}_i \mmult \chi_{i} =
  {\bf  r}' = (\sigma_1,{\rho_2}', \ldots, {\rho_k}')$, $\chi_i \in \rr$,
 and we get
 $(\sigma_1,{\rho_2}', \ldots, {\rho_k}')\equiv_{\cal S} (\sigma_1, \ldots, \sigma_k)$.
Hence $(\sigma_1,{\rho_2}', \ldots, {\rho_k}') - (\sigma_1, \ldots, \sigma_k)
  = (0, {\rho_2}' - \sigma_2, \ldots, {\rho_k}' - \sigma_k) \in {\cal S}$ implying
  $({\rho_2}' - \sigma_2, \ldots, {\rho_k}' - \sigma_k) \in {\cal M}$.
Now we have to be more careful since we cannot conclude that
  $({\rho_2}', \ldots, {\rho_k}'),(\sigma_2, \ldots, \sigma_k) \in {\cal M}$.
But we know $(\sigma_1,{\rho_2}', \ldots, {\rho_k}') = 
 (\sigma_1, \ldots, \sigma_k) + (0, {\rho_2}' - \sigma_2, \ldots, {\rho_k}' - \sigma_k) = 
 (\sigma_1, \ldots, \sigma_k) + \sum_{i=1}^w  {\bf d}_i \mmult \eta_{i}$ 
 where $(0, {\rho_2}' - \sigma_2, \ldots, {\rho_k}' - \sigma_k) = 
 \sum_{i=1}^w {\bf d}_i \mmult \eta_{i}$
 for  $\eta_{i} \in \rr$, i.e., $(\sigma_1,{\rho_2}', \ldots, {\rho_k}') 
 \equiv_{\langle {\bf d}_1, \ldots, {\bf d}_w\rangle} (\sigma_1, \ldots, \sigma_k)$.
Hence, as $\{ (\delta_2^i, \ldots, \delta^i_{k})\mid 1 \leq i \leq w\}$
 is a right Gr\"obner basis of ${\cal M}$ both vectors $(\sigma_1,{\rho_2}', \ldots, {\rho_k}')$ and
 $(\sigma_1, \ldots, \sigma_k)$ must have a common normal form using
 $\{ {\bf d}_i = (0,\delta_2^i, \ldots, \delta^i_{k})
 \mid  1 \leq i \leq w \}$ for reduction\footnote{The elements in this set cannot influence the
 first coordinate which is $\sigma_1$ for both vectors.} and we are done.
\\
The same argument applies to show local confluence. 
Let us assume there are ${\bf r}$, ${\bf s}_1$, ${\bf s}_2 \in \rr^k$
 such that ${\bf r} \myr_{G} {\bf s}_1$ and ${\bf r} \myr_{G} {\bf s}_2$.
Then by the definition of $G$ the first coordinates $\sigma^1_1$ and $\sigma_1^2$ 
 of ${\bf s}_1$ respectively ${\bf s}_2$ are joinable
 by $\{\gamma_1, \ldots, \gamma_s \}$ to some element say $\sigma$ giving rise to elements
 ${\bf r}_1 = {\bf s}_1 + \sum_{i=1}^s  {\bf c}_i \mmult \chi_{i}$ and
 ${\bf r}_2 = {\bf s}_2 + \sum_{i=1}^s {\bf c}_i \mmult \psi_{i}$ with first
 coordinate $\sigma$.
Again  we know $(\sigma,{\rho^1_2}, \ldots, {\rho^1_k}) = 
 (\sigma,{\rho^2_2}, \ldots, {\rho^2_k}) + (0,{\rho^1_2} - {\rho^2_2}, \ldots, {\rho^1_k} - {\rho^2_k})$
 with $({\rho^1_2} - {\rho^2_2}, \ldots, {\rho^1_k} - {\rho^2_k}) \in {\cal M}$.
Hence $(\sigma,{\rho^1_2}, \ldots, {\rho^1_k}) = 
 (\sigma,{\rho^2_2}, \ldots, {\rho^2_k})
  + \sum_{i=1}^w {\bf d}_i \mmult \eta_{i}$
 for $\eta_{i} \in \rr$, i.e., $(\sigma_1,{\rho_2}', \ldots, {\rho_k}') 
 \equiv_{\langle {\bf d}_1, \ldots, {\bf d}_w\rangle} (\sigma_1, \ldots, \sigma_k)$.
As again  $\{ (\delta_2^i, \ldots, \delta^i_{k})\mid 1 \leq i \leq w\}$
 is a right Gr\"obner basis of ${\cal M}$
 both vectors must have a common normal with respect to reduction using
 $\{ {\bf d}_i = (0,\delta_2^i, \ldots, \delta^i_{k})
 \mid  1 \leq i \leq w \}$.
\\
\qed

The task of describing two-sided syzygy modules is much more complicated.
We follow the ideas given by Apel in his habilitation  \cite{Ap98}.

Let ${\cal R}$ be the free Abelian group with basis elements $\alpha \otimes \beta$ where
  $\alpha,\beta \in \rr$.
We define a new vector space ${\cal S}$ with formal sums as elements 
 $\sum_{i=1}^n \gamma_i \skm \alpha_i \otimes \beta_i \skm \delta_i$  where $\gamma_i, \delta_i \in \rr$
 and $\alpha_i \otimes \beta_i \in {\cal R}$.
Let ${\cal U}$ be the subspace of ${\cal S}$ generated by the  vectors 
$$\alpha \otimes (\beta_1+\beta_2) - \alpha \otimes \beta_1 - \alpha \otimes \beta_2$$
$$(\alpha_1 + \alpha_2) \otimes \beta - \alpha_1 \otimes \beta - \alpha_2 \otimes \beta$$
$$\alpha \otimes (\gamma \skm \beta) - \gamma \skm (\alpha \otimes \beta)$$
$$(\gamma \skm \alpha) \otimes \beta - \gamma \skm (\alpha \otimes \beta)$$
$$\alpha \otimes (\beta \skm \gamma) -  (\alpha \otimes \beta) \skm \gamma$$
$$(\alpha \skm \gamma) \otimes \beta - (\alpha \otimes \beta) \skm \gamma$$
where $\alpha, \alpha_i, \beta, \beta_i, \gamma \in \rr$.
Then the quotient ${\cal S} / {\cal U}$ is called the \betonen{tensor product} denoted
 by $\rr \otimes \rr$.

The sets we are interested in can be defined as follows:
Let $R$ be some subset of $\rr$.
Syzygies of $R$ are solutions of the equations 
$\sum_{i=1}^n \sum_{j=1}^{n_i} \alpha_{i,j} \skm \rho_i \skm \beta_{i,j} = 0,
 \alpha_{i,j}, \beta_{i,j} \in \rr, \rho_i \in R$.
The set containing all such solutions is called the syzygy module of $R$.
We can now describe these sets  using objects of the ``polynomial'' structure 
 ${\cal S}[\rr]$ which contains formal sums of the form
 $\sum_{i=1}^n \sum_{j=1}^{n_i}(\alpha_{i,j} \otimes \beta_{i,j}) \skm \gamma_i$, $\alpha_i, \beta_i, \gamma_i \in \rr$.
We can associate a mapping $\phi: {\cal S}[\rr] \myr \rr$ by
 $\sum_{i=1}^n \sum_{j=1}^{n_i}(\alpha_{i,j} \otimes \beta_{i,j}) \skm \gamma_i \mapsto \sum_{i=1}^n \sum_{j=1}^{n_i}\alpha_{i,j} \skm \gamma_i \skm \beta_{i,j}$.
Then for the set $R$ we are interested in, the set of ``solutions'' is
 $\bigcup_{\rho_1, \ldots, \rho_k \in R, k \in \n} S_{\rho_1, \ldots, \rho_k}$
 with ordered lists of not necessarily different elements from $R$ such that 
 $S_{\rho_1, \ldots, \rho_k} = \{ (\sum_{j=1}^{n_1}\alpha_{1,j} \otimes \beta_{1,j}, \ldots, \sum_{j=1}^{n_k}\alpha_{k,j} \otimes \beta_{k,j}) 
 \mid \phi(\sum_{i=1}^k \sum_{j=1}^{n_i} (\alpha_{i,j} \otimes \beta_{i,j}) \skm \rho_i) = 0,
  \alpha_{i,j}, \beta_{i,j} \in \rr \}$.
Then these sets $S_{\rho_1, \ldots, \rho_k}$ are in fact modules
\begin{enumerate}
\item $S_{\rho_1, \ldots, \rho_k}$ is closed under scalar multiplication,
       i.e., $(\sum_{j=1}^{n_1}\alpha_{1,j} \otimes \beta_{1,j}, \ldots, \sum_{j=1}^{n_k}\alpha_{k,j} \otimes \beta_{k,j}) 
        \in S_{\rho_1, \ldots, \rho_k}$ and $\gamma \in \rr$ implies
        $\gamma \skm (\sum_{j=1}^{n_1}\alpha_{1,j} \otimes \beta_{1,j}, \ldots, \sum_{j=1}^{n_k}\alpha_{k,j} \otimes \beta_{k,j})  = 
        (\gamma \skm(\sum_{j=1}^{n_1}\alpha_{1,j} \otimes \beta_{1,j}), \ldots,
         \gamma \skm ( \sum_{j=1}^{n_k}\alpha_{k,j} \otimes \beta_{k,j})) 
        \in S_{\rho_1, \ldots, \rho_k}$:\\
      $\phi(\sum_{i=1}^k \sum_{j=1}^{n_i}(\alpha_{i,j} \otimes \beta_{i,j}) \skm \rho_i) = 0$ implies
      $\phi(\gamma \skm (\sum_{i=1}^k \sum_{j=1}^{n_i}(\alpha_{i,j} \otimes \beta_{i,j}) \skm \rho_i)) = 0$
      as $\gamma \skm (\alpha_{i,j} \otimes \beta_{i,j}) =
      (\gamma \skm \alpha_{i,j}) \otimes \beta_{i,j}$ and hence 
      $\phi(\gamma \skm (\sum_{i=1}^k \sum_{j=1}^{n_i}(\alpha_{i,j} \otimes \beta_{i,j}) \skm \rho_i))
      = \sum_{i=1}^k \sum_{j=1}^{n_i} \gamma \skm \alpha_{i,j} \skm \rho_i \skm \beta_{i,j} = 0$.
      Multiplication from the right can be treated similarly.
\item $S_{\rho_1, \ldots, \rho_k}$ is closed under addition, i.e., 
      $(\sum_{j=1}^{n_1}\alpha_{1,j} \otimes \beta_{1,j}, 
        \ldots, \sum_{j=1}^{n_k}\alpha_{k,j} \otimes \beta_{k,j})$, 
      $(\sum_{j=1}^{\tilde{n_1}}\tilde{\alpha}_{1,j} \otimes \tilde{\beta}_{1,j}, 
        \ldots, 
        \sum_{j=1}^{\tilde{n_k}}\tilde{\alpha}_{k,j} \otimes \tilde{\beta}_{k,j})
       \in S_{\rho_1, \ldots, \rho_k}$ implies
       $(\sum_{j=1}^{n_1}\alpha_{1,j} \otimes \beta_{1,j} + 
         \sum_{j=1}^{\tilde{n_1}}\tilde{\alpha}_{1,j} \otimes \tilde{\beta}_{1,j}, 
         \ldots, 
         \sum_{j=1}^{n_k}\alpha_{k,j} \otimes \beta_{k,j} +
         \sum_{j=1}^{\tilde{n_k}}\tilde{\alpha}_{k,j} \otimes \tilde{\beta}_{k,j})
       \in S_{\rho_1, \ldots, \rho_k}$:%\\
%      Using the defining equations for the quotient we get that
%      $\delta_1 \otimes \gamma_1 + \delta_2 \otimes \gamma_2 =  
%       - (\delta_1 + \delta_2) \otimes \gamma_1 + \delta_2 \otimes (\gamma_1 - \gamma_2)$ and
%      applying $\phi$ to the sum then gives us
%      $\sum_{j=1}^k (\alpha_{i} \skm \rho_i \skm \beta_{i} + \underbrace{\tilde{\alpha}_{i} \skm \rho_i \skm \beta_{i} - \tilde{\alpha}_{i} \skm \rho_i \skm \beta_{i}}_{=0} + \tilde{\alpha}_{i} \skm \rho_i \skm \tilde{\beta}_{i}) = 0$.
\end{enumerate}
The question arises when such modules have useful bases for characterizing syzygy modules in
 non-commutative reduction rings.
This would mean the existence of sets $B_{\rho_1, \ldots, \rho_k} = \{ B_i \in (\rr\otimes\rr)^k \mid i \in I \}$
 such that for each $(\sum_{j=1}^{n_1}\alpha_{1,j} \otimes \beta_{1,j}, 
        \ldots, \sum_{j=1}^{n_k}\alpha_{k,j} \otimes \beta_{k,j})
        \in S_{\rho_1, \ldots, \rho_k}$ there exist $\gamma_{ij}, \delta_{ij} \in \rr$
 with $(\sum_{j=1}^{n_1}\alpha_{1,j} \otimes \beta_{1,j}, 
        \ldots, \sum_{j=1}^{n_k}\alpha_{k,j} \otimes \beta_{k,j}) = 
 \sum_{i\in I} \sum_{j=1}^{n_i} \gamma_{ij} \skm B_i \skm \delta_{ij}$.
But even if this is possible it still remains the problem that we have to handle infinitely many
 sets of solutions associated to ordered subsets of a set admitting elements to occur more than once.
This problem arises from the fact that in contrary to one-sided syzygy modules or
 syzygy modules in commutative structures the summands in the representations
 cannot be ``collected'' and ``combined'' in such a way that for a set $R$ the 
 sums can be written as a $\sum_{\rho \in R} \alpha_{\rho} \skm \rho \skm \beta_{\rho}$. 

Let us close this section by illustrating the situation with two examples.
\begin{example}~\\
{\rm
Let $\Sigma = \{ a,b\}$ and $\Sigma^*$ the free monoid on the alphabet $\Sigma$.
Further let $\rr = \q [\Sigma^*]$ the monoid ring over $\Sigma^*$ and $\q$.
Let us look at the syzygy module of the set $\{a,b\} \subset \rr$, i.e.~the
 set of solutions of the equations $\sum_{j=1}^{n_1} \alpha_{1,j} \skm a \skm \beta_{1,j} 
 + \sum_{j=1}^{n_2} \alpha_{2,j} \skm b \skm \beta_{2,j} = 0,
 \alpha_{i,j}, \beta_{i,j} \in \rr$.
Then we find $\{ ( -1 \otimes b, a \otimes 1), (-b \otimes 1, 1 \otimes a)
 \} \subseteq S_{a,b}$
 and this set is a finite basis for $S_{a,b}$.
\exaend
}
\end{example}
\begin{example}~\\
{\rm
Let $\m$ be the monoid presented by $( \{ a,b, c\}; \{ ab=a, ac=a, bc=b \})$
 and $\rr = \q[\m]$ the monoid ring over $\m$ and $\q$.
Let us look at the syzygy module of the set $\{a,b\} \subset \rr$.
Then we find
 $\{ 
%(1 \otimes 1, a \otimes -1), (-1 \otimes 1, a \otimes 1),
% (1 \otimes b, a \otimes -1), (-1 \otimes b, a \otimes 1),
% (1 \otimes b^i, a \otimes b^j), (-1 \otimes b^i, a \otimes b^j),
 (1 \otimes 1, -a \otimes c^ib^j)\mid i,j \in \n \} \subseteq S_{a,b}$ and hence
 $S_{a,b}$ has no finite basis. 
\exaend
}
\end{example}
Hence the task of two-sided syzygies is much more complicated than the one-sided case.
This was also observed by Apel for graded structures where we have more structural information
 \cite{Ap98}.

%%%%%%%%%%%%%%%%%%%%%%%%%%%%%%%%%%
\section{Polynomial Rings over Reduction Rings}\label{section.polyrings}
For a  ring $\rr$ with a reduction relation $\R$  fulfilling (A1) -- (A3) we
 adopt the usual notations in $\rr[X]$ the polynomial ring in one variable $X$
 where multiplication is denoted by $\rmult$.
Notice that for scalar multiplication with $\alpha \in \rr$
 we assume $\alpha \skm X = X \skm \alpha$  (see \cite{Pe97} for other possibilities).
We specify an ordering on the set of terms in one variable by defining
 that if $X^i$ divides $X^j$, i.e.~$0 \leq i \leq j$, then  $X^i \preceq X^j$.
Using this ordering, the head term $\hterm(p)$, the head monomial $\hm(p)$ and
 the head coefficient $\hc(p)$ of
 a polynomial $p \in \rr[X]$ are defined as usual,
 and $\reductum(p) = p - \hm(p)$.
We extend the function $\hterm$  to sets of polynomials $F \subseteq\rr[X]$ by
 $\hterm(F)  = \{ \hterm(f) \mid f \in F \}$.

Let ${\mathfrak i} \subseteq \rr[X]$ be a finitely
 generated ideal in $\rr[X]$.
It is easy to see that given a term $t$ the set 
 $C(t,{\mathfrak i}) = \{ \hc(f) \mid f \in {\mathfrak i}, \hterm(f) = t \} \cup \{ 0 \}$
 is an ideal in $\rr$.
In order to guarantee that these ideals are also finitely generated we will
 assume that $\rr$ is a  Noetherian ring\footnote{We run into similar problems
 as in the module case in Section \ref{section.modules} as we cannot conclude that
 the ideal $C(t,{\mathfrak i})$ is finitely generated from the fact that ${\mathfrak i}$ is.}.
Note that for any two terms $t$ and $s$ such that $t$ divides $s$ we have
 $C(t,{\mathfrak i}) \subseteq C(s,{\mathfrak i})$.
This follows, as for $s = t \rmult u$, $u \in \{ X^i \mid i \in \n \}$, we find that $\hc(f) \in C(t,{\mathfrak i})$
 implies $\hc(f \rmult u) = \hc(f) \in C(s,{\mathfrak i})$ since $f \in {\mathfrak i}$ implies $f \rmult u \in {\mathfrak i}$.

We additionally define a partial ordering on $\rr$ by setting for 
 $\alpha, \beta \in \rr$, $\alpha >_{\rr} \beta$ if and only if
 there exists a finite set $B \subseteq \rr$ such that $\alpha
 \red{+}{\Longrightarrow}{}{B} \beta$.
Then we can define an ordering on
 $\rr[X]$ as follows:
For $f,g \in \rr[X]$, $f > g$ if and only if either
 $\hterm(f) \succ \hterm(g)$ or $(\hterm(f) = \hterm(g)$ and 
 $\hc(f) >_{\rr} \hc(g))$ or $(\hm(f) = \hm(g)$ and 
 $\reductum(f) > \reductum(g))$. 
Notice that this ordering in general is neither total nor Noetherian 
 on $\rr[X]$.
\begin{definition}\label{def.red}~\\
{\rm
Let $p, f$ be two non-zero polynomials in $\rr[X]$. 
We say $f$ \betonen{reduces} $p$ to $q$ at a monomial 
 $\alpha \skm X^i$ in one step, denoted by $p \red{}{\myr}{}{f} q$, if
\begin{enumerate}
\item[(a)] $\hterm(f)$ divides $X^i$, i.e.~$\hterm(f) \rmult X^j  = X^i$ for some term $X^j$,
\item[(b)] $\alpha \Longrightarrow_{\hc(f)} \beta$, with  $\alpha = \beta + \sum_{i=1}^k \gamma_i \skm \hc(f) \skm \delta_i$
 for some $\beta, \gamma_i, \delta_i \in \rr$, $1 \leq i \leq k$, and
\item[(c)] $q = p - \sum_{i=1}^k (\gamma_i \skm f \skm \delta_i) \rmult X^j$.
\dend
\end{enumerate}
}
\end{definition}
Notice that if $f$ reduces $p$ to $q$ at a monomial $\alpha
 \skm t$ the term $t$ can still occur in the resulting polynomial $q$.
Hence termination of this reduction cannot be shown by arguments involving
 terms only as in the case of polynomial rings over fields.
But when using a {\em finite} set of polynomials  for reduction we know by (A1)
 that reducing $\alpha$ in $\rr$ with respect to the finite set of head coefficients of
 the applicable polynomials must terminate and then either the monomial containing the term $t$
 disappears or is irreducible.  
Hence the reduction relation as defined in Definition \ref{def.red} is Noetherian when using  {\em finite} sets of
 polynomials. Therefore it fulfills Axiom (A1).
It is easy to see that (A2) and (A3) are also true and if the reduction relation $\R$
 satisfies (A4) this 
 is inherited by the reduction relation $\myr$ in $\rr[X]$.
\begin{theorem}\label{theo.poly}~\\
{\sl
If $(\rr, \R)$ is a Noetherian reduction ring, then $(\rr[X], \myr)$
 is a Noetherian reduction ring.
\theoend
}
\end{theorem}
\Ba{}~\\
By Hilbert's basis theorem $\rr[X]$ is Noetherian as
 $\rr$ is Noetherian.
We only have to prove that every ideal ${\mathfrak i} \neq \{ 0 \}$ 
 in $\rr[X]$ has a finite Gr\"obner basis.
\\
A finite basis $G$ of ${\mathfrak i}$ will be defined in stages
 according to the degree of terms occurring as head terms among
 the polynomials in ${\mathfrak i}$ and then we will show that
 $G$ is in fact a Gr\"obner basis.
\\
Let $G_0$ be a finite Gr\"obner basis of the ideal $C(X^0,{\mathfrak i})$
 in $\rr$, which must exist since $\rr$ is supposed to be Noetherian
 and a reduction ring.
Further, at stage $i>0$, if for each $X^j$ with $j < i$ 
 we have $C(X^j,{\mathfrak i}) \subsetneqq C(X^{i},{\mathfrak i})$, include 
 for each $\alpha$ in {\sc Gb}$(C(X^{i}, {\mathfrak i}))$ 
 (a finite Gr\"obner basis of $C(X^i, {\mathfrak i})$) a polynomial
 $p_{\alpha}$ from ${\mathfrak i}$ in $G_{i}$
 such that $\hm(p) = \alpha \skm X^i$.
Notice that in this construction we use the axiom of choice, when
 choosing the $p_{\alpha}$ from the infinite set ${\mathfrak i}$, and hence
 the construction is non-constructive.
At each stage only a finite number of polynomials can be added since  the respective
 Gr\"obner bases {\sc Gb}$(C(X^i, {\mathfrak i}))$
 are always finite, and at most one polynomial
 from ${\mathfrak i}$ is included for each element in
 {\sc Gb}$(C(X^i, {\mathfrak i}))$.
\\
If  a polynomial
 with head term $X^i$ is included, then
 $C(X^j,{\mathfrak i}) \subsetneqq C(X^i,{\mathfrak i})$ for every $j < i$.
So if $X^i \in HT({\mathfrak i})$ is not included as a head term of a
 polynomial in $G_i$, then there is a term $X^j$ occurring as a head term
 in some set $G_j$, $j<i$,
 $C(X^i,{\mathfrak i}) = C(X^j,{\mathfrak i})$ and $C(X^j,G_j)$ is a Gr\"obner basis
 for the ideal $C(X^j,{\mathfrak i})=C(X^i,{\mathfrak i})$ in $\rr$.
\\
We claim that the set $G = \bigcup_{i\geq 0} G_i$ is a finite Gr\"obner basis of
 ${\mathfrak i}$.
\\
To show that $G$ is finite it suffices to prove that the set $\hterm(G)$
 is finite, since in every stage only finitely many polynomials all
 having {\em new} head terms are added.
Assuming that $\hterm(G)$ is infinite, there is
 a sequence $X^{n_i}$, $i \in \n$ of different terms such that $n_i < n_{i+1}$.
But then by construction there is an ascending sequence of ideals in $\rr$,
 namely $C(X^{n_0},{\mathfrak i}) \subsetneqq C(X^{n_1},{\mathfrak i}) \subsetneqq \ldots$
 which contradicts the fact that
 $\rr$ is supposed to be Noetherian.
\\
So after some step $m$ no more polynomials $p$ from ${\mathfrak i}$ can
 be found such that for $\hterm(p) = X^i$ the set $C(X^i,{\mathfrak i})$ is
 different from all $C(X^j,{\mathfrak i})$, $j<i$.
\\
Notice that for all $p \in {\mathfrak i}$ we have
 $p \red{*}{\myr}{}{G} 0$ and $G$ generates ${\mathfrak i}$.
This follows immediately from the construction of $G$.
Hence $G$ is at least a wesk Gr\"obner basis.
\\
To see that $\red{}{\myr}{}{G}$ is confluent, let $p$ be a
 polynomial which has two distinct normal forms with respect to $G$,
 say $p_1$ and $p_2$.
Let $t$ be the largest term on which $p_1$ and $p_2$ differ and let
 $\alpha_1$ and $\alpha_2$ be the respective coefficients of $t$
 in $p_1$ and $p_2$.
Since $p_1 - p_2 \in {\mathfrak i}$ this polynomial
 reduces to $0$ using $G$ and
 without loss of generality we can assume that these reductions
 always take place at the respective head terms of the polynomials in the
 reduction sequence.
Let $s \in \hterm(G)$ be the head term of the polynomial in $G$
 which reduces  $\hterm(p_1 - p_2)$, i.e., $s$ divides $t$,
 $\alpha_1 - \alpha_2 \in C(s,{\mathfrak i})$, and hence
 $\alpha_1 \equiv_{{\mathfrak i}} \alpha_2$.
Therefore, not both $\alpha_1$ and $\alpha_2$ can be in normal form with 
 respect to any Gr\"obner basis of $C(s,{\mathfrak i})$ and hence with respect to 
 the set
 of head coefficients of polynomials in $G$ with head term $s$.
So both, $\alpha_1 \skm t$ and $\alpha_2 \skm t$ cannot be in
 normal form with respect to $G$, which is a contradiction to the fact 
 that $p_1$ and $p_2$ are supposed to be in normal form with respect to $G$.
\\
Finally we have to prove  $\equiv_{{\mathfrak i}}\; = \red{*}{\lr}{}{G}$.
Let $p \equiv_{{\mathfrak i}} q$ both be in normal form with respect to
 $G$.
Then as before $p-q \red{*}{\myr}{}{G} 0$ implies $p=q$.
Hence we have shown that $G$ is in fact a finite Gr\"obner basis of ${\mathfrak i}$. 
\\ \qed
This theorem of course can be applied to $\rr[X]$ and a new variable $X_2$ and
 by iteration we immediately get the following:
\begin{corollary}~\\
{\sl
If $(\rr, \R)$ is a Noetherian reduction ring, then $\rr[X_1, \ldots, X_n]$
 is a Noetherian reduction ring with the respective extended
 reduction relation.
\corend
}
\end{corollary}
Notice that other definitions of reduction relations in $\rr[X_1, \ldots, X_n]$ are known
 in the literature.
These are usually based on divisibility of terms and admissible term orderings
 on the set of terms to distinguish the head terms.
The proof of Theorem \ref{theo.poly} can be generalized for these cases.

Moreover, these results also hold for weak reduction rings.
\begin{corollary}~\\
{\sl
If $(\rr, \R)$ is a Noetherian weak reduction ring, then $\rr[X_1, \ldots, X_n]$
 is a Noetherian weak reduction ring with the respective extended reduction relation.
\corend
}
\end{corollary}
\Ba{}~\\
This follows immediately by using weak Gr\"obner bases $G_i$ for the definition of
 $G$ in the proof of Theorem \ref{theo.poly}. 
As before the property
 that for all $p \in {\mathfrak i}$ we have
 $p \red{*}{\myr}{}{G} 0$ and $G$ generates ${\mathfrak i}$ follows immediately
 from the construction of $G$.
Hence the result holds for $\rr[X_1]$ and can be extended to $\rr[X_1, \ldots, X_n]$. 
\\
\qed
Now if $(\rr,\R)$ is an effective reduction ring, 
 then  addition and  multiplication in $\rr[X]$ as well as reduction as defined
 in Definition \ref{def.red} are computable operations.
However, the proof of Theorem \ref{theo.poly} does not specify
 how Gr\"obner bases for finitely generated ideals in $\rr[X]$ can be constructed using
 Gr\"obner basis methods for $\rr$.
So we cannot conclude that for effective reduction rings the polynomial ring again will
 be effective.
A more suitable characterization of Gr\"obner bases requiring $\rr$ to fulfill additional conditions 
 is needed.

In order to provide completion procedures to compute Gr\"obner bases,
 various characterizations of Gr\"obner bases by finite test sets of special
 polynomials in certain commutative
 reduction rings (e.g.~the integers and Euclidean domains)
 can be found in the literature (see e.g.~\cite{KaNa85,KaKa84,Mo88}).
A general approach to characterize commutative reduction rings allowing
 the computation of Gr\"obner bases  using Buchberger's approach was presented
 by Stifter in \cite{St87}.

Let us close this section by providing similar characterizations for polynomial rings
 over non-commutative reduction rings and outlining the arising problems.
For simplicity we restrict ourselves to the case of $\rr[X]$ but this is no general restriction.
Given a generating set $F \subseteq \rr[X]$ the key idea is  to distinguish special
 elements of $\ideal{}{}(F)$ which have representations $\sum_{i=1}^n g_i \rmult f_i \rmult h_i$,
 $g_i,h_i \in \rr[X]$, $f_i \in F$ such that the head terms $\hterm(g_i \rmult f_i \rmult h_i)$
 are all the same within the representation.
Then on one hand the respective coefficients $\hc(g_i \rmult f_i \rmult h_i)$ can add up to zero which
 in the commutative case 
 means that the sum of the head coefficients is in an appropriate module generated by the coefficients  $\hc(f_i)$ ---
 m(odule)-polynomials
 are related to these situations.
If the result is not zero the sum of the  coefficients $\hc(g_i \rmult f_i \rmult h_i)$ as in the commutative case
 can be described in terms
 of a Gr\"obner basis of the coefficients  $\hc(f_i)$ --- g(r\"obner)-polynomials are related to these situations.
Zero divisors in the reduction ring occur as a special instance of m-polynomials
 where $F = \{ f \}$ and $\alpha \rmult f \rmult \beta$, $\alpha, \beta \in \rr$ are considered.

In case $\rr$ is a commutative or one-sided reduction ring the first problem is related
 to solving linear homogeneous equations in $\rr$ and to the existence of
 finite bases of the respective modules.

Let us become more precise and look into the definitions of m- and g-polynomials
 for the special case of rings with {\em right} reduction relations.
\begin{definition}\label{def.rr.one-sided.gpol}~\\
{\rm
Let $P = \{ p_1, \ldots, p_k \}$ be a finite set of  polynomials in $\rr[X]$,
 $u_1, \ldots, u_k$ terms in $\{ X^j \mid j \in \n \}$ such that for the term
 $t = \max \{ \hterm(p_i) \mid 1 \leq i \leq k \}$ we have $t  = \hterm(p_i) \rmult u_i$
 and $\gamma_i = \hc(p_i)$ for  $1 \leq i \leq k$.
\\
Let $G$ be a right
  Gr\"obner basis of the right ideal generated by $\{ \gamma_i \mid 1 \leq i \leq k \}$ in $\rr$ and
$$\alpha = \sum_{i=1}^k \gamma_i \skm \beta_{i}^{\alpha}$$
 for $\alpha \in G$, $\beta_{i}^{\alpha} \in \rr$.
Then we define the  \betonen{g-polynomials (Gr\"obner polynomials)}
 corresponding to $P$ and $t$ by setting
$$ g_{\alpha} =  \sum_{i=1}^k  p_i \rmult  u_i \skm  \beta_{i}^{\alpha}$$
 where $\hterm(p_i) \rmult u_i = t$.
Notice that $\hm(g_{\alpha})= \alpha \skm t$.
\\
For the right module $M = \{ (\delta_1, \ldots, \delta_k) \mid  
 \sum_{i=1}^k  \gamma_i \skm \delta_i = 0 \}$, let the set
 $\{B_j \mid j \in I_M \}$ be a basis with
 $B_j = (\beta_{j, 1}, \ldots, \beta_{j, k})$ for $\beta_{j,l} \in \rr$
  and $1 \leq l \leq k$.
We define the 
 \betonen{m-polynomials (module polynomials)}
 corresponding to $P$ and $t$ by setting
$$ h_j = \sum_{i=1}^k p_i  \rmult u_i \skm  \beta_{j, i}
 \mbox{ for each } j \in I_M$$
 where $\hterm(p_i) \rmult u_i = t$.
Notice that $\hterm(h_j) \prec t$ for each $j \in I_M$.
\dend
}
\end{definition} 

Given a set of polynomials $F$ the corresponding m- and g-polynomials are those
 resulting for every subset $P \subseteq F$ according to this definition.

In case we want effectiveness, we have to require that the bases in this definition are  computable.
Of course for commutative reduction rings the definition extends to characterize two-sided ideals.
However, the whole situation
 becomes more complicated for non-commutative two-sided reduction rings, as the
 equations are no longer linear and we have to distinguish right and left
 multipliers  simultaneously.
Moreover the set of m-polynomials is a much more complicated structure.
In some cases the problem for two-sided ideals can be translated into the
 one-sided case and hence solved via one-sided reduction techniques \cite{KaWe90}. 
But the general case is much more involved, see Definition
 \ref{def.rr.two-sided.gpol} below.

The g-polynomials corresponding to right Gr\"obner bases of
 right ideals in $\rr$
 can successfully be treated whenever finite right Gr\"obner bases exist.
Here, if we want effectiveness, we have to require that a right Gr\"obner basis as well as
  representations for its elements in terms of the generating set are computable.

Using m- and g-polynomials, right Gr\"obner bases can be characterized
 similar to the characterizations in terms of syzygies (a direct generalization of 
 the approaches by Kapur and Narendran in \cite{KaNa85} respectively M\"oller in \cite{Mo88}):
In case for the respective subsets $P \subseteq F$ the respective terms $t
 = \max \{ \hterm(p) \mid p \in P \}$ only give rise
 to finitely many m- and g-polynomials, these situations can be
 localized to finitely many terms.
One can provide a completion procedure based on this characterization which will
 indeed compute a finite right Gr\"obner basis if $\rr$ is Noetherian.
In principal ideal rings, where the function ${\sf gcd}$ (greatest common divisor)
 is defined it is sufficient to consider
 subsets $P \subseteq F$ of size $2$ (compare  \cite{KaNa85}).

Now let us look at two-sided ideals and two-sided reduction relations.

\begin{definition}\label{def.rr.two-sided.gpol}~\\
{\rm
Let $P = \{ p_1, \ldots, p_k \}$ be a finite set of  polynomials in $\rr[X]$,
 $u_1, \ldots, u_k$ terms in $\{ X^j \mid j \in \n \}$ such that for the term
 $t = \max \{ \hterm(p_i) \mid 1 \leq i \leq k \}$ we have $t  = \hterm(p_i) \rmult u_i$
 and $\gamma_i = \hc(p_i)$ for  $1 \leq i \leq k$.
\\
Let $G$ be a 
  Gr\"obner basis of the ideal generated by $\{ \gamma_i \mid 1 \leq i \leq k \}$ in $\rr$ and
$$\alpha = \sum_{i=1}^k \sum_{j = 1}^{n_i} \beta_{i,j}^{\alpha} \skm \gamma_i \skm \delta_{i,j}^{\alpha}$$
 for $\alpha \in G$, $\beta_{i,j}^{\alpha}, \delta_{i,j}^{\alpha} \in \rr$,  $1 \leq i \leq k$,$1 \leq j \leq n_i$.
Then we define the  \betonen{g-polynomials (Gr\"obner polynomials)}
 corresponding to $P$ and $t$ by setting
$$ g_{\alpha} =  \sum_{i=1}^k \sum_{j = 1}^{n_i}  \beta_{i,j}^{\alpha} \skm p_i\rmult  u_i \skm  \delta_{i,j}^{\alpha} $$
where $\hterm(p_i) \rmult u_i = t$.
Notice that $\hm(g_{\alpha})= \alpha \skm t$.
\\
We define the 
 \betonen{m-polynomials (module polynomials)}
 corresponding to $P$ and $t$ as 
$$h =  \sum_{i=1}^k \sum_{j = 1}^{n_i}  \beta_{i,j} \skm p_i \rmult  u_i \skm  \delta_{i,j}$$
where $\sum_{i=1}^k \sum_{j = 1}^{n_i} \beta_{i,j} \skm \gamma_i \skm \delta_{i,j}=0$.
Notice that $\hterm(h) \prec t$.
\dend
}
\end{definition}
Given a set of polynomials $F$, the set of corresponding
 g- and m-polynomials contains those which are specified by
 Definition \ref{def.rr.two-sided.gpol} for each subset $P \subseteq F$
 fulfilling the respective conditions.
For a set consisting of one polynomial the corresponding
m-polynomials also reflect the multiplication of the polynomial with
zero-divisors of the head coefficient, i.e., by a basis of the annihilator
 of the head coefficient.
Notice that given a finite set of polynomials the corresponding sets of
 g- and m-polynomials in general can be infinite.

We can use g- and m-polynomials to characterize finite weak Gr\"obner bases.
Notice that this characterization does not require
 $\rr$ to be Noetherian.
In order to characterize Gr\"obner bases in this fashion the Translation Lemma must hold for
 the reduction ring. 

\begin{theorem}\label{theo.rr.cp}~\\
{\sl
Let $F$ be a finite set of polynomials in $\rr[X]\backslash \{ 0 \}$.
Then $F$ is a  weak Gr\"obner basis of the ideal it generates if and only if
all g-polynomials and all m-polynomials corresponding to $F$ as specified in
       Definition \ref{def.rr.two-sided.gpol} reduce to zero.
\theoend
}
\end{theorem}
\Ba{}~\\
First let $F$ be a weak Gr\"obner basis.
By Definition \ref{def.rr.two-sided.gpol} the g- and m-polynomials are elements of the
 ideal generated by $F$ and hence reduce to zero using $F$.
\\
It remains to show that every 
 $g \in \ideal{}{}(F) \backslash \{ 0 \}$ reduces to zero by $F$.
Remember that for
 $g \in \ideal{}{}(F)$, $ g \red{}{\myr}{}{F} g'$ implies $g' \in \ideal{}{}(F)$.
As  $\red{}{\myr}{}{F}$ is Noetherian\footnote{To achieve this we have demanded
 that $F$ is finite.}, thus
 it suffices to show that every  $g \in \ideal{}{}(F) \backslash \{ 0 \}$ 
 is $\red{}{\myr}{}{F}$-reducible.
%\\
Let $g = \sum_{i=1}^m \alpha_i \skm f_{i} \rmult u_{i} \skm \beta_i$ be an arbitrary
 representation of $g$ with $\alpha_i, \beta_i \in \rr$, $u_i \in \{ X^j \mid j \in \n \}$, and  $f_i \in F$
 (not necessarily different polynomials).
%\\
Depending on this representation of $g$ and the degree  ordering 
 $\succeq$ on $\{ X^j \mid j \in \n \}$ we define the maximal occurring term of 
 this representation of $g$ to be
 $t = \max \{ \hterm(f_{i} \rmult u_{i}) \mid 1 \leq i \leq m \}$ and
 $K$ is the number of polynomials $f_i \rmult u_i$ containing $t$ as a term.
%\\
Then $t \succeq \hterm(g)$. 
We will show that $G$ is reducible by induction
 on $(t,K)$, where
 $(t',K')<(t,K)$ if and only if $t' \prec t$ or $(t'=t$ and
 $K'<K)$\footnote{Note that this ordering is well-founded since $\succ$
                  is well-founded on $\{ X^j \mid j \in \n \}$ and $K \in\n$.}.
Without loss of generality let the first $K$ multiples occurring in our representation
 of $g$ be those with head term $t$, i.e., for $\sum_{i=1}^K \alpha_i \skm f_{i} \rmult u_{i} \skm \beta_i$
 we have $\hterm(f_{i} \rmult u_{i}) = t$ for $1 \leq i \leq K$,
 and $\hterm( \alpha_i \skm f_{i} \rmult u_{i} \skm \beta_i) \pred t$ for $K<i\leq m$.
In case $t \succ \hterm(g)$ there is an m-polynomial corresponding to the set of polynomials
 $P = \{ f_1, \ldots, f_K \}$ and by our assumption this polynomial is reducible
 to zero using $F$ hence yielding the  existence of a representation 
 $\sum_{i=1}^n \gamma_i \skm f_i \rmult v_i \skm \delta_i$ with 
 $t \succ \tilde{t} = \max \{\hterm(f_{i} \rmult v_{i}) \mid i \in \{ 1, \ldots n \}  \}$.
%\\
We can then change the original representation of $g$ by substituting this sum for
 $\sum_{i=1}^K \alpha_i \skm f_{i} \rmult u_{i} \skm \beta_i$ yielding a new representation
 with smaller maximal term than $t$.
\\
On the other hand, if $t = \hterm(g)$ then again we can assume that the first $K$
 multiples have head term $t$.
In this case there exists a g-polynomial corresponding to the set  of polynomials
 $P = \{ f_1, \ldots, f_K \}$ and by our assumption this polynomial is reducible
 to zero using $F$.
Now as the head monomial of the g-polynomial and the head monomial of $g$ are equal, then
 $g$ must be reducible by $F$ as well.
\\
\qed

In order to characterize infinite sets $F$ as weak
 Gr\"obner bases we have to be more careful since we can no longer assume that $\red{}{\myr}{}{F}$ is terminating\footnote{This can of course be achieved by requiring
 the stronger axiom (A1') to hold for the reduction relation.}.
But inspecting the proof of the previous theorem closely we see that this is not necessary.
Under the stronger assumption that the g-polynomial reduces to zero using reduction at head monomials
 only, i.e., we have a terminating reduction sequence using finitely many polynomials in $F$ only, we can
 conclude that the polynomials used to extinguish the term $t$ in the g-polynomial can equally be applied to
 extinguish the head monomial of $g$.
Since there cannot be an infinite sequence of decreasing terms $t$ one can show that $g$ reduces to zero
 by iterating arguments involving g- and m-polynomials.

\begin{corollary}~\\
{\sl
Let $F$ be a  set of polynomials in $\rr[X]\backslash \{ 0 \}$.
Then $F$ is a  weak Gr\"obner basis of the ideal it generates if and only if
 all g-polynomials and all m-polynomials corresponding to $F$ as specified in
 Definition \ref{def.rr.two-sided.gpol} reduce to zero using reduction
 at head monomials only.
\corend
}
\end{corollary}
\begin{corollary}~\\
{\sl
Let $F$ be a  set of polynomials in $\rr[X]\backslash \{ 0 \}$.
Additionally let the Translation Lemma hold in $\rr$.
Then $F$ is a  Gr\"obner basis of the ideal it generates  if and only if
 all g-polynomials and all m-polynomials corresponding to $F$ as specified in
 Definition \ref{def.rr.two-sided.gpol} reduce to zero using reduction
 at head monomials only.
\corend
}
\end{corollary}

Still the problem remains that the set of m-polynomials does not have a nice characterization 
 as an algebraic structure.
Remember that in the one-sided case or the case of commutative reduction rings the m-polynomials for 
 a finite set of polynomials $P$ correspond to submodules of $\rr^{|P|}$, as they correspond to
 solutions of linear equations.
When attempting to describe the setting for two-sided ideals in non-commutative reduction rings 
 one runs into the same problems as in the previous section on modules.

%% file: generalization_k.tex
In the literature Gr\"obner bases and  reduction relations have been
 introduced to various algebraic structures such as the classical
 commutative polynomial rings over fields, non-commutative polynomial
 rings over fields, commutative polynomial rings over reduction rings,
 skew polynomial rings, Lie algebras, monoid and group rings and many more. 
This chapter is intended to give a generalized setting subsuming these
 approaches and outlining a framework for introducing reduction relations and
 Gr\"obner bases to other structures fitting the appropriate requirements.
An additional aim was to work out what conditions are necessary at what
 point in order to give more insight into the  ideas  behind
 algebraic characterizations such as specialized standard representations
 for ideal elements as well as into the idea of using rewriting techniques 
 for achieving confluent reduction relations describing the ideal congruence.

This chapter is organized as follows:
Section \ref{section.general} introduces the general structure we are looking into called
 function rings.
Section \ref{section.right.standard} gives the algebraic characterization for the case of right
 ideals in form of right standard representations.
To work out the difficulties involved by our notion of terms and coefficients
 separately, Section \ref{section.right.field} first treats the easier case of function rings over
 fields while Section \ref{section.right.rr} then goes into the details when taking  a reduction
 ring as introduced in Chapter \ref{chapter.reduction.rings} as coefficient domain.
Since for function rings over general reduction rings only a feasible
 characterization of weak Gr\"obner bases is possible, we show that this
 situation can be improved when looking at the special case of function
 rings over the integers in Section \ref{section.right.integers}.
Section \ref{section.right.module} is dedicated to the study of a generalization of the concept of
 right ideals -- right modules.
The remaining Sections \ref{section.ideal.standard} -- \ref{section.two-sided.module} 
 then treat the same concepts and problems now in the
 more complex setting of two-sided ideals. 
%%%%%%%%%%%%%%%%%%%%%%%%%%%%%%%%%%%%%%%%%%%%%%%%%%%%%%%%%%%%%%%%%%%%%%%%%%%
\section{The General Setting}\label{section.general}
Let $\myt$ be a set
%\footnote{We require this set to be linearly independent,
% i.e.~the zero function introduced cannot be represented by a 
% non-trivial sum of of the form $\sum_{i=1}^n t_i \radd \sum_{j=1}^n -t_j$, in
% order to have unique representations of elements of the function ring as formal 
% sums of monomials.} 
 and
 let $\rr = (\rr,+,\skm,0,1)$ be an associative ring with $1$.
By $\f^{\myt}_{\rr}$ we will denote the set of all functions
 $f : \myt \myr \rr$
 with finite \betonen{support} $\supp(f) = \{ t \mid t \in \myt, f(t) \neq 0 \}$.
We will simply write $\f$ if the context is clear.
By $\zero$ we will denote the function with empty support, i.e.,
 $\supp(\zero) = \emptyset$.
This function will be called the \betonen{zero function}.
Two elements of $\f$ are equal if they are equal as functions, i.e.,
 they have the same support
 and coincide in their respective values.
We require the set $\myt$ to be independent in the sense that a
 function $f$ has unique support.

$\f$ can be viewed as a group with respect to a binary
 operation
 $$\radd : \f \times \f \myr \f$$
 called \betonen{addition}
 by associating to $f,g$ in $\f$  the function in $\f$, denoted by $f \radd g$,
 which has support $\supp(f \radd g) \subseteq \supp(f) \cup \supp(g)$ and
 values $(f \radd g)(t) = f(t) + g(t)$ for $t \in \supp(f) \cup \supp(g)$.
The zero function $\zero$ fulfills $\zero \radd f = f \radd \zero = f$,
 hence is neutral with respect to $\radd$.
For an element $f \in \f$ we define the element $- f$ with
 $\supp(-f) = \supp(f)$ and for all $t \in \supp(f)$ the value of
 $(-f)(t)$ is the inverse of the element $f(t)$ with respect to $+$
 in $\rr$ denoted by $-f(t)$.
Notice that since in $\rr$ every element has such an inverse the inverse
 of an element in $\f \backslash \{ \zero \}$ is always defined.
Then $-f$ is the (left and right) inverse of $f$, since $f \radd (- f)$ as well
 as  $(- f) \radd f$ equals $\zero$, i.e., has empty support.
This follows as for all $t \in \supp(f)$ we have
 $(f \radd (- f))(t) = f(t) + (- f)(t) =f(t) - f(t) = 0 = -f(t) + f(t) = 
  (-f)(t) + f(t) =  ((- f) \radd f)(t)$.
We will write  $f-g$ to
 abbreviate $f \radd (- g)$ for $f,g$ in $\f$.
If the context is clear we will also write $f + g$ instead of $f \radd g$.
Notice that $(\f, \radd, \zero)$ is an Abelian group
 since $(\rr,+,0)$ is  Abelian.
Sums of functions $f_1, \ldots, f_m$ will be abbreviated by
 $f_1 \radd \ldots \radd f_m = \sum_{i=1}^m f_i$ as usual.
Now if $\rr$ is a computable ring\footnote{A ring $\rr$ is called computable, 
 if the ring operations $+$ and $\skm$ are computable, i.e.~for
 $\alpha, \beta \in \rr$ we can compute $\alpha + \beta$ and $\alpha \skm \beta$.},
 then $(\f, \radd)$ is a computable group. 
 
In the next lemma we provide a syntactical representation for elements
 of the function ring.
\begin{lemma}\label{lem.monomrep}~\\
{\sl
Every $f \in \f \backslash \{ \zero \}$ has a
 finite representation of the form
$$f = \sum_{t \in \supp(f)} m_t$$
where $m_t \in \f$ such that $\supp(m_t) = \{ t \}$
 and 
$f(t) = m_t(t)$.
The representation of $\zero$ is the empty sum.
\lemend
}
\end{lemma}
\Ba{}~\\
This can be shown by induction on $n = |\supp(f)|$.
For $n = 0$ we have the empty sum which is the zero function $\zero$ and are done.
Hence let $\supp(f) = \{ t_1, \ldots, t_n \}$ and $n > 0$.
Furthermore let $f(t_1) = \alpha \in \rr$ and $m \in \f$ be the unique
 function with  $\supp(m) = \{ t_1 \}$ and $m(t_1) = \alpha$.
Then there exists an inverse function $- m$ and a function
 $(-m) \radd f \in \f$ such that
 $$ f = (m \radd (- m)) \radd f = m \radd ((- m) \radd f)$$
 and $\supp((- m) \radd f) = \{ t_2, \ldots t_n \}$.
Hence by our induction hypothesis $\supp((- m) \radd f)$ has a
 representation $\sum_{t \in \{ t_2, \ldots t_n \}} m_t$
 yielding
 $$ f = m \radd ((- m) \radd f) 
      = m \radd \sum_{t \in \{ t_2, \ldots t_n \}} m_t
      = \sum_{t \in \supp(f)} m_t$$
 with $m_{t_1} = m$.
\\
\qed
This presentation is unique up to permutations.
We will call such a representation of an element
 as a formal sum of special functions
 a \betonen{polynomial representation} or a \betonen{polynomial}
 to stress the similarity with the objects known as polynomials
 in other fields of mathematics.
Polynomial representations in terms of these functions are unique up to
 permutations of the respective elements of their support.
Since these special functions are of interest we define the
 following subsets of $\f$:
 $$\monoms(\f) = \{ f \in \f \mid |\supp(f)| = 1 \}$$
will be called the set of \betonen{monomial functions} or 
 \betonen{monomials} in $\f$.
Monomials will often be denoted by $m_t$ where the suffix $t$ is
 the element of the support, i.e., $\supp(m_t) = \{ t \}$.
A subset of
 this set, namely
 $$\terms(\f) = \{ m_t \in \monoms(\f) \mid m_t(t) = 1\}$$ where $1$ denotes the
 unit in $\rr$
 will be called the set of \betonen{term functions} or \betonen{terms} of $\f$.
Notice that this set can be viewed as an embedding of $\myt$ in $\f$ via
 the mapping $t \longmapsto f$ with $\supp(f)= \{ t \}$ and $f(t) = 1$.

Further we assume the existence of a second binary operation called
 \betonen{multiplication}
 $$\rmult :  \f \times \f \myr \f$$
 such that $(\f,\radd,\rmult,\zero)$ is a ring.
In particular we have
 $\zero \rmult f = f \rmult \zero = \zero$ for all $f$ in $\f$.
This ring is called a \betonen{function ring}\footnote{Notice that in the
 literature the term function ring is usually restricted to those rings
 where the multiplication is defined pointwise as in Example \ref{exa.punktweise}.
 Here we want to allow more interpretations for $\rmult$.}.
In case $\rmult$ is a computable operation, $\f$ is a computable function ring.

\begin{definition}~\\
{\rm
An element $\one^r_{\f} \in \f$ is called a \betonen{right unit} of $\f$ if for all $f\in \f$ we have
 $f \rmult \one^r_{\f} = f$. 
Similarly $\one^{\ell}_{\f} \in \f$ is called a \betonen{left unit} of $\f$ if for all $f\in \f$ we have
 $ \one^{\ell}_{\f} \rmult f = f$.
An element $\one_{\f} \in \f$ is called a \betonen{unit} if for all $f \in \f$ we have $\one_{\f} \rmult f =
 f \rmult \one_{\f} = f$.
\dend
}
\end{definition}
In general $\f$ need not have a left or right unit.
If $\f$ does not have a unit this can be achieved by enlarging
 the set $\myt$ by a new element, say $\Lambda$, and associating to $\Lambda$
 a function $f_{\Lambda}$ with support $\{ \Lambda \}$ and $f_{\Lambda}(\Lambda) = 1$.
The definition of $\rmult$ must be extended such that
 for all $f \in \f$ we have $f \rmult f_{\Lambda} = f_{\Lambda} \rmult f = f$.
Similarly we could add a left or right unit by requiring $f \rmult f_{\Lambda}^r = f$
 respectively $f_{\Lambda}^{\ell} \rmult f = f$.
When adding a new element $f_{\Lambda}$ as a unit
 to $\f$ we have
 $f_{\Lambda} \in \terms(\f) \subseteq \monoms(\f)$.

We will not specify our ring multiplication $\rmult$ further at the moment
 except for giving some examples.

Our first example outlines the situation for multiplying two elements
 by multiplying the respective values of the support.
This is the definition of multiplication normally associated to function rings in
 the mathematical literature.
\begin{example}\label{exa.punktweise}~\\
{\rm
Let us specify our multiplication $\rmult$ by associating to $f,g$ in $\f$ 
 the function in $\f$, denoted by $f \rmult g$,
 which has support $\supp(f \rmult g) \subseteq \supp(f) \cap \supp(g)$ and 
 values $(f \rmult g)(t) := f(t) \skm g(t)$ for $t \in \supp(f) \cap \supp(g)$.
Notice that in this case $\f$ can only contain a (right, left) unit if $\myt$ is finite,
 since otherwise a unit function would have infinite support and hence be no
 element of $\f$.
But the set of special functions $u_S = \sum_{t \in S} u_t$ where $S \subseteq \myt$ 
 finite, $\supp(u_t) = \{ t \}$ and $u_t(t) =1$ is an approximation of a unit, since
 for every function $f$ in $\f$ and all functions $u_S$ with $\supp(f) \subseteq S$
 we have $f \rmult u_S = u_S \rmult f = f$.
However, if we want a real unit, adding a new symbol $\Lambda$ to
 $\myt$ and $f_{\Lambda}$ with $f_{\Lambda}(\Lambda) = 1$ to $\f$ together with an extension of the definition
 of $\rmult$ by $f_{\Lambda} \rmult f = f \rmult f_{\Lambda} = f$ for all
 $f \in \f$ will do the trick.
\exaend
}
\end{example}
Remember that by Lemma \ref{lem.monomrep} polynomials have representations of
 the form $f = \sum_{t \in \supp(f)} m_t$ and $g = \sum_{s \in \supp(g)} n_s$ 
 yielding
 $$f \rmult g  = (\sum_{t \in \supp(f)} m_t) \rmult (\sum_{s \in \supp(g)} n_s)
 = \sum_{t \in \supp(f), s \in \supp(g)} m_t \rmult n_s$$
 since the multiplication $\rmult$ must satisfy the
 distributivity law of the ring axioms.
Hence knowing the behaviour of the multiplication for monomials, 
 i.e.~$\rmult : \monoms(\f) \times \monoms(\f) \myr \f$, is enough to
 characterize the multiplication $\rmult$.

For all examples from the literature mentioned in this work, we can even
 state that the multiplication can be defined by specifying
 $\rmult : \myt \times \myt \myr \f$, and then lifting it to
 $\monoms(\f)$ and $\f$.
This is done by defining $m_t \rmult n_s = (m_t(t) \skm n_s(s)) \skm (t \rmult s)$
 and extending this to the formal sums of monomials\footnote{Notice that this lifting requires that when writing a monomial $m_t$ as $m_t(t) \skm t$ we have $m_t(t) \skm t = t \skm m_t(t)$.}.

A well-known example for the special instance
 $\rmult : \myt \times \myt \myr \myt$ are the polynomial rings from Section
 \ref{section.buchberger}.
\begin{example}\label{exa.polyring}~\\
{\rm
For a set of variables $X_1, \ldots, X_n$ let us define the set of commutative terms
 $\myt = \{ X_1^{i_1} \ldots X_n^{i_n} \mid i_1, \ldots i_n \in \n \}$
 and let $\f^{\myt}_{\q}$ be the set of all functions $f : \myt \myr \q$ with finite support,
 where $\q$ are the rational numbers.
Multiplication $\rmult : \myt \times \myt \myr \myt$ is specified as
 $X_1^{i_1} \ldots X_n^{i_n} \rmult X_1^{j_1} \ldots X_n^{j_n} = X_1^{i_1 + j_1} \ldots X_n^{i_n + j_n}$.
Hence here we have an example where the set $\myt$ is a monoid
 with unit element $X_1^{0} \ldots X_n^{0}$.
Then $\f$ can be interpreted as the ordinary polynomial ring $\q[X_1, \ldots, X_n]$
  with the usual multiplication $(\alpha \skm t) \rmult (\beta \skm s) = (\alpha \skm \beta) \skm (t \rmult s)$ where $\alpha , \beta \in \q, s,t \in \myt$.
\exaend
}
\end{example}
Notice that in this example the unit element is an element of the set $\myt$ embedded in $\f$.
This does not have to be the case as the next example shows.
\begin{example}\label{exa.matrix}~\\
{\rm
Let us fix a finite set $\myt = \{ e_{11}, e_{12}, e_{21}, e_{22} \}$ and
 let $\f^{\myt}_{\q}$ be the set of all functions $f : \myt \myr \q$,
 where $\q$ are the rational numbers.
We specify the multiplication $\rmult$ on $\f^{\myt}_{\q}$ by the
 action on $\myt$ as follows: 
 $e_{ij} \rmult e_{kl} = \zero$ in case $j \neq k$ and
 $e_{ij} \rmult e_{jl} = e_{il}$ for $i,j,l,k \in \{ 1,2 \}$.
Then multiplication is not Abelian since
 $e_{11} \rmult e_{12} = e_{12}$ whereas $e_{12} \rmult e_{11} = \zero$.
$(\f^{\myt}_{\q},\radd,\rmult,\zero)$ is a ring, in fact isomorphic to the
 ring of $2 \times 2$ rational matrices\footnote{This interpretation can be extended to arbitrary rings of $n \times n$ matrices over a field $\myk$ by setting $\myt = \{ e_{ij} \mid 1 \leq i,j \leq n \}$, $e_{ij} \rmult e_{kl} = \zero$ in case $j \neq k$ and
 $e_{ij} \rmult e_{jl} = e_{il}$ else. The unit element then is $e_{11} + \ldots + e_{nn}$.}
It contains a unit element, namely $e_{11} + e_{22}$.
\exaend
}
\end{example}
Notice that in this example the unit element is {\em not}
 an element of the set $\myt$ embedded in $\f$.
Moreover, the multiplication here arises from the situation
 $\rmult : \myt \times \myt \myr \myt \cup \{ \zero \}$.
The next example even allows multiplications of terms to result in polynomials,
 i.e., $\rmult : \myt \times \myt \myr \f$.
\begin{example}\label{exa.skew}~\\
{\rm
For a set of variables $X_1, X_2, X_3$ let us define the set of commutative terms
 $\myt = \{ X_1^{i_1} X_2^{i_2} X_3^{i_3} \mid i_1, i_2, i_3 \in \n \}$
 and let $\f^{\myt}_{\q}$ be the set of all functions $f : \myt \myr \q$ with finite support,
 where $\q$ are the rational numbers.
Multiplication $\rmult : \myt \times \myt \myr \f$ is lifted from the following
 multiplication of the variables:
 $X_2 \rmult X_1 = X_2 + X_3$, $X_3 \rmult X_1 = X_1X_3$, $X_3 \rmult X_2 = X_2X_3$
 and $X_i \rmult X_j = X_iX_j$ for $i<j$.
Then $\f$ can be interpreted as a skew-polynomial ring $\q[X_1, X_2,X_3]$
 with unit element $X_1^{0} X_2^0 X_3^{0} \in \f^{\myt}_{\q}$.
\exaend
}
\end{example}
%
%\begin{example}\label{exa.quaternionen}~\\
%{\rm
%Let $\myt = \{ e,i,j,k \}$ be a finite set.
%Let $\f^{\myt}_{\q}$ be the set of all functions $f : \myt \myr \q$.
%We specify the multiplication $\rmult$ on $\f^{\myt}_{\q}$ by the
% action on $\myt$ as follows: 
% $i \rmult i = j \rmult j = k \rmult k = i \rmult j \rmult k = -e$, 
% $e \rmult e = e$, $e \rmult i = i \rmult e = i$,
% $e \rmult j = j \rmult e = j$, $e \rmult k = k \rmult e = k$,
% $i \rmult j = - j \rmult i = k$, $j \rmult k = -k \rmult j = i$,
% $k \rmult i = -i \rmult k = j$.
%The elements $\{ e,-e,i,-i,j,-j,k,-k \}$ form a non-Abelian group of order
% 8 under the product $\rmult$.
%One can also show that the non-zero elements of this function
% ring form a non-Abelian
% group under multiplication, i.e.~$\f^{\myt}_{\q}$ is in fact a division ring
% or a skew field.
%This ring is known as the ring of quaternions introduced by Hamilton.
%Its unit element is $e$.
%\exaend
%}
%\end{example}
%
Finally, many examples for function rings will be taken from monoid rings
 and hence we close this subsection by giving an example of a monoid ring.
\begin{example}\label{exa.free.group.ring}~\\
{\rm
Let $\myt = \{ a^i, b^i, 1 \mid i \in \n^+\}$,
 where $1$ is the empty
 word in $\{a,b\}^*$, and let the multiplication $\rmult$ be
 defined by  the following multiplication table: \\[1ex]
\renewcommand{\baselinestretch}{1.5}\small\normalsize
\begin{tabular}{c||c|c|c}
         & $1$ & $a^j$ & $b^{j}$ \\
\hline
\hline
$1$ & $1$ & $a^j$ & $b^{j}$ \\
\hline
$a^{i}$  & $a^{i}$ &   $a^{i+j}$ & $a^{i \mbox { {\tiny \sf monus} } j}b^{j \mbox { {\tiny \sf monus} } i}$ \\
\hline
$b^i$    & $b^i$ &  $a^{j \mbox { {\tiny \sf monus} } i}b^{i \mbox { {\tiny \sf monus} } j}$ & $b^{i+j}$ \\
\end{tabular}
\renewcommand{\baselinestretch}{1}\small\normalsize
\\[1ex]
where $i, j \in \n^+$ and $i \mbox{ {\tiny \sf monus} } j = i - j$ if $i \geq j$ and $0$ else.
In fact $\myt$ is the free group on one generator which can be presented as a monoid 
 by $(\{ a,b\}; \{ ab = ba = 1 \})$.
Let $\f^{\myt}_{\q}$ be the set of all functions $f : \myt \myr \q$ with finite support.
Then $\f^{\myt}_{\q}$ is a ring and is known as a special case of 
 the free group ring.
Its unit element is $1 \in \f^{\myt}_{\q}$.
\exaend
}
\end{example}
For the special case that we have $\rmult:  \myt \times \myt \myr \myt$, and
 some subring $\rr' \subseteq \rr$ we get that the function ring
 $\f^{\myt}_{\rr'}$ is a subring of $\f^{\myt}_{\rr}$.
This follows directly as then  for $f,g \in \f^{\myt}_{\rr'}$ we have
 $f + (-g),f \rmult g \in \f^{\myt}_{\rr'}$.
This is no longer true if $\rmult:  \myt \times \myt \myr \f^{\myt}_{\rr}$.
Let $\rr = \q$, $\rr'=\z$ and
 $\myt = \{ X_1^iX_2^j \mid i,j \in \n \}$ with $\rmult$ induced by
 $X_2 \rmult X_1 = \frac{1}{2}\skm X_1X_2$, $X_1 \rmult X_2 = X_1X_2$.
Then for $X_2,X_1 \in \f^{\myt}_{\z}$ we get $X_2 \rmult X_1 = \frac{1}{2}\skm X_1X_2
 \in \f^{\myt}_{\q}$.

Similarly, if we have $\myt' \subseteq \myt$ and $\rmult: \myt' \times \myt'
 \myr \f^{\myt'}_{\rr}$, then $\f^{\myt'}_{\rr}$ is a subring of
 $\f^{\myt}_{\rr}$. 
Again this follows as for $f,g \in \f^{\myt'}_{\rr}$ we have
 $f + (-g),f \rmult g \in \f^{\myt'}_{\rr}$.
Let us review Example \ref{exa.skew}:
There we have $\myt = \{ X_1^{i_1} X_2^{i_2} X_3^{i_3} \mid i_1, i_2, i_3 \in \n \}$
 and the multiplication $\rmult : \myt \times \myt \myr \f^{\myt}_{\q}$
 is lifted from the following
 multiplication of the variables:
 $X_2 \rmult X_1 = X_2 + X_3$, $X_3 \rmult X_1 = X_1X_3$, $X_3 \rmult X_2 = X_2X_3$
 and $X_i \rmult X_j = X_iX_j$ for $i<j$.
Then for $\myt' = \{ X_2^{i_2} X_3^{i_3} \mid i_2, i_3 \in \n \}$
 we have $\rmult : \myt' \times \myt' \myr \f^{\myt'}_{\q}$ and hence
 $\f^{\myt'}_{\q}$ is a subring of $\f^{\myt}_{\q}$.
%%%%%%%%%%%%%%%%%%%%%%%%%%%%%%%%%%%%%%%%%%%%%%%%%%%%%%%%%%%%%%%%%%%%%%
\section{Right Ideals and Right Standard Representations}\label{section.right.standard}
Since $\f$ is a ring, we can define right, left or two-sided ideals.
In this section in a first step we will restrict our attention to one-sided ideals,
 in particular to right ideals since left ideals in general
 can be treated in a symmetrical manner.

A subset $\mathfrak{i} \subseteq \f$ is called a
 \betonen{right ideal}, if
\begin{enumerate}
\item $\zero \in \mathfrak{i}$,
\item for $f,g \in \mathfrak{i}$ we have $f \radd g \in \mathfrak{i}$, and
\item for $f \in \mathfrak{i}$, $g \in \f$ we have
          $f \rmult g \in \mathfrak{i}$.
\end{enumerate}
Right ideals can also be specified in terms of 
 generating sets.
For $F \subseteq \f \backslash \{ \zero\}$ let
 $\ideal{r}{}(F) = \{ \sum_{i=1}^n f_i \rmult g_i \mid
  f_i \in F, g_i \in \f, n \in \n \} = \{ \sum_{i=1}^n f_i \rmult m_i \mid
 f_i \in F, m_i \in \monoms(\f),  n \in \n \}$.
These generated sets are subsets of $\f$ since for
 $f,g \in \f$ $f \rmult g$
 as well as $f \radd g$ are again elements of $\f$, and it is easily checked
 that they are in fact  right ideals:
 \begin{enumerate}
 \item $\zero \in \ideal{r}{}(F)$ since $\zero$ can be written as
        the empty sum.
 \item For two elements $\sum_{i=1}^n f_i \rmult g_i$ and 
        $\sum_{j=1}^m f_j \rmult h_j$ in $\ideal{r}{}(F)$,
        the resulting sum $\sum_{i=1}^n f_i \rmult g_i \radd
        \sum_{j=1}^m f_j \rmult h_j$ is again
        an element in $\ideal{r}{}(F)$.
 \item For an element $\sum_{i=1}^n f_i \rmult g_i$ in $\ideal{r}{}(F)$
        and a polynomial $h$ in $\f$, the product
        $(\sum_{i=1}^n f_i \rmult g_i) \rmult h =
          \sum_{i=1}^n f_i \rmult (g_i \rmult h)$  is again
        an element in $\ideal{r}{}(F)$.
 \end{enumerate}
Given a right ideal 
 $\mathfrak{i} \subseteq \f$ we
 call a set $F \subseteq \f \backslash \{ \zero\}$ a
 \betonen{basis} or a \betonen{generating set} of $\mathfrak{i}$ if
 $\mathfrak{i} = \ideal{r}{}(F)$.
Then every element $g \in \ideal{r}{}(F) \backslash \{ \zero\}$
 has different representations
 of the form
 $$g  =  \sum_{i=1}^n f_i \rmult h_i,  
   f_i \in F, h_i \in \f, n \in \n.$$
Of course the distributivity law in $\f$ then allows to convert any such
 representation into one of the form 
$$g  =  \sum_{j=1}^m f_i \rmult m_i,  
   f_i \in F, m_i \in \monoms(\f), m \in \n.$$
As we have seen in Section \ref{section.intro.applications}, it is not obvious
 whether some polynomial belongs to an ideal.
Let again $f_1 = X_1^{2} + X_2$ and $f_2 = X_1^{2} + X_3$ be two polynomials
 in the polynomial ring $\q[X_1,X_2,X_3]$ and
 ${\mathfrak i} = \{ f_1 \mrm g_1 + f_2
 \mrm g_2 \mid g_1, g_2 \in \q[X_1,X_2,X_3] \}$
 the (right) ideal generated by them. 
It is not hard to see that the polynomial $X_2 - X_3$ belongs to
 ${\mathfrak i}$ since $X_2 - X_3 = f_1 - f_2$ is a representation of
 $X_2 - X_3$ in terms of $f_1$ and $f_2$.
The same is true for the polynomial $X_2^2 - X_2X_3$ where now we have
 to use multiples of $f_1$ and $f_2$, namely $X_2^2 - X_2X_3 =
 f_1 \rmult X_2 - f_2 \rmult X_2$.
However, when looking at the polynomial $X_3^{3} + X_1 + X_3$
 we find that there is no obvious algorithm to find such appropriate multiples.
The problem is that for an arbitrary generating set for an ideal we have to
 look at arbitrary polynomial multiples with no boundary.
One first improvement for the situation can be achieved if we can 
 represent ideal elements by special representations in terms of the given
 generating set.
In polynomial rings such representations are studied as variations of the term
 \betonen{standard representations} in the literature
 (see also Section \ref{section.buchberger}).
They will also be introduced in this setting.
Since  standard representations are in general  distinguished by
 conditions involving an ordering on
 the set of polynomials, we will start by introducing the notion of
 an ordering to
 $\f$.

Let $\succeq$ be a total well-founded ordering on the set $\myt$. 
This enables us to make our polynomial representations of
 functions  unique by using the ordering $\succeq$ to arrange the elements
 of the support:
$$f = \sum_{i=1}^k m_{t_i}, \mbox{ where }
 \supp(f) = \{ t_1 , \ldots, t_k \}, t_1 \succ \ldots \succ t_k.$$

Using the ordering $\succeq$ on $\myt$ we are now able to give some notions for polynomials which
 are essential in introducing standard representations, standard bases
 and Gr\"obner bases in the classical approach.
We call  the monomial
 with the largest term according to $\succeq$ the \betonen{head monomial}
 of $f$ denoted by $\hm(f)$, consisting of  the \betonen{head term} 
 denoted by $\hterm(f)$ and the
 \betonen{head coefficient} denoted by  $\hc(f)= f(\hterm(f))$.
$f - \hm(f)$ is called the \betonen{reductum} of $f$ denoted by $\reductum(f)$.
Note that $\hm(f) \in \monoms(\f)$, $\hterm(f) \in \myt$ and $\hc(f) \in \rr$.
These notions can be extended to sets of functions 
$F \subseteq \f \backslash \{ \zero\}$ by
 setting $\hm(F) = \{ \hm(f) \mid f \in F \}$, $\hterm(F) = \{
 \hterm(f) \mid f \in F \}$ and $\hc(F) = \{ \hc(f) \mid f \in F \}$.

Notice that for some polynomial $f = \sum_{i=1}^k m_{t_i}\in\f$, 
 and some term $t \in \myt$
 we cannot conclude
 that for the terms occurring in the multiple
 $f \rmult t = \sum_{i=1}^k m_{t_i} \rmult t$
 we have $t_1 \rmult t \succ \ldots \succ t_k \rmult t$ (in case the
 multiplication of terms again results in terms)
 or $\hterm(t_1 \rmult t) \succ \ldots \succ \hterm(t_k \rmult t)$ as the ordering
 need not be compatible with multiplication in $\f$. 

\begin{example}\label{exa.not.stable}~\\ 
{\rm
Let $\myt = \{ x, 1 \}$ and $\rmult$ induced by
 the following multiplication on $\myt$: $x \rmult x = 1 \rmult 1 = 1$, 
 $x \rmult 1 = 1 \rmult x = x$.
Then assuming $x \succ 1$, after multiplying both sides of the equation
 with $x$, we get
 $x \rmult x  = 1 \pred 1 \rmult x = x$.
On the other hand, assuming the precedence
 $1 \succ x$ similarly we get $x = 1 \rmult x \pred
 1 = x \rmult x$.
Hence the ordering is not compatible with multiplication using elements in  $\myt$.
\exaend
}
\end{example}

We will later on see that this lack of compatibility leads to additional 
 requirements when defining standard representations, standard bases and 
 Gr\"obner bases.
Since the elements of $\myt$ can be identified with the terms in
 $\terms(\f)$, the ordering $\succeq$ can be extended as a total well-founded\footnote{An ordering $\succeq$ on a set ${\cal M}$
 will be called well-founded if its strict part $\succ$ is well-founded, i.e.,
 does not allow infinite descending chains of the form
 $m_1 \succ m_2 \succ \ldots$.}
 ordering on $\terms(\f)$.
Additionally we can provide orderings on $\monoms(\f)$ and $\f$ as follows.

\begin{definition}~\\
{\rm
Let $\succeq$ be a total well-founded ordering on $\myt$.
Let $>_{\rr}$ be a (not necessarily total) well-founded ordering  on $\rr$.
We define an
 ordering on $\monoms(\f)$
 by $m_{t_1} \succ m_{t_2}$
 if $t_1 \succ t_2$ or ($t_1 = t_2$ and $m_{t_1}(t_1) >_{\rr} m_{t_2}(t_2)$).
\\
For two elements $f,g$ in $\f$ we define
 $f \succ g$ iff $\hm(f) \succ \hm(g)$ or $(\hm(f) = \hm(g)$ and
 $\reductum(f) \succ \reductum(g))$.
We further define $f \succ \zero$ for all $f \in\f \backslash \{ \zero\}$.
\dend
}
\end{definition}

Notice that the total well-founded ordering on $\terms(\f)$ extends to
 a well-founded ordering on $\monoms(\f)$.

For a field $\myk$ we have the trivial ordering $>_{\myk}$ where
 $\alpha >_{\myk} 0$ for all $\alpha \in \myk\backslash{\{0\}}$ and no other
 elements are comparable.
Then the resulting ordering on the respective function ring corresponds to
 the one given in Definition \ref{def.ordering.polyring} for polynomial rings over
 fields.
\begin{lemma}\label{lem.wellfounded}~\\
{\sl
The ordering $\succ$ on $\f$ is well-founded.
\lemend
}
\end{lemma}
\Ba{}~\\
The proof of this lemma will use a method known as Cantor's
second diagonal argument (compare e.g.~\cite{BeWe92} Chapter 4).
Let us assume that $\succ$ is not well-founded on $\f$.
We will show that this gives us a contradiction to the fact that the
ordering $\succeq$ on $\monoms(\f)$ inducing $\succ$ is well-founded.
%\\
Hence, let us suppose $f_0 \succ f_1 \succ \ldots \succ f_k \succ \ldots\;$,
 $k \in \n$ is a strictly descending chain in $\f$.
%\\
Then we can construct a sequence of sets of pairs 
 $\{ \{ (m_{t_k}, g_{kn}) \mid n \in \n \} \mid k \in \n \}$ 
 recursively as follows:
%\\
For $k=0$ let $m_{t_0} = \min_{\succeq} \{ \hm(f_i) \mid i \in \n \}$
 which is well-defined since $\succeq$ is well-founded on $\monoms(\f)$.
Now let $j \in \n$ be the least index such that we have $m_{t_0} =\hm(f_j)$.
Then $m_{t_0} = \hm(f_{j + n})$ holds for all $n \in \n$ and we can set
$g_{0n} = f_{j+n} - \hm(f_{j+n})$, i.e., $m_{t_0} \succ \hm(g_{0n})$
for all $n \in \n$.
%\\
For $k+1$ we let $m_{t_{k+1}} = \min_{\succeq} \{  \hm(g_{ki}) \mid i \in \n \}$ and
again let $j \in \n$ be the least index such that $m_{t_{k+1}} =
\hm(g_{kj})$ holds, i.e., $m_{t_{k+1}} = \hm(g_{k(j+n)})$ for all $n \in \n$.
Again we set $g_{(k+1)n} = g_{k(j+n)} - \hm(g_{k(j+n)})$.
\\
Then the following statements hold for every $k \in \n$:
\begin{enumerate}
\item For all monomials $m$ occuring in the polynomials $g_{kn}$, $n \in \n$,
       we have $m_{t_k} \succ m$.
\item $g_{k0} \succ g_{k1} \succ \ldots\;$ is a strictly descending chain in $\f$.
\end{enumerate}
Hence we get that $m_{t_0} \succ m_{t_1} \succ \ldots \;$ is a strictly
descending chain in $\monoms(\f)$ contradicting the
fact that $\succeq$ is supposed to be well-founded on this set.
\\
\qed

Characterizations of ideal bases in terms of special standard representations
 they allow are mainly provided 
 for polynomial rings over {\em fields} in the literature
 (compare \cite{BeWe92} and Section \ref{section.buchberger}).
Hence we will first take a closer look at possible generalizations
 of these concepts to function rings over fields.

%%%%%%%%%%%%%%%%%%%%%%%%%%%%%%%%%%%%%%%%%%%%%%%%%%%%%%%%%%%%%%%%%%%
\subsection{The Special Case of Function Rings over Fields}\label{section.right.field}
Let $\f_{\myk}$ be a function ring over a field $\myk$.
Remember that for a set $F$ of polynomials in $\f_{\myk}$ every polynomial
 $g \in \ideal{r}{}(F)$ has a representation of the form
 $g  =  \sum_{i=1}^n f_i \rmult h_i,  
   f_i \in F, h_i \in \f_{\myk}, n \in \n.$
However, such an arbitrary representation can contain monomials
 larger than $\hm(g)$ which are cancelled in the sum.
A first idea of standard representations in the literature now is to represent
 $g$ as a sum of polynomial multiples $f_i \rmult h_i$ such that no
 cancellation of monomials larger than $\hm(g)$ takes place, i.e.
 $\hm(g) \succeq \hm(f_i \rmult h_i)$.
Hence in a first step we look at the following analogon of a definition
 of standard representations (compare \cite{BeWe92}, page 218):
\begin{definition}\label{def.general.standard.rep}~\\
{\rm
Let $F$ be a set of polynomials in $\f_{\myk}$
 and $g$ a non-zero polynomial in $\ideal{r}{}(F)$.
A representation of the form
 \begin{eqnarray}
  g & = & \sum_{i=1}^n f_i \rmult h_i, 
   f_i \in F, h_i \in \f_{\myk}, n \in \n\label{eqn.general.standard}
 \end{eqnarray}
 where additionally $\hterm(g) \succeq \hterm(f_i \rmult h_i)$ holds for
 $1 \leq i \leq n$ is called a \betonen{(general) right standard representation}
 of $g$ in terms of $F$.
If every $g \in \ideal{r}{}(F) \backslash \{ \zero\}$
 has such a representation in terms of $F$,
 then $F$ is called a \betonen{(general) right standard basis} of $\ideal{r}{}(F)$.
\dend
}
\end{definition}
What distinguishes an arbitrary representation from a (general) right
 standard representation is the fact that the former may contain polynomial
 multiples $f_i \rmult h_i$ with head terms $\hterm(f_i \rmult h_i)$
 larger than the head term of the represented polynomial $g$.
Therefore, in order to change an arbitrary representation into one fulfilling our
 additional condition (\ref{eqn.general.standard})
 we have to deal with special sums of polynomials.
\begin{definition}\label{def.general.critical.situations}~\\
{\rm
Let $F$ be a set of polynomials in $\f_{\myk}$ and $t$ an element in $\myt$.
Then we define the \betonen{critical set}
 ${\cal C}_{gr}(t,F)$ to contain all tuples of the form
 $(t, f_1, \ldots, f_k, h_1, \ldots, h_k)$, $k \in \n$, $f_1, \ldots, f_k \in F$\footnote{As in 
 the case of commutative polynomials, $f_1, \ldots, f_k$ are not
 necessarily different polynomials from $F$.},
 $h_1, \ldots, h_k \in \f_{\myk}$ such that
 \begin{enumerate}
 \item $\hterm(f_i \rmult h_i) = t$, $1 \leq i \leq k$, and
 \item $\sum_{i=1}^k \hm(f_i \rmult h_i) = \zero$.
 \end{enumerate}
We set ${\cal C}_{gr}(F) = \bigcup_{t \in \myt} {\cal C}_{gr}(t,F)$.
\dend
}
\end{definition}
Notice that for the sums of polynomial multiples in this definition we
 get $\hterm(\sum_{i=1}^k f_i \rmult h_i) \pred t$.
This definition is motivated by the definition
 of syzygies of polynomials
 in commutative polynomial rings over rings.
However, it differs from the original definition insofar as we need not have 
 $\hterm(f \rmult h) = \hterm(\hterm(f) \rmult \hterm(h))$, i.e.,
 we cannot localize the definition 
 to the head monomials of the polynomials in $F$.
Still we can characterize (general) right standard bases using this concept.
\begin{theorem}\label{theo.general.standard.basis}~\\
{\sl
Let $F$ be a set of polynomials in $\f_{\myk} \backslash \{ \zero\}$.
Then $F$ is a (general) right standard basis of $\ideal{r}{}(F)$ if and only if
 for every tuple 
 $(t, f_1, \ldots, f_k, h_1, \ldots, h_k)$ in ${\cal C}_{gr}(F)$
 the polynomial $\sum_{i=1}^k f_i \rmult h_i$ (i.e., the element in $\f_{\myk}$
 corresponding to this sum) has a  (general) right standard representation
 with respect to $F$.
\theoend
}
\end{theorem}
\Ba{}~\\
In case $F$ is a (general) right standard basis,
 since these polynomials are all elements of $\ideal{r}{}(F)$, they must
 have  (general) right standard representations with respect to $F$.
\\
To prove the converse, it remains to show that every element in
 $\ideal{r}{}(F)$ has a (general) right standard representation
 with respect to $F$.
Hence, let $g = \sum_{j=1}^m f_{j} \rmult h_{j}$ be an arbitrary
 representation of a non-zero  polynomial $g\in \ideal{r}{}(F)$ such that
 $f_j \in F$, $h_{j} \in \f_{\myk}$, $m \in \n$.
Depending on this  representation of $g$ and the
 well-founded total ordering $\succeq$ on $\myt$ we define
 $t = \max_{\succeq} \{ \hterm(f_{j} \rmult h_{j}) \mid 1\leq j \leq m \}$ and
 $K$ as the number of polynomials $f_j \rmult h_j$ with head term $t$.
%\\
Then $t \succeq \hterm(g)$ and 
 in case $\hterm(g) = t$ this immediately implies that this representation is
 already a (general) right standard  one. 
%\\
Else we proceed by induction
 on $t$.
%\\
Without loss of generality let $f_1, \ldots, f_K$ be the polynomials
 in the corresponding representation
 such that  $t=\hterm(f_i \rmult h_i)$, $1 \leq i \leq K$.
Then the tuple $(t, f_1, \ldots, f_K, h_1, \ldots, h_K)$
 is in ${\cal C}_{gr}(F)$ and let $h = \sum_{i=1}^K f_i \rmult h_i$.
%\\
We will now change our representation of $g$ in such a way that for the new
 representation of $g$ we have a smaller maximal term.
%\\
Let us assume $h$ is not $\zero$\footnote{In case  $h =\zero$,
 just substitute the empty sum for the representation of $h$
 in the equations below.}. 
%\\
By our assumption, $h$ has a (general)  right standard representation
 with respect to $F$, say $\sum_{j=1}^n p_j \rmult q_j$, 
 where $p_j \in F$, $q_j \in \f_{\myk}$, $n \in \n$ 
 and  all terms occurring in the sum are bounded by
 $t \succ \hterm(h)$ as $\sum_{i=1}^K \hm(f_i \rmult h_i) = \zero$.
%\\
This gives us: 
\begin{eqnarray}
  g   & = & \sum_{i=1}^K f_i \rmult h_i + \sum_{i=K+1}^m f_i \rmult h_i
             \nonumber\\                                                           
  & = & \sum_{j=1}^n p_j \rmult q_j + \sum_{i=K+1}^m f_i \rmult h_i
              \nonumber
\end{eqnarray}
which is a representation of $g$ where the maximal term of
 the involved polynomial multiples is smaller than $t$.
\\
\qed
Remember that by the distributivity law in $\f_{\myk}$
 any representation of a polynomial $g$ of the form
 $g  =  \sum_{i=1}^n f_i \rmult h_i,  
   f_i \in F, h_i \in \f_{\myk}, n \in \n$
can be converted into one of the form 
 $g  =  \sum_{j=1}^m f_j \rmult m_j,  
   f_j \in F, m_j \in \monoms(\f_{\myk}), m \in \n.$
Now for polynomial rings the conversion of a (general right) standard representation
 from a sum of polynomial multiples into a sum of monomial multiples
 again results in a standard representation.
This is due to the fact that the orderings used for the polynomial rings are
 compatible with multiplication.
Now let us look at a second analogon to this kind of
 standard representations in our setting.
\begin{definition}\label{def.standard.rep}~\\
{\rm
Let $F$ be a set of polynomials in $\f_{\myk}$
 and $g$ a non-zero polynomial in $\ideal{r}{}(F)$.
A representation of the form
 \begin{eqnarray}
  g & = & \sum_{i=1}^n f_i \rmult m_i, 
   f_i \in F, m_i \in \monoms(\f_{\myk}), n \in \n\label{eqn.standard}
 \end{eqnarray}
 where additionally $\hterm(g) \succeq \hterm(f_i \rmult m_i)$ holds for
 $1 \leq i \leq n$ is called a \betonen{right standard representation} of $g$ in terms of $F$.
If every $g \in \ideal{r}{}(F) \backslash \{ \zero\}$ has such a representation in terms of $F$,
 then $F$ is called a \betonen{right standard basis} of $\ideal{r}{}(F)$.
\dend
}
\end{definition}
If our ordering $\succ$ on $\f_{\myk}$ is compatible with $\rmult$
 we can conclude that the conversion of a general right 
 standard representation into a sum involving only monomial multiples again
 results in a right standard representation as defined in Definition \ref{def.standard.rep}.
But since in general the ordering and the multiplication are not compatible
 (review Example \ref{exa.not.stable}) a polynomial multiple $f \rmult h$ can
 contain monomials $m, m' \in \monoms(f \rmult m_j)$ where
 $h = \sum_{j=1}^n m_j$ such that $m$ and $m'$ are
 larger than $\hm(f \rmult h)$ and $m = m'$. 
Hence just applying the distributivity to a sum of polynomial multiples
 no longer changes a standard representation as defined in
 Definition \ref{def.general.standard.rep} into
 one as defined in Definition \ref{def.standard.rep}.
Remember that this was true for polynomial rings over fields where both definitions
 are equivalent.
Let us look at the monoid ring $\q[\m]$ where $\m$ is the monoid presented
 by $(\{a,b,c\};ab=a)$.
Moreover, let $\succ$ be the length-lexicographical ordering induced
 by the precedence $c \succ b \succ a$.
Then for the polynomials $f = ca +1$, $h = b^2 -b \in \q[\m]$ we get
 $\hterm(f \rmult b^2) = \hterm(ca + b^2) = ca$ and
 $\hterm(f \rmult b) = \hterm(ca + b) = ca$.
On the other hand $\hterm(f \rmult h) = \hterm(ca + b^2-ca-b) =
 \hterm( b^2 -b) = b^2$.
Hence for the polynomial $g = b^2-b$ the polynomial multiple $f \rmult h$
 is a general right standard representation as defined in Definition 
 \ref{def.general.standard.rep} while the sum of monomial multiples
 $f \rmult b^2 - f \rmult b$ is no right standard representation as
 defined in Definition \ref{def.standard.rep}.
We can even state that $g$ has no right standard representation in terms of
 the polynomial $f$.

Now as our aim is to link standard representations of polynomials to
 reduction relations, a closer inspection of the concept of general right
 standard representations shows that a reduction relation related to them
 has to involve polynomial multiples for defining the reduction steps. 
Right standard representations can also be linked to special instances
 of such reduction relations but are traditionally linked to reduction relations
 involving monomial multiples.
There is no example known from the literature where reduction relations involving
 polynomial multiples gain real advantages over reduction relations involving
 monomial multiples only\footnote{Examples where reduction relations involving polynomial multiples are studied for the original case of Gr\"obner bases in commutative polynomial rings can be found in \cite{Tr78,Za78}.}.
Therefore we will restrict our attention to right standard representations
 as presented in Definition \ref{def.standard.rep}.

Again, in order to change an arbitrary representation into one fulfilling our
 additional condition (\ref{eqn.standard}) of Definition \ref{def.standard.rep}
 we have to deal with special sums of polynomials.
\begin{definition}\label{def.critical.situations}~\\
{\rm
Let $F$ be a set of polynomials in $\f_{\myk}$ and $t$ an element in $\myt$.
Then we define the \betonen{critical set}
 ${\cal C}_{r}(t,F)$ to contain all tuples of the form
 $(t, f_1, \ldots, f_k, m_1, \ldots, m_k)$, $k \in \n$, $f_1, \ldots, f_k \in F$\footnote{As in 
 the case of commutative polynomials, $f_1, \ldots, f_k$ are not
 necessarily different polynomials from $F$.},
 $m_1, \ldots, m_k \in \monoms(\f)$ such that
 \begin{enumerate}
 \item $\hterm(f_i \rmult m_i) = t$, $1 \leq i \leq k$, and
 \item $\sum_{i=1}^k \hm(f_i \rmult m_i) = \zero$.
 \end{enumerate}
We set ${\cal C}_{r}(F) = \bigcup_{t \in \myt} {\cal C}_{r}(t,F)$.
\dend
}
\end{definition}
As before, we can characterize right standard bases using this concept.
\begin{theorem}\label{theo.standard.basis}~\\
{\sl
Let $F$ be a set of polynomials in $\f_{\myk} \backslash \{ \zero\}$.
Then $F$ is a  right standard basis of $\ideal{r}{}(F)$ if and only if
 for every tuple 
 $(t, f_1, \ldots, f_k, m_1, \ldots, m_k)$ in ${\cal C}_{r}(F)$
 the polynomial $\sum_{i=1}^k f_i \rmult m_i$ (i.e., the element in $\f$
 corresponding to this sum) has a  right standard representation
 with respect to $F$.
\theoend
}
\end{theorem}
\Ba{}~\\
In case $F$ is a  right standard basis,
 since these polynomials are all elements of $\ideal{r}{}(F)$, they must
 have  right standard representations with respect to $F$.
\\
To prove the converse, it remains to show that every element in
 $\ideal{r}{}(F)$ has a  right standard representation
 with respect to $F$.
Hence, let $g = \sum_{j=1}^m f_{j} \rmult m_{j}$ be an arbitrary
 representation of a non-zero  polynomial $g\in \ideal{r}{}(F)$ such that
 $f_j \in F$, $m_{j} \in \monoms(\f_{\myk})$, $m \in \n$.
Depending on this  representation of $g$ and the
 well-founded total ordering $\succeq$ on $\myt$ we define
 $t = \max_{\succeq} \{ \hterm(f_{j} \rmult m_{j}) \mid 1\leq j \leq m \}$ and
 $K$ as the number of polynomials $f_j \rmult m_j$ with head term $t$.
%\\
Then $t \succeq \hterm(g)$ and 
 in case $\hterm(g) = t$ this immediately implies that this representation is
 already a  right standard  one. 
%\\
Else we proceed by induction
 on $t$.
%\\
Without loss of generality let $f_1, \ldots, f_K$ be the polynomials
 in the corresponding representation
 such that  $t=\hterm(f_i \rmult m_i)$, $1 \leq i \leq K$.
Then the tuple $(t, f_1, \ldots, f_K, m_1, \ldots, m_K)$
 is in ${\cal C}_{r}(F)$ and let $h = \sum_{i=1}^K f_i \rmult m_i$.
%\\
We will now change our representation of $g$ in such a way that for the new
 representation of $g$ we have a smaller maximal term.
%\\
Let us assume $h$ is not $\zero$\footnote{In case  $h =\zero$,
 just substitute the empty sum for the representation of $h$
 in the equations below.}. 
%\\
By our assumption, $h$ has a  right standard representation
 with respect to $F$, say $\sum_{j=1}^n h_j \rmult l_j$, 
 where $h_j \in F$, $l_j \in \monoms(\f_{\myk})$, $n \in \n$ 
 and  all terms occurring in the sum are bounded by
 $t \succ \hterm(h)$ as $\sum_{i=1}^K \hm(f_i \rmult m_i) = \zero$.
%\\
This gives us: 
\begin{eqnarray}
  g   & = & \sum_{i=1}^K f_i \rmult m_i + \sum_{i=K+1}^m f_i \rmult m_i
             \nonumber\\                                                           
  & = & \sum_{j=1}^n h_j \rmult l_j + \sum_{i=K+1}^m f_i \rmult m_i
              \nonumber
\end{eqnarray}
which is a representation of $g$ where the maximal term of
 the involved monomial multiples is smaller than $t$.
\\
\qed
For commutative polynomial rings over fields standard bases are in fact
 Gr\"obner bases.
Remember that in algebraic terms a set $F$ is a Gr\"obner basis of
 the ideal $\ideal{}{}(F)$ it generates if and only if 
 $\hterm(\ideal{}{}(F)) = \{ t \rmult w \mid
 t \in \hterm(F), w \mbox{ a term} \}$
 (compare Definition \ref{def.buchberger.gb.algebraisch}).
The localization to the set of head terms only is possible as
 the ordering and multiplication are compatible, i.e.~$\hterm(f \rmult w)
 = \hterm(f) \rmult w$ for any $f \in F$ and any term $w$.
Then of course if every $g \in \ideal{}{}(F)$ has a standard representation
 in terms of $F$ we immediately get that $\hterm(g) = \hterm(f \rmult w)
 = \hterm(f) \rmult w$ for some $f \in F$ and some term $w$.
Moreover, for any reduction relation based on divisibility of terms we get that
 $g$ is reducible at its head monomial by this polynomial $f$.
This of course corresponds to the second definition of Gr\"obner bases in
 rewriting terms -- a set $F$ is a Gr\"obner basis of the ideal it generates
 if and only if the reduction relation $\red{}{\myr}{b}{F}$ associated
 to the polynomials in $F$ is confluent\footnote{The additional properties of capturing the ideal congruence
 and being terminating required by Definition \ref{def.gb} trivially hold for polynomial rings over fields.}
(compare Definition \ref{def.buchberger.gb}).
Central in both definitions of Gr\"obner bases is the idea of ``dividing'' terms.
Important in this context is the fact that divisors are smaller
 than the terms they divide with respect to term orderings
 and moreover the ordering on the terms is stable under multiplication
 with monomials.
The algebraic definition states that every head term of a polynomial in $\ideal{}{}(G)$
 has a head term of a polynomial in $G$ as a divisor\footnote{
When generalizing this definition to our setting of function rings
 we have to be very careful
 as in reality this implies that every polynomial in the ideal is reducible
 to zero which is the definition of a weak Gr\"obner basis (compare
 Definition \ref{def.weak.gb.reduction.ring}).
Gr\"obner bases and weak Gr\"obner bases coincide in  polynomial rings over
 fields due to the Translation Lemma (compare Lemma \ref{lem.buchberger.confluence} (2)).}. 
Similarly the reduction relation is based on divisibility of
 terms (compare Definition \ref{def.buchberger.red}).
The stability of the ordering under multiplication is important
 for the correctness of these characterizations of Gr\"obner bases since it allows
 finite localizations for the test sets to s-polynomials (Lemma \ref{lem.buchberger.confluence} is
 central in this context).

In our context now the ordering $\succ$ and the multiplication $\rmult$ on
 $\f_{\myk}$ in general are not compatible.
Hence, a possible algebraic definition of Gr\"obner bases and a definition
 of a reduction relation related to right standard representations must involve
 the whole polynomials and not only their head terms.

\begin{definition}\label{def.weak.gb}~\\
{\rm
A subset $F$ of $\f_{\myk} \backslash \{ \zero\}$ is called a 
 \betonen{weak right Gr\"obner basis} of $\ideal{r}{}(F)$ if
$\hterm(\ideal{r}{}(F) \backslash \{ \zero \}) =
 \hterm( \{ f \rmult m \mid f \in F, m \in \monoms(\f_{\myk}) \}
 \backslash \{ \zero \} )$.
\dend
}
\end{definition}
Instead of considering multiples of head terms of the generating set $F$
 we look at head terms of monomial multiples of polynomials in $F$.

In the setting of function rings over fields,
 in order to localize the definitions of standard representations and
 weak Gr\"obner bases to head terms instead of head monomials
 and show their equivalence
 we have to view $\f$ as a vector space
 with scalars from $\myk$.
We define a natural left \betonen{scalar multiplication} \label{page.scalar}
 $\skm : \myk \times \f \myr \f$
 by associating to $\alpha \in \myk$ and $f \in \f$  the function
 in $\f$, denoted by  $\alpha \skm f$, which has support
 $\supp(\alpha \skm f) \subseteq \supp(f)$
 and values $(\alpha \skm f)(t) = \alpha \skm f(t)$ for $t \in \supp(f)$.
Notice that if $\alpha \neq 0$ we have
 $\supp(\alpha \skm f) = \supp(f)$.
Similarly, we can define a natural right scalar multiplication
 $\skm :\f \times \myk \myr \f$
 by associating to $\alpha \in \myk$ and $f \in \f$  the function
 in $\f$, denoted by  $f \skm \alpha$, which has support
 $\supp(f \skm \alpha) \subseteq \supp(f)$
 and values $(f \skm \alpha)(t) = f(t) \skm \alpha$ for $t \in \supp(f)$. 
Since $\myk$ is associative we have 
\begin{eqnarray*}
((\alpha \skm f) \skm \beta) (t) & = & (\alpha \skm f)(t) \skm \beta \\
                                 & = & (\alpha \skm f(t)) \skm \beta \\
                                 & = & \alpha \skm (f(t) \skm \beta) \\
                                 & = & \alpha \skm ((f \skm \beta)(t)) \\
                                 & = & (\alpha \skm ( f \skm \beta))(t)
\end{eqnarray*}
and we will write $\alpha \skm f \skm \beta$.
Monomials can be represented as $m = \alpha \skm t$ where
 $\supp(m) = \{ t \}$ and $m(t) = \alpha$.
 
Additionally we have to state how scalar multiplication and ring multiplication are
 compatible.
Remember that we have introduced the elements of our function rings as formal sums of monomials.
We want to treat these objects similar to those occurring in the examples known from the literature.
In particular we want to achieve that multiplication in $\f_{\myk}$ can be specified by defining a multiplication on the terms and lifting it to the monomials.
Hence we 
 require the following equations
 $(\alpha \skm f) \rmult g = \alpha \skm (f \rmult g)$
 and $f \rmult (g \skm \alpha) = (f \rmult g) \skm \alpha$ to hold\footnote{Then of course since $\myk$ is Abelian we have
 $(\alpha \skm f) \rmult g = \alpha \skm (f \rmult g) = f \rmult (\alpha \skm g)
 = f \rmult (g \skm \alpha) = (f \rmult g) \skm \alpha$.}.
These equations are valid in the examples from the literature studied here.
The condition of course then implies that multiplication in $\f_{\myk}$
 can be specified by knowing $\rmult : \myt \times \myt \myr \f_{\myk}$.
This follows as for $\alpha, \beta \in \myk$ and $t,s \in \myt$ we have
\begin{eqnarray*}
(\alpha \skm t) \rmult (\beta \skm s) & = & \alpha \skm (t \rmult 
(\beta \skm  s)) \\
 & = & \alpha \skm( t \rmult ( s \skm \beta)) \\
 & = & \alpha \skm ( t \rmult s ) \skm \beta \\
 & = & (\alpha \skm \beta) \skm (t \rmult s).
\end{eqnarray*}
If $\f$ contains a unit element $\one$ the field can be embedded into
 $\f$ by $\alpha \longmapsto \alpha \skm \one$. 
Then for $\alpha \in \myk$ and $f \in \f_{\myk}$
 the equations $\alpha \skm f = (\alpha \skm \one) \rmult f$
 and $f \skm \alpha = f \rmult (\alpha \skm \one)$ hold.
Moreover, as $\myk$ is Abelian
 $\alpha \skm f \skm \beta = \alpha \skm \beta \skm f$ for any
 $\alpha, \beta \in \myk$, $f \in \f_{\myk}$.

In the next lemma we show that in fact both characterizations of
 special bases, right standard bases and weak Gr\"obner bases,
 coincide as in the case of polynomial
 rings over fields.
\begin{lemma}\label{lem.rsb=gb}~\\
{\sl
Let $F$ be a subset of $\f_{\myk} \backslash \{ \zero\}$.
Then $F$ is a right standard basis if and only
 if it is a weak right Gr\"obner basis.
\lemend
}
\end{lemma}
\Ba{}~\\
Let us first assume that $F$ is a right standard basis, i.e., every polynomial
 $g$ in $\ideal{r}{}(F)$ has a right standard  representation with respect
 to $F$.
In case $g \neq \zero$ this implies the existence of a polynomial $f \in F$ and a monomial
 $m \in \monoms(\f_{\myk})$ such that $\hterm(g) = \hterm(f \rmult m)$.
Hence $\hterm(g) \in \hterm( \{  f \rmult m \mid m \in \monoms(\f_{\myk}),
 f \in F \}\backslash \{ \zero\})$.
As the converse, namely $\hterm( \{  f \rmult m \mid m \in \monoms(\f_{\myk}),
 f \in F \}\backslash \{ \zero\})
 \subseteq \hterm(\ideal{r}{}(F) \backslash \{ \zero \})$ trivially holds,
 $F$ then is a weak right Gr\"obner basis.
\\
Now suppose that $F$ is a weak right Gr\"obner basis and again
 let $g \in \ideal{r}{}(F)$.
We have to show that $g$ has a right standard  representation with respect to $F$.
This will be done by induction on $\hterm(g)$.
In case $g = \zero$ the empty sum is our required right standard  representation.
Hence let us assume $g \neq \zero$.
Since then $\hterm(g) \in \hterm(\ideal{r}{}(F)\backslash \{ \zero\})$
 by the definition of weak
 right Gr\"obner bases we know there exists
 a polynomial $f \in F$ and a monomial $m \in \monoms(\f_{\myk})$ such that
 $\hterm(g) = \hterm(f \rmult m)$.
Then there exists a monomial $\tilde{m} \in \monoms(\f_{\myk})$ such that
 $\hm(g)=\hm(f \rmult \tilde{m})$,
 namely\footnote{Notice that this step requires that we can view $\f_{\myk}$
 as a vector space. In order to get a similar result without introducing
 vector spaces we would have to use a different definition of weak right
 Gr\"obner bases. E.g.~requiring that $\hm(\ideal{r}{}(F) \backslash \{ \zero \})
 = \hm(\{ f \rmult m \mid f \in F, m \in \monoms(\f_{\myk})\}\backslash \{ \zero \}\})$
 would be a possibility. However, then no localization of critical situations
 to head terms is possible, which is {\em the} advantage of having a field as 
 coefficient domain.} $\tilde{m}=(\hc(g) \skm \hc(f \rmult m)^{-1})\skm m)$.
Let $g_1 = g - f \rmult \tilde{m}$. 
Then $\hterm(g) \succ \hterm(g_1)$ implies the
 existence of a right standard  representation for $g_1$ which can be
 added to the multiple $f \rmult \tilde{m}$ to give
 the desired right standard  representation of $g$.
\\ 
\qed
Inspecting this proof closer we get the following corollary.
\begin{corollary}\label{cor.right.rep}~\\
{\sl
Let a subset $F$  of $\f_{\myk} \backslash \{ \zero\}$ be a 
 weak right Gr\"obner basis.
Then every $g \in \ideal{r}{}(F)$ has a right standard representation
       in terms of $F$ of the form
        $g = \sum_{i=1}^n f_i \rmult m_i,
         f_i \in F, m_i \in \monoms(\f_{\myk}), n \in \n$
         such that
         $\hm(g) = \hm(f_1 \rmult m_1)$ and
         $\hterm(f_1 \rmult m_1) \succ \hterm(f_2 \rmult m_2) \succ \ldots \succ
          \hterm(f_n \rmult m_n)$.
\corend
}
\end{corollary}
Notice that 
 we hence get stronger representations as
 specified in Definition \ref{def.standard.rep} for the case that the set $F$
 is a weak right Gr\"obner basis or a right standard basis.

In the literature Gr\"obner bases are linked to reduction relations.
These reduction relations in general then correspond to the respective
 standard representations as follows: if $g \red{*}{\myr}{}{F} \zero$, 
 then the monomial multiples involved in the respective reduction steps
 add up to a standard representation of $g$ in terms of $F$.
One possible reduction relation related to right standard representations
 as defined in Definition \ref{def.standard.rep}
 is called \betonen{strong reduction}\label{page.strong.reduction}\footnote{Strong
 reduction has been studied extensively
 for monoid rings in \cite{Re95}.} where a monomial $m_1$ is reducible
 by some polynomial $f$, if there exists some monomial $m_2$ such that
 $m_1 = \hm(f \rmult m_2)$.
Notice that such a reduction step eliminates the occurence of the term $\hterm(m_1)$
 in the resulting reductum $m_1 - f \rmult m_2$.
When generalizing this reduction relation to function rings we can no longer
 localize the reduction step to checking whether $\hm(f)$ divides $m_1$, as
 now the whole polynomial is involved in the reduction step.
We can no longer conclude that $\hm(f)$ divides $m_1$ but only that
 $m_1 = \hm(f \rmult m_2)$.

Our definition of weak right Gr\"obner bases using the condition
 $\hterm(\ideal{r}{}(F) \backslash \{ \zero \})$ $=
 \hterm( \{ f \rmult m \mid f \in F, m \in \monoms(\f_{\myk}) \}
 \backslash \{ \zero \} )$ in Definition \ref{def.weak.gb}
 corresponds to this problem that in many cases 
 orderings on $\myt$ are not compatible with the multiplication $\rmult$.
Let us review Example \ref{exa.not.stable}
  where the ordering $\succeq$ induced by $x \succ 1$ on
 terms respectively monomials
 is well-founded but in general not
 compatible with multiplication, due to the
 algebraic structure of $\myt$.
There for the polynomial $f = x + 1$ and the term $x$ we get
 $\hm(f \rmult x) = x$ while $\hm(f) \rmult x = 1$.\label{page.not.stable}

Behind this phenomenon lies the fact that the  definition
 of ``divisors'' arising from the algebraic characterization of weak Gr\"obner bases
 in the context of function rings does not have the same
 properties as divisors in polynomial rings.
One such important property is that divisors are smaller with respect to
 the ordering on terms and that this ordering is transitive.
Hence if $t_1$ is a divisor of $t_2$ and $t_2$ is a divisor of $t_3$
 then $t_1$ is also a divisor of $t_3$.
This is the basis of localizations when checking for the Gr\"obner basis
 property in polynomial rings over fields (compare Lemma \ref{lem.buchberger.confluence}).
Unfortunately this is no longer true for function rings in general.
Now $m_1 \in \hm(\ideal{r}{}(G))$ implies the existence of $m_2 \in \monoms(\f_{\myk})$
 such that $\hm(f \rmult m_2) = m_1$.
Reviewing the previous example we see that for $f = x + 1$, $m_2 = x$ and $m_1 = \hm(f) = x$ we get $\hm(f \rmult m_2) = \hm((x+1) \rmult x) = x$, i.e.~$\hm(f\rmult m_2)$
 divides $m_1$.
On the other hand $m_1=x$ divides $1$ as $x \rmult x = 1$.
But $\hm(\hm(f \rmult m_2) \rmult x) = 1$ while
 $\hm(f \rmult m_2 \rmult x) = x$, i.e.~the head monomial of the multiple
 involving the polynomial $f \rmult m_2$ does not divide $1$.

Notice that even if we restrict the concept of right divisors to monomials only we do not
 get transitivity.
We are interested when for some monomials $m_1,m_2,m_3 \in \monoms(\f_{\myk})$
 the facts that $m_1$ divides $m_2$ and $m_2$ divides $m_3$ imply that
 $m_1$ divides $m_3$.
Let $m, m' \in \monoms(\f_{\myk})$ such that
 $\hm(m_1 \rmult m) = m_2$ and $\hm(m_2 \rmult m') = m_3$.
Then $m_3 = \hm(m_2 \rmult m') = \hm( \hm(m_1 \rmult m) \rmult m')$.
When does this equal $\hm(m_1 \rmult m \rmult m')$ or even
 $\hm(m_1 \rmult \hm(m \rmult m'))$?\label{page.divisor}
Obviously if we have $\rmult : \monoms(\f_{\myk}) \times \monoms(\f_{\myk}) \mapsto \monoms(\f_{\myk})$, which is true for the Examples \ref{exa.punktweise},
 \ref{exa.polyring}
 and \ref{exa.matrix},
 this is true.
However if multiplication of monomials results in polynomials we are in trouble.
Let us look at the skew-polynomial ring $\q[X_1,X_2,X_3]$,
 $X_1 \succ X_2 \succ X_3$, defined in Example \ref{exa.skew},
 i.e.~$X_2 \rmult X_1 = X_2 + X_3$,$X_3 \rmult X_1 = X_1X_3$, $X_3 \rmult X_2 = X_2X_3$
 and $X_i \rmult X_j = X_iX_j$ for $i<j$.
Then from the fact that $X_2$ divides $X_2$ we get
 $\hm(X_2 \rmult X_1) = X_2$ and since again $X_2$ divides $X_2$,
 $\hm(\hm(X_2 \rmult X_1) \rmult X_1) = \hm(X_2 \rmult X_1) = X_2$.
But $\hm(X_2 \rmult X_1 \rmult X_1) = \hm(X_1X_3 + X_2 + X_3) = X_1X_3$.
Next we will show how using a restricted set of divisors only
 will enable some sort of transitivity.

To establish a certain kind of compatibility for the ordering $\succeq$ and
 the multiplication
 $\rmult$, additional requirements can be added.
One way to do this is by giving an additional ordering on $\myt$
 which is in some sense weaker than $\succeq$ but adds more information
 on compatibility with right multiplication.
Examples from the literature, where this technique is successfully applied,
 include special monoid and group rings (see e.g.~\cite{Re95,MaRe97a,MaRe95}).
There restrictions of the respective orderings on the monoid or group
 elements are of syntactical nature involving the presentation of the monoid
 or group (e.g.~prefix orderings of various kinds for commutative
 monoids and groups, free groups and polycyclic groups).
\begin{definition}\label{def.refined.ordering}~\\
{\rm
We will call an ordering $\geq$ on $\myt$ a 
 \betonen{right reductive restriction}
 of the ordering $\succeq$ or simply \betonen{right reductive}, 
 if the following hold:
\begin{enumerate}
\item $t \geq s$ implies $t \succeq s$ for $t, s \in \myt$.
\item $\geq$ is a partial ordering on $\myt$ which is
 compatible with multiplication $\rmult$ from the right in the following sense:
 if for $t, t_1, t_2, w \in \myt$,
 $t_2 \geq t_1$, $t_1 \succ t$ and
 $t_2 = \hterm(t_1 \rmult w)$ hold, then $t_2 \succ t \rmult w$.
\dend
\end{enumerate}
}
\end{definition}
Notice that if $\succeq$ is a partial well-founded ordering on $\myt$ so is $\geq$.

We can now distinguish special ``divisors'' of
 monomials:
For $m_1, m_2 \in \monoms(\f_{\myk})$ we call $m_1$ a \betonen{stable left  divisor} of $m_2$ if and only if 
 $\hterm(m_2) \geq \hterm(m_1)$
 and there exists $m \in \monoms(\f_{\myk})$ such that $m_2 = \hm(m_1 \rmult m)$.
Then $m$ is called a \betonen{stable right multiplier} of $m_1$.

If $\myt$ contains a unit element\footnote{I.e.~$\one \rmult t = t \rmult \one = t$ for all $t \in \myt$.} $\one$ and $\one \predeq t$
 for all terms $t \in \myt$ this immediately\footnote{As there are
 no terms smaller than $\one$ the second condition of Definition \ref{def.refined.ordering} trivially holds.} implies $\one \leq t$ and
 hence $\one$ is a stable divisor of any monomial $m$.
It remains to show that stable division is also transitive.
For three monomials $m_1, m_2, m_3 \in \monoms(\f)$ let
 $m_1$ be a stable divisor of $m_2$ and $m_2$ a stable divisor of $m_3$.
Then there exist monomials $m, m' \in \monoms(\f)$ such that
 $m_2 = \hm(m_1 \rmult m)$ with $\hterm(m_2) \geq \hterm(m_1)$ and
 $m_3 = \hm(m_2 \rmult m')$ with $\hterm(m_3) \geq \hterm(m_2)$.
Let us have a look at the monomial $\hm(\hm(m_1 \rmult m) \rmult m')$.
Remember how on page \pageref{page.divisor} we have seen that the case
 $m_1 \rmult m \in \monoms(\f)$ is not critical as then we immediately
 have that this monomial equals $\hm(m_1 \rmult m \rmult m') = 
 \hm(m_1 \rmult \hm(m \rmult m'))$.
Hence let us assume that $m_1 \rmult m \not\in \monoms(\f)$.
Then for all terms $s \in \terms(m_1 \rmult m) \backslash \hterm(m_1 \rmult m)$
 we know $s \pred \hterm(m_1 \rmult m)= \hterm(m_2)$.
Moreover $\hterm(m_3) \geq \hterm(m_2)$ and $\hterm(m_3) =
 \hterm(\hterm(m_2) \rmult \hterm(m'))$ then implies
 $\hterm(m_3) \succ \hterm(s \rmult \hterm(m'))$ and hence
 $\hm(\hm(m_1 \rmult m) \rmult m') = \hm(m_1 \rmult m \rmult m')$.
In both cases now $\hterm(m_3) \geq \hterm(m_1)$.
However, we cannot conclude that $\hm(m_1 \rmult m \rmult m') =
 \hm(m_1 \rmult \hm(m \rmult m'))$.
Still $m_1$ is a stable
 right divisor of $m_3$ as in case $m \rmult m'$ is a polynomial
 there exists some monomial $\tilde{m}$ in this polynomial such
 that $\hm(m_1 \rmult m \rmult m') = \hm(m_1 \rmult \tilde{m})$. 

The intention of restricting the ordering
 is that now, if $\hterm(m_2) \geq \hterm(m_1)$
 and $m_2 = m_1 \rmult m$, then for all terms $t$ with $\hterm(m_1) \succ t$ we
 then can conclude $\hterm(m_2) \succ \hterm(t \rmult m)$,
 which will be used to localize the
 multiple $\hterm(m_1 \rmult m)$ to $\hterm(m_1)$ achieving an equivalent to the properties
 of ``divisors'' in the case of commutative polynomial rings.
Under certain conditions
 reduction relations based on this divisibility property for
 terms will have the stability
 properties we desire.
On the other hand, restricting the choice of divisors in this way will lead
 to reduction relations which in general no longer capture the
 respective right ideal congruences\footnote{Prefix reduction for monoid rings
 is an example where the right ideal congruence is lost. See e.g.~\cite{MaRe95}
 for more on this topic.}.
\begin{example}\label{exa.skew2}~\\
{\rm
In Example \ref{exa.polyring} of a commutative polynomial ring 
 we can state a reductive restriction of any term ordering by
 $t \geq s$ for two terms $t$ and $s$ if and only if $s$ divides $t$ as 
 a term, i.e.~for $t=X_1^{i_1} \ldots X_n^{i_n}$, $s=X_1^{j_1} \ldots X_n^{j_n}$
 we have $j_l \leq i_l$, $1 \leq l \leq n$.
The same is true for skew-polynomial rings as defined by Kredel in his
 PhD thesis \cite{Kr93}.
The situation changes if for the defining equations of skew-polynomial rings,
 $X_j \rmult X_i = c_{ij} \skm X_iX_j + p_{ij}$ where
 $i<j$, $p_{ij} \pred X_iX_j$, we allow $c_{ij} = 0$.
Then other restrictions of the ordinary term orderings have to be
 considered due to the possible vanishing of head terms.
Let $X_2 \rmult X_1 = X_1, X_3 \rmult X_1 = X_1X_3,
 X_3 \rmult X_2 = X_2X_3$ and $\succ$ a term ordering
 with precedence $X_3 \succ X_2 \succ X_1$.
Then, although $X_2 \succ X_1$, as
 $X_2 \rmult (X_1X_2) = X_1X_2$ and $X_1 \rmult (X_1X_2) = X_1^2X_2 
 \succ X_1X_2$, we get
 $X_2 \rmult (X_1X_2) \pred X_1 \rmult (X_1X_2)$.
Hence, since $X_2$ is a divisor of $X_1X_2$ as a term, the classical restriction
 for polynomial rings no longer holds as $X_2$ is no stable divisor
 of $X_1X_2$.
For these cases the  restriction to $u < v$ if and only if $u$ is a prefix of $v$ as a word will work.
Then we know that for the respective term $w$ with $u \rmult w = v$ multiplication
 is just concatenation of $u$ and $w$ as words
 and hence for all $t \pred u$ the result of $t \rmult w$
 is again smaller than $u \rmult w$.
\exaend
}
\end{example}
Let us continue with algebraic consequences related to
 the right reductive restriction of our ordering by distinguishing special 
 standard representations.
Notice that for standard representations in commutative polynomial rings
 we already have that $\hterm(g) = \hterm(f_i \rmult m_i)$ implies
 $\hterm(g) = \hterm(f_i) \rmult \hterm(m_i)$ and for all $t \pred \hterm(f_i)$
 we have $t \rmult w \pred \hterm(f_i) \rmult w$ for any term $w$.
In the setting of function rings an analogon to the
 latter property now can be achieved by restricting the monomial multiples in
 the representation to stable ones.
Herefore we have different possibilities to incorporate these restrictions
 into the condition  $\hterm(g) \succeq \hterm(f_i \rmult m_i)$ of Definition
 \ref{def.standard.rep.polyring} and
 Definition \ref{def.standard.rep}.
The most general one is to require 
 $\hterm(g) = \hterm(f_1 \rmult m_1)= \hterm(\hterm(f_1) \rmult m_1)
 \geq \hterm(f_1)$ and $\hterm(g) \succeq \hterm(f_i \rmult m_i)$ 
 for all $2 \leq i \leq n$.
Then a representation of $g$ can contain further monomial multiples
 $f_j \rmult m_j$, $2 \leq j \leq n$ with $\hterm(g) = \hterm(f_j \rmult m_j)$
 not fullfilling the restriction on the first multiple of $f_1$.
Hence when defining critical situations we have to look at the same set
 as in Definition \ref{def.critical.situations}.
Another generalization is to demand
 $\hterm(g) = \hterm(f_1 \rmult m_1)= \hterm(\hterm(f_1) \rmult m_1)
 \geq \hterm(f_1)$ and $\hterm(g) \succeq \hterm(f_i \rmult m_i)=
 \hterm(\hterm(f_i) \rmult m_i) \geq \hterm(f_i)$ 
 for all $2 \leq i \leq n$.
Then critical situations can be localized to stable multiplers.
But we can also give a weaker analogon as follows:
\begin{definition}\label{def.right_reductive}~\\
{\rm
Let $F$ be a set of polynomials in $\f_{\myk}$ 
 and $g$ a non-zero polynomial in $\ideal{r}{}(F)$.
A representation of the form 
$$g = \sum_{i=1}^n f_i \rmult m_i,
 f_i \in F, m_i \in \monoms(\f_{\myk}), n \in \n$$
 such that
 $\hterm(g)  =
 \hterm(f_i \rmult m_i)= \hterm(\hterm(f_i) \rmult m_i) \geq \hterm(f_i)$ 
 for $1 \leq i \leq k$, for some $k \geq 1$, and
 $\hterm(g) \succ \hterm(f_i \rmult m_i)$ 
 for $k < i \leq n$ is called a
 \betonen{right reductive standard  representation} in terms of $F$.
\dend
}
\end{definition}
Notice that we restrict the possible multipliers to stable ones if the monomial
 multiple has the same head term as $g$, i.e.~contributes to the head term of $g$.  
For definitions sake we will let the empty sum be the right reductive standard 
 representation of $\zero$.
The idea behind right reductive standard  representations is that for an appropriate
 definition of a reduction relation based now on stable divisors
 such representations will again allow a reduction step to take
 place at the head monomial.

In case we have $\rmult : \myt \times \myt \myr \myt$ we can rephrase
 the condition in Definition \ref{def.right_reductive} to
 $\hterm(g) = 
 \hterm(f_i \rmult m_i) = \hterm(f_i) \rmult \hterm(m_i) \geq \hterm(f_i)$, $1 \leq i \leq k$.
\begin{definition}~\\
{\rm
A set $F \subseteq \f_{\myk} \backslash \{ \zero\}$ is called a
 \betonen{right reductive standard basis} (with respect to the reductive ordering $\geq$)
 of $\ideal{r}{}(F)$ if every polynomial $f \in \ideal{r}{}(F)$
 has a right reductive standard  representation in terms of $F$. 
\dend
}
\end{definition}
Again, in order to change an arbitrary representation into one fulfilling our
 additional condition of Definition \ref{def.right_reductive}
 we have to deal with special sums of polynomials.
\begin{definition}\label{def.right_reductive.critical.situations}~\\
{\rm
Let $F$ be a set of polynomials in $\f_{\myk}$ and $t$ an element in $\myt$.
Then we define the \betonen{critical set}
 ${\cal C}_{rr}(t,F)$ to contain all tuples of the form
 $(t, f_1, \ldots, f_k, m_1, \ldots, m_k)$, $k \in \n$, $f_1, \ldots, f_k \in F$\footnote{As in 
 the case of commutative polynomials, $f_1, \ldots, f_k$ are not
 necessarily different polynomials from $F$.},
 $m_1, \ldots, m_k \in \monoms(\f)$ such that
 \begin{enumerate}
 \item $\hterm(f_i \rmult m_i) =\hterm(\hterm(f_i) \rmult m_i)= t$,
       $1 \leq i \leq k$, 
 \item $\hterm(f_i \rmult m_i) \geq \hterm(f_i)$, $1 \leq i \leq k$, and
 \item $\sum_{i=1}^k \hm(f_i \rmult m_i) = \zero$.
 \end{enumerate}
We set ${\cal C}_{rr}(F) = \bigcup_{t \in \myt} {\cal C}_{rr}(t,F)$.
\dend
}
\end{definition}
Unfortunately, in contrary to the characterization of right standard
 bases in Theorem \ref{theo.standard.basis} these critical situations
 will not be sufficient to characterize right reductive standard bases.
To see this let us consider the following example:

\begin{example}\label{exa.free.group}~\\
{\rm
Let us recall the description of the free group ring
 in Example \ref{exa.free.group.ring} with $\myt = \{a^i,b^i,1 \mid i \in \n^+\}$
 and let $\succeq$ be the ordering induced by the length-lexicographical odering 
 on $\myt$ resulting from the precedence $a \succ b$.
\\
Then the set consisting of the polynomial $a + 1$ does not give rise to non-trivial
 critical situations,
 but still is no right reductive standard basis as the polynomial
 $b + 1 \in \ideal{r}{}(\{ a+ 1 \})$ has
 no right reductive standard representation with respect to $a + 1$.
\exaend
}
\end{example}

However, the failing situation $b+1 = (a+1)\rmult b$
 described in Example \ref{exa.free.group} 
 describes the only kind of additional critical
 situations which have to be resolved in order to characterize 
 right reductive standard bases.
\begin{theorem}\label{theo.right_reductive.standard.basis}~\\
{\sl
Let $F$ be a set of polynomials in $\f_{\myk} \backslash \{ \zero\}$.
Then $F$ is a  right reductive standard basis of $\ideal{r}{}(F)$ if and only if
\begin{enumerate}
 \item for every $f \in F$ and every $m \in \monoms(\f_{\myk})$ the multiple
 $f \rmult m$ has a right reductive standard representation in terms of $F$,
 \item
 for every tuple 
 $(t, f_1, \ldots, f_k, m_1, \ldots, m_k)$ in ${\cal C}_{rr} (F)$
 the polynomial $\sum_{i=1}^k f_i \rmult m_i$ (i.e., the element in $\f$
 corresponding to this sum) has a  right reductive standard representation
 with respect to $F$.
\end{enumerate}
\theoend
}
\end{theorem}
\Ba{}~\\
In case $F$ is a  right reductive standard basis,
 since these polynomials are all elements of $\ideal{r}{}(F)$, they must
 have  right reductive standard representations with respect to $F$.
\\
To prove the converse, it remains to show that every element in
 $\ideal{r}{}(F)$ has a  right reductive standard representation
 with respect to $F$.
Hence, let $g = \sum_{j=1}^m f_{j} \rmult m_{j}$ be an arbitrary
 representation of a non-zero  polynomial $g\in \ideal{r}{}(F)$ such that
 $f_j \in F$, and $m_{j} \in \monoms(\f_{\myk})$.
By our first statement every such monomial multiple $f_j \rmult m_j$ has
 a right reductive standard representation in terms of $F$ and we can
 assume that all multiples are replaced by them.
Depending on this  representation of $g$ and the
 well-founded total ordering $\succeq$ on $\myt$ we define
 $t = \max_{\succeq} \{ \hterm(f_{j} \rmult m_{j}) \mid 1\leq j \leq m \}$ and
 $K$ as the number of polynomials $f_j \rmult m_j$ with head term $t$.
Then for each multiple
 $f_j \rmult m_j$ with $\hterm(f_j \rmult m_j) =t$ we know
 that $\hterm(f_j \rmult m_j) = \hterm(\hterm(f_j) \rmult m_j)
 \geq \hterm(f_j)$ holds.
%\\
Then $t \succeq \hterm(g)$ and 
 in case $\hterm(g) = t$ this immediately implies that this representation is
 already a  right reductive standard  one. 
%\\
Else we proceed by induction
 on $t$.
%\\
Without loss of generality let $f_1, \ldots, f_K$ be the polynomials
 in the corresponding representation
 such that  $t=\hterm(f_i \rmult m_i)$, $1 \leq i \leq K$.
Then the tuple $(t, f_1, \ldots, f_K, m_1, \ldots, m_K)$
 is in ${\cal C}_{rr}(F)$ and let $h = \sum_{i=1}^K f_i \rmult m_i$.
%\\
We will now change our representation of $g$ in such a way that for the new
 representation of $g$ we have a smaller maximal term.
%\\
Let us assume $h$ is not $\zero$\footnote{In case  $h =\zero$,
 just substitute the empty sum for the representation of $h$
 in the equations below.}. 
%\\
By our assumption, $h$ has a  right reductive standard representation
 with respect to $F$, say $\sum_{j=1}^n h_j \rmult l_j$, 
 where $h_j \in F$, and $l_j \in \monoms(\f_{\myk})$ 
 and  all terms occurring in the sum are bounded by
 $t \succ \hterm(h)$ as $\sum_{i=1}^K \hm(f_i \rmult m_i) = \zero$.
%\\
This gives us: 
\begin{eqnarray}
  g   & = & \sum_{i=1}^K f_i \rmult m_i + \sum_{i=K+1}^m f_i \rmult m_i
             \nonumber\\                                                           
  & = & \sum_{j=1}^n h_j \rmult l_j + \sum_{i=K+1}^m f_i \rmult m_i
              \nonumber
\end{eqnarray}
which is a representation of $g$ where the maximal term is smaller than $t$.
\\
\qed
We can similarly refine Definition \ref{def.weak.gb} with respect to
 a reductive restriction $\geq$ of the ordering $\succeq$.
\begin{definition}\label{def.gb.fr}~\\
{\rm
A set $F \subseteq \f_{\myk} \backslash \{ \zero\}$ is called a 
 \betonen{weak right reductive Gr\"obner basis} (with respect to
 the reductive ordering $\geq$) of $\ideal{r}{}(F)$ if
$\hterm(\ideal{r}{}(F) \backslash \{ \zero \}) =
 \hterm( \{ f \rmult m \mid f \in F, m \in \monoms(\f_{\myk}),
  \hterm(f \rmult m) =  \hterm(\hterm(f) \rmult m) \geq \hterm(f) \}
 \backslash \{ \zero \} )$.
\dend
}
\end{definition}
This definition now localizes the characterization of the Gr\"obner basis
 to the head terms of the generating set of polynomials.

The next lemma states that in fact both characterizations of
 special bases, right reductive standard bases and weak right reductive
 Gr\"obner bases,
 coincide as in the case of polynomial
 rings over fields.
\begin{lemma}\label{lem.sb=gb}~\\
{\sl
Let $F$ be a subset of $\f_{\myk} \backslash \{ \zero\}$.
Then $F$ is a right reductive standard basis if and only
 if it is a weak right reductive Gr\"obner basis.
\lemend
}
\end{lemma}
\Ba{}~\\
Let us first assume that $F$ is a right reductive standard basis, i.e., every polynomial
 $g$ in $\ideal{r}{}(F)$ has a right reductive standard  representation with respect
 to $F$.
In case $g \neq \zero$ this implies the existence of a polynomial $f \in F$ and a monomial
 $m \in \monoms(\f_{\myk})$ such that $\hterm(g) = 
 \hterm(f \rmult m) = \hterm(\hterm(f) \rmult m) \geq \hterm(f)$.
Hence $\hterm(g) \in \hterm( \{  f \rmult m \mid m \in \monoms(\f_{\myk}),
 f \in F, \hterm(f \rmult m) = \hterm(\hterm(f) \rmult m) \geq \hterm(f) \}\backslash \{ \zero\})$.
As the converse, namely $\hterm( \{  f \rmult m \mid m \in \monoms(\f_{\myk}),
 f \in F, \hterm(f \rmult m) = \hterm(\hterm(f) \rmult m) \geq \hterm(f) \}\backslash \{ \zero\})
 \subseteq \hterm(\ideal{r}{}(F) \backslash \{ \zero \})$ trivially holds,
 $F$ is then a weak right reductive Gr\"obner basis.
\\
Now suppose that $F$ is a weak right reductive Gr\"obner basis and again let $g \in \ideal{r}{}(F)$.
We have to show that $g$ has a right reductive standard  representation with respect to $F$.
This will be done by induction on $\hterm(g)$.
In case $g = \zero$ the empty sum is our required right reductive standard  representation.
Hence let us assume $g \neq \zero$.
Since then $\hterm(g) \in \hterm(\ideal{r}{}(F)\backslash \{ \zero\})$
 by the definition of weak
 right  reductive Gr\"obner bases we know there exists
 a polynomial $f \in F$ and a monomial $m \in \monoms(\f_{\myk})$ such that
 $\hterm(g) = \hterm(f \rmult m) = \hterm(\hterm(f) \rmult m)\geq \hterm(f)$.
Then there exists a monomial $\tilde{m} \in \monoms(\f)$ such that
 $\hm(g)=\hm(f \rmult \tilde{m})$,
 namely\footnote{Notice that this step again requires that we can view $\f$
 as a vector space.} $\tilde{m}=(\hc(g) \skm \hc(f \rmult m)^{-1})\skm m)$.
Let $g_1 = g - f \rmult \tilde{m}$. 
Then $\hterm(g) \succ \hterm(g_1)$ implies the
 existence of a right reductive standard  representation for $g_1$ which can be
 added to the multiple $f \rmult \tilde{m}$ to give
 the desired right reductive standard  representation of $g$.
\\ 
\qed
An inspection of the proof shows that in fact we can require a stronger condition
 for the head terms of the monomial multiples involved in right
 reductive standard representations in terms of right reductive Gr\"obner bases.
\begin{corollary}\label{cor.right_reductive}~\\
{\sl
Let a subset $F$  of $\f_{\myk} \backslash \{ \zero\}$ be a 
 weak right reductive Gr\"obner basis.
Then every $g \in \ideal{r}{}(F)$ has a right reductive standard representation
       in terms of $F$ of the form
        $g = \sum_{i=1}^n f_i \rmult m_i,
         f_i \in F, m_i \in \monoms(\f), n \in \n$
         such that
         $\hterm(g) = \hterm(f_1 \rmult m_1) 
                    \succ \hterm(f_2 \rmult m_2) \succ \ldots \succ
                          \hterm(f_n \rmult m_n)$, and
         $\hterm(f_i \rmult m_i) =\hterm(\hterm(f_i) \rmult m_i)\geq \hterm(f_i)$
          for all $1 \leq i \leq n$.
\corend
}
\end{corollary}
The importance of Gr\"obner bases in commutative polynomial
 rings stems from the fact that they can be characterized by
 special polynomials, the so-called s-polynomials, and that only finitely many
 such polynomials have to be checked in order to decide whether a set is a Gr\"obner basis.
This test can be combined with adding ideal elements to the generating set leading to an algorithm
 which computes finite Gr\"obner bases by means of completion.
These finite sets then can be used to solve many problems related to the ideals
 they generate.

Given a field as coefficient domain
 the critical situations for function rings now
 lead to s-polynomials as in the original case
 and can be identified by studying term multiples of polynomials.
Let $p$ and $q$ be two non-zero polynomials in $\f_{\myk}$.
We are interested in terms $t, u_1, u_2$ such that
 $ \hterm(p \rmult u_1) = \hterm(\hterm(p) \rmult u_1) = t =
   \hterm(q \rmult u_2) = \hterm(\hterm(q) \rmult u_2)$
 and $\hterm(p) \leq t$, $\hterm(q) \leq t$.
Let ${\cal C}_{s}(p,q)$ (this is a specialization of Definition
 \ref{def.right_reductive.critical.situations})
 be the critical set containing all such tuples $(t,u_1,u_2)$
 (as a short hand for $(t, p, q, u_1, u_2)$).
We call the polynomial $\hc(p \rmult u_1)^{-1} \skm p \rmult u_1 -
 \hc(q \rmult u_2)^{-1} \skm q \rmult u_2 = \spol{r}(p,q,t,u_1,u_2)$
 the \betonen{s-polynomial}
 of $p$ and $q$ related to the tuple $(t,u_1,u_2)$.

\begin{theorem}\label{theo.rrsb=rgb.k}~\\
{\sl
Let $F$ be a set of polynomials in $\f_{\myk} \backslash \{ \zero\}$.
Then $F$ is a weak right reductive Gr\"obner basis of $\ideal{r}{}(F)$ if and only if
\begin{enumerate}
\item for all $f$ in $F$ and for $m \in \monoms(\f_{\myk})$ the multiple
      $f \rmult m$ has a right reductive standard  representation in terms of
      $F$, and
\item for all $p$ and $q$ in $F$ and every tuple $(t,u_1,u_2)$
      in ${\cal C}_{s}(p,q)$ the respective s-polynomial $\spol{r}(p,q,t,u_1,u_2)$
      has a right reductive standard 
      representation in terms of $F$.
\end{enumerate}
}
\end{theorem}
\Ba{}~\\
In case $F$ is a weak right reductive Gr\"obner basis it is also a right reductive standard basis, and
 since all multiples $f \rmult m$ and s-polynomials $\spol{r}(p,q,t,u_1,u_2)$ stated above
 are elements of $\ideal{r}{}(F)$, they must
 have right reductive standard  representations in terms of $F$.
\\
The converse will be proven by showing that every element in
 $\ideal{r}{}(F)$ has a right reductive standard  representation in terms of $F$.
Now, let $g = \sum_{j=1}^m f_{j} \rmult m_{j}$ be an arbitrary
 representation of a non-zero  polynomial $g\in \ideal{r}{}(F)$ such that
 $f_j \in F$, $m_{j} \in \monoms(\f)$, $m \in \n$.
By our first assumption every multiple $f_j \rmult m_j$ in this
 sum has a right reductive representation.
Hence without loss of generaltity we can assume that
 $\hterm(\hterm(f_j) \rmult m_j) = \hterm(f_j \rmult m_j) \geq
 \hterm(f_j)$ holds.
\\
Depending on this  representation of $g$ and the
 well-founded total ordering $\succeq$ on $\myt$ we define
 $t = \max_{\succeq} \{ \hterm(f_{j} \rmult m_{j}) \mid 1\leq j \leq m \}$ and
 $K$ as the number of polynomials $f_j \rmult m_j$ with head term $t$.
Without loss of generality we can assume that the multiples
 with head term $t$ are just $f_1  \rmult m_1, \ldots , f_K \rmult m_K$.
We proceed by induction
 on $(t,K)$, where
 $(t',K')<(t,K)$ if and only if $t' \prec t$ or $(t'=t$ and
 $K'<K)$\footnote{Note that this ordering is well-founded since $\succ$
                  is well-founded on $\myt$ and $K \in\n$.}.
\\
Obviously, $t \succeq \hterm(g)$ holds. 
If $K = 1$ this gives us $t = \hterm(g)$ and by our assumptions our
 representation is already of the required form.
Hence let us assume $K > 1$, then there are  two not
 necessarily different polynomials $f_1,f_2$ and corresponding monomials
 $m_1 = \alpha_1 \skm w_1$, $m_2 = \alpha_2 \skm w_2$ with $\alpha_1,\alpha_2
 \in \myk$, $w_1,w_2 \in \myt$
 in the corresponding   representation
      such that  $t=\hterm(\hterm(f_1) \rmult w_1) =\hterm(f_1 \rmult w_1) =
       \hterm(f_2 \rmult w_2) =  \hterm(\hterm(f_2) \rmult w_2)$ and
       $t \geq \hterm(f_1)$, $t \geq \hterm(f_2)$.
     Then the tuple $(t, w_1,w_2)$ is in ${\cal C}_{s}(f_1,f_2)$ and
      we have an s-polynomial
      $h= \hc(f_1 \rmult w_1)^{-1} \skm  f_1
      \rmult w_1 -\hc(f_2 \rmult w_2)^{-1}\skm f_2 \rmult w_2$ corresponding
      to this tuple.
%\\
We will now change our representation of $g$ by using the additional
 information on this s-polynomial in such a way that for the new
 representation of $g$ we either have a smaller maximal term or
 the occurrences of the term $t$
 are decreased by at least 1.
%\\
Let us assume the s-polynomial is not $\zero$\footnote{In case 
               $h =\zero$,
               just substitute the empty sum
               for the right reductive representation of $h$
               in the equations below.}. 
%\\
By our assumption, $h$ has a right reductive standard  representation in terms of $F$, say
  $\sum_{i=1}^n h_i \rmult l_i$, 
  where $h_i \in F$, and $l_i \in \monoms(\f_{\myk})$ 
  and  all terms occurring in the sum are bounded by
  $t \succ \hterm(h)$.
%\\
This gives us: 
     \begin{eqnarray}
       &   & f_1 \rmult m_1 + f_2 \rmult m_2      \nonumber\\  
       &   &  \nonumber\\ 
       & = & \alpha_1 \skm f_1 \rmult w_1 + \alpha_2 \skm f_2 \rmult w_2
             \nonumber\\  
       &   &  \nonumber\\                                                         
       & = &  \alpha_1 \skm f_1 \rmult w_1 +
              \underbrace{ \alpha'_2 \skm \beta_1 \skm f_1 \rmult w_1
                   - \alpha'_2 \skm \beta_1 \skm f_1 \rmult w_1}_{=\, 0} 
                   + \underbrace{\alpha'_2\skm \beta_2 }_{= \alpha_2} \skm f_2 \rmult w_2 \nonumber\\
       &   & \nonumber\\ 
       & = & (\alpha_1 + \alpha'_2 \skm \beta_1) \skm f_1 \rmult w_1 - \alpha'_2 \skm 
               \underbrace{(\beta_1 \skm f_1 \rmult w_1
             -  \beta_2 \skm f_2 \rmult w_2)}_{=\,
             h} 
             \nonumber\\
       &   & \nonumber\\ 
       & = & (\alpha_1 + \alpha'_2 \skm \beta_1) \skm f_1 \rmult w_1 - \alpha'_2 \skm
                   (\sum_{i=1}^n h_{i} \rmult l_{i}) \label{s}
     \end{eqnarray}
     where $\beta_1=\hc(f_1 \rmult w_1)^{-1}$, $\beta_2=\hc(f_2 \rmult w_2)^{-1}$
      and  $\alpha'_2 \skm \beta_2 = \alpha_2$.
     By substituting (\ref{s}) in our representation of $g$ it becomes smaller.
\\
\qed
Notice that both test sets in this characterization in general are not finite.

Remember that in commutative polynomial rings over fields we can
 restrict these critical situations to one s-polynomial arising from
 the least common multiple of the head terms $\hterm(p)$ and $\hterm(q)$.
Here we can introduce a similar concept of least common multiples,
 but now two terms can have no, one, finitely many and even infinitely
 many such multiples.

Given two non-zero polynomials $p$ and $q$ in $\f_{\myk}$ let
 $S(p,q) = \{ t \mid \mbox{ there exist } u_1,u_2 \in \myt \mbox{ such that }
   \hterm(p \rmult u_1) = \hterm(\hterm(p) \rmult u_1) = t =
   \hterm(q \rmult u_2) = \hterm(\hterm(q) \rmult u_2)
   \mbox{ and } \hterm(p) \leq t, \hterm(q) \leq t \}$.
A subset $LCM(p,q)$ of $S(p,q)$ is called a set of least common multiples for
 \label{page.lcm}
 $p$ and $q$ if for any $t \in S(p,q)$ there exists $t' \in LCM(p,q)$
 such that $t' \leq t$ and all other $s \in LCM(p,q)$ are not comparable with $t'$
 with respect to the reductive ordering $\leq$.

For polynomial rings over fields  a term $t$
 is smaller than another term $s$ with respect to the reductive
 ordering if $t$ is a divisor of $s$
 and $LCM(p,q)$ consists of the least common multiple
 of the head terms $\hterm(p)$ and $\hterm(q)$.
But for function rings in general other situations are possible.
Two polynomials do not have to give rise to
 any s-polynomial.
Just take $\myt$ to be the free monoid on $\{ a,b\}$ and $\myk = \q$.
Then for the two polynomials $p = a+1$ and $q = b+1$ we have
 $S(p,q) = \emptyset$ as there are no terms $u_1,u_2$ in $\myt$ such that
 $a \rmult u_1 = b \rmult u_2$.

Next we give an example where the set $LCM(p,q)$ is finite but larger that one
 element.
\begin{example}~\\
{\rm
Let our set of terms $\myt$ be presented as a monoid by
 $(\{a,b,c,d_1,d_2,x_1,x_2\};\{ ax_i = cx_i, bx_i = cx_i, 
  d_jx_i = x_id_j \mid i,j \in \{1,2\}\})$, 
 $\succeq$ is the length-lexicographical ordering induced by the precedence 
 $x_2 \succ x_1 \succ a \succ b \succ c \succ d_1 \succ d_2$ and
 the reductive ordering $\geq$ is the prefix ordering.
Then for the two polynomials $p = a + d_1$ and $q = b + d_2$ we get the
 respective sets $S(p,q) = \{ cx_1w, cx_2w \mid w \in \myt \}$
 and $LCM(p,q) = \{ cx_1,cx_2 \}$ with resulting s-polynomials
 $\spol{r}(p,q,cx_1,x_1,x_1) = x_1d_1 - x_1d_2$ and
 $\spol{r}(p,q,cx_2,x_2,x_2) = x_2d_1 - x_2d_2$.
\exaend
}
\end{example}
It is also possible to have infinitely many least common multiples.
\begin{example}~\\
{\rm
Let our set of terms $\myt$ be presented as a monoid by
 $(\{a,b,c,d_1,d_2,x_i \mid i \in \n\};\{ ax_i = cx_i, bx_i = cx_i, 
  d_jx_i = x_id_j \mid i \in \n,j \in \{1,2\}\})$, 
 $\succeq$ is the length-lexicographical ordering induced by the precedence 
 $\ldots \succ x_n \succ \ldots \succ x_1 \succ a \succ b \succ c \succ d_1 \succ d_2$ and
 the reductive ordering $\geq$ is the prefix ordering.
Then for the two polynomials $p = a + d_1$ and $q = b + d_2$ we get the
 respective set $S(p,q) = \{ cx_iw \mid i \in \n, w \in \myt \}$
 and the infinite set
 $LCM(p,q) = \{ cx_i \mid i \in \n \}$ with infinitely many resulting s-polynomials
 $\spol{r}(p,q,cx_i,x_i,x_i) = x_id_1 - x_id_2$.
\exaend
}
\end{example}
However, we have to show that we can restrict the set ${\cal C}_{s}(p,q)$
 to those tuples corresponding to terms in $LCM(p,q)$.

Remember that one
 problem which is related to the fact that the ordering $\succeq$
 and the multiplication $\rmult$ in general are not compatible
 is that an important property fulfilled for representations of 
 polynomials in commutative polynomial rings over fields no longer holds.
This property in fact underlies Lemma \ref{lem.buchberger.confluence} (4),
 which is essential in Buchberger's characterization of Gr\"obner bases in polynomial rings:
 $p \red{*}{\myr}{b}{F} 0$ implies $p \rmult m \red{*}{\myr}{b}{F} 0$ for
 any monomial $m$.
Notice that $p \red{*}{\myr}{b}{F} 0$ implies that $p$ has a
 standard representation with  respect to $F$,
 say $\sum_{i=1}^n f_i \rmult m_i$, and
 it is easy to see that then $\sum_{i=1}^n f_i \rmult m_i \rmult m$ is a
 standard representation of $p \rmult m$ with respect to $F$.
This lemma is central in localizing {\em all} the critical situations related to two
 polynomials to the {\em one} s-polynomial resulting from the
 least common multiple of the respective head terms.

Unfortunately, neither the lemma nor its implication for the existence of
 the respective standard representations holds in our more general setting.
There, if $g \in \ideal{r}{}(F)$ has a right reductive standard  representation
 $g = \sum_{i=1}^n f_i \rmult m_i$, then the sum 
 $\sum_{i=1}^n f_i \rmult m_i \rmult m$ in general is {\em no} right reductive
  standard representation not even a right standard representation
 of the multiple $g \rmult m$ for $m \in \monoms(\f_{\myk})$.
Even while $g \in \ideal{r}{}(\{g\})$ has the trivial right reductive standard 
 representation $g = g$, the multiple $g \rmult m$ is in general {\em no}
 right reductive standard  representation of the function $g \rmult m$ for 
 $m \in \monoms(\f_{\myk})$.
Recall the example on page \pageref{page.not.stable} where for $g = x+1$
 we have $\hm(g \rmult x) = x$ while $\hm(g) \rmult x = 1$ as
 $x \rmult x = 1$.
Similarly, while $g \red{}{\myr}{}{g} 0$ must hold for any
 reduction relation, this no longer will imply $g \rmult m \red{*}{\myr}{}{g} 0$.
To see this let us review Example \ref{exa.free.group}:
 For $g = a+ 1$ and $m=b$ we get the multiple $g \rmult m = 
 (a+ 1) \rmult b = 1 + b$, but 
 $\hterm(g  \rmult m) = b \neq 1 = \hterm(\hterm(g) \rmult m)$.
Moreover, $b+ 1$ is not reducible by $a+ 1$ for any reduction
 relation based on head monomial divisibility. 

In order to give localizations of the test sets from Theorem \ref{theo.rrsb=rgb.k}
 it is important to study under which conditions
 the stability of right reductive standard  representations
 with respect to multiplication by monomials can be restored.
The next lemma provides a sufficient condition.

\begin{lemma}\label{lem.red.reps}~\\
{\sl
Let  $F \subseteq \f_{\myk} \backslash \{ \zero\}$ and
 $p$ a non-zero polynomial in $\f_{\myk}$.
Moreover, we assume that $p$ has a right reductive standard  representation in terms of $F$ and
 $m$ is a monomial such that $ \hterm(p \rmult m) =
 \hterm(\hterm(p) \rmult m) \geq \hterm(p)$. 
Then $p \rmult m$ again has a right reductive standard representation in terms of $F$.
\lemend
}
\end{lemma}
\Ba{}~\\
Let $p = \sum_{i=1}^n f_i \rmult m_i$ with $n \in \n$, $f_i \in F$,
 $m_i \in \monoms(\f_{\myk})$
 be a right reductive standard  representation of $p$ in terms of $F$, i.e.,
 $\hterm(p) = \hterm(f_i \rmult m_i) =
 \hterm(\hterm(f_i) \rmult m_i) \geq \hterm(f_i)$, $1 \leq i \leq k$ and
 $\hterm(p) \succ  \hterm(f_i \rmult m_i)$
 for all $k+1 \leq i \leq n$. 
\\
Let us first analyze $f_j \rmult m_j \rmult m$ for $1 \leq j \leq k$:\\
Let $\terms(f_j \rmult m_j) = \{ s_1, \ldots, s_l \}$ with $s_1 \succ s_i$, $2 \leq i \leq l$,
 i.e.~$s_1 = \hterm(f_j \rmult m_j) = \hterm(\hterm(f_j) \rmult m_j) = \hterm(p)$.
Hence $\hterm(\hterm(p) \rmult m) = \hterm(s_1 \rmult m) \geq \hterm(p) = s_1$ and as $s_1 \succ s_i$,
 $2 \leq i \leq l$, by Definition \ref{def.refined.ordering} we can conclude
 $\hterm(\hterm(p) \rmult m) 
 = \hterm(s_1 \rmult m) \succ s_i \rmult m \succeq \hterm(s_i \rmult m)$
 for $2 \leq i \leq l$.
This implies $\hterm(\hterm(f_j \rmult m_j) \rmult m) = \hterm(f_j \rmult m_j \rmult m)$.
Hence we get
\begin{eqnarray*}
 \hterm(p \rmult m) & = & \hterm(\hterm(p) \rmult m) \\
                    & = & \hterm(\hterm(f_j \rmult m_j) \rmult m), \mbox{ as } 
                          \hterm(p) = \hterm(f_j \rmult m_j) \\
                    & = & \hterm(f_j \rmult m_j \rmult m)
\end{eqnarray*}
and
since $\hterm(p \rmult m) \geq \hterm(p) \geq \hterm(f_j)$ we can conclude
 $\hterm(f_j \rmult m_j \rmult m) \geq \hterm(f_j)$.
It remains to show that $f_j \rmult m_j \rmult m$ has
 a right reductive standard representation in terms of $F$.
First we show that $\hterm(\hterm(f_j) \rmult m_j \rmult m) \geq \hterm(f_j)$:
We know $\hterm(f_j) \rmult m_j \succeq 
 \hterm(\hterm(f_j) \rmult m_j) = \hterm(f_j \rmult m_j)$
%\footnote{Notice that
% $\hterm(f_j) \rmult m_j$ can be a polynomial and hence we cannot conclude 
% $\hterm(f_j) \rmult m_j = \hterm(\hterm(f_j) \rmult m_j)$. But for 
% $s \in \terms(\hterm(f_j) \rmult m_j) \backslash
%  \{ \hterm(\hterm(f_j) \rmult m_j)\}$ we know $\hterm(\hterm(f_j) \rmult m_j)
%  \succ s$. Then $\hterm(\hterm(\hterm(f_j) \rmult m_j) \rmult m) \geq
%  \hterm(\hterm(f_j) \rmult m_j)$ implies $\hterm(f_j) \rmult m_j \rmult m
%  \succeq \hterm(\hterm(f_j) \rmult m_j) \rmult m
%  \succ s \rmult m$.}
 and hence
 $\hterm(\hterm(f_j) \rmult m_j \rmult m) = \hterm(\hterm(f_j \rmult m_j) \rmult m) =
  \hterm(f_j \rmult m_j \rmult m) \geq \hterm(f_j)$.
\\
%Weiter
Now in case $m_j \rmult m \in \monoms(\f_{\myk})$ we are done as then $f_j \rmult ( m_j \rmult m)$
 is a right reductive standard representation in terms of $F$.
\\
Hence let us assume
 $m_j \rmult m = \sum_{i=1}^{k} \tilde{m}_i$, $\tilde{m}_i \in \monoms(\f_{\myk})$.
Let $\terms (f_j) = \{ t_1, \ldots, t_s \}$ with $t_1 \succ t_p$,  $2 \leq p \leq s$,
 i.e.~$t_1 = \hterm(f_j)$.
As $\hterm( \hterm(f_j) \rmult m_j) \geq \hterm(f_j) \succ t_p$,$2 \leq p \leq s$,
 again by Definition \ref{def.refined.ordering}
 we can conclude $\hterm(\hterm(f_j) \rmult m_j) \succ t_p \rmult m_j \succeq \hterm(t_p \rmult m_j)$,
 and $\hterm(f_j) \rmult m_j \succ \sum_{p=2}^s t_p \rmult m_1$.
Then for each $s_i$, $2 \leq i \leq l$
 there exists $t_p \in \terms(f_j)$ such that
 $s_i \in \supp(t_p \rmult m_j)$.
Since $\hterm(p) \succ s_i$ and even\footnote{$\hterm(p) \succ t_p \rmult m_j$
 if $\hterm(f_j \rmult m_j) \not\in \supp(t_p \rmult m_j)$.}
 $\hterm(p) \succeq t_p \rmult m_j$ we find that
 either
 $\hterm(p \rmult m) \succeq \hterm((t_p \rmult m_j) \rmult m) = \hterm(t_p \rmult (m_j \rmult m))$
 in case $\hterm(t_p \rmult m_j ) = \hterm(f_j \rmult m_j)$
 or 
 $\hterm(p \rmult m) \succ \hterm((t_p \rmult m_j) \rmult m) = \hterm(t_p \rmult (m_j \rmult m))$.
Hence we can conclude $f_j \rmult \tilde{m}_i \predeq \hterm(p \rmult m)$, $1 \leq i \leq l$ and
 for at least one $\tilde{m}_i$ we get
 $\hterm(f_j \rmult \tilde{m}_i) = \hterm(f_j \rmult m_j \rmult m) \geq \hterm(f_j)$.
\\
It remains to analyze the situation for the function
 $(\sum_{i=k+1}^n f_i \rmult m_i) \rmult m$.
Again we find that for all terms $s$ in the $f_i \rmult m_i$, $k+1 \leq i \leq n$,
 we have $\hterm(p) \succ s$ and we get 
 $\hterm(p \rmult m) \succ \hterm(s \rmult m)$.
Hence all polynomial multiples of the $f_i$  in the representation 
 $\sum_{i=k+1}^n \sum_{j=1}^{k_i} f_i \rmult \tilde{m}^i_j$, where
 $m_i \rmult m = \sum_{j=1}^{k_i}\tilde{m}^i_j$, are bounded by 
 $\hterm(p \rmult m)$.
\\
\qed

Notice that these observations are no longer true in case we only require
 $\hterm(p \rmult m) =
 \hterm(\hterm(p) \rmult m) \succeq \hterm(p)$, as 
 then $\hterm(p) \succ s$ no longer implies that
 $\hterm(p \rmult m) \succ \hterm(s \rmult m)$ will hold.

Of course this lemma now implies that if for two polynomials $p$ and $q$
 in $\f_{\myk}$ all s-polynomials related to the set $LCM(p,q)$ have 
 right reductive standard representations so have all s-polynomials related
 to any tuple in ${\cal C}_{s}(p,q)$.

So far we have characterized weak right reductive 
 Gr\"obner bases as special right ideal bases
 providing right reductive standard  representations for the right ideal elements.
In the literature the existence of such representations is
 normally established by means of reduction relations.
The special representations presented here can be related to a reduction relation
 based on the divisibility of terms as defined in the context of 
 right reductive restrictions of our ordering following  Definition \ref{def.refined.ordering}.
Let $\geq$ be such a right reductive restriction of the ordering $\succeq$.
\begin{definition}\label{def.rred}~\\
{\rm
Let $f,p$ be two non-zero polynomials in $\f_{\myk}$.
We say $f$ \betonen{right reduces} $p$ \betonen{to} $q$ \betonen{at a monomial}
 $\alpha \skm t$ \betonen{in one step}, denoted by $p \red{}{\myr}{r}{f} q$, if
 there exists $m \in \monoms(\f_{\myk})$ such that
\begin{enumerate}
\item $t \in \supp(p)$ and $p(t) = \alpha$,
\item $\hterm(f\rmult m) = \hterm(\hterm(f) \rmult m) = t \geq \hterm(f)$,
\item $\hm(f \rmult m) = \alpha \skm t$, and
\item $q = p - f \rmult m$.
\end{enumerate}
We write $p \red{}{\myr}{r}{f}$ if there is a polynomial $q$ as defined
above and $p$ is then called right reducible by $f$. 
%\\
Further, we can define $\red{*}{\myr}{r}{}, \red{+}{\myr}{r}{}$ and
 $\red{n}{\myr}{r}{}$ as usual.
%\\
Right reduction by a set $F \subseteq \f_{\myk}$ is denoted by
 $p \red{}{\myr}{r}{F} q$ and abbreviates $p \red{}{\myr}{r}{f} q$
 for some $f \in F$.
\dend
}
\end{definition}
Notice that if $f$ right reduces $p$ to $q$ at $ \alpha \skm t$ then
 $t \not\in \supp(q)$.
If for some $w \in \myt$ we have $\hterm(f \rmult w) =
 \hterm(\hterm(f) \rmult w) = t \geq \hterm(f)$ we can always reduce $\alpha \skm t$
 in $p$ by $f$ using the monomial $m = (\alpha \skm \hc(f \rmult w)^{-1}) \skm w$. 
Other definitions of reduction relations are possible, e.g.~substituting item 2 by the condition
 $\hterm(\hterm(f)\rmult m) = \hterm(f \rmult m)$ (called right reduction in the context 
 of monoid rings in \cite{Re95}; such a reduction relation would be 
 connected to standard representations as defined in
 Definition \ref{def.standard.rep})
or by the condition
 $\hterm(f \rmult m) = t$ (called strong reduction in the context of monoid rings in \cite{Re95} and for function rings on page \pageref{page.strong.reduction}).
We have chosen this particular reduction relation as it provides the necessary information to apply
 Lemma \ref{lem.red.reps} to give localizations for the test sets in Theorem \ref{theo.rrsb=rgb.k}
 later on.
Let us continue by studying some of the properties of our reduction relation.
\begin{lemma}\label{lem.sred}~\\ 
{\sl 
Let $F$ be a set of polynomials in $\f_{\myk} \backslash \{ \zero\}$.
\begin{enumerate}
\item For $p,q \in \f_{\myk}$, $p \red{}{\myr}{r}{f \in F} q$ implies $p \succ q$, in particular $\hterm(p)
  \succeq \hterm(q)$.
\item $\red{}{\myr}{r}{F}$ is Noetherian.
\lemend
\end{enumerate}
}
\end{lemma}
\Ba{}
\begin{enumerate}
\item Assuming that the reduction step takes place at a monomial $\alpha \skm t$,
      by Definition \ref{def.rred} we know $\hm(f \rmult m) = \alpha \skm t$ which yields
      $p \succ p  -  f \rmult m$
      since $\hm(f \rmult m) \succ \reductum(f \rmult m)$.
\item This follows directly from 1. as the ordering $\succeq$ on $\myt$ is well-founded
 (compare Theorem \ref{lem.wellfounded}).
\end{enumerate}\renewcommand{\baselinestretch}{1}\small\normalsize
\qed
The next lemma shows how reduction sequences and right reductive standard
 representations are related.
\begin{lemma}\label{lem.red.rep}~\\
{\sl
Let $F \subseteq \f_{\myk} \backslash \{ \zero\}$ and $p \in \f_{\myk} \backslash \{ \zero\}$.
Then $p \red{*}{\myr}{r}{F} \zero$ implies that $p$ has a  right reductive standard
 representation in terms of $F$.
\lemend
}
\end{lemma}
\Ba{}~\\
This follows directly by adding up the polynomials used in the
 reduction steps occurring in the reduction sequence $p \red{*}{\myr}{r}{F} \zero$,
 say $p \red{}{\myr}{r}{f_1} p_1 \red{}{\myr}{r}{f_2} \ldots \red{}{\myr}{r}{f_n}
 \zero$.
If the reduction steps take place at the respective head monomials only, we can
 additionally state that $ p = \sum_{i=1}^n f_i \rmult m_i$,
 $\hterm(f_i \rmult m_i) = \hterm(\hterm(f_i) \rmult m_i) \geq \hterm(f_i)$,
 $1 \leq i \leq n$, and even $\hterm(f_1 \rmult m_1) \succ 
 \hterm(f_2 \rmult m_2) \succ \ldots \hterm(f_n \rmult m_n)$.
\\
\qed
If $p \red{*}{\myr}{r}{F} q$, then $p$ has a right reductive standard  representation in 
 terms of $F \cup \{ q \}$, respectively $p-q$ has  a right reductive standard  representation in 
 terms of $F$.
On the other hand, if a polynomial $g$ has a right reductive standard representation
 in terms of some set $F$ it is reducible by a polynomial in $F$.
To see this let $g = \sum_{i=1}^n f_i \rmult m_i,
 f_i \in F, m_i \in \monoms(\f_{\myk}), n \in \n$ be a right reductive standard representation of $g$ in
 terms of $F$.
Then $\hterm(g) = \hterm(f_1 \rmult m_1) =
 \hterm(\hterm(f_1) \rmult m_1) \geq \hterm(f_1)$ 
 and by Definition \ref{def.rred} this implies
 that $g \red{}{\myr}{r}{f_1} g - \alpha \skm f_1 \rmult m_1 = g'$ where 
 $\alpha \in \myk$ such that $\alpha \skm \hc(f_1 \rmult m_1) = \hc(g)$.

So far we have given an algebraic characterization of {\em weak}
 right reductive Gr\"obner bases in Definition
 \ref{def.gb.fr} and a characterization of them as right reductive standard bases in
 Lemma \ref{lem.sb=gb}.
Another characterization known from the literature is that for a Gr\"obner basis in 
 a polynomial ring every element of the ideal it generates reduces to zero using the
 Gr\"obner basis.
Reviewing Definition \ref{def.weak.gb.reduction.ring} we find that this is in fact only the
 definition of a weak Gr\"obner basis.
However in polynomial rings over fields and many other structures in the literature
 the definitions of weak Gr\"obner bases and Gr\"obner bases coincide as
 the Translation Lemma holds (see Lemma \ref{lem.buchberger.confluence} (2)).
This is also true for function rings over fields.

The first part of the following lemma is only needed for the proof of the second
 part which is an analogon of the Translation Lemma for function rings over fields.

\begin{lemma}\label{lem.trans.fr}~\\
{\sl
Let $F$ be a set of polynomials in $\f_{\myk}$ and $p,q,h$ polynomials in $\f_{\myk}$.
\begin{enumerate}
\item
Let $p-q \red{}{\myr}{r}{F} h$.
Then there exist $p', q' \in \f_{\myk}$ such that $p \red{*}{\myr}{r}{F} p'$
 and $q \red{*}{\myr}{r}{F} q'$ and $h = p'-q'$.
\item
Let $\zero$ be a normal form of $p-q$ with respect to $F$.
Then there exists $g \in \f_{\myk}$ such that $p \red{*}{\myr}{r}{F} g$ and
 $q \red{*}{\myr}{r}{F} g$.
\lemend
\end{enumerate}
}
\end{lemma}
\Ba{}~\\
\begin{enumerate}
\item
Let $p-q \red{}{\myr}{r}{F} h$ at the monomial $\alpha \skm t$, i.e.,
 $h = p-q - f \rmult m$ for some $f \in F$,$m \in \monoms(\f_{\myk})$ such that
 $\hterm(\hterm(f)\rmult m) = \hterm(f \rmult m) = t \geq \hterm(f)$
 and $\hm(f \rmult m) = \alpha \skm t$, i.e., $\alpha$ is the coefficient
 of $t$ in $p-q$.
We have to distinguish three cases:
\begin{enumerate}
\item $t \in \supp(p)$ and $t \in \supp(q)$:
      Then we can eliminate the occurrence of $t$ in the respective polynomials
        by right reduction and get $p \red{}{\myr}{r}{f} p - \alpha_1 \skm f \rmult m
        = p'$, $q \red{}{\myr}{r}{f} q - \alpha_2 \skm f \rmult m = q'$, 
        where $\alpha_1 \skm \hc(f \rmult m)$
        and $\alpha_2 \skm \hc(f \rmult m)$ are the coefficients of $t$ in
        $p$ respectively $q$.
      Moreover, $\alpha_1 \skm \hc(f \rmult m) - \alpha_2 \skm \hc(f \rmult m) = 
       \alpha$ and hence $\alpha_1 - \alpha_2 = 1$, as $\hc(f \rmult m) = \alpha$.
      This gives us $p' - q' = p - \alpha_1 \skm f \rmult m - 
        q + \alpha_2 \skm f \rmult m =
        p - q - (\alpha_1 - \alpha_2) \skm f \rmult m =
        p-q - f \rmult m = h$.
\item $t \in \supp(p)$ and $t \not\in \supp(q)$:
      Then we can eliminate the term $t$ in the polynomial $p$ by right reduction
        and get $p \red{}{\myr}{r}{f} p - f \rmult m = p'$, $q=q'$,
        and, therefore, $p' - q' = p - f \rmult m - q =h$.
\item $t \in \supp(q)$ and $t \not\in \supp(p)$:
      Then we can eliminate the term $t$ in the polynomial $q$ by right reduction
        and get $q \red{}{\myr}{r}{f} q + f \rmult m = q'$, $p=p'$, and, therefore,
        $p' - q' = p - (q + f \rmult m) =h$.
\end{enumerate}
\item We show our claim by induction on $k$, where $p-q \red{k}{\myr}{r}{F} \zero$.
In the base case $k=0$ there is nothing to show as then $p=q$.
Hence, let $p-q \red{}{\myr}{r}{F} h \red{k}{\myr}{r}{F} \zero$.
Then by 1.~there are polynomials $p',q' \in \f_{\myk}$ such that 
 $p \red{*}{\myr}{r}{F} p'$
 and $q \red{*}{\myr}{r}{F} q'$ and $h = p'-q'$.
Now the induction hypothesis for $p'-q' \red{k}{\myr}{r}{F} \zero$ yields the
 existence of a polynomial $g \in \f_{\myk}$ such that
 $p \red{*}{\myr}{r}{F} g$ and $q \red{*}{\myr}{r}{F} g$.
\end{enumerate}
\qed 

The essential part of the proof is that right reducibility is connected to
 stable divisors of terms.
We will later see that for function rings over arbitrary reduction rings, when the coefficient
 is also involved in the reduction step, this lemma no longer holds.
\begin{definition}\label{def.gb.reduction}~\\
{\rm
A subset $G$ of  $\f_{\myk}$ is called a 
 \betonen{right  Gr\"obner basis} (with respect to the reduction relation
 $\red{}{\myr}{r}{}$) of the right 
 ideal ${\mathfrak i}= \ideal{r}{}(G)$ it generates, if
 $\red{*}{\longleftrightarrow}{r}{G} = \;\;\equiv_{{\mathfrak i}}$
 and $\red{}{\myr}{r}{G}$ is confluent.
}
\end{definition}
Recall the free group ring in Example \ref{exa.free.group}.
There the polynomial $b + 1$ lies
 in the right ideal generated by the polynomial $a + 1$.
Unlike in the case of polynomial rings over fields where for any set
 of polynomials $F$ we have
 $\red{*}{\longleftrightarrow}{b}{F} = \;\;\equiv_{\ideal{}{}(F)}$,
 here we have
 $b + 1 \equiv_{\ideal{r}{}(\{a + 1 \})} \zero$ but
 $b + 1 \nred{*}{\longleftrightarrow}{r}{a + 1} \zero$.
Hence the first condition of Definition \ref{def.gb.reduction} now
 becomes necessary while it can be omitted in the definition
 of Gr\"obner bases for ordinary polynomial rings.

Now by Lemma \ref{lem.trans.fr} and Theorem \ref{theo.weak+translation}
 weak right reductive Gr\"obner bases are right Gr\"obner bases
 and can be characterized as follows:

\begin{corollary}\label{lem.gb=red.to.zero}~\\
{\sl
Let $G$ be a set of polynomials in $\f_{\myk} \backslash \{ \zero\}$.
$G$ is a right Gr\"obner basis if and only if for every $g \in \ideal{r}{}(G)$
 we have $g\red{*}{\myr}{r}{G} \zero$.  
}
\end{corollary}

Finally we can characterize right Gr\"obner bases similar to Theorem 
 \ref{theo.buchberger.completion}.

\begin{theorem}\label{theo.s-pol}~\\
{\sl
Let $F$ be a set of polynomials in $\f_{\myk} \backslash \{ \zero \}$.
Then $F$ is a right Gr\"obner basis if and only if
\begin{enumerate}
\item for all $f$ in $F$ and for all $m \in \monoms(\f_{\myk})$ we have
       $f \rmult m \red{*}{\myr}{r}{F} \zero$, and
\item for all  $p$ and $q$ in $F$ and every tuple $(t,u_1,u_2)$
      in ${\cal C}_{s}(p,q)$ and the respective s-polynomial $\spol{r}(p,q,t,u_1,u_2)$
       we have $\spol{r}(p,q,t,u_1,u_2) \red{*}{\myr}{r}{F} \zero$.
\end{enumerate}
\theoend
}
\end{theorem}
However, the importance of Gr\"obner bases in the classical case stems
 from the fact that we only have to check a finite set of s-polynomials for $F$
 in order to decide, whether $F$ is a Gr\"obner basis.
Hence, we are interested in localizing the test sets in Theorem \ref{theo.s-pol} --
 if possible to finite ones.
\begin{definition}\label{def.weakly.saturated}~\\
{\rm
A set of polynomials $F \subseteq \f_{\myk} \backslash \{ \zero\}$ is called
 \betonen{weakly saturated}, if for every monomial $m \in \monoms(\f_{\myk})$
 and every polynomial
 $f$ in $F$ we have
 $f \rmult m \red{*}{\myr}{r}{F} \zero$.
\dend
}
\end{definition}
Then for a weakly saturated set $F$ and any monomial $m \in \monoms(\f_{\myt})$,
 $f \in F$, the multiple $f \rmult m$ has a right reductive standard representation
 in terms of $F$.
Notice that since the coefficient domain is a field and $\f$ a vector
 space we can even restrict
 ourselves to multiples with elements of $\myt$.
However, for reduction rings as coefficient
 domains, we will need monomial multiples and hence we give the more general definition.
For the free group ring in Example \ref{exa.free.group} the set $\{ a+ 1, b+ 1 \}$ is
 weakly saturated.
\begin{definition}\label{def.saturator}~\\
{\rm
Let $F$ be a set of polynomials in $\f_{\myk} \backslash \{ 0 \}$.
A set
 $\sat(F) \subseteq \{ f \rmult m \mid f \in F, m \in \monoms(\f_{\myk}) \}$
 is called a \betonen{stable saturator} for $F$ if for any $f \in F$, $m \in \monoms(\f_{\myk})$
 there exist $s \in \sat(F)$, $m' \in\monoms(\f_{\myk})$ such that
 $f \rmult m = s \rmult m'$ and
 $\hterm(f \rmult m) = \hterm(\hterm(s) \rmult m')
 \geq \hterm(s)$. 
}
\end{definition}
Notice that a stable saturator need not be weakly saturated.
Let $s \in \sat(F) \subseteq \{ f \rmult m \mid f \in F, m \in \monoms(\f_{\myk}) \}$ and $m'\in \monoms(\f_{\myk})$.
For $\sat(F)$ to be weakly saturated then $s \rmult m' \red{*}{\myr}{}{\sat(F)} \zero$ must hold.
We know that $s = f \rmult m$ for some $f \in F, m \in \monoms(\f_{\myk})$.
In case $m \rmult m' \in \monoms(\f_{\myk})$ we are done.
But this is no longer true if the monomial multiple results in a polynomial.
Let our set of terms consist of words on the alphabet $\{a,b,c\}$ with
 multiplication $\rmult$ deduced form the
 equations $a\rmult b=a, b\rmult a=b^2-b,a \rmult a=\zero$.
As ordering on $\myt$ we take the length lexicographical ordering with precedence
 $a \succ b \succ c$ and as reductive restriction the prefix ordering.
For the polynomial $f = ca + 1$ we get a stable saturator
 $\sat(\{f\})= \{ ca + 1, ca + b, ca + b^2, b^3 + ca, a \}$.
Then the polynomial multiple $(f \rmult b) \rmult a = f \rmult (b \rmult a) =
 f \rmult (b^2-b) = ca + b^2 - (ca +b) = b^2 -b$
 is not reducible by $\sat(\{f\})$ while $f \rmult b = ca + b \in \sat(\{f\})$.
%
%f * b = ca +b
%f * b^2 = ca + b^2
%f * b^3 = ca + b^3 nicht reduzibel!
%Weitere Vielfache mit Potenzen von b mit b^3 + ca reduzibel
%(ca + b^3) \rmult a = b^4-b^3 -> mit b^3+ca nach b^4 +ca -> mit (b^3+ca)*b nach 0
%Eigentlich auch nicht zu betrachten, da kein Vielfaches der Form f * m sondern
%der Form f * p, p Polynom
%gleiches gilt fuer f * b * a = f * (b^2-b), welches nicht nach 0 geht!!!!!!
%f * a = a
%
\begin{corollary}~\\
{\sl
Let $\sat(F)$ be a stable saturator of a set $F \subseteq \f_{\myk}$.
Then for any $f \in F$, $m \in \monoms(\f_{\myk})$ there exists $s \in \sat(F)$
 such that $f \rmult m \red{}{\myr}{r}{s} \zero$.
}
\end{corollary}
\begin{lemma}~\\
{\sl
Let $F$ be a set of polynomials in $\f_{\myk} \backslash \{ 0 \}$.
If for all $s$ in a stable saturator $\sat(F)$ we have $ s \red{*}{\myr}{r}{F} \zero$, then
 for every $m$ in $\monoms(\f_{\myk})$
 and every polynomial $f$ in $F$ the right multiple
 $f \rmult m$ has a right reductive standard  representation in terms of $F$.
}
\end{lemma}
\Ba{}~\\
This follows immediately from Lemma \ref{lem.red.rep} and
 Lemma \ref{lem.red.reps}.
\\
\qed
\begin{definition}\label{def.s-poly.localization}~\\
{\rm
Let $p$ and $q$ be two non-zero polynomials in $\f_{\myk}$.
Then a subset $C \subseteq \{ \spol{r}(p,q,t,u_1,u_2) \mid
 (t,u_1,u_2) \in {\cal C}_{s}(p,q) \}$ is called a \betonen{stable
 localization} for the critical situations if for every s-polynomial
 $\spol{r}(p,q,t,u_1,u_2)$ related to a tuple $(t,u_1,u_2)$
 in ${\cal C}_{s}(p,q)$ there exists a polynomial $h \in C$ and
 a monomial $m \in \monoms(\f_{\myk})$ such that
 \begin{enumerate}
 \item $\hterm(h) \leq \hterm(\spol{r}(p,q,t,u_1,u_2))$,
 \item $\hterm(h\rmult m) = \hterm(\hterm(h)\rmult m) =
       \hterm(\spol{r}(p,q,t,u_1,u_2))$,
 \item $\spol{r}(p,q,t,u_1,u_2) = h \rmult m$.
\dend
 \end{enumerate}
}
\end{definition}
The set $LCM(p,q)$ (see page \ref{page.lcm}) allows a stable localization
 as follows: $C = \{ \spol{r}(p,q,t,u_1,u_2) \mid t \in LCM(p,q), (t, u_1,u_2)
 \in {\cal C}_{s}(p,q)\}$.
\begin{corollary}~\\
{\sl
Let $C \subseteq \{ \spol{r}(p,q,t,u_1,u_2) \mid
 (t,u_1,u_2) \in {\cal C}_{s}(p,q) \}$ be a stable localization for
 two polynomials $p,q \in \f_{\myk}$.
Then for any s-polynomial $\spol{r}(p,q,t,u_1,u_2)$ there exists $h \in C$
 such that $\spol{r}(p,q,t,u_1,u_2) \red{}{\myr}{r}{h} \zero$.
}
\end{corollary}
\begin{lemma}~\\
{\sl
Let $F$ be a set of polynomials in $\f_{\myk} \backslash \{ 0 \}$.
If for all $h$ in 
      a stable localization  $C \subseteq \{ \spol{r}(p,q,t,u_1,u_2) \mid
    (t,u_1,u_2) \in {\cal C}_{s}(p,q) \}$,
  we have $h \red{*}{\myr}{r}{F} \zero$, then
 for every $(t,u_1,u_2)$
 in ${\cal C}_{s}(p,q)$ the s-polynomial  $\spol{r}(p,q,t,u_1,u_2)$
 has a right reductive standard  representation in terms of $F$.
}
\end{lemma}
\Ba{}~\\
This follows immediately from Lemma \ref{lem.red.rep} and
 Lemma \ref{lem.red.reps}.
\\
\qed
So far we have seen that basically the theory for right Gr\"obner bases and the refined
 notion of right reductive standard bases (for right ideals of course)
 carries over similar from the case of polynomial rings over fields.
Now Lemma \ref{lem.red.reps} and
 Lemma \ref{lem.red.rep} allow a localization of the test situations
 from Theorem \ref{theo.s-pol}.

\begin{theorem}\label{theo.s-pol.2}~\\
{\sl
Let $F$ be a set of polynomials in $\f_{\myk} \backslash \{ 0 \}$.
Then $F$ is a right Gr\"obner basis if and only if
\begin{enumerate}
\item for all $s$ in a stable saturator $\sat(F)$  we have
       $s \red{*}{\myr}{r}{F} \zero$, and
\item for all  $p$ and $q$ in $F$, and every polynomial $h$ in 
      a stable localization  $C \subseteq \{ \spol{r}(p,q,t,u_1,u_2) \mid
      (t,u_1,u_2) \in {\cal C}_{s}(p,q) \}$,
      we have $h \red{*}{\myr}{r}{F} \zero$.
%      for all $h \in \sat(s)$
%      we have $h \red{*}{\myr}{r}{F} \zero$.
\end{enumerate}
\theoend
}
\end{theorem}
\Ba{}~\\
In case $F$ is a right Gr\"obner basis by Lemma \ref{lem.gb=red.to.zero}
  all elements of $\ideal{r}{}(F)$ must right reduce to zero by $F$.
Since the polynomials in question all belong to the right ideal generated
 by $F$ we are done.
\\
The converse will be proven by showing that every element in
 $\ideal{r}{}(F)$ has a right reductive representation in terms of $F$.
Now, let $g = \sum_{j=1}^m f_{j} \rmult m_{j}$ be an arbitrary
 representation of a non-zero  polynomial $g\in \ideal{r}{}(F)$ such that
 $f_j \in F$, and $m_{j} \in \monoms(\f_{\myk})$.
\\
By our first assumption for every multiple $f_j \rmult m_j$ in this
 sum we have some $s \in \sat(F)$, $m \in \monoms(\f_{\myk})$ such that $f_{j} \rmult m_{j} = s \rmult m$
 and $\hterm(f_{j} \rmult m_{j}) = \hterm(s \rmult m) = \hterm(\hterm(s) \rmult m) \geq \hterm(s)$.
Since we have $s \red{*}{\myr}{r}{F} \zero$, by Lemma
 \ref{lem.red.reps} we can conclude that each $f_{j} \rmult m_{j}$ has a right reductive
 standard representation in terms of $F$.
Therefore, we can assume that
 $\hterm(\hterm(f_j) \rmult m_j) = \hterm(f_j \rmult m_j) \geq
 \hterm(f_j)$ holds.
\\
Depending on this  representation of $g$ and the
 well-founded total ordering $\succeq$ on $\myt$ we define
 $t = \max_{\succeq} \{ \hterm(f_{j} \rmult m_{j}) \mid 1\leq j \leq m \}$ and
 $K$ as the number of polynomials $f_j \rmult m_j$ with head term $t$.
\\
Without loss of generality we can assume that the polynomial multiples
 with head term $t$ are just $f_1  \rmult m_1, \ldots , f_K \rmult m_K$.
We proceed by induction
 on $(t,K)$, where
 $(t',K')<(t,K)$ if and only if $t' \prec t$ or $(t'=t$ and
 $K'<K)$\footnote{Note that this ordering is well-founded since $\succ$
                  is well-founded on $\myt$ and $K \in\n$.}.
Obviously, $t \succeq \hterm(g)$ must hold. 
If $K = 1$ this gives us $t = \hterm(g)$ and by our assumption our
 representation is already of the required form.
\\
Hence let us assume $K > 1$, then for the  two
 not necessarily different polynomials $f_1,f_2$
 and corresponding monomials $m_1 = \alpha_1 \skm w_1$, $m_2 = \alpha_2 \skm w_2$,
 $\alpha_1,\alpha_2 \in \myk$, $w_1,w_2 \in \myt$,
 in the corresponding   representation we have
       $t=\hterm(\hterm(f_1) \rmult w_1) =\hterm(f_1 \rmult w_1) =
       \hterm(f_2 \rmult w_2) =  \hterm(\hterm(f_2) \rmult w_2)$ and
       $t \geq \hterm(f_1)$, $t \geq \hterm(f_2)$.
Then the tuple $(t, w_1,w_2)$ is in ${\cal C}_{s}(f_1,f_2)$ and
      we have a polynomial $h$ in a stable localization 
      $C \subseteq \{ \spol{r}(f_1,f_2,t,w_1,w_2)
       \mid (t, w_1, w_2) \in {\cal C}_{s}(f_1,f_2) \}$
      and $m \in \monoms(\f_{\myk})$ such that 
      $\spol{r}(f_1,f_2,t,w_1,w_2)= \hc(f_1 \rmult w_1)^{-1} \skm  f_1
      \rmult w_1 -\hc(f_2 \rmult w_2)^{-1}\skm f_2 \rmult w_2
      = h \rmult m$ and 
       $\hterm(\spol{r}(f_1,f_2,t,w_1,w_2)) =
        \hterm(h\rmult m) = \hterm(\hterm(h)\rmult m) \geq \hterm(h)$.
%      we have a polynomial $s$ in a stable localization 
%      $C \subseteq \{ \spol{r}(f_1,f_2,t,w_1,w_2) \mid (t, w_1, w_2) \in {\cal C}(f_1,f_2) \}$
%      and $\alpha \skm w \in \monoms(\f_{\k})$ such that 
%      $\spol{r}(f_1,f_2,t,w_1,w_2)= \hc(f_1 \rmult w_1)^{-1} \skm  f_1
%      \rmult w_1 -\hc(f_2 \rmult w_2)^{-1}\skm f_2 \rmult w_2
%      = s \rmult (\alpha \skm w)$.
%Hence there is $h \in \sat(s)$ such that
% $s \rmult (\alpha \skm w) = h \rmult l$ for some $l \in \monoms(\f_{\myk})$ and
% $\hterm(s \rmult (\alpha \skm w)) = \hterm(h \rmult l) =
% \hterm(\hterm(h) \rmult l) \geq \hterm(h)$.
\\
We will now change our representation of $g$ by using the additional
 information on this situation in such a way that for the new
 representation of $g$ we either have a smaller maximal term or
 the occurrences of the term $t$
 are decreased by at least 1.
%\\
Let us assume the s-polynomial is not $\zero$\footnote{In case 
               $h =\zero$,
               just substitute the empty sum
               for the right reductive representation of $h$
               in the equations below.}. 
%\\
By our assumption, $h \red{*}{\myr}{r}{F} \zero$ 
 and by Lemma \ref{lem.red.rep} $h$ then 
 has a right reductive standard  representation in terms of $F$. 
Then by Lemma \ref{lem.red.reps} the multiple $h \rmult m$ again
 has a right reductive standard representation in terms of $F$, say
  $\sum_{i=1}^n h_i \rmult l_i$, 
  where $h_i \in F$, and $l_i \in \monoms(\f_{\myk})$ and all terms occurring in this sum are bounded by
  $t \succ \hterm(h \rmult m)$.
%\\
This gives us: 
     \begin{eqnarray}
       &   & \alpha_1 \skm f_1 \rmult w_1 + \alpha_2 \skm f_2 \rmult w_2
             \nonumber\\  
       &   &  \nonumber\\                                                         
       & = &  \alpha_1 \skm f_1 \rmult w_1 +
              \underbrace{ \alpha'_2 \skm \beta_1 \skm f_1 \rmult w_1
                   - \alpha'_2 \skm \beta_1 \skm f_1 \rmult w_1}_{=\, 0} 
                   + \underbrace{\alpha'_2\skm \beta_2 }_{= \alpha_2} \skm f_2 \rmult w_2 \nonumber\\
       &   & \nonumber\\ 
       & = & (\alpha_1 + \alpha'_2 \skm \beta_1) \skm f_1 \rmult w_1 - \alpha'_2 \skm 
               \underbrace{(\beta_1 \skm f_1 \rmult w_1
             -  \beta_2 \skm f_2 \rmult w_2)}_{=\,
             h \rmult m} 
             \nonumber\\
       &   & \nonumber\\ 
       & = & (\alpha_1 + \alpha'_2 \skm \beta_1) \skm f_1 \rmult w_1 - \alpha'_2 \skm
                   (\sum_{i=1}^n h_{i} \rmult l_{i}) \label{s1}
     \end{eqnarray}
     where $\beta_1=\hc(f_1 \rmult w_1)^{-1}$, $\beta_2=\hc(f_2 \rmult w_2)^{-1}$
      and  $\alpha'_2 \skm \beta_2 = \alpha_2$.
     By substituting (\ref{s1}) our representation of $g$ becomes smaller.
\\
\qed
Obviously we now have criteria for when a set is a right Gr\"obner basis.
As in the case of completion procedures such as the Knuth-Bendix procedure or 
 the Buchberger algorithm,
 elements from these test sets which do not reduce to zero can be added to the set being tested,
 to gradually describe a not necessarily finite right Gr\"obner basis.
Of course in order to get a computable completion procedure certain assumptions on the test
 sets have to be made, e.g.~they should themselves be recursively enumerable,
 and normal forms with respect 
 to finite sets have to be computable.
Then provided such enumeration procedures for stable saturators and critical
 situations, an enumeration procedure for a respective right Gr\"obner basis
 has to ensure that all necessary candidates are enumerated and
 tested for reducibility to $\zero$.
If this is not the case they are added to the right Gr\"obner basis, have
 to be added to the enumeration of the stable saturator candidates and the new
 arising critical situations have to be added to the respective enumeration process.

%%%%%%%%%%%%%%%%%%%%%%%%%%Examples%%%%%%%%%%%%%%%%%%%%%%%%%%%%%%%%%%%%%%%%%%%%%%%%%%
We close this subsection by outlining how different structures known
 to allow finite Gr\"obner bases can be interpreted as function rings.
Using the respective interpretations the terminology can be
 adapted at once to the respective structures and in general the resulting
 characterizations of Gr\"obner bases coincide with the
 results known from literature.
%%%%%%%%%%%%%%%%%%%%%%%
\subsubsection{Polynomial Rings}\label{section.polynomials}
A commutative polynomial ring $\myk[x_1, \ldots, x_n]$ is a function
 ring according to the following interpretation:
 \begin{itemize}
 \item $\myt$ is the set of terms $\{ x_1^{i_1} \ldots x_n^{i_n} \mid i_1, \ldots, i_n \in \n \}$.
 \item $\succ$ can be any admissible term ordering on $\myt$.
       For the reductive ordering $\geq$ we have $t \geq s$ if $s$ divides $t$ as
        as term\footnote{Apel has studied another possible reductive ordering $\geq$
        where we have $t \geq s$ if $s$ is a prefix of $t$. This ordering gives
        rise to Janet bases.}.
 \item Multiplication $\rmult$ is specified by the action on terms, i.e.~$\rmult: \myt \times \myt \myr \myt$ where 
       $x_1^{i_1} \ldots x_n^{i_n} \rmult x_1^{j_1} \ldots x_n^{j_n} =
        x_1^{i_1+j_1} \ldots x_n^{i_n+j_n}$.
 \end{itemize}
We do not need the concept of weak saturation.
A stable localization of ${\cal C}_s(p,q)$ is already provided by
 the tuple corresponding to the least
 common multiple of the terms $\hterm(p)$ and $\hterm(q)$.

Since this structure is Abelian, one-sided and two-sided ideals coincide.
Buchberger's Algorithm provides an effictive procedure to compute finite
 Gr\"obner bases.
%%%%%%%%%%%%%%%%%%%%%%%%%%%%%%%%%%%%%%%%%%%%
\subsubsection{Solvable Polynomial Rings}\label{section.skewpolynomials}
According to \cite{KaWe90,Kr93}, a solvable polynomial ring 
 $\myk\{x_1, \ldots, x_n; p_{ij};c_{ij}\}$ with
 $1 \leq j < i \leq n$, $p_{ij} \in \myk[x_1, \ldots, x_n]$, $c_{ij} \in \myk^*$ is a function
 ring according to the following interpretation:
 \begin{itemize}
 \item $\myt$ is the set of terms $\{ x_1^{i_1} \ldots x_n^{i_n} \mid i_1, \ldots, i_n \in \n \}$.
 \item $\succ$ can be any admissible term ordering on $\myt$ for which
       $x_jx_i \succ p_{ij}$, $j<i$, must hold.
       For the reductive ordering $\geq$ we have $t \geq s$ if $s$ divides $t$ as
        as term.
 \item Multiplication $\rmult$ is specified by lifting the following action
       on the variables:
       $x_i \rmult x_j = x_ix_j$ if $i \leq j$ and
       $x_i \rmult x_j = c_{ij} \skm x_jx_i + p_{ij}$ if $i > j$.
 \end{itemize}
We do not need the concept of weak saturation except in case we also allow
 $c_{ij} = 0$. Then appropriate term multiples which ``delete'' head
 terms have to be taken into account. This critical set can be described
 in a finitary manner.
For the reductive ordering $\geq$ then  we can chose
  $t \geq s$ if $s$ is a prefix of $t$
 (compare Example \ref{exa.skew2}).

The set ${\cal C}_s(p,q)$ again contains as a stable localization the
 tuple corresponding to the least
 common multiple of the terms $\hterm(p)$ and $\hterm(q)$.

This structure is no longer Abelian, but finite Gr\"obner bases
 can be computed for one- and two-sided ideals (see \cite{KaWe90,Kr93}).
%%%%%%%%%%%%%%%%%%%%%%%%%%%%%%%%%%%%%%%%%%%%
%\subsubsection{Apel's Graded Structures}
%%%%%%%%%%%%%%%%%%%%%%%%%%%%%%%%%%%%%%%%%%%%
\subsubsection{Non-commutative Polynomial Rings}
A non-commutative polynomial ring $\myk[\{x_1, \ldots, x_n\}^*]$ is a function
 ring according to the following interpretation:
 \begin{itemize}
 \item $\myt$ is the set of words on $\{ x_1, \ldots, x_n\}$.
 \item $\succ$ can be any admissible ordering on $\myt$.
       For the reductive ordering $\geq$ we can chose
         $t \geq s$ if $s$ is a subword 
        of $t$.
 \item Multiplication $\rmult$ is specified by the action on words which is just concatenation.
 \end{itemize}
We do not need the concept of weak saturation.
A stable localization of ${\cal C}_s(p,q)$ is already provided by
 the tuples corresponding to word overlaps resulting from the equations $u_1\hterm(p)v_1 = \hterm(q)$, 
 $u_2\hterm(q)v_2 = \hterm(p)$, $u_3\hterm(p) = \hterm(q)v_3$ respectively $u_4\hterm(q) = \hterm(p)v_4$
 with the restriction that $|u_3| < |\hterm(q)|$ and $|u_4| < | \hterm(p) |$,
 $u_i,v_i \in \myt$.

This structure is not Abelian.
For the case of one-sided ideals finite Gr\"obner bases can be computed (see e.g.~\cite{Mo94}).
The case of two-sided ideals only allows an enumerating procedure.
This is not surprising as the word problem for monoids can be reduced to the problem
 of computing the respective Gr\"obner bases (see e.g.~\cite{Mo87,MaRe95}).
%%%%%%%%%%%%%%%%%%%%%%%%%%%%%%%%%%%%%%%%%%%%
%\subsubsection{Path Algebras}\label{section.pathalgebras}
%%%%%%%%%%%%%%%%%%%%%%%%%%%%%%%%%%%%%%%%%%%%
\subsubsection{Monoid and Group Rings}\label{section.monoidring}
A monoid or group ring $\myk[\m]$ is a function
 ring according to the following interpretation:
 \begin{itemize}
 \item $\myt$ is the monoid or group $\m$. In the cases studied by us
       as well as in \cite{Ro93,Lo96}, it is assumed that the
       elements of the monoid or group have a certain form.
       This presentation is essential in the approach.
       We will assume that the given monoid or group is presented
       by a convergent semi-Thue system.
 \item $\succ$ will be the completion ordering induced from the
       presentation of $\m$ to $\m$ and hence to $\myt$.
       The reductive ordering $\geq$ depends on the choice of the
       presentation.
 \item Multiplication $\rmult$ is specified by lifting the monoid or
       group operation.
 \end{itemize}
The concept of weak saturation and the choice of stable localizations of
 ${\cal C}_s(p,q)$ again depend on the choice of the presentation.
We will close this section by listing monoids and groups which allow
 finite Gr\"obner bases for the respective monoid or group ring
 and pointers to the literature where the appropriate solutions can be
 found.

\vspace*{1cm}
\begin{tabular}{l|l|l}
Structure & Ideals & Quote \\
\hline
\hline
Finite monoid & one- and two-sided &\cite{Re96,MaRe96a} \\
Free monoid & one-sided &\cite{Mo94,MaRe96a} \\
\hline
Finite group & one- and two-sided & \cite{Re95,MaRe96a} \\
Free group & one-sided & \cite{MaRe93b,Ro93,Re95,MaRe96a} \\
Plain group & one-sided &\cite{MaRe93b,Re95,MaRe96a} \\
Context-free group & one-sided & \cite{Re95,MaRe96a} \\
Nilpotent group &  one- and two-sided & \cite{Re95,MaRe96b} \\
Polycyclic group &  one- and two-sided &  \cite{Lo96,Re96} \\
\hline
\end{tabular}

%%% Local Variables: 
%%% mode: latex
%%% TeX-master: "testlauf"
%%% TeX-master: t
%%% TeX-master: t
%%% End: 

%% file: generalization_rr.tex
%%%%%%%%%%%%%%%%%%%%%%%%%%%%%%%%%%%%%%%%%%%%%%%%%%%%%%%%%%%%%%%%%%%%%%%%%%%%%%%%%%%%%%%%
\subsection{Function Rings over Reduction Rings}\label{section.right.rr}

The situation becomes more complicated for a function ring $\f_{\rr}$ where
 $\rr$ is not a field.
We will abbreviate $\f_{\rr}$ by $\f$.

Notice that similar to the previous section it is possible to study generalizations
 of standard representations for function rings over reduction rings with respect to the orderings $\succeq$ and $\geq$ on $\myt$.
General right standard representations as defined in Definition 
 \ref{def.general.standard.rep}, as well as the corresponding critical
 situations from Definition \ref{def.general.critical.situations}
 and the characterization of general right standard bases as 
 in Theorem \ref{theo.general.standard.basis} carry over to our function ring
 $\f$.
The same is true for right standard representations as defined in Definition
 \ref{def.standard.rep}, the corresponding critical
 situations from Definition \ref{def.critical.situations}
 and the characterization of right standard bases as 
 in Theorem \ref{theo.general.standard.basis}.
However, these standard representations can no longer be linked to weak right
 Gr\"obner bases as defined in Definition \ref{def.weak.gb}.
This is of course obvious as for function rings over fields we have a 
 characterization of such Gr\"obner bases by head terms which is no longer 
 possible for function rings over reduction rings.
This is already the case for polynomial rings over the integers.
For example take the polynomial $3 \skm X$ in $\q[X]$.
Then obviously for $F_1 = \{ 3 \skm X\}$ and $F_2 = \{ X \}$ we get that
 $\hterm(\ideal{r}{}(F_1)\backslash \{0\}) =
 \hterm(\{3 \skm X \rmult X^i \mid i \in \n \}) =
 \hterm(\{X \rmult X^i \mid i \in \n \}) =
 \hterm(\ideal{r}{}(F_2)\backslash \{0\})$
 while of course $F_1$ is no right Gr\"obner basis of $\ideal{r}{}(F_2)$
 and $F_2$ is no right Gr\"obner basis of $\ideal{r}{}(F_1)$.
One possible generalizing of Definition \ref{def.weak.gb} is as follows:
 $F$ is a weak right Gr\"obner basis of $\ideal{r}{}(F)$ if
 $\hm(\ideal{r}{}(F)\backslash \{0\}) = \hm(\{ f \rmult m \mid f \in F,
 m \in \monoms(\f) \})$.
But this does not solve the problem as there is no
 equivalent to Lemma \ref{lem.rsb=gb} to link these right Gr\"obner bases
 to the respective standard bases.
The reason for this is that the definitions of standard representations as 
 provided by Definition \ref{def.general.standard.rep} and \ref{def.standard.rep}
 are no longer related to reduction relations corresponding to Gr\"obner bases.
Of course it is possible to study other generalizations of these definitions,
 e.g.~involving the ordering on the coefficients, but we take a
 different approach.
Our studies of standard representations for function rings
 over fields revealed that in fact we can identify stronger conditions for such 
 representations in terms of weak right Gr\"obner bases (review e.g.~Corollary
 \ref{cor.right.rep} and \ref{cor.right_reductive}).
These special represenations arise from reduction sequences.
Hence we will proceed by studying such standard representations
 which can be directly related to reduction relations in our function ring.

%Let $\rr$ be a reduction ring as described in Section \ref{section.reductionrings}.
%We denote the redution relation by $\Longrightarrow$.
%Since our reduction relations for $\f$ will involve $\Longrightarrow$,
% in order to get terminating reduction relations we require for
% the ordering $>_{\rr}$ on $\rr$  that $\alpha \Longrightarrow_{B} \beta$
% for some $\alpha, \beta \in \rr$, $B \subseteq \rr$ implies $\alpha >_{\rr} \beta$.

Similar to function rings over fields we need to view $\f$ as a vector
 space now over $\rr$, a reduction ring as described in Section \ref{section.reductionrings}.
In general $\rr$ is not Abelian and hence we have to distinguish right and
 left scalar multiplication as defined on page \pageref{page.scalar}.
However, since $\rr$ is associative as in the case of fields we can write
 $\alpha \skm f \skm \beta$.
 
Notice that for $f,g$ in $\f$ and $\alpha,\beta \in \rr$ we have
\begin{enumerate}
\item $\alpha \skm (f \radd g) = \alpha \skm f \radd \alpha \skm g$
\item $\alpha \skm (\beta \skm f) = (\alpha \skm \beta) \skm f$
\item $(\alpha + \beta) \skm f = \alpha \skm f \radd \beta \skm f$,
\end{enumerate}
i.e., $\f$ is a left $\rr$-module.
Similarly we have
\begin{enumerate}
\item $(f \radd g)\skm \alpha =  f\skm \alpha \radd g \skm \alpha$
\item $(f \skm \alpha) \skm \beta = f \skm (\alpha \skm \beta)$
\item $f \skm (\alpha + \beta)  = f \skm \alpha \radd  f \skm \beta$,
\end{enumerate}
i.e., $\f$ is a right $\rr$-module as well.
Moreover, as $(\alpha \skm f) \skm \beta = \alpha \skm (f \skm \beta)$
 for all $f \in \f$, $\alpha, \beta \in \rr$, $\f$ is an 
 $\rr$-$\rr$ bimodule.

In order to state how scalar multiplication and ring multiplication are
 compatible,
we again require $(\alpha \skm f) \rmult g = \alpha \skm (f \rmult g)$
 and $f \rmult (g \skm \alpha) = (f \rmult g) \skm \alpha$ to hold.
This is true for all examples we know from the literature.

If we additionally require that for
 $\alpha, \beta \in \rr$ and $t,s \in \myt$ we have
 $(\alpha \skm t) \rmult (\beta \skm s) = (\alpha \skm \beta) \skm (t \rmult s)$,
then the multiplication in $\f$
 can be specified by knowing $\rmult : \myt \times \myt \myr \f$.

If $\f$ contains a unit $\one$, $\rr$ can be embedded into $\f$ via the mapping
 $\alpha \longmapsto \alpha \skm \one$.
Then for $\alpha \in \rr$ and $f \in \f$
 the equations $\alpha \skm f = (\alpha \skm \one) \rmult f$
 and $f \skm \alpha = f \rmult (\alpha \skm \one)$ hold.
Since for $\alpha \in \rr$ and $t \in \myt$ we have $\alpha \skm t = t \skm \alpha$
 this implies $(\alpha \skm t) \rmult (\beta \skm s) = (\alpha \skm \beta) \skm (t \rmult s)$\footnote{$(\alpha \skm t) \rmult (\beta \skm s) = 
 (\alpha \skm t) \rmult ((\beta \skm \one) \rmult s) =
 ((\alpha \skm t) \rmult (\beta \skm \one)) \rmult s =
 (\alpha \skm (t \rmult (\beta \skm \one)) \rmult s =
 (\alpha \skm (t \skm \beta)) \rmult s =
 (\alpha \skm (\beta \skm t)) \rmult s = (\alpha \skm \beta) \skm (t \rmult s)$.
  }.

Moreover, if $\rr$ is Abelian,
 we get $\alpha \skm (f \rmult g) = f \rmult (\alpha \skm g)$
 and $\f$ is an algebra.

Remember that we want to study standard representations directly related to
 reduction relations on $\f$.
Since we have a function ring over a reduction ring such a reduction relation 
 originates from the reduction relation on the reduction ring $\rr$.
Here we want to distinguish two such reduction relations on $\f$.

One possible generalization in the spirit of these ideas 
 for function rings over reduction rings is as follows:
\begin{definition}\label{def.right_reductive_rr}~\\
{\rm
Let $F$ be a set of polynomials in $\f$
 and $g$ a non-zero polynomial in $\ideal{r}{}(F)$.
A representation of the form 
$$g = \sum_{i=1}^n f_i \rmult m_i,
 f_i \in F, m_i \in \monoms(\f), n \in \n$$
 such that
 $\hterm(g) = \hterm(\hterm(f_1) \rmult m_1) =
 \hterm(f_1 \rmult m_1) \geq \hterm(f_1)$ and
 $\hterm(g) \succ \hterm(f_i \rmult m_i)$ 
 for all $2 \leq i \leq n$ is called a
 \betonen{right reductive standard  representation} in terms of $F$.
A set $F \subseteq \f \backslash \{ \zero\}$ is called a
 \betonen{right reductive standard basis} of $\ideal{r}{}(F)$
 if every polynomial $f \in \ideal{r}{}(F)$
 has a right reductive standard representation in terms of $F$. 
\dend
}
\end{definition}
Notice that that this definition
 differs from Definition \ref{def.right_reductive} insofar as we
 demand $\hterm(g) \succ \hterm(f_i \rmult m_i)$ 
 for all $2 \leq i \leq n$.
In fact we use those special standard representations which arise
 in the case of function rings for $g \in \ideal{r}{}(F)$ when $F$ 
 already is a right reductive standard basis 
 (compare Corollary \ref{cor.right_reductive}).
This definition is directly related to the reduction relation presented
 in Definition \ref{def.rred} for $\f_{\myk}$ generalized to $\f$. 
A possible definition of reduction can be given in the following fashion
 where we require that the reduction step eliminates the respective monomial
 it is applied to.
\begin{definition}\label{def.rred_rr}~\\
{\rm
Let $f,p$ be two non-zero polynomials in $\f$.
We say $f$ \betonen{right reduces} $p$ \betonen{to} $q$ \betonen{eliminating the monomial}
 $\alpha \skm t$ \betonen{in one step}, denoted by $p \red{}{\myr}{r,e}{f} q$, if
 there exists $m \in \monoms(\f)$ such that
\begin{enumerate}
\item $t \in \supp(p)$ and $p(t) = \alpha$,
\item $\hterm(\hterm(f)\rmult m) = \hterm(f \rmult m) = t \geq \hterm(f)$,
\item $\hm(f \rmult m) = \alpha \skm t$, such that
      $\alpha \R_{\hc(f \rmult m)} 0$, and
\item $q = p - f \rmult m$.
\end{enumerate}
We write $p \red{}{\myr}{r,e}{f}$ if there is a polynomial $q$ as defined
above and $p$ is then called right reducible by $f$. 
%\\
Further, we can define $\red{*}{\myr}{r,e}{}, \red{+}{\myr}{r,e}{}$ and
 $\red{n}{\myr}{r,e}{}$ as usual.
%\\
Right reduction by a set $F \subseteq \f$ is denoted by
 $p \red{}{\myr}{r,e}{F} q$ and abbreviates $p \red{}{\myr}{r,e}{f} q$
 for some $f \in F$.
\dend
}
\end{definition}
This reduction relation is related to a special instance\footnote{Compare
 Pan's reduction relation for the integers as defined in Example
 \ref{exa.Z}.} of
 the reduction relation $\R$. Notice that by Axiom (A2)
 $\alpha \R_{\hc(f \rmult m)} 0$ implies $\alpha \in \ideal{r}{\rr}(\hc(f \rmult m))$ and hence $\alpha = \hc(f \rmult m) \skm \beta$ for some $\beta \in \rr$.

Notice that in contrary to $\f_{\myk}$ now for $g, f \in \f$ and $m \in \monoms(\f)$
 the situation $\hterm(g) = \hterm(f \rmult m) = \hterm(\hterm(f) \rmult m) \geq \hterm(f)$ alone no
 longer implies that $\hm(g)$ is right reducible by $f$.
This is due to the fact that we can no longer modify the coefficients involved in the reduction step
 in the appropriate manner since reduction rings in general will not contain inverse elements.

Let us continue by studying our reduction relation.
\begin{lemma}\label{lem.rred.rr}~\\ 
{\sl 
Let $F$ be a set of polynomials in $\f \backslash \{ \zero\}$.
\begin{enumerate}
\item For $p,q \in \f$ $p \red{}{\myr}{r,e}{F} q$ implies $p \succ q$, in particular $\hterm(p)
  \succeq \hterm(q)$.
\item $\red{}{\myr}{r,e}{F}$ is Noetherian.
\lemend
\end{enumerate}
}
\end{lemma}
\Ba{}
\begin{enumerate}
\item Assuming that the reduction step takes place at a monomial $\alpha \skm t$,
      by Definition \ref{def.rred_rr} we know $\hm(f \rmult m) = \alpha \skm t$ 
      which yields
      $p \succ p  -  f \rmult m $
      since $\hm(f \rmult m) \succ \reductum(f \rmult m)$.
\item This follows from 1. 
\end{enumerate}\renewcommand{\baselinestretch}{1}\small\normalsize
\qed
The Translation Lemma no longer holds for this reduction relation.
This is already so for polynomial rings over the integers.
\begin{example}~\\
{\rm
Let $\z[X]$ be the polynomial ring in one indeterminant over $\z$.
Moreover, let $\R$ be the reduction relation on $\z$ where for
 $\alpha, \beta \in \z$,
 $\alpha \R_{\beta}$ if and only if there exists $\gamma \in \z$ such that
 $\alpha = \beta \skm \gamma$ (compare Example \ref{exa.Z}).
Let $p = 2 \skm x$, $q = -3 \skm X$ and $f = 5 \skm X$.
Then $p-q = 5 \skm X \red{}{\myr}{r,e}{f} 0$ while $p \nred{}{\myr}{r,e}{f}$
 and $q \nred{}{\myr}{r,e}{f}$.
\exaend
}
\end{example}
The reduction relation $\red{}{\myr}{r,e}{}$ in polynomial rings over the integers is
 known as Pan's reduction in the literature. 
The generalization of Gr\"obner bases then 
 are weak Gr\"obner bases as by completion one can achieve
 that all ideal elements reduce to zero.
Next we present a proper algebraic characterization of weak right Gr\"obner bases
 related to right reductive standard representations
 and the reduction relation defined in Definition \ref{def.rred_rr}.
Notice that it differs from Definition \ref{def.gb.fr} for function rings over fields
 insofar as we now have to look at the head monomials of the right ideal instead of the head terms only.
\begin{definition}\label{def.gb.rr}~\\
{\rm
A set $F \subseteq \f \backslash \{ \zero\}$ is called a 
 \betonen{weak right reductive Gr\"obner basis} of $\ideal{r}{}(F)$ if
$\hm(\ideal{r}{}(F) \backslash \{ \zero \}) =
 \hm( \{ f \rmult m \mid f \in F, m \in \monoms(\f),
 \hterm(\hterm(f) \rmult m) = \hterm(f \rmult m) \geq \hterm(f) \}
 \backslash \{ \zero \} )$.
\dend
}
\end{definition}
Similar to Lemma \ref{lem.sb=gb} right reductive standard bases and weak right reductive  Gr\"obner bases
 coincide.
\begin{lemma}\label{lem.sb=gb_rr}~\\
{\sl
Let $F$ be a subset of $\f \backslash \{ \zero\}$.
Then $F$ is a right reductive standard basis if and only
 if it is a weak right reductive  Gr\"obner basis.
\lemend
}
\end{lemma}
\Ba{}~\\
Let us first assume that $F$ is a right reductive standard basis, i.e., every polynomial
 $g$ in $\ideal{r}{}(F)$ has a right reductive standard representation with respect 
 to $F$.
In case $g \neq \zero$ this implies the existence of a polynomial $f \in F$ and a monomial 
 $m \in \monoms(\f)$ such that $\hterm(g) = \hterm(\hterm(f) \rmult m) 
 = \hterm(f \rmult m) \geq \hterm(f)$ and $\hm(g) = \hm(f \rmult m)$\footnote{Notice
 that if we had generalized the original Definition \ref{def.right_reductive} this
 would not holds.}.  
Hence $\hm(g) \in \hm( \{  f \rmult m \mid m \in \monoms(\f),  
 f \in F, \hterm(\hterm(f) \rmult m) = \hterm(f \rmult m) \geq \hterm(f) \}\backslash \{ \zero\})$.
As the converse, namely $\hm( \{  f \rmult m \mid m \in \monoms(\f), 
 f \in F, \hterm(\hterm(f) \rmult m) = \hterm(f \rmult m) \geq \hterm(f) \}\backslash \{ \zero\})
 \subseteq \hm(\ideal{r}{}(F) \backslash \{ \zero \})$ trivially holds, 
 $F$ is a weak right reductive  Gr\"obner basis.
\\
Now suppose that $F$ is a weak right reductive  Gr\"obner basis and again let $g \in \ideal{r}{}(F)$.
We have to show that $g$ has a right reductive standard representation with respect to $F$.
This will be done by induction on $\hterm(g)$.
In case $g = \zero$ the empty sum is our required right reductive standard representation.
Hence let us assume $g \neq \zero$.
Since then $\hm(g) \in \hm(\ideal{r}{}(F)\backslash \{ \zero\})$
 by the definition of weak
 right reductive  Gr\"obner bases we know there exists
 a polynomial $f \in F$ and a monomial $m \in \monoms(\f)$ such that
 $\hterm(\hterm(f) \rmult m) = \hterm(f \rmult m) \geq \hterm(f)$ and
 $\hm(g)=\hm(f \rmult m)$.
Let $g_1 = g - f \rmult m$. 
Then $\hterm(g) \succ \hterm(g_1)$ implies the 
 existence of a right reductive standard representation for $g_1$ which can be
 added to the multiple $f \rmult m$ to give
 the desired right reductive standard representation of $g$.
\\ 
\qed 
\begin{corollary}~\\
{\sl
Let $F$ a subset of $\f \backslash \{ \zero\}$ be a 
 weak right reductive  Gr\"obner basis.
Then every $g \in \ideal{r}{}(F)$ has a right reductive standard representation
       in terms of $F$ of the form
        $g = \sum_{i=1}^n f_i \rmult m_i,
         f_i \in F, m_i \in \monoms(\f), n \in \n$
         such that
         $\hterm(g) = \hterm(\hterm(f_1) \rmult m_1) =
         \hterm(f_1 \rmult m_1) \geq \hterm(f_1)$ and
         $\hterm(f_1 \rmult m_1) \succ \hterm(f_2 \rmult m_2) \succ \ldots \succ
          \hterm(f_n \rmult m_n)$.
}
\end{corollary}
\Ba{}~\\
This follows from inspecting the proof of Lemma \ref{lem.sb=gb_rr}.
\\
\qed
Another consequence of Lemma \ref{lem.sb=gb_rr} is the characterization of
 weak right reductive Gr\"obner bases in rewriting terms.
\begin{lemma}\label{lem.weak.gb}~\\
{\sl
A subset $F$ of $\f \backslash \{ \zero\}$ is a 
 weak right reductive Gr\"obner basis if for all $g \in \ideal{r}{}(F)$ 
 we have $g \red{*}{\myr}{r}{F} \zero$.
}
\end{lemma}

Now to find some analogon to s-polynomials in $\f$ we again study what polynomial
 multiples occur when changing arbitrary representations of right ideal elements
 into right reductive standard representations.

Given a generating set $F \subseteq \f$ of a right ideal in $\f$ the key idea 
 in order to characterize weak right Gr\"obner bases is  to distinguish special
 elements of $\ideal{r}{}(F)$ which have representations $\sum_{i=1}^n f_i \rmult m_i$,
 $f_i \in F$, $m_i \in \monoms(\f)$ such that the head terms $\hterm(f_i \rmult m_i)$
 are all the same within the representation.
Then on one hand the respective coefficients $\hc(f_i \rmult m_i)$ can add up to zero which
means that the sum of the head coefficients is in an appropriate module in $\rr$ ---
 m-polynomials
 are related to these situations (see also Definition \ref{def.critical.situations}).
If the result is not zero the sum of the coefficients $\hc(f_i \rmult m_i)$ can be described in terms
 of a (weak) right Gr\"obner basis in $\rr$ --- g-polynomials are related to these situations.
Zero divisors in the reduction ring eliminating the head
 monomial of a polynomial occur as a special instance of m-polynomials
 where $F = \{ f \}$ and $f \skm \alpha$, $\alpha \in \rr$ are considered.

The first problem is related
 to solving linear homogeneous equations in $\rr$ and to the existence of
 possibly finite bases of the respective modules.
In case we want effectiveness, we have to require that these bases are  computable.

The g-polynomials can successfully be treated when possibly 
 finite (weak) right Gr\"obner bases exist for finitely
 generated right ideals in $\rr$.
Here, in case we want effectiveness, we have to require that the (weak) right Gr\"obner bases as well as
  representations for their elements in terms of the generating set are computable.

Using m- and g-polynomials, weak right Gr\"obner bases can again be characterized
 as in Section \ref{section.polyrings}.
The definition of m- and g-polynomials is inspired by Definition \ref{def.rr.one-sided.gpol}.
One main difference however is that in function rings multiples of one polynomial by different
 terms can result in the same head terms for the multiples while the multiples themselves are different.
These multiples have to be treated as different ones contributing to the same overlap although they
 arise from the same polynomial.
Hence when looking at sets of polynomials we now have to assume that we have multisets which can
 contain polynomials more than once.
Additionally, while in Definition \ref{def.rr.one-sided.gpol} we can restrict our attention to
 overlaps equal to the maximal head term of the polynomials involved now we have to introduce
 the overlapping term as an additional variable. 

\begin{definition}\label{def.gpol}~\\
{\rm
Let $P = \{p_1, \ldots, p_k\}$ be a multiset of not necessarily different 
 polynomials in $\f$ and
 $t$ an element in $\myt$
 such that there are  $w_1, \ldots, w_k \in \myt$ with
 $\hterm(p_i \rmult w_i) = \hterm(\hterm(p_i) \rmult w_i) = t \geq \hterm(p_i)$, for all $1 \leq i \leq k$.
Further let $\gamma_i = \hc(p_i \rmult w_i)$ for  $1 \leq i \leq k$.
\\
Let $G$ be a (weak)
 right Gr\"obner basis of $\{ \gamma_1, \ldots, \gamma_k \}$ in $\rr$ with respect to $\R$.
Additionally let 
$$\alpha = \sum_{i=1}^k \gamma_i \skm \beta_{i}^{\alpha}$$
 for $\alpha \in G$, $\beta^{\alpha}_i \in \rr$,  $1 \leq i \leq  k$.
Then we define the  \betonen{g-polynomials (Gr\"obner polynomials)}
 corresponding to $p_1, \ldots, p_k$ and $t$ by setting
$$ g_{\alpha} = \sum_{i=1}^k  p_i \rmult  w_i \skm  \beta^{\alpha}_i.$$
Notice that $\hm(g_{\alpha})= \alpha \skm t$.
\\
For the right module $M = \{ (\delta_1, \ldots, \delta_k) \mid  
 \sum_{i=1}^k  \gamma_i \skm \delta_i = 0 \}$, let the set
 $\{B_j \mid j \in I_M \}$ be a basis with
 $B_j = (\beta_{j,1}, \ldots, \beta_{j,k})$ for $\beta_{j,l} \in \rr$
  and $1 \leq l \leq k$.
%\\
Then we define the 
 \betonen{m-polynomials (module polynomials)}
 corresponding to $P$ and $t$ by setting
$$ h_j = \sum_{i=1}^k p_i \rmult w_i \skm  \beta_{j,i}
 \mbox{ for each } j \in I_M.$$
Notice that $\hterm(h_j) \prec t$ for each $j \in I_M$.
\dend
}
\end{definition}

Given a set of polynomials $F$, the set of corresponding
 g- and m-polynomials contains those which are specified by
 Definition \ref{def.gpol} for each term $t \in \myt$ fulfilling the respective conditions.
For a set consisting of one polynomial the corresponding
m-polynomials reflect the multiplication of the polynomial with
zero-divisors of the head monomial, i.e., by a basis of the annihilator
 of the head monomial.
Notice that given a finite set of polynomials the corresponding sets of
 g- and m-polynomials in general can be infinite.

As in Theorem \ref{theo.rrsb=rgb.k} we can use g- and m-polynomials instead of s-polynomials
 to characterize special bases
 in function rings.
As before we also have to take into account right multiples of the generating set
 as Example \ref{exa.free.group} does not require a field as coefficient domain.
\begin{theorem}\label{theo.rr.cp.i}~\\
{\sl
Let $F$ be a set of polynomials in $\f \backslash \{ \zero \}$.
Then $F$ is a weak right Gr\"obner basis of $\ideal{r}{}(F)$ if and only if
\begin{enumerate}
\item for all $f$ in $F$ and for all $m$ in $\monoms(\f)$, 
       $f \rmult m$ has a right reductive standard representation in terms of $F$, and
\item all g- and m-polynomials corresponding to $F$ as specified in
       Definition \ref{def.gpol}
       have right reductive standard representations in terms of $F$.
\end{enumerate}
\theoend
}
\end{theorem}
\Ba{}~\\
In case $F$ is a weak right Gr\"obner basis it is also a right reductive standard basis, and
 since the multiples $f \rmult m$  and the respective g- and m-polynomials
 are all elements of $\ideal{r}{}(F)$ they must
 have right reductive standard representations.
\\
The converse will be proven by showing that every element in
 $\ideal{r}{}(F)$ has a right reductive standard representation in terms of $F$.
Let $g \in \ideal{r}{}(F)$ have a representation in terms of $F$ of
 the following form:
$g = \sum_{j=1}^m f_{j} \rmult (w_{j} \skm \alpha_j)$  such that
 $f_j \in F$, $w_j \in \myt$ and $\alpha_{j} \in \rr$.
Depending on this  representation of $g$ and the
 well-founded total ordering $\succeq$ on $\myt$ we define
 $t = \max_{\succeq} \{ \hterm(f_{j} \rmult (w_{j} \skm \alpha_j)) \mid 1\leq j \leq m \}$ and
 $K$ as the number of polynomials $f_j \rmult (w_{j} \skm \alpha_j)$ with head term $t$.
We show our claim by induction on $(t,K)$, where
 $(t',K')<(t,K)$ if and only if $t' \prec t$ or $(t'=t$ and
 $K'<K)$.
\\
Since by our first assumption every multiple $f_j \rmult (w_{j} \skm \alpha_j)$ in this
 sum has a right reductive standard representation in terms of $F$,
 we can assume that
 $\hterm(\hterm(f_j) \rmult w_j) = \hterm(f_j \rmult w_j) \geq
 \hterm(f_j)$ holds.
Moreover, without loss of generality we can assume that the polynomial multiples
 with head term $t$ are just $f_1  \rmult w_1, \ldots , f_K \rmult w_K$.
Notice that these assumptions on the representation of $g$ neither change
 $t$ nor $K$.
\\
Obviously, $t \succeq \hterm(g)$ must hold. 
If $K=1$ this gives us $t = \hterm(g)$ and by our assumptions our
 representation is already a right reductive one and we are done.
\\
Hence let us assume $K>1$.
\\
First let $\sum_{j=1}^K \hm(f_j \rmult (w_{j} \skm \alpha_j)) = \zero$.
Then by Definition \ref{def.gpol} there exists a tuple
 $(\alpha_1, \ldots, \alpha_K) \in M$, as 
 $\sum_{j=1}^K \hc(f_j \rmult w_{j}) \skm \alpha_j = 0$.
Hence there are $\delta_1, \ldots, \delta_K \in \rr$ such that
 $\sum_{i=1}^l A_i \skm \delta_i = (\alpha_1, \ldots, \alpha_K)$ for some $l \in \n$,
 $A_i=(\alpha_{i,1}, \ldots, \alpha_{i,K}) \in \{ A_j \mid j \in I_M \}$,
 and $\alpha_j = \sum_{i=1}^l \alpha_{i,j} \skm \delta_i$, $1 \leq j \leq K$.
By our assumption there are module polynomials $h_i = \sum_{j=1}^K f_j \rmult w_j \skm \alpha_{i,j}$,$1 \leq i \leq l$,
 all having right reductive standard representations in terms of $F$.
\\
Then since
\begin{eqnarray}
\sum_{j=1}^K f_j \rmult (w_{j} \skm \alpha_j) & = & \sum_{j=1}^K f_j \rmult w_{j} \skm (\sum_{i=1}^l \alpha_{i,j} \skm \delta_i) \nonumber\\
         &=& \sum_{j=1}^K \sum_{i=1}^l (f_j \rmult w_{j} \skm \alpha_{i,j}) \skm \delta_i \nonumber\\
         &=& \sum_{i=1}^l (\sum_{j=1}^K f_j \rmult w_{j} \skm \alpha_{i,j}) \skm \delta_i \nonumber\\
         &=& \sum_{i=1}^l h_i \skm \delta_i \nonumber
\end{eqnarray}
 we can change the representation of $g$ to
 $\sum_{i=1}^l h_i \skm \delta_i + \sum_{j=K+1}^m f_{j} \rmult (w_{j} \skm \alpha_j)$
 and replace each $h_i$ by its right reductive standard representation in terms of $F$.
Remember that for all $h_i$, $1 \leq i \leq l$ we have $\hterm(h_i) \pred t$.
Hence, for this new representation we now have maximal term smaller than $t$
 and by our induction hypothesis we have a right reductive standard representation
 for $g$ in terms of $F$ and are done.
\\
It remains to study the case where $\sum_{j=1}^K \hm(f_j \rmult (w_j \skm \alpha_j)) \neq 0$.
Then
 we have $\hterm(f_1 \rmult (w_1 \skm \alpha_1) + \dots + f_K \rmult (w_K \skm \alpha_K)) = t = \hterm(g)$,
 $\hc(g) = \hc(f_1 \rmult (w_1 \skm \alpha_1) + \dots + f_K \rmult (w_K \skm \alpha_K)) \in \ideal{r}{}(\{ \hc(f_1 \rmult w_1), \ldots, 
 \hc(f_K \rmult w_K) \})$ and $\hm(f_1 \rmult (w_1 \skm \alpha_1) + \dots + f_K \rmult (w_K \skm \alpha_K)) = \hm(g)$.
Hence $\hc(g) = \alpha \skm \delta$ with $\delta \in \rr$ and $\alpha \in G$\footnote{Remember
 that we assume the reduction relation $\R$ on $\rr$ based on division, see the remark after
 Definition \ref{def.rred_rr}.}, $G$ being a
 (weak) right Gr\"obner basis of $\ideal{r}{}(\{ \hc(f_1 \rmult w_1), \ldots, 
 \hc(f_K \rmult w_K) \})$
 (compare Definition \ref{def.gpol}).
Let $g_{\alpha}$ be the respective g-polynomial corresponding to $\alpha$.
Then the polynomial $g' = g - g_{\alpha} \skm \delta$ lies in
 $\ideal{r}{}(F)$.
Since the multiple\footnote{Note that right reductive standard representations are stable
 under multiplication with coefficients which are no zero-divisors of the head coefficient.} $g_{\alpha} \skm \delta$ has a right reductive standard representation in terms of $F$,
 say $\sum_{j=1}^l f_j \rmult m_j$,
 for the situation $\sum_{j=1}^K f_{j} \rmult (w_{j} \skm \alpha_j)- f_1 \rmult m_1$
 all polynomial multiples involved in this sum have head term $t$ and their head monomials add up to $\zero$.
Therefore, this situation again corresponds to an m-polynomial of $F$.
Hence we can apply our results from above and get that the polynomial $g'$ 
 has a smaller representation than $g$,
 especially the maximal term $t'$ is smaller.
Moreover, we
 can assume that $g'$ has a right reductive standard representation in terms
 of $F$, say $g' = \sum_{i=1}^n f_i \rmult \tilde{m}_i$.
Then $g = \sum_{i=1}^n f_i \rmult \tilde{m}_i + g_{\alpha} \skm \delta =
 \sum_{i=1}^n f_i \rmult \tilde{m}_i + \sum_{j=1}^l f_j \rmult m_j$ is a right
 reductive standard representation
 of $g$ in terms of $F$ and we are done.
\\
\qed
Since in general we will have infinitely many g- and m-polynomials related to $F$,
 it is important to look for possible localizations of these situations.
We are looking for concepts similar to those of
 weak saturation and stable localizations in the previous section.
Remember that Lemma \ref{lem.red.reps} is central there.
It describes when the existence of a right reductive standard representation for some
 polynomial implies the existence of a right reductive standard representation
 for a multiple of the polynomial.
Unfortunately we cannot establish an analogon to this lemma for
 right reductive standard representations in $\f$ as defined in Definition \ref{def.right_reductive_rr}.
\begin{example}~\\
{\rm
Let $\f$ be a function ring over the integers with $\myt = \{X_1, \ldots, X_7 \}$
 and multiplication $\rmult : \myt \times \myt \mapsto \f$ defined by the
 following equations:
$X_1 \rmult X_2 = X_4$, $X_4 \rmult X_3 = X_5$, $X_2 \rmult X_3 = X_6 + X_7$,
 $X_1 \rmult X_6 = 3 \skm X_5$, $X_1 \rmult X_7 = -2 \skm X_5$ and else
 $X_i \rmult X_j = \zero$.
Additionally let $X_5 > X_4 > X_1 \succ X_2 \succ X_3 \succ X_6 \succ X_7$.
\\
Then for $p = X_4$, $f = X_1$ and $m = X_3$ we find that
\begin{enumerate}
\item $p$ has a right reductive standard representation in terms of $\{ f \}$, namely $p = f \rmult X_2$.
\item $\hterm( p \rmult m) = \hterm( \hterm(p) \rmult m ) \geq \hterm(p)$ as
      $X_5 = X_4 \rmult X_3 > X_4$ and for all $X_i \pred X_4$ we have $X_i \rmult X_3 \pred X_5$.
\item $p \rmult m = X_5$ has no right reductive standard representation in terms of $\{ f \}$ as only
      $X_1 \rmult X_j \neq \zero$ for $j = \{ 2,6,7 \}$, 
        namely $X_1 \rmult X_2 = X_4$, $X_1 \rmult X_6 = 3 \skm X_5$,
        $X_1 \rmult X_7 = -2 \skm X_5$,  and
      $X_1 \rmult (X_j \skm \alpha) \neq X_5$ for all $j \in \{ 2,6,7 \}$, $\alpha \in \z$.
\end{enumerate}
Notice that these problems are due to the fact that while $(X_1 \rmult X_2) \rmult X_3 =
 X_1 \rmult (X_2 \rmult X_3) = X_5$, $X_1 \rmult (X_2 \rmult X_3) = X_1 \rmult (X_6 + X_7) =
 X_1 \rmult X_6 + X_1 \rmult X_7$ does not give us a right reductive standard representation in terms of $X_1$ as $\hterm(X_1 \rmult X_6) = X_5$ {\em and} $\hterm(X_1 \rmult X_7) = X_5$ (compare Definition \ref{def.right_reductive_rr}).
This was the crucial point in the proof of Lemma \ref{lem.red.reps} and it is only fulfilled for the weaker form 
 of right reductive standard representations in $\f_{\myk}$ as defined in Definition \ref{def.right_reductive}.
\exaend
}
\end{example}
As this example shows an analogon to Lemma \ref{lem.red.reps} does not hold in our general case.
Note that the trouble arises from the fact that we allow multiplication of two terms to result
 in a polynomial.
If we restrict ourselves to multiplications where multiples of monomials are again monomials,
 the proof of Lemma \ref{lem.red.reps} carries over and we can look for appropriate localizations.

However, the reduction relation defined in Definition \ref{def.rred_rr} is only one way of
 defining a reduction relation in $\f$ and we stated that the main motivation behind it
 is to link the reduction relation with special standard representations as it is done in the case of $\f_{\myk}$.
The question now arises whether this motivation is as appropriate for $\f$ as it was for $\f_{\myk}$.
In $\f_{\myk}$ any reduction relation based on stable divisibility of terms can be linked
 to right reductive standard representations as defined in Definition \ref{def.right_reductive}
 and hence the approach is very powerful.
It turns out that for different reduction relations in $\f$ based on stable right divisibility
 this is no longer so.
Let us look at another familiar way of generalizing a reduction relation for $\f$ from one 
 defined in the reduction ring.
From now on we require a (not necessarily Noetherian) partial ordering on $\rr$:
 for $\alpha, \beta \in \rr$, $\alpha >_{\rr} \beta$ if and only if
 there exists a finite set $B \subseteq \rr$ such that $\alpha
 \red{+}{\Longrightarrow}{}{B} \beta$.
This ordering ensures that reduction in $\f$ is terminating when using a finite set 
 of polynomials.
\begin{definition}\label{def.rred.rr}~\\
{\rm
Let $f,p$ be two non-zero polynomials in $\f$.
We say $f$ \betonen{right reduces} $p$ \betonen{to} $q$ \betonen{at a monomial}
 $\alpha \skm t$ \betonen{in one step}, denoted by $p \red{}{\myr}{r}{f} q$, if
 there exists $ m \in \monoms(\f)$ such that
\begin{enumerate}
\item $t \in \supp(p)$ and $p(t) = \alpha$,
\item $\hterm(\hterm(f)\rmult m) = \hterm(f \rmult m) = t \geq \hterm(f)$, 
\item $\alpha \Longrightarrow_{\hc(f \rmult m)} \beta$, with  
      $\alpha = \hc(f \rmult m) + \beta$
      for some $\beta\in \rr$, and
\item $q = p -  f \rmult m$.
\end{enumerate}
We write $p \red{}{\myr}{r}{f}$ if there is a polynomial $q$ as defined
above and $p$ is then called right reducible by $f$. 
%\\
Further, we can define $\red{*}{\myr}{r}{}, \red{+}{\myr}{r}{}$ and
 $\red{n}{\myr}{r}{}$ as usual.
%\\
Right reduction by a set $F \subseteq \f \backslash \{ \zero \}$ is denoted by
 $p \red{}{\myr}{r}{F} q$ and abbreviates $p \red{}{\myr}{r}{f} q$
 for some $f \in F$.
\dend
}
\end{definition}
Notice that in specifying this reduction relation we use a special instance of
 $\alpha \Longrightarrow_{\hc(f \rmult m)} \beta$,
 namely the case that $\alpha = \hc(f \rmult m) + \beta$
 for some $\beta\in \rr$.
Moreover, for this reduction relation we can still have $t \in \supp(q)$.
Hence other arguments than used in the proof of Lemma \ref{lem.rred.rr} have
 to be provided to show termination.
It turns out that for infinite subsets of polynomials $F$ in $\f$ the reduction
 relation $\red{}{\myr}{r}{F}$ need not terminate.
\begin{example}\label{exa.reduction.not.terminating}~\\
{\rm
Let $\rr = \q[\{ X_i \mid i \in \n \}]$ with $X_1 \succ X_2 \succ \ldots$
 be the polynomial ring over the rationals with infinitely many indeterminates. 
We associate this ring with the reduction relation based on divisibility of terms.
Let $\f = \rr[Y]$ be our function ring.
Elements of $\f$ are polynomials in $Y^i$, $i \in \n$ with coefficients in $\rr$.
Then for $p = X_1 \skm Y$ and the infinite set $F = \{ f_i = (X_i - X_{i+1}) \skm Y
 \mid i \in \n \}$ we get the infinite reduction sequence
 $p \red{}{\myr}{r}{f_1} X_2 \skm Y \red{}{\myr}{r}{f_2} 
  X_3 \skm Y \red{}{\myr}{r}{f_3} \ldots$
\exaend
}
\end{example}
However, if we restrict ourselves to finite sets of polynomials the reduction
 relation is Noetherian.
\begin{lemma}\label{lem.sred.rr}~\\ 
{\sl 
Let $F$ be a finite set of polynomials in $\f \backslash \{ \zero\}$.
\begin{enumerate}
\item For $p,q \in \f$ $p \red{}{\myr}{r}{F} q$ implies $p \succ q$, in particular $\hterm(p)
  \succeq \hterm(q)$.
\item $\red{}{\myr}{r}{F}$ is Noetherian.
\lemend
\end{enumerate}
}
\end{lemma}
\Ba{}
\begin{enumerate}
\item Assuming that the reduction step takes place at a monomial $\alpha \skm t$,
      by Definition \ref{def.rred.rr} we know 
        $\hm(\alpha \skm t -  f \rmult m) = \beta \skm t$ 
      which yields
      $p \succ p  -  f \rmult m$
      since $\alpha >_{\rr} \beta$.
\item This follows from 1.~and Axiom (A1) as long as only finite sets of polynomials are involved.
      Since we have $\hterm(f \rmult m) = \hterm( \hterm(f) \rmult m) \geq \hterm(f)$ we get
      $\hc(f \rmult m) = \hc(f) \skm \hc(\hterm(f) \rmult m)$.
      Then $\alpha \Longrightarrow_{\hc(f \rmult m)} \beta$ implies
      $\alpha \R_{\hc(f)}$.
      Hence an infinite reduction sequence would give rise to an infinite
       reduction sequence in $\rr$ with respect to 
       the finite set of head coefficients $\{ \hc(f) \mid f \in F \}$ 
       contradicting our assumption.
\end{enumerate}\renewcommand{\baselinestretch}{1}\small\normalsize
\qed
Now if we try to link the reduction relation in Definition 
 \ref{def.rred.rr} to special standard representations, we find
 that this is no longer as natural as in the cases studied before,
 where for $\f_{\myk}$ we linked the reduction relation from Definition \ref{def.rred}
 to the right reductive standard representations in Definition
 \ref{def.right_reductive} respectively for $\f$ the right reduction relation
 from Definition \ref{def.rred_rr} to right reductive standard 
 representations as defined in Definition \ref{def.right_reductive_rr}.
Hence we claim that for generalizing Gr\"obner bases to $\f$, the rewriting
 approach is more suitable.
Hence we use the following definition of weak right Gr\"obner bases
 in terms of our reduction relation.
\begin{definition}~\\
{\rm
A set $F \subseteq \f \backslash \{ \zero \}$ is called a weak
 right Gr\"obner basis (with respect to $\red{}{\myr}{r}{}$)
 of $\ideal{r}{}(F)$ if for all
 $g \in \ideal{r}{}(F)$ we have $g \red{*}{\myr}{r}{F} \zero$. 
\dend
}
\end{definition}
Every reduction sequence $g \red{*}{\myr}{r}{F} \zero$ gives rise to a special
 representation of $g$ in terms of $F$ which could be taken as a new definition
 of standard representations.
\begin{corollary}\label{cor.representation.rr}~\\
{\sl
Let $F$ be a set of polynomials in $\f$ 
 and $g$ a non-zero polynomial in $\ideal{r}{}(F)$
 such that $g \red{*}{\myr}{r}{F} \zero$.
Then $g$ has a representation of the form 
$$g = \sum_{i=1}^n f_i \rmult m_i,
 f_i \in F, m_i \in \monoms(\f), n \in \n$$
 such that
 $\hterm(g) = \hterm(\hterm(f_i) \rmult m_i) =
 \hterm(f_i \rmult m_i) \geq \hterm(f_i)$ for $1 \leq i \leq k$, and
 $\hterm(g) \succ \hterm(f_i \rmult m_i)$ 
 for all $k+1 \leq i \leq n$.
\corend
}
\end{corollary}
\Ba{}~\\
We show our claim by induction on $n$ where $g \red{n}{\myr}{r}{F} \zero$.
If $n=0$ we are done.
Else let $g \red{1}{\myr}{r}{F} g_1 \red{n}{\myr}{r}{F} \zero$.
In case the reduction step takes place at the head monomial,
 there exists a polynomial $f \in F$
 and a monomial $m \in \monoms(\f)$
 such that $\hterm(\hterm(f) \rmult m) = \hterm(f \rmult m) = \hterm(g)
 \geq \hterm(f)$ and $\hc(g) \R_{\hc(f \rmult m)} \beta$ with
 $\hc(g) = \hc(f\rmult m) + \beta$ for some $\beta \in \rr$.
Moreover the induction hypothesis then is applied to  $g_1 = 
 g - f \rmult m \skm \beta$.
If the reduction step takes place at a monomial with term smaller $\hterm(g)$
 for the respective monomial multiple $f \rmult m$ we immediately get
 $\hterm(g) \succ \hterm(f \rmult m)$ and we can apply our induction hypothesis
 to the resulting polynomial $g_1$.
In both cases we can arrange the monomial multiples $f \rmult m$ arising from
 the reduction steps in such a way that gives us th desired representation.
\\
\qed 
Notice that on the other hand the existence of such a representation for a
 polynomial does not imply reducibility.
For example take the polynomial ring $\z[X]$ with Pan's reduction. Then with respect to
 the polynomials $F = \{2 \skm X, 3 \skm X \}$ the polynomial $g = 5 \skm X$ has a
 representation $5 \skm X = 2 \skm X + 3 \skm X$ of the desired form but is neither
 reducible by $2 \skm X$ nor $3 \skm X$.
This is of course a consequence of the fact that $\{2,3\}$ is no Gr\"obner basis
 in $\z$ with respect to Pan's reduction.

In fact Corollary \ref{cor.representation.rr} provides additional
 information for the head coefficient of $g$, namely
 $\hc(g) = \sum_{i=1}^k \hc(f_i) \skm \hc(m_i)$ and this is a
 standard representation
 of $\hc(g)$ in terms of $\{\hc(f_i) \mid 1 \leq i \leq k \}$ in the reduction
 ring $\rr$.

We can characterize weak right Gr\"obner bases similar to
 Theorem \ref{theo.rr.cp.i}.
Of course the g-polynomials in Definition \ref{def.gpol} depend on the reduction
 relation $\R$ in $\rr$ which now is defined according to Definition \ref{def.rred.rr}.
Notice that the characterization will only hold for finite
 sets as the proof requires the reduction relation to be Noetherian.
Additionally we need that the reduction ring fulfills Axiom (A4),
 i.e., for $\alpha, \beta, \gamma, \delta \in \rr$,
             $\alpha \Longrightarrow_{\beta}$ and
             $\beta \Longrightarrow_{\gamma} \delta$ 
             imply $\alpha \Longrightarrow_{\gamma}$ or 
             $\alpha \Longrightarrow_{\delta}$\footnote{Notice that
 (A4) is no basis for localizing test sets, as this would require that
 $\alpha \Longrightarrow_{\beta}$ and $\beta \Longrightarrow_{\gamma} \delta$ 
 imply $\alpha \Longrightarrow_{\gamma}$. Hence even if the reduction relation
 in $\f$ satisfies (A4), this does not substitute Lemma \ref{lem.red.reps}
 or its variants.}.
\begin{theorem}\label{theo.gb.reduction}~\\
{\sl
Let $F$ be a finite set of polynomials in $\f \backslash \{ \zero \}$
 where the reduction ring satisfies (A4).
Then $F$ is a weak right Gr\"obner basis  of $\ideal{r}{}(F)$ if and only if
\begin{enumerate}
\item for all $f$ in $F$ and for all $m$ in $\monoms(\f)$ we have 
       $f \rmult m \red{*}{\myr}{r}{F} \zero$, and
\item all g- and m-polynomials corresponding to $F$ as specified in
       Definition \ref{def.gpol}
       reduce to $\zero$ using $F$.
\theoend
\end{enumerate}

}
\end{theorem}
\Ba{}~\\
In case $F$ is a weak right Gr\"obner basis,
 since the multiples $f \rmult m$ and the respective g- and m-polynomials
 are all elements of $\ideal{r}{}(F)$ they must
 reduce to zero using $F$.
\\
The converse will be proven by showing that every element in
 $\ideal{r}{}(F)$ is reducible  by $F$.
Then as $g \in \ideal{r}{}(F)$ and $g \red{}{\myr}{r}{F} g'$ implies
 $g' \in \ideal{r}{}(F)$ we have $g\red{*}{\myr}{r}{F} \zero$.
Notice that this only holds in case the reduction relation $\red{}{\myr}{r}{F}$
 is Noetherian.
This follows as by our assumption $F$ is finite (Lemma \ref{lem.sred.rr}).
\\
Let $g \in \ideal{r}{}(F)$ have a representation in terms of $F$ of
 the following form:
$g = \sum_{j=1}^m f_{j} \rmult (w_{j} \skm \alpha_j)$ such that
 $f_j \in F$, $w_{j} \in \myt$, $\alpha_j \in \rr$.
Depending on this  representation of $g$ and the
 well-founded total ordering $\succeq$ on $\myt$ we define
 $t = \max_{\succeq} \{ \hterm(f_{j} \rmult (w_{j} \skm \alpha_j)) \mid 1\leq j \leq m \}$ and
 $K$ as the number of polynomials $f_j \rmult (w_{j} \skm \alpha_j)$ with head term $t$.
We show our claim by induction on $(t,K)$, where
 $(t',K')<(t,K)$ if and only if $t' \prec t$ or $(t'=t$ and
 $K'<K)$.
\\
Since by our first assumption every multiple $f_j \rmult (w_{j} \skm \alpha_j)$ in this
 sum reduces to zero using $F$ and hence
 has a right  representation as defined
 in Corollary \ref{cor.representation.rr}, we can assume that
 $\hterm(\hterm(f_j) \rmult w_j) = \hterm(f_j \rmult w_j) \geq
 \hterm(f_j)$ holds.
Moreover, without loss of generality we can assume that the polynomial multiples
 with head term $t$ are just $f_1  \rmult(w_{1} \skm \alpha_1) , \ldots , f_K \rmult (w_K \skm \alpha_K)$.
Notice that these assumptions neither change $t$ nor $K$ for our representation
 of $g$.
\\
Obviously, $t \succeq \hterm(g)$ must hold. 
If $K=1$ this gives us $t = \hterm(g)$ and even
 $\hm(g) = \hm(f_1 \rmult (w_{1} \skm \alpha_1))$, implying that $g$ is right
 reducible at $\hm(g)$ by $f_1$.
\\
Hence let us assume $K>1$.
\\
First let $\sum_{j=1}^K \hm(f_j \rmult (w_{j} \skm \alpha_j)) = \zero$.
Then by Definition \ref{def.gpol} we know $(\alpha_1, \ldots, \alpha_K) \in M$, as 
 $\sum_{j=1}^K \hc(f_j \rmult w_{j}) \skm \alpha_j = 0$.
Hence there are $\delta_1, \ldots, \delta_K \in \rr$ such that
 $\sum_{i=1}^l A_i \skm \delta_i = (\alpha_1, \ldots, \alpha_K)$ for some $l \in \n$,
 $A_i=(\alpha_{i,1}, \ldots, \alpha_{i,K}) \in \{ A_j \mid j \in I_M \}$,
 and $\alpha_j = \sum_{i=1}^l \alpha_{i,j} \skm \delta_i$, $1 \leq j \leq K$.
By our assumption there are module polynomials $h_i = \sum_{j=1}^K f_j \rmult w_j \skm \alpha_{i,j}$,$1 \leq i \leq l$,
 all having representations in terms of $F$ as defined
 in Corollary \ref{cor.representation.rr}.
\\
Then since
\begin{eqnarray}
\sum_{j=1}^K f_j \rmult (w_{j} \skm \alpha_j) & = & \sum_{j=1}^K f_j \rmult w_{j} \skm (\sum_{i=1}^l \alpha_{i,j} \skm \delta_i) \nonumber\\
         &=& \sum_{j=1}^K \sum_{i=1}^l (f_j \rmult w_{j} \skm \alpha_{i,j}) \skm \delta_i \nonumber\\
         &=& \sum_{i=1}^l (\sum_{j=1}^K f_j \rmult w_{j} \skm \alpha_{i,j}) \skm \delta_i \nonumber\\
         &=& \sum_{i=1}^l h_i \skm \delta_i \nonumber
\end{eqnarray}
 we can change the representation of $g$ to
 $\sum_{i=1}^l h_i \skm \delta_i + \sum_{j=K+1}^m f_{j} \rmult (w_{j} \skm \alpha_j)$
 and replace each $h_i$ by its respective representation in terms of $F$.
Remember that for all $h_i$, $1 \leq i \leq l$ we have $\hterm(h_i) \pred t$.
Hence, for this new representation we now have maximal term smaller than $t$
 and by our induction hypothesis $g$ is reducible by $F$ and we are done.
\\
It remains to study the case where $\sum_{j=1}^K \hm(f_j \rmult (w_j \skm \alpha_j)) \neq 0$.
Then
 we have $\hterm(f_1 \rmult (w_1 \skm \alpha_1) + \dots + f_K \rmult (w_K \skm \alpha_K)) = t = \hterm(g)$,
 $\hc(g) = \hc(f_1 \rmult (w_1 \skm \alpha_1) + \dots + f_K \rmult (w_K \skm \alpha_K)) \in \ideal{r}{}(\{ \hc(f_1 \rmult w_1), \ldots, 
 \hc(f_K \rmult w_K) \})$ and 
 even $\hm(f_1 \rmult (w_1 \skm \alpha_1) + \dots + f_K \rmult (w_K \skm \alpha_K)) = \hm(g)$.
Hence $\hc(g)$ is $\R$-reducible  by some $\alpha$, $\alpha \in G$, a
 (weak) right Gr\"obner basis of $\ideal{r}{}(\{ \hc(f_1 \rmult w_1), \ldots, 
 \hc(f_K \rmult w_K) \})$ in $\rr$ with respect to the reduction relation $\R$.
Let $g_{\alpha}$ be the respective g-polynomial corresponding to $\alpha$ and $t$.
Then we know that $g_{\alpha}\red{*}{\myr}{r}{F} \zero$.
Moreover, we know that the head monomial of $g_{\alpha}$ is reducible
 by some polynomial $f \in F$ and we assume
 $\hterm(g_{\alpha}) = \hterm(\hterm(f) \rmult m) = \hterm(f \rmult m) \geq \hterm(f)$
 and $\hc(g_{\alpha}) \R_{\hc(f \rmult m)}$.
Then, as $\hc(g)$ is $\R$-reducible by $\hc(g_{\alpha})$, 
 $\hc(g_{\alpha})$ is $\R$-reducible  and (A4) holds,
 the head monomial of $g$ is also reducible
%\footnote{Remember that by (A4)
% for $\alpha, \beta_1, \ldots, \beta_k \in \rr$, $\alpha \red{*}{\R}{}{\{
% \beta_1,\ldots, \beta_k \}} 0$ implies $\alpha \R_{\beta_i}$ for some
% $1 \leq i \leq k$.} 
by some $f' \in F$ and we are done. 
\\
\qed
Of course this theorem is also true for infinite $F$ if we can show
 that for the respective function ring the reduction relation is terminating.

Now the question arises when the critical situations in this characterization
 can be localized to subsets of the respective sets
 as in Theorem \ref{theo.s-pol.2}.
Reviewing the Proof of Theorem \ref{theo.s-pol.2} we find that 
 Lemma \ref{lem.red.reps} is central as it describes when multiples of
 polynomials which have a right reductive standard representation in terms
 of some set $F$ again have such a representation.
As we have seen above, this will not hold for function rings over reduction
 rings in general.
Now one way to introduce localizations would be to restrict the attention
 to those $\f$ satisfying Lemma \ref{lem.red.reps}.
Then appropriate adaptions of Definition \ref{def.weakly.saturated},
 \ref{def.saturator} and \ref{def.s-poly.localization} would allow a localization
 of the critical situations.
However, we have stated that it is not natural to link right reduction as defined
 in Definition \ref{def.rred_rr} to special standard representations.
Hence, to give localizations of Theorem \ref{theo.gb.reduction} another
 property for $\f$ is sufficient:
\begin{definition}\label{def.stable.loc}~\\
{\rm
A set $C \subset S \subseteq \f$ is called a \betonen{stable localization} of
 $S$ if for every $g \in S$ there exists $f \in C$
 such that $g \red{}{\myr}{r}{f} \zero$.
% and if $\hm(f)$ is reducible by $F$, so is $\hm(g)$.
\dend
}
\end{definition}
In case $\f$ and $\red{}{\myr}{r}{}$ allow such stable localizations,
 we can rephrase Theorem \ref{theo.gb.reduction} as follows:
\begin{theorem}\label{theo.loc}~\\
{\sl
Let $F$ be a finite set of polynomials in $\f \backslash \{ \zero \}$
 where the reduction ring satisfies (A4).
Then $F$ is a weak right Gr\"obner basis of $\ideal{r}{}(F)$ if and only if
\begin{enumerate}
\item for all $s$ in a stable localization of 
       $\{ f \rmult m \mid f \in \f, m \in \monoms(\f) \}$ we have 
       $s \red{*}{\myr}{r}{F} \zero$, and
\item for all $h$ in a stable localization of the g- and m-polynomials corresponding to $F$ as specified in
       Definition \ref{def.gpol} we have 
       $h \red{*}{\myr}{r}{F} \zero$.
\end{enumerate}
\theoend
}
\end{theorem} 
We have stated that for arbitrary reduction relations in $\f$ it is not natural to
 link them to special standard representations.
Still, when proving Theorem \ref{theo.loc}, we will find that in order to change the
 representation of an arbitrary right ideal element,
 Definition \ref{def.stable.loc} is not enough to ensure reducibility.
However, we can substitute the critical situation using an analogon of
 Lemma \ref{lem.red.reps}, which, while not related to reducibility, in this case will
 still be sufficient to make the representation smaller.
\begin{lemma}\label{lem.red.reps_rr}~\\
{\sl
Let  $F$ be a subset of polynomials in $\f \backslash \{ \zero\}$ and
 $f$, $p$  non-zero polynomials in $\f$.
If $p \red{}{\myr}{r}{f} \zero$ and $f \red{*}{\myr}{r}{F} \zero$,
 then $p$ has a standard representation of the form
$$p = \sum_{i=1}^n f_i \rmult l_i,
 f_i \in F, l_i \in \monoms(\f), n \in \n$$
 such that
 $\hterm(p) = \hterm(\hterm(f_i) \rmult l_i) =
 \hterm(f_i \rmult l_i) \geq \hterm(f_i)$ for $1 \leq i \leq k$ and
 $\hterm(p) \succ \hterm(f_i \rmult l_i)$ 
 for all $k+1 \leq i \leq n$ (compare Definition \ref{def.right_reductive}).
\lemend
}
\end{lemma}
\Ba{}~\\
If $p \red{}{\myr}{r}{f} \zero$ then $p = f \rmult m$ with $m \in \monoms(\f)$
 and $\hterm(p) = \hterm(\hterm(f) \rmult m) =
 \hterm(f \rmult m) \geq \hterm(f)$.
Similarly $f \red{*}{\myr}{r}{F} \zero$ implies
 $f = \sum_{i=1}^n f_i \rmult m_i,
 f_i \in F, m_i \in \monoms(\f), n \in \n$
 such that
 $\hterm(f) = \hterm(\hterm(f_i) \rmult m_i) =
 \hterm(f_i \rmult m_i) \geq \hterm(f_1)$, $1 \leq i \leq k$, and
 $\hterm(f) \succ \hterm(f_i \rmult m_i)$ 
 for all $k+1 \leq i \leq n$ (compare Corollary \ref{cor.representation.rr}). 
\\
Let us first analyze $f_i \rmult m_i \rmult m$ with $\hterm(f_i \rmult m_i) = \hterm(f)$,
 $1 \leq i \leq k$.\\
Let $\terms(f_i \rmult m_i) = \{ s_1^i, \ldots, s_{k_i}^i \}$ with
 $s_1^i \succ s_j^i$,
 $2 \leq j \leq k_i$, i.e., $s_1^i = \hterm(f_i \rmult m_i) = \hterm(\hterm(f_i) \rmult m_i) =\hterm(f)$.
Hence $\hterm(f) \rmult m = s_1^i \rmult m \geq \hterm(f) = s_1^i$ and
 as $s_1^i \succ s_j^i$, $2 \leq j \leq k_i$, by Definition \ref{def.refined.ordering}
 we can conclude that $\hterm(\hterm(f) \rmult m) = \hterm(s_1^i \rmult m)
 \succ s_j^i \rmult m \succeq \hterm(s_j^i \rmult m)$ for $2 \leq j \leq k_i$.
This implies $\hterm(\hterm(f_i \rmult m_i) \rmult m) = \hterm(f_i \rmult m_i \rmult m)$.
Hence we get
\begin{eqnarray*}
 \hterm(f \rmult m) & = & \hterm(\hterm(f) \rmult m) \\
                    & = & \hterm(\hterm(f_i \rmult m_i) \rmult m), \mbox{ as } 
                          \hterm(f) = \hterm(f_i \rmult m_i) \\
                    & = & \hterm(f_i \rmult m_i \rmult m)
\end{eqnarray*}
 and since $\hterm(f \rmult m) \geq \hterm(f) \geq \hterm(f_i)$ we can conclude
 $\hterm(f_i \rmult m_i \rmult m) \geq \hterm(f_i)$.
It remains to show that the $f_i \rmult m_i \rmult m$ have
 representations of the desired form in terms of $F$.
First we show that $\hterm(\hterm(f_i) \rmult m_i \rmult m) \geq \hterm(f_i)$.
We know $\hterm(f_i) \rmult m_i \succeq 
 \hterm(\hterm(f_i) \rmult m_i) = \hterm(f_i \rmult m_i)$\footnote{Notice that
 $\hterm(f_i) \rmult m_i$ can be a polynomial and hence we cannot conclude 
 $\hterm(f_i) \rmult m_i = \hterm(\hterm(f_i) \rmult m_i)$.}
 and hence
 $\hterm(\hterm(f_i) \rmult m_i \rmult m) = \hterm(\hterm(f_i \rmult m_i) \rmult m) =
  \hterm(f_i \rmult m_i \rmult m) \geq \hterm(f_i)$.
Then in case $m_i \rmult m \in \monoms(\f)$ we are done as then $f_i \rmult ( m_i \rmult m)$
 is a representation of the desired form.
\\
Hence let us assume
 $m_i \rmult m = \sum_{r=1}^{k_i} \tilde{m}^i_r$, $\tilde{m}^i_r \in \monoms(\f)$.
Let $\terms(f_i) = \{ t^i_1, \ldots, t^i_{w_i} \}$ with $t^i_1 \succ t^i_l$,  $2 \leq l \leq w_i$,
 i.e., $t^i_1 = \hterm(f_i)$.
As $\hterm( \hterm(f_i) \rmult m_i) \geq \hterm(f_i) \succ t_l^i$, $2 \leq l \leq w_i$,
 again by Definition \ref{def.refined.ordering} we can conclude that
 $\hterm(\hterm(f_i) \rmult m_i) \succ t^i_l \rmult m_i \succeq 
  \hterm(t^i_l \rmult m_i)$, $2 \leq l \leq w_i$, and
 $\hterm(f_i) \rmult m_i \succ \sum_{l=2}^{w_i} t_l^i \rmult m_i$.
Then for each $s_j^i$, $2 \leq j \leq k_i$,
 there exists $t_l^i \in \terms(f_i)$ such that
 $s \in \supp(t_l^i \rmult m_i)$.
Since $\hterm(f) \succ s_j^i$ and even $\hterm(f) \succ t_l^i \rmult m_i$ we find that
 either
 $\hterm(f \rmult m) \succeq \hterm((t_l^i \rmult m_i) \rmult m) = \hterm(t_l^i \rmult (m_i \rmult m))$
 in case $\hterm(t_l^i \rmult m_i) = \hterm(f_i \rmult m_i)$
 or 
 $\hterm(f \rmult m) \succ (t_l^i \rmult m_i) \rmult m = t_l^i \rmult (m_i \rmult m)$.
% as $\f$ is associative.
Hence we can conclude $f_i \rmult \tilde{m}^i_r \predeq \hterm(f \rmult m)$, $1 \leq r \leq k_i$ and
 for at least one $\tilde{m}^i_r$ we get
 $\hterm(f_i \rmult \tilde{m}^i_r) = \hterm(f_i \rmult m_i \rmult m) \geq \hterm(f_i)$.
\\
It remains to analyze the situation for the function
 $(\sum_{i=k+1}^n f_i \rmult m_i) \rmult m$.
Again we find that for all terms $s$ in the $f_i \rmult m_i$, $k+1 \leq i \leq n$,
 we have $\hterm(f) \succ s$ and we get 
 $\hterm(f \rmult m) \succ\hterm(s \rmult m)$.
Hence all polynomial multiples of the $f_i$  in the representation 
 $\sum_{i=k+1}^n \sum_{j=1}^{k_i} f_i \rmult \tilde{m}^i_j$, where
 $m_i \rmult m = \sum_{j=1}^{k_i}\tilde{m}^i_j$, are bounded by 
 $\hterm(f \rmult m)$.
\\
\qed

Now we are able to prove Theorem \ref{theo.loc}.

\Ba{of Theorem \ref{theo.loc}}~\\
The proof is basically the same as for Theorem \ref{theo.gb.reduction}.
Due to Lemma \ref{lem.red.reps_rr} we can substitute the multiples $f_j \rmult m_j$
 by appropriate representations without changing $(t,K)$.
Hence, we only have to ensure that despite testing less polynomials we are able to
 apply our induction hypothesis.
Taking the notations from the proof of Theorem \ref{theo.gb.reduction},
let us first check the situation for m-polynomials.
\\
Let $\sum_{j=1}^K \hm(f_j \rmult (w_{j} \skm \alpha_j)) = \zero$.
Then by Definition \ref{def.gpol} we know
 $(\alpha_1, \ldots, \alpha_K) \in M$, as 
 $\sum_{j=1}^K \hc(f_j \rmult w_{j}) \skm \alpha_j = 0$.
Hence there are $\delta_1, \ldots, \delta_K \in \rr$ such that
 $\sum_{i=1}^l A_i \skm \delta_i = (\alpha_1, \ldots, \alpha_K)$ for some $l \in \n$,
 $A_i=(\alpha_{i,1}, \ldots, \alpha_{i,K}) \in \{ A_j \mid j \in I_M \}$,
 and $\alpha_j = \sum_{i=1}^l \alpha_{i,j} \skm \delta_i$, $1 \leq j \leq K$.
There are module polynomials
 $h_i = \sum_{j=1}^K f_j \rmult w_j \skm \alpha_{i,j}$,$1 \leq i \leq l$
 and by our assumption there are polynomials $h_i'$
 in the stable localization such that $h_i \red{}{\myr}{r}{h_i'} \zero$.
Moreover, $h_i' \red{*}{\myr}{r}{F} \zero$.
Then by Lemma \ref{lem.red.reps_rr} the m-polynomials $h_i$ all have representations
 bounded by $t$.
Again we get
\begin{eqnarray}
\sum_{j=1}^K f_j \rmult (w_{j} \skm \alpha_j) & = & \sum_{j=1}^K f_j \rmult w_{j} \skm (\sum_{i=1}^l \alpha_{i,j} \skm \delta_i) \nonumber\\
         &=& \sum_{j=1}^K \sum_{i=1}^l (f_j \rmult w_{j} \skm \alpha_{i,j}) \skm \delta_i \nonumber\\
         &=& \sum_{i=1}^l (\sum_{j=1}^K f_j \rmult w_{j} \skm \alpha_{i,j}) \skm \delta_i \nonumber\\
         &=& \sum_{i=1}^l h_i \skm \delta_i \nonumber
\end{eqnarray}
 and we can change the representation of $g$ to
 $\sum_{i=1}^l h_i \skm \delta_i + \sum_{j=K+1}^m f_{j} \rmult (w_{j} \skm \alpha_j)$
 and replace each $h_i$ by the respective special standard representation in terms of $F$.
Remember that for all $h_i$, $1 \leq i \leq l$ we have $\hterm(h_i) \pred t$.
Hence, for this new representation we now have maximal term smaller than $t$
 and by our induction hypothesis $g$ is reducible by $F$ and we are done.
\\
It remains to study the case where $\sum_{j=1}^K \hm(f_j \rmult (w_j \skm \alpha_j)) \neq 0$.
Then
 we have $\hterm(f_1 \rmult (w_1 \skm \alpha_1) + \dots + f_K \rmult (w_K \skm \alpha_K)) = t = \hterm(g)$,
 $\hc(g) = \hc(f_1 \rmult (w_1 \skm \alpha_1) + \dots + f_K \rmult (w_K \skm \alpha_K)) \in \ideal{r}{}(\{ \hc(f_1 \rmult w_1), \ldots, 
 \hc(f_K \rmult w_K) \})$ and 
 even $\hm(f_1 \rmult (w_1 \skm \alpha_1) + \dots + f_K \rmult (w_K \skm \alpha_K)) = \hm(g)$.
Hence $\hc(g)$ is $\R$-reducible  by some $\alpha$, $\alpha \in G$, $G$ being a
 (weak) right Gr\"obner basis of $\ideal{r}{}(\{ \hc(f_1 \rmult w_1), \ldots, 
 \hc(f_K \rmult w_K) \})$ in $\rr$ with respect to the reduction relation $\R$.
Let $g_{\alpha}$ be the respective g-polynomial corresponding to $\alpha$ and $t$.
Then we know that $g_{\alpha}\red{}{\myr}{r}{g_{\alpha}'} \zero$ for some
 $g_{\alpha}'$ in the stable localization and $g_{\alpha}'\red{*}{\myr}{r}{F} \zero$.
Moreover, we know that the head monomial of $g_{\alpha}'$ is reducible
 by some polynomial $f \in F$ and we assume
 $\hterm(g_{\alpha}) = \hterm(\hterm(f) \rmult m) = \hterm(f \rmult m) \geq \hterm(f)$
 and $\hc(g_{\alpha}) \R_{\hc(f \rmult m)}$.
Then, as $\hc(g)$ is $\R$-reducible by $\hc(g_{\alpha})$, 
 $\hc(g_{\alpha})$ is $\R$-reducible by $\hc(g_{\alpha}')$, 
 $\hc(g_{\alpha}')$ is $\R$-reducible to zero and (A4) holds,
 the head monomial of $g$ is also reducible by some $f' \in F$ and we are done. 
\\
\qed
Again, if for infinite $F$ we can assure that the reduction relation is Noetherian,
 the proof still holds.

%%% Local Variables: 
%%% mode: latex
%%% TeX-master: "testlauf"
%%% End: 

%% file: generalization_z.tex
%%%%%%%%%%%%%%%%%%%%%%%%%%%%%%%%%%%%%%%%%%%%%%%%%%%%%%%%%
\subsection{Function Rings over the Integers}\label{section.right.integers}
In the previous section we have seen that for the reduction relations for
 $\f$ as defined in Definition \ref{def.rred_rr} and \ref{def.rred.rr}
 the Translation Lemma no longer holds.
This is due to the fact that the first definition is based on divisibility
 in $\rr$ and hence too weak and the second definition is based on the abstract
 notion of the reduction relation $\R$ and hence there is not
 enough information on the reduction step involving the coefficient.

When studying special reduction rings where we have more information
 on the specific reduction relation $\R$ the situation often can
 be improved.
Here we want to go into the details for the case that
 $\rr$ is the ring of the integers $\z$.
Remember that there are various ways of defining a reduction relation
 for the integers.
In Example \ref{exa.Z} two possibilities are presented.
Here we want to use the second one based on division with remainders
 in order to introduce a reduction
 relation to $\f_{\z}$.
We follow the ideas presented in \cite{MaRe93a} for characterizing 
 prefix Gr\"obner
 bases in monoid rings $\z[\m]$ where $\m$ is presented by a
 finite convergent string rewriting system.

In order to use elements of $\f_{\z}$ as rules for a reduction 
 relation we need an ordering
 on $\z$.
We  specify a total well--founded ordering on $\z$ as
 follows\footnote{If not stated otherwise $<$ is the usual ordering on $\z$, i.e.~$\ldots < -3 < -2 < -1 < 0 <1 <2<3 \ldots$.}:
% $a<_{Z} b$ iff $|a|<|b|$ or $(|a|=|b|, a>0, b<0)$, i.e.~$0<_Z 1<_Z -1<_Z  \ldots$.
\[ \alpha <_Z \beta \mbox{ iff }
  \left\{ \begin{array}{l} \alpha \geq 0 \mbox{ and } \beta<0\\ 
                           \alpha \geq 0 , \beta >0 \mbox{ and } \alpha < \beta \\
                           \alpha <0, \beta <0 \mbox{ and } \alpha > \beta \end{array} \right. \]
and $\alpha \leq_Z \beta$ iff $\alpha=\beta$ or $\alpha <_Z \beta$.
Hence we get $0 \leq_Z 1 \leq_Z 2 \leq_Z 3 \leq_Z \ldots \leq_Z -1 \leq_Z -2 \leq_Z -3 \leq_Z \ldots$.
Then we can make the following important observation:
Let $\gamma \in \n$.
We call the positive numbers $0, \ldots, \gamma -1$ the remainders of $\gamma$.
Then for each $\delta \in \z$ there are unique $\alpha,\beta \in \z$
 such that $\delta=\alpha \cd \gamma + \beta$ and $\beta$
 is a remainder of $\gamma$.
We get $\beta<\gamma$ and in case $\delta>0$ and $\alpha \not = 0$ even $\gamma \leq \delta$.
Further $\gamma$ does not divide $\beta_1 - \beta_2$, if $\beta_1, \beta_2$ are different remainders of $\gamma$.

As we will later on only use polynomials with head coefficients in $\n$ for
 reduction, we will mainly require the part of the ordering on $\n$ which then
 coincides with the natural ordering on this set.
Then we will drop the suffix\footnote{In the literature other orderings on 
 the integers are used by Buchberger and Stifter \cite{St87} and
 Kapur and Kandri-Rody \cite{KaKa88}. They then have to consider
 s- and t-polynomials as 
 critical situations.}.

This ordering $<_{Z}$ can be used to induce an ordering on $\f_{\z}$ as follows:
for two elements $f,g$ in $\f$ we define
 $f \succ g$ iff $\hterm(f) \succ \hterm(g)$ or ($(\hterm(f) = \hterm(g)$ and
 $\hc(f) >_Z \hc(g)$)
 or ($(\hm(f) = \hm(g)$ and $\reductum(f) \succ \reductum(g)$).

The reduction relation presented in Definition \ref{def.rred.rr}
 now can be adapted to this special case:
Let $\R$ be our reduction relation on $\z$ where $\alpha \R_{\gamma} \beta$, if
 $\gamma > 0$ and for some $\delta \in \z$ we 
 have $\alpha = \gamma \skm \delta + \beta$ with $0 \leq \beta < \gamma$,
 i.e.~$\beta$ is the remainder of $\alpha$ modulo $\gamma$.

\begin{definition}\label{def.rred.z}~\\
{\rm
Let $p$, $f$ be two non-zero polynomials in $\f_{\z}$. 
We say $f$ \betonen{right reduces} $p$ \betonen{to} $q$ \betonen{at a monomial} $\alpha \skm t$
 \betonen{in one step},
 i.e. $p \red{}{\myr}{r}{f} q$, if there exists $ s \in \terms(\f_{\z})$ such that
\begin{enumerate}
\item $t \in \supp(p)$ and $p(t) = \alpha$,
\item $\hterm(\hterm(f)\rmult s) = \hterm(f \rmult s) = t \geq \hterm(f)$, 
\item $\alpha \geq_{\z} \hc(f \rmult m) > 0$ and
 	 $\alpha \R_{\hc(f \rmult s)} \delta$ where
	 $\alpha = \hc(f\rmult s)\skm \beta + \delta$ with 
         $\beta, \delta \in \z$,
         $0 \leq \delta < \hc(f \rmult s)$, and
\item $q = p -  f \rmult m$ where $m = \beta \skm s$.
\end{enumerate}
We write $p \red{}{\myr}{r}{f}$ if there is a polynomial $q$ as defined
above and $p$ is then called right reducible by $f$. 
%\\
Further, we can define $\red{*}{\myr}{r}{}, \red{+}{\myr}{r}{}$ and
 $\red{n}{\myr}{r}{}$ as usual.
%\\
Right reduction by a set $F \subseteq \f \backslash \{ \zero \}$ is denoted by
 $p \red{}{\myr}{r}{F} q$ and abbreviates $p \red{}{\myr}{r}{f} q$
 for some $f \in F$.
\dend
}
\end{definition}
As before, for this reduction relation we can still have $t \in \supp(q)$.
Hence other arguments than those used in the proof of Lemma \ref{lem.rred.rr} have
 to be used to show termination.
The important part now is that if we still have $t \in \supp(q)$
 then its coefficient will be smaller according to our ordering $<_{\z}$ chosen
 for $\z$ and since this ordering is well-founded we are done. 
Notice that in contrary to Lemma \ref{lem.sred.rr} we do not have to restrict
 ourselves to finite sets of polynomials in order to ensure termination.

The additional information we have on the coefficients before
 and after the reduction step now enables us to prove
 an analogon of the Translation Lemma for function rings over the integers.
The first and second part of the lemma are only needed to prove the essential
 third part.

\begin{lemma}\label{lem.trans.integers}~\\
{\sl
Let $F$ be a set of polynomials in $\f_{\z}$ and $p,q,h$ polynomials in $\f_{\z}$.
\begin{enumerate}
\item
Let $p-q \red{}{\myr}{r}{F} h$ such that the reduction step takes place
 at the monomial $\alpha \skm t$ and we additionally have $t \not\in \supp(h)$.
Then there exist $p', q' \in \f_{\z}$ such that $p \red{*}{\myr}{r}{F} p'$
 and $q \red{*}{\myr}{r}{F} q'$ and $h = p'-q'$.
\item 
Let $\zero$ be the unique normal form of $p$ with respect to $F$ and $t = \hterm(p)$.
Then there exists a polynomial $f \in F$ such that $p \red{}{\myr}{r}{f} p'$
 and $t \not \in \supp(p')$.
\item
Let $\zero$ be the unique normal form of $p-q$ with respect to $F$.
Then there exists $g \in \f_{\z}$ such that $p \red{*}{\myr}{r}{F} g$ and
 $q \red{*}{\myr}{r}{F} g$.
\end{enumerate}
\lemend}
\end{lemma}
\Ba{}~\\
\begin{enumerate}
\item
Let $p-q \red{}{\myr}{r}{F} h$ at the monomial $\alpha \skm t$, i.e.,
 $h = p-q - f \rmult m$ for some $m= \beta \skm s \in \monoms(\f_{\z})$ such that
 $\hterm(\hterm(f)\rmult s) = \hterm(f \rmult s) = t \geq \hterm(f)$
 and $\hc(f \rmult s) > 0$.
Remember that $\alpha$ is the coefficient
 of $t$ in $p-q$.
Then as $t \not\in \supp(h)$ we know $\alpha = \hc(f \rmult m)$.
Let $\alpha_1$ respectively $\alpha_2$ be the coefficients of $t$ in
 $p$ respectively $q$ and $\alpha_1 = \hc(f \rmult m) \skm \beta_1 + \gamma_1$
 respectively $\alpha_2 = \hc(f \rmult m) \skm \beta_2 + \gamma_2$ for some
 $\beta_1, \beta_2, \gamma_1, \gamma_2 \in \z$ where $0 \leq \gamma_1, \gamma_2
 < \hc(f \rmult s) \leq \hc(f \rmult m)$.
Then $\alpha = \hc(f \rmult m) = \alpha_1 - \alpha_2 =
 \hc(f \rmult m) \skm (\beta_1 - \beta_2) + (\gamma_1 - \gamma_2)$, and
 as $\gamma_1 - \gamma_2$ is no multiple of $\hc(f \rmult m)$ we have
 $\gamma_1 - \gamma_2 = 0$ and hence $\beta_1 - \beta_2 = 1$.
We have to distinguish two cases:
\begin{enumerate}
\item $\beta_1 \neq 0$ and $\beta_2 \neq 0$:
      Then $p \red{}{\myr}{r}{F} p -  f \rmult m \skm \beta_1 = p'$,
        $q \red{}{\myr}{r}{F} q - f \rmult m \skm \beta_2= q'$ and
        $p'-q' = p - f \rmult m \skm \beta_1 - q + f \rmult m \skm \beta_2=
        p-q - f \rmult m = h$.
\item $\beta_1 = 0$ and $\beta_2 = - 1$ (the case $\beta_2 = 0$
       and $\beta_1 = 1$ being symmetric):
      Then $p' = p$, 
        $q \red{}{\myr}{r}{F} q - f \rmult m \skm \beta_2 = 
         q + f \rmult m \skm \beta = q'$ and
        $p'-q' = p-q - f \rmult m = h$.
\end{enumerate}
\item Since $p \red{*}{\myr}{r}{F} \zero$, $\hm(p) = \alpha \skm t$
 must be $F$-reducible.
Let $f_1, \ldots, f_k \in F$ be all polynomials in $F$
 such that $\alpha \skm t$ is reducible by them.
Let $m_1, \ldots m_k$ be the respective monomials involved in possible reduction
 steps.
Moreover, let $\gamma = \min_{1 \leq i \leq k} \{ \hc(f_i \rmult m_i) \}$
 and without loss of generality $\hm(f \rmult m) = \gamma \skm t$ for some
 $f \in F$, $\hterm(\hterm(f) \rmult m) = \hterm(f \rmult m) \geq \hterm(f)$.
We claim that for $p \red{}{\myr}{r}{f_1} p - f \rmult m = p'$
 we have $t \not \in \supp(p')$.
Suppose $\hterm (p') = t$.
Then by our definition of reduction we must have $0 < \hc(p') < \hc(f \rmult m)$.
But then $p'$ would no longer be $F$-reducible contradicting our assumption
 that $\zero$ is the unique normal form of $p$.
\item Since $\zero$ is the unique normal form of $p-q$ by 2.~there exists
 a reduction sequence $p-q \red{}{\myr}{r}{f_{i_1}} h_1 \red{}{\myr}{r}{f_{i_2}}
 \ldots \red{}{\myr}{r}{f_{i_k}} \zero$ such that for the head terms we get $\hterm(p-q) \succ \hterm(h_1)
 \succ \ldots$.
We show our claim by induction on $k$, where $p-q \red{k}{\myr}{r}{F} \zero$
 is such a reduction sequence.
In the base case $k=0$ there is nothing to show as then $p=q$.
Hence, let $p-q \red{}{\myr}{r}{F} h \red{k}{\myr}{r}{F} \zero$.
Then by 1.~there are polynomials $p',q' \in \f_{\z}$ such that 
 $p \red{*}{\myr}{r}{F} p'$
 and $q \red{*}{\myr}{r}{F} q'$ and $h = p'-q'$.
Now the induction hypothesis for $p'-q' \red{k}{\myr}{r}{F} \zero$ yields the
 existence of a polynomial $g \in \f_{\z}$ such that
 $p \red{*}{\myr}{r}{F} g$ and $q \red{*}{\myr}{r}{F} g$.
\end{enumerate}
\qed 

Hence weak Gr\"obner bases are in fact Gr\"obner bases and can be characterized as
 follows:

\begin{definition}~\\
{\rm
A set $F \subseteq \f_{\z} \backslash \{ \zero \}$ is called a (weak)
 right Gr\"obner basis of $\ideal{r}{}(F)$ if for all
 $g \in \ideal{r}{}(F)$ we have $g \red{*}{\myr}{r}{F} \zero$. 
\dend
}
\end{definition}
\begin{corollary}\label{cor.standard.rep}~\\
{\sl
Let $F$ be a set of polynomials in $\f_{\z}$ 
 and $g$ a non-zero polynomial in $\ideal{r}{}(F)$
 such that $g \red{*}{\myr}{r}{F} \zero$.
Then $g$ has a representation of the form 
$$g = \sum_{i=1}^n f_i \rmult m_i,
 f_i \in F, m_i \in \monoms(\f_{\z}), n \in \n$$
 such that
 $\hterm(g) = \hterm(\hterm(f_i) \rmult m_i) =
 \hterm(f_i \rmult m_i) \geq \hterm(f_i)$, $1 \leq i \leq k$, and
 $\hterm(g) \succ \hterm(f_i \rmult m_i) = 
 \hterm(\hterm(f_i) \rmult m_i)$ 
 for all $k+1 \leq i \leq n$.
}
\end{corollary}
\Ba{}~\\
We show our claim by induction on $n$ where $g \red{n}{\myr}{r}{F} \zero$.
If $n=0$ we are done.
Else let $g \red{1}{\myr}{r}{F} g_1 \red{n}{\myr}{r}{F} \zero$.
In case the reduction step takes place at the head monomial,
 there exists a polynomial $f \in F$
 and a monomial $m= \beta \skm s \in \monoms(\f)$
 such that $\hterm(\hterm(f) \rmult s) = \hterm(f \rmult s) = \hterm(g)
 \geq \hterm(f)$ and $\hc(g) \R_{\hc(f \rmult s)} \delta$ with
 $\hc(g) = \hc(f\rmult s) \skm \beta + \delta$ for some $\beta, \delta \in \z$,
 $0 \leq \delta < \hc(f \rmult s)$.
Moreover the induction hypothesis then is applied to  $g_1 = 
 g - f \rmult m$.
If the reduction step takes place at a monomial with term smaller $\hterm(g)$
 for the respective monomial multiple $f \rmult m$ we immediately get
 $\hterm(g) \succ \hterm(f \rmult m)$ and we can apply our induction hypothesis
 to the resulting polynomial $g_1$.
In both cases we can arrange the monomial multiples $f \rmult m$ arising from
 the reduction steps in such a way that gives us the desired representation.
\\
\qed 
We can even state that $\hc(g) \red{*}{\R}{}{\{\hc(f_i \rmult m_i)\mid 1 \leq i \leq k\}}0$.
Now right Gr\"obner bases
 can be characterized using the concept of 
 s-polynomials combined with the technique of saturation
 which is necessary in order to describe the
 whole right ideal congruence by the reduction relation.
\begin{definition}\label{def.s-poly.z}~\\
{\rm
Let $p_{1}, p_{2}$ be  two polynomials in $\f_{\z}$.
If there are respective terms
 $t,u_1, u_2 \in \myt$ such that
 $\hterm(\hterm(p_i) \rmult u_i) = \hterm(p_i \rmult u_i)=t
 \geq \hterm(p_i)$ let $HC(p_{i} \rmult u_i) = \gamma_i$.\\
Assuming $\gamma_1 \geq \gamma_2 > 0$\footnote{Notice that $\gamma_i > 0$
 can always be achieved by studying the situation for $- p_i$ in case 
 we have $HC(p_{i} \rmult u_i) < 0$.},
 there are $\beta, \delta \in \z$ such that
 $\gamma_1 = \gamma_2 \skm \beta + \delta$ and $0 \leq \delta < \gamma_2$
 and we get the following s-polynomial
      $$\spol{r}(p_{1}, p_{2},t,u_1, u_2) =  p_2 \rmult u_2 \skm \beta -  p_1 \rmult u_1.$$
The set $\spols(\{p_1,p_2\})$ then is the set of all such
 s-polynomials corresponding to $p_1$ and $p_2$.
\dend
}
\end{definition}
These sets can be infinite\footnote{This is due to the fact that in general we cannot always find finite locations for $t$. One well-studied field are monoid rings.}.
\begin{theorem}\label{theo.gb.reduction.z}~\\
{\sl
Let $F$ be a  set of polynomials in $\f_{\z} \backslash \{ \zero \}$.
Then $F$ is a right Gr\"obner basis of $\ideal{r}{}(F)$ if and only if
\begin{enumerate}
\item for all $f$ in $F$ and for all $m$ in $\monoms(\f_{\z})$ we have 
       $f \rmult m \red{*}{\myr}{r}{F} \zero$, and
\item all s-polynomials corresponding to $F$ as specified in
       Definition \ref{def.s-poly.z}
       reduce to $\zero$ using $F$.
\end{enumerate}
\theoend
}
\end{theorem}
\Ba{}~\\
In case $F$ is a right Gr\"obner basis,
 since the multiples  $f \rmult m$ and the respective s-polynomials
 are all elements of $\ideal{r}{}(F)$ they must
 reduce to zero using $F$.
\\
The converse will be proven by showing that every element in
 $\ideal{r}{}(F)$ is reducible  by $F$.
Then as $g \in \ideal{r}{}(F)$ and $g \red{}{\myr}{r}{F} g'$ implies
 $g' \in \ideal{r}{}(F)$ we have $g\red{*}{\myr}{r}{F} \zero$.
Notice that this is sufficient as the reduction relation $\red{}{\myr}{r}{F}$
 is Noetherian.
\\
Let $g \in \ideal{r}{}(F)$ have a representation in terms of $F$ of
 the following form:
$g = \sum_{j=1}^m f_{j} \rmult w_{j} \skm \alpha_j$ such that
 $f_j \in F$, $w_{j} \in \myt$ and $\alpha_j \in \z$.
Depending on this  representation of $g$ and the
 well-founded total ordering $\succeq$ on $\myt$ we define
 $t = \max_{\succeq} \{ \hterm(f_{j} \rmult w_{j}) \mid 1\leq j \leq m \}$, 
 $K$  as the number of polynomials $f_j \rmult w_j$ with head term $t$, and
 $M = \{\{ \hc(f_i \rmult w_i) \mid \hterm(f_j \rmult w_j) = t \}\}$ a multiset
 in $\z$.
We show our claim by induction on $(t,M)$, where
 $(t',M')<(t,M)$ if and only if $t' \prec t$ or $(t'=t$ and $M' \ll M)$\footnote{We 
 define $M' \ll M$ if $M$ can be transformed into $M'$ by substituting
 elements in $M$ with sets of smaller elements (with respect to our ordering on the integers.}.
\\
Since by our first assumption every multiple $f_j \rmult w_j$ in this
 sum reduces to zero using $F$ and hence
 has a  representation as specified
 in Corollary \ref{cor.standard.rep}, we can assume that
 $\hterm(\hterm(f_j) \rmult w_j) = \hterm(f_j \rmult w_j) \geq
 \hterm(f_j)$ holds.
Moreover, without loss of generality we can assume that the polynomial multiples
 with head term $t$ are just $f_1  \rmult w_1, \ldots , f_K \rmult w_K$ and
 additionally we can assume $\hc(f_j \rmult w_j) > 0$\footnote{This can easily be achieved
 by adding $-f$ to $F$ for all $f \in F$ and using $(-f_j) \rmult w_j$
 in case $\hc(f_j \rmult w_j) < 0$.}.
\\
Obviously, $t \succeq \hterm(g)$ must hold. 
If $K = 1$ this gives us $t = \hterm(g)$ and even
 $\hm(g) = \hm(f_1 \rmult w_1 \skm \alpha_1)$, implying that $g$ is right
 reducible at $\hm(g)$ by $f_1$.
\\
Hence let us assume $K>1$.
\\
Without loss of generality we can assume that $\hc(f_1 \rmult w_1) \geq
 \hc(f_2 \rmult w_2) > 0$ and there are $\alpha, \beta \in \z$ such that
 $\hc(f_2 \rmult w_2) \skm \alpha + \beta = \hc(f_1 \rmult w_1)$ and
 $\hc(f_2 \rmult w_2) > \beta \geq 0$.
Since $t = \hterm(f_1 \rmult w_1) =  \hterm(f_2 \rmult w_2)$ by Definition
 \ref{def.s-poly.z} we have an s-polynomial $\spol{r}(f_1,f_2,t,w_1,w_2)
 = f_2 \rmult w_2 \skm \alpha - f_1 \rmult w_1$.
If $\spol{r}(f_1,f_2,t,w_1,w_2) \neq \zero$\footnote{In case $\spol{r}(f_1,f_2,t,w_1,w_2) = \zero$
 the proof is similar. We just have to substitute $\zero$ in the equations below which
 immediately gives us a smaller representation of $g$.}
 then $\spol{r}(f_1,f_2,t,w_1,w_2)\red{*}{\myr}{r}{F} \zero$ implies
 $\spol{r}(f_1,f_2,t,w_1,w_2) = \sum_{i=1}^k \delta_i \skm h_i \rmult v_i$,
 $\delta_i \in \z$, $h_i \in F$, $v_i \in \myt$ where this
 sum is a representation in the sense of Corollary \ref{cor.standard.rep}
 with terms bounded by $\hterm(\spol{r}(f_1,f_2,t,w_1,w_2)) \leq t$.
This gives us
\begin{eqnarray}
g & = & f_1 \rmult w_1 \skm \alpha_1 + f_2 \rmult w_2 \skm \alpha_2 + \sum_{j=3}^m f_j \rmult w_j \skm \alpha_j \\ \nonumber
  & = & f_1 \rmult w_1 \skm \alpha_1 + \underbrace{f_2 \rmult w_2 \skm \alpha_1 \skm \alpha - f_2 \rmult w_2 \skm \alpha_1 \skm \alpha}_{=\zero} +
        f_2 \rmult w_2 \skm \alpha_2 + \sum_{j=3}^m f_j \rmult w_j \skm \alpha_j \\ \nonumber
  & = & f_2 \rmult w_2 \skm (\alpha_1 \skm \alpha + \alpha_2) -  \underbrace{(f_2 \rmult  w_2 \skm \alpha  - f_1 \rmult w_1}_{=\spol{r}(f_1,f_2,t,w_1,w_2)}
              \skm \alpha_1 + \sum_{j=3}^m f_j \rmult w_j \skm \alpha_j \\ \nonumber
  & = & f_2 \rmult w_2 \skm (\alpha_1 \skm \alpha + \alpha_2) - (\sum_{i=1}^k \delta_i \skm h_i \rmult v_i) \skm \alpha_1
         + \sum_{j=3}^m f_j \rmult w_j \skm \alpha_j \\ \nonumber
\end{eqnarray}
and depending on this new representation of $g$ we define
 $t' = \max_{\succeq} \{ \hterm(f_{j} \rmult w_{j}),
  \hterm(h_{j} \rmult v_{j}) \mid f_j, h_j \mbox{ appearing in the new representation } \}$, and
 $M' = \{\{ \hc(f_i \rmult w_i), \hc(h_{j} \rmult v_{j}) \mid \hterm(f_j \rmult w_j) = \hterm(h_{j} \rmult v_{j}) = t' \}\}$
 and we either get $t' \pred t$ and have a smaller representation for $g$
 or in case $t'=t$ we have to distinguish two cases
 \begin{enumerate}
 \item $\alpha_1 \skm \alpha   + \alpha_2= 0$. \\
       Then $M' = M - \{\{ \hc(f_1 \rmult w_1), \hc(f_2 \rmult w_2) \}\} \cup 
        \{\{ \hc(h_j \rmult v_j) \mid \hterm(h_j \rmult v_j) = t \}\}$.
       As those polynomials $h_j$ with $\hterm(h_j \rmult v_j) = t$ are used to right reduce the
       monomial $\beta \skm t = \hm(\spol{r}(f_1,f_2,t,w_1,w_2))$ we know that for them we have
       $0 < \hc(h_j \rmult v_j) \leq \beta < \hc(f_2 \rmult w_2) \leq 
     \hc(f_1 \rmult w_1)$. Hence $M' \ll M$ and we have a smaller representation for $g$.
 \item $\alpha_1 \skm \alpha   + \alpha_2 \neq 0$. \\
       Then $M' = (M - \{\{ \hc(f_1 \rmult w_1) \}\}) \cup 
     \{ \{ \hc(h_j \rmult v_j) \mid \hterm(h_j \rmult v_j) = t \}\}$.
       Again  $M' \ll M$ and we have a smaller representation for $g$.
 \end{enumerate}
Notice that the case $t'=t$ and $M'\ll M$ cannot occur infinitely often but has to result
 in either $t' <t$ or will lead to $t'=t$ and $K=1$ and hence to reducibility by $\red{}{\myr}{r}{F}$.
\\
\qed
Now the question arises when the critical situations in this characterization
 can be localized to subsets of the respective sets
 as in Theorem \ref{theo.s-pol.2}.
Reviewing the Proof of Theorem \ref{theo.s-pol.2} we find that 
 Lemma \ref{lem.red.reps} is central as it describes when multiples of
 polynomials which have a right reductive standard representation in terms
 of some set $F$ again have such a representation.
As we have seen before, this will not hold for function rings over reduction
 rings in general.
As in Section \ref{section.right.rr}, to give localizations of Theorem \ref{theo.gb.reduction.z} the
 concept of stable subsets is sufficient:
\begin{definition}\label{def.stable.loc.z}~\\
{\rm
A set $C \subset S \subseteq \f_{\z}$ is called a \betonen{stable localization} of
 $S$ if for every $g \in S$ there exists $f \in C$
 such that $g \red{}{\myr}{r}{f} \zero$.
% and if $\hm(f)$ is reducible by $F$, so is $\hm(g)$.
\dend
}
\end{definition}
Stable localizations for the sets of s-polynomials again arise from the appropriate
 sets of least common multiples as presented on page \ref{page.lcm}.
In case $\f_{\z}$ and $\red{}{\myr}{r}{}$ allow such stable localizations,
 we can rephrase Theorem \ref{theo.gb.reduction.z} as follows:
\begin{theorem}\label{theo.loc.z}~\\
{\sl
Let $F$ be a set of polynomials in $\f_{\z} \backslash \{ \zero \}$.
Then $F$ is a right Gr\"obner basis of $\ideal{r}{}(F)$ if and only if
\begin{enumerate}
\item for all $s$ in a stable localization of 
       $\{ f \rmult m \mid f \in \f_{\z}, m \in \monoms(\f_{\z}) \}$ we have 
       $s \red{*}{\myr}{r}{F} \zero$, and
\item for all $h$ in a stable localization of the s-polynomials corresponding to $F$ as specified in
       Definition \ref{def.s-poly.z} we have 
       $h \red{*}{\myr}{r}{F} \zero$.
\end{enumerate}
\theoend
}
\end{theorem} 
When proving Theorem \ref{theo.loc.z}, we can substitute the critical situation using an analogon of
 Lemma \ref{lem.red.reps}, which  will be sufficient to make the representation used in the proof smaller.
It is a direct consequence of Lemma \ref{lem.red.reps_rr}.
\begin{corollary}\label{cor.red.reps.z}~\\
{\sl
Let  $F \subseteq \f_{\z} \backslash \{ \zero\}$ and
 $f$, $p$  non-zero polynomials in $\f_{\z}$.
If $p \red{}{\myr}{r}{f} \zero$ and $f \red{*}{\myr}{r}{F} \zero$,
 then $p$ has a representation of the form
$$p = \sum_{i=1}^n f_i \rmult l_i,
 f_i \in F, l_i \in \monoms(\f_{\z}), n \in \n$$
 such that
 $\hterm(p) = \hterm(\hterm(f_i) \rmult l_i) =
 \hterm(f_i \rmult l_i) \geq \hterm(f_i)$ for $1 \leq i \leq k$ and
 $\hterm(p) \succ \hterm(f_i \rmult l_i)$ 
 for all $k+1 \leq i \leq n$ (compare Definition \ref{def.right_reductive}).
\corend
}
\end{corollary}

\Ba{Theorem \ref{theo.loc.z}}~\\
The proof is basically the same as for Theorem \ref{theo.gb.reduction.z}.
Due to Corollary \ref{cor.red.reps.z} we can substitute the multiples $f_j \rmult w_j$
 by appropriate representations.
Hence, we only have to ensure that despite testing less polynomials we are able to
 apply our induction hypothesis.
Taking the notations from the proof of Theorem \ref{theo.gb.reduction.z},
let us  check the situation for $K>1$.
\\
Without loss of generality we can assume that $\hc(f_1 \rmult w_1) \geq
 \hc(f_2 \rmult w_2) > 0$ and there are $\alpha, \beta \in \z$ such that
 $\hc(f_2 \rmult w_2)\skm \alpha  + \beta = \hc(f_1 \rmult w_1)$ and
 $\hc(f_2 \rmult w_2) > \beta \geq 0$.
Since $t = \hterm(f_1 \rmult w_1) =  \hterm(f_2 \rmult w_2)$ by Definition
 \ref{def.s-poly.z} we have an s-polynomial $h \in \spols(f_1,f_2)$ and $m \in \monoms(\f_{\z})$ such that
 $h \rmult m  = \alpha \skm f_2 \rmult w_2 - f_1 \rmult w_1$.
If $h \neq \zero$\footnote{In case $h = \zero$
 the proof is similar. We just have to substitute $\zero$ in the equations below which
 immediately gives us a smaller representation of $g$.}
 then by Corollary \ref{cor.red.reps.z}
 $f_2 \rmult w_2 \skm \alpha  - f_1 \rmult w_1 \red{}{\myr}{r}{h} \zero$ and
 $h\red{*}{\myr}{r}{F} \zero$ imply
 $f_2 \rmult w_2 \skm \alpha - f_1 \rmult w_1 = \sum_{i=1}^k  h_i \rmult v_i \skm \delta_i$,
 $\delta_i \in \z$, $h_i \in F$, $v_i \in \myt$ where this
 sum is a representation in the sense of Corollary \ref{cor.standard.rep}
 with terms bounded by $\hterm(h \rmult m) \leq t$.
As in the proof of Theorem \ref{theo.gb.reduction.z} we now can use this bounded
 representation to get a smaller representation of $g$ and are done.
\\
\qed

%%%%%%%%%%%%%%%%%%%%%%%%%%Examples%%%%%%%%%%%%%%%%%%%%%%%%%%%%%%%%%%%%%%%%%%%%%%%%%%
We close this subsection by outlining how different structures known
 to allow finite Gr\"obner bases can be interpreted as function rings.
Using the respective interpretations the terminology can be
 adapted at once to the respective structures and in general the resulting
 characterizations of Gr\"obner bases coincide with the
 results known from literature.
%%%%%%%%%%%%%%%%%%%%%%%
\subsubsection{Polynomial Rings}
A commutative polynomial ring $\z[x_1, \ldots, x_n]$ is a function
 ring according to the following interpretation:
 \begin{itemize}
 \item $\myt$ is the set of terms $\{ x_1^{i_1} \ldots x_n^{i_n} \mid i_1, \ldots, i_n \in \n \}$.
 \item $\succ$ can be any admissible term ordering on $\myt$.
       For the reductive ordering $\geq$ we have $t \geq s$ if $s$ divides $t$ as
 	as term.
 \item Multiplication $\rmult$ is specified by the action on terms, i.e.~$\rmult: \myt \times \myt \myr \myt$ where 
       $x_1^{i_1} \ldots x_n^{i_n} \rmult x_1^{j_1} \ldots x_n^{j_n} =
        x_1^{i_1+j_1} \ldots x_n^{i_n+j_n}$.
 \end{itemize}
We do not need the concept of weak saturation.

Since the integers are an instance of euclidean domains,
 similar reductions to those given by Kandri-Rodi and Kapur in \cite{KaKa88} arise. 
A stable localization of ${\cal C}_s(p,q)$ is already provided by
 the tuple corresponding to the least
 common multiple of the terms $\hterm(p)$ and $\hterm(q)$.
In contrast to the s- and t-polynomials studied by Kandri-Rodi and Kapur,
 we restrict ourselves to s-polynomials as described in Definition \ref{def.s-poly.z}.

Since this structure is Abelian, one-sided and two-sided ideals coincide.
Buchberger's Algorithm provides an effictive procedure to compute finite
 Gr\"obner bases.
%%%%%%%%%%%%%%%%%%%%%%%%%%%%%%%%%%%%%%%%%%%%
%%%%%%%%%%%%%%%%%%%%%%%%%%%%%%%%%%%%%%%%%%%%
\subsubsection{Non-commutative Polynomial Rings}
A non-commutative polynomial ring $\z[\{x_1, \ldots, x_n\}^*]$ is a function
 ring according to the following interpretation:
 \begin{itemize}
 \item $\myt$ is the set of words on $\{ x_1, \ldots, x_n\}$.
 \item $\succ$ can be any admissible ordering on $\myt$.
       For the reductive ordering $\geq$ we can chose
	 $t \geq s$ if $s$ is a subword 
        of $t$.
 \item Multiplication $\rmult$ is specified by the action on words which is just concatenation.
 \end{itemize}
We do not need the concept of weak saturation.
A stable localization of ${\cal C}_s(p,q)$ is already provided by
 the tuples corresponding to word overlaps resulting from the equations $u_1\hterm(p)v_1 = \hterm(q)$, 
 $u_2\hterm(q)v_2 = \hterm(p)$, $u_3\hterm(p) = \hterm(q)v_3$ respectively $u_4\hterm(q) = \hterm(p)v_4$
 with the restriction that $|u_3| < |\hterm(q)|$ and $|u_4| < | \hterm(p) |$,
 $u_i,v_i \in \myt$.
The coefficients arise as described in Definition \ref{def.s-poly.z}.

This structure is not Abelian.
For the case of one-sided ideals finite Gr\"obner bases can be computed.
The case of two-sided ideals only allows an enumerating procedure.
This is not surprising as the word problem for monoids can be reduced to the problem
 of computing the respective Gr\"obner bases (see e.g.~\cite{Mo87,MaRe95}).
%%%%%%%%%%%%%%%%%%%%%%%%%%%%%%%%%%%%%%%%%%%%
%\subsubsection{Path Algebras}\label{section.pathalgebras}
%%%%%%%%%%%%%%%%%%%%%%%%%%%%%%%%%%%%%%%%%%%%
\subsubsection{Monoid and Group Rings}
A monoid or group ring $\z[\m]$ is a function
 ring according to the following interpretation:
 \begin{itemize}
 \item $\myt$ is the monoid or group $\m$. In the cases studied by us
       as well as in \cite{Ro93,Lo96}, it is assumed that the
       elements of the monoid or group have a certain form.
       This presentation is essential in the approach.
       We will assume that the given monoid or group is presented
       by a convergent semi-Thue system.
 \item $\succ$ will be the completion ordering induced from the
       presentation of $\m$ to $\m$ and hence to $\myt$.
       The reductive ordering $\geq$ depends on the choice of the
       presentation.
 \item Multiplication $\rmult$ is specified by lifting the monoid or
       group operation.
 \end{itemize}
The concept of weak saturation and the choice of stable localizations of
 ${\cal C}_s(p,q)$ again depend on the choice of the presentation.
More on this topic can be found in \cite{Re95}.
%%% Local Variables: 
%%% mode: latex
%%% TeX-master: "testlauf"
%%% End: 

%% file: rmodule.tex
%%%%%%%%%%%%%%%%%%%%%%%%%%%%%%%%%%%%%%%%%%%%%%%%%%%%%%%%%%%%%%%%%%%%%%%%%%%%%%%%%%%%
\section{Right $\f$-Modules}\label{section.right.module}
The concept of  modules arises naturally as a generalization of the concept of
 an  ideal in a ring:
Remember that an  ideal of a ring is an additive subgroup of the ring which is 
 additionally closed under  multiplication with ring elements.
Extending this idea to arbitrary additive groups then gives us the concept of 
  modules.

In this section we turn our attention to right modules, but left modules can be defined similarly
 and all results carry over (with the respective modifications of the terms ``right'' and ``left'').
Let $\f$ be a function ring with unit $\one$.

\begin{example}~\\
{\rm
Let us provide some examples for right $\f$-modules.
\begin{enumerate}
\item
Any right ideal in $\f$ is of course a right $\f$-module.
\item
The set ${\cal M} = \{ {\bf 0 } \}$ with right scalar multiplication ${\bf 0} \rmult f = {\bf 0}$
 is a right $\f$-module called the trivial right $\f$-module.
\item
Given a function ring $\f$  and a natural number $k$, 
 let $\f^k = \{ (f_1, \ldots, f_k) \mid f_i \in \f \}$
 be the set of all vectors of length $k$ with coordinates in $\f$.
Obviously $\f^k$ is an additive commutative group with respect to ordinary vector addition.
Moreover, $\f^k$ is a right $\f$-module
 with right scalar multiplication $\rmult : \f^k \times \f \myr \f^k$ defined by
 $(f_1, \ldots, f_k) \rmult f = (f_1 \rmult f , \ldots, f_k \rmult f)$.
\exaend
\end{enumerate}
}
\end{example}

\begin{definition}~\\
{\rm
A subset of a right $\f$-module ${\cal M}$ which is again a right $\f$-module is
 called a \betonen{right submodule} of ${\cal M}$.
\dend
}
\end{definition}

For example any right ideal of $\f$ is a right  submodule of the right $\f$-module $\f^1$.
Provided a set of vectors $S \subset {\cal M}$ the set
 $\{ \sum_{i=1}^s  {\bf m}_i \rmult {g_i} \mid s \in \n, g_i \in \f, {\bf m}_i \in S \}$ is a right submodule of ${\cal M}$.
This set is denoted as $\langle S \rangle_r$ and $S$ is called its generating set.
If $\langle S \rangle_r = {\cal M}$ then $S$ is a generating set of the right module itself.
If $S$ is finite then ${\cal M}$ is said to be finitely generated.
A generating set is called linearly independent or a basis if for all $s \in \n$, pairwise different
 ${\bf m}_1, \ldots, {\bf m}_s \in S$ and $g_1, \ldots, g_s \in \f$, $\sum_{i=1}^s {\bf m}_i \rmult g_i = {\bf 0}$
 implies $g_1 = \ldots = g_s = \zero$. 
A right $\f$-module is called \betonen{free} if it has a basis.
The right $\f$-module $\f^k$ is free and one
 such basis is the set of unit vectors ${\bf e}_1 = (\one, \zero, \ldots, \zero), {\bf e}_2 = (\zero, \one, \zero, \ldots, \zero),
 \ldots, {\bf e}_k = (\zero, \ldots, \zero, \one)$.
Using this basis the elements of $\f^k$ can be written uniquely as ${\bf f} = \sum_{i=1}^k {\bf e}_i \rmult f_i$
 where ${\bf f} =  (f_1, \ldots, f_k)$.
Moreover, $\f^k$ has special properties similar to the special case of 
 $\myk[x_1, \ldots, x_n]$ and we will continue to state some of them.

\begin{theorem}\label{theo.right.submodule}~\\
{\sl
Let $\f$ be right Noetherian.
Then every right submodule of $\f^k$ is finitely generated.
}
\end{theorem}
\Ba{}~\\
Let ${\cal S}$ be a right submodule of $\f^k$.
We show our claim by induction on $k$.
For $k=1$ we find that ${\cal S}$ is in fact a right ideal in $\f$ and hence by our hypothesis finitely
 generated.
For $k>1$ let us look at the set $I = \{ f_1 \mid (f_1, \ldots, f_k) \in {\cal S} \}$.
Then again $I$ is a right ideal in $\f$ and hence finitely generated.
Let $\{g_1, \ldots, g_s \mid g_i \in \f \}$ be a generating set of $I$.
Choose ${\bf g}_1, \ldots, {\bf g}_s \in {\cal S}$ such that the first coordinate of ${\bf g}_i$ is $g_i$.
Similarly, the set $\{ (f_2, \ldots, f_k) \mid (\zero,f_2, \ldots, f_k) \in  {\cal S}\}$ is a submodule
 of $\f^{k-1}$ and hence finitely generated by some set $\{(n_2^i, \ldots, n_k^i), 1 \leq i \leq w  \}$.
Then the set $\{{\bf g}_1, \ldots, {\bf g}_s \} \cup \{ {\bf n}_i=(\zero, n_2^i, \ldots, n_k^i) \mid 
 1 \leq i \leq w \}$ is a generating set for ${\cal S}$.
To see this assume ${\bf m} = (m_1, \ldots, m_k)\in {\cal S}$.   
Then $m_1 = \sum_{i=1}^s g_i \rmult h_i$ for some $h_i \in \f$ and
 ${\bf m'} = {\bf m} - \sum_{i=1}^s {\bf g}_i \rmult h_i \in {\cal S}$ with first coordinate $\zero$.
Hence ${\bf m'} = \sum_{i=1}^w  {\bf n}_i \rmult l_i$ for some $l_i \in \f$ giving rise to
$${\bf m} = {\bf m'} + \sum_{i=1}^s  {\bf g}_i \rmult h_i = 
\sum_{i=1}^w {\bf n}_i \rmult l_i + \sum_{i=1}^s  {\bf g}_i \rmult h_i.$$
\qed
$\f^k$ is called right Noetherian if and only if all its right submodules are finitely generated.

If $\f$ is a right reduction ring, results on the existence of right Gr\"obner bases for the 
 right submodules carry over from modifications of the proofs in Section \ref{section.right.module}.

A natural reduction relation using the right reduction relation in $\f$ denoted by $\R$
 can be defined using the representation as (module) polynomials with respect to
 the basis of unit vectors as follows:
\begin{definition}\label{def.rred.module}~\\
{\rm 
Let ${\bf f} = \sum_{i=1}^k  {\bf e}_i \rmult f_i$, ${\bf p} = \sum_{i=1}^k  {\bf e}_i \rmult p_i\in \f^k$.
We say that
 ${\bf f}$ \betonen{reduces} ${\bf p}$ to ${\bf q}$ at ${\bf e}_s \rmult p_s$ in one step, 
 denoted by ${\bf p} \myr_{\bf f} {\bf q}$, if 
\begin{enumerate}
\item $p_j = \zero$ for $1 \leq j < s$,
\item $p_s \R_{f_s} q_s$, % with $p_s = q_s + f_s \skm d$, $d \in \monoms(\f)$, and
\item \begin{tabbing}
       ${\bf q}$ \= = \= ${\bf p} -  \sum_{i=1}^n  {\bf f} \skm d_i$ \\
               \> = \> $(0, \ldots , 0, q_s, 
             p_{s+1} - \sum_{i=1}^n f_{s+1} \skm d, \ldots , 
             p_k - \sum_{i=1}^n f_k \skm d)$. 
     \end{tabbing}
\dend
\end{enumerate}
}
\end{definition}
Notice that item 2 of this definition is dependant on the definition
 of the reduction relation $\R$ in $\f$.
If we assume that the reduction relation is the one specified in
 Definition \ref{def.rred_rr} we get
           $p_s = q_s + f_s \skm d$, 
           $d \in \monoms(\f)$, but there are other possibilities.
Reviewing the introduction of right modules to reduction rings we could
 substitute 2.~by $p_s = q_s + f_s \skm d$, 
           $d \in \f$ as well (compare Definition \ref{def.module.red.right}).

To show that our reduction relation is terminating we have to extend the ordering from
 $\f$ to $\f^k$.
For two elements ${\bf p} = (p_1, \ldots, p_k)$, ${\bf q} = (q_1, \ldots, q_k) \in \f^k$ we define
 ${\bf p} \succ {\bf q}$ if and only if there exists $1 \leq s \leq k$ such that
 $p_i = q_i$, $1 \leq i < s$, and $p_s \succ q_s$.
\begin{lemma}~\\ 
{\sl 
Let $F$ be a finite set of module polynomials in $\f^k$.
\begin{enumerate}
\item For ${\bf p},{\bf q} \in \f^k$ ${\bf p} \red{}{\myr}{}{F} {\bf q}$
       implies ${\bf p} \succ {\bf q}$.
\item $\red{}{\myr}{}{F}$ is Noetherian in case $\R_{F_i}$ is for $1 \leq i \leq k$ and $F_i = \{ f_i \mid f= (f_1, \ldots , f_k) \in F \}$.\footnote{Notice
 that $F_i\subseteq \f$.}.
\lemend
\end{enumerate}
}
\end{lemma}
\Ba{}
\begin{enumerate}
\item Assuming that the reduction step takes place at ${\bf e}_s \rmult p_s$,
      by Definition \ref{def.rred.module} we know $p_s \R_{f_s} q_s$ and $p_s > q_s$
      implying ${\bf p} \succ {\bf q}$.      
\item This follows from 1.~and Axiom (A1).
\end{enumerate}\renewcommand{\baselinestretch}{1}\small\normalsize
\qed
\begin{definition}~\\
{\rm
A subset $B$ of  $\f^k$ is called a 
 \betonen{right Gr\"obner basis} of the right submodule ${\cal S} = \langle B \rangle_r$, if
 $\red{*}{\lr}{}{B} = \;\;\equiv_{{\cal S}}$
 and $\red{}{\myr}{}{B}$ is convergent. 
\dend
}
\end{definition}
For any reduction relation in $\f$ fulfilling the Axioms
 (A1)--(A3), the following theorem holds.
\begin{theorem}~\\
{\sl
If in $(\f,\red{}{\R}{}{})$ every finitely generated right ideal has a finite right Gr\"obner basis, then 
 the same holds for finitely generated right submodules in $(\f^k, \myr)$.
\theoend
}
\end{theorem}
\Ba{}~\\
Let ${\cal S} = \langle \{ {\bf s}_1, \ldots, {\bf s}_n \}\rangle$ 
 be a finitely generated right submodule of $\f^k$.
We show our claim by induction on $k$.
For $k=1$ we find that ${\cal S}$ is in fact a finitely generated right ideal in $\f$ and hence by our 
 hypothesis must have a finite right Gr\"obner basis.
For $k>1$ let us look at the set $I = \{ f_1 \mid (f_1, \ldots, f_k) \in {\cal S} \}$ which is in fact
 the right ideal generated by $\{ s_1^i \mid {\bf s}_i = (s_1^i, \ldots, s_k^i), 1 \leq i \leq n\}$.
Hence $I$ must have a finite right Gr\"obner basis $H = \{g_1, \ldots, g_s \mid g_i \in \f \}$.
Choose ${\bf g}_1, \ldots, {\bf g}_s \in  {\cal S}$ such that the first coordinate of 
 ${\bf g}_i$ is $g_i$.
Similarly the set ${\cal S}' = \{(f_2, \ldots, f_k) \mid  (\zero, f_2, \ldots, f_k) \in {\cal S} \}$ is a submodule in
 $\f^{k-1}$ which by our induction hypothesis then must have a finite right Gr\"obner basis 
 $\{(\tilde{g}_2^i, \ldots, \tilde{g}_k^i), 1 \leq i \leq w \}$.
Then the set $G = \{{\bf g}_1, \ldots, {\bf g}_s \} \cup 
 \{ {\bf \tilde{g}}_i = (\zero, \tilde{g}_2^i, \ldots, \tilde{g}_k^i) \mid 
  1 \leq i \leq w \}$ is a right Gr\"obner basis for ${\cal S}$.
As shown in the proof of Theorem \ref{theo.right.submodule},
 $G$ is a generating set for ${\cal S}$.
It remains to show that $G$ is in fact a right Gr\"obner basis.
\\
First we have to show $\red{*}{\lr}{}{G} = \;\;\equiv_{{\cal S}}$.
By the definition of the reduction relation in $\f^k$ we immediately find
 $\red{*}{\lr}{}{G} \subseteq \;\;\equiv_{{\cal S}}$.
To see the converse let ${\bf p} = (p_1, \ldots , p_k) \equiv_{{\cal S}}
 {\bf q} = (q_1, \ldots, q_k)$.
Then $p_1 \equiv_{ \langle \{ s^1_i \mid {\bf s}_i = 
 (s_1^i, \ldots, s_k^i), 1 \leq i \leq n\} \rangle_r} q_1$ 
 and hence by the definition of $G$ we get 
 $p_1 \red{*}{\lr}{}{\{ g_1^i \mid {\bf g}_i = (g_1^i, \ldots, g_k^i), 1 \leq i \leq s \}} q_1$.
But this gives us ${\bf p} \red{*}{\lr}{}{H} {\bf p} + \sum_{i=1}^s {\bf g}_i \rmult r_i = {\bf  p}' = (q_1, {p_2}', \ldots, {p_k}')$, $r_i \in \f$, and we get
 $(q_1,{p_2}', \ldots, {p_k}')\equiv_{{\cal S}} (q_1,q_2, \ldots, q_k)$ and hence 
 $(q_1,{p_2}', \ldots, {p_k}')- (q_1,q_2, \ldots, q_k) = (\zero, {p_2}' - q_2, \ldots, {p_k}' - q_k) \in {\cal S}$
 implying $({p_2}' - q_2, \ldots, {p_k}' - q_k) \in {\cal S}'$ and
 $(\zero,{p_2}' - q_2, \ldots, {p_k}' - q_k) = \sum_{i=1}^w{\bf \tilde{g}}_i \rmult {\eta_{i}}$ for  $\eta_{i} \in \f$.
Hence $(q_1,{p_2}', \ldots, {p_k}')$ and $(q_1,q_2, \ldots, q_k) = (q_1,{p_2}', \ldots, {p_k}') - (\zero,{p_2}' - q_2, \ldots, {p_k}' - q_k) = (q_1,{p_2}', \ldots, {p_k}') - \sum_{i=1}^w{\bf \tilde{g}}_i \rmult {\eta_{i}}$ must be
 joinable by $\{ {\bf \tilde{g}}_i \mid 1 \leq i \leq w \}$ as the restriction of this set 
 without the first coordinate is a right Gr\"obner basis of ${\cal S}'$.
\\
Since the reduction relation using the finite set $G$ is terminating we only have to show local confluence. 
Let us assume there are ${\bf p}$, ${\bf q}_1$, ${\bf q}_2 \in \f^k$
 such that ${\bf p} \myr_{G} {\bf q}_1$ and ${\bf p} \myr_{G} {\bf q}_2$.
Then by the definition of $G$ the first coordinates $q^1_1$ and $q^2_1$ are joinable to some element say $s$
 by $H = \{g_1, \ldots, g_s \}$ giving rise to the elements
 ${\bf p}_1 = {\bf q}_1 + \sum_{i=1}^s {\bf g}_i \rmult h_{i}$ and
 ${\bf p}_2 = {\bf q}_2 + \sum_{i=1}^s {\bf g}_i \rmult \tilde{h}_{i}$ with first
 coordinate $s$.
As before, ${\bf p}_1 = {\bf p}_2 + \sum_{i=1}^w{\bf \tilde{g}}_i \rmult {\eta_{i}}$ and 
 hence ${\bf p}_1$ and ${\bf p}_2$ must be joinable by
 $\{ {\bf \tilde{g}}_i \mid 1 \leq i \leq w \}$.
\\
\qed

Now given a right submodule ${\cal S}$ of ${\cal M}$, we can define 
 ${\cal M}/{\cal S} = \{ {\bf f} + {\cal S} \mid {\bf f} \in {\cal M} \}$.
Then with addition defined as $({\bf f} + {\cal S}) + ({\bf g} + {\cal S}) = ({\bf f}+ {\bf g}) + {\cal S}$
 the set ${\cal M}/{\cal S}$ is an Abelian group and can be turned into
 a right $\f$-module by the action $({\bf f} + {\cal S}) \rmult g = {\bf f} \rmult g + {\cal S}$.
${\cal M}/{\cal S}$ is called the \betonen{right quotient module} of ${\cal M}$ by ${\cal S}$.

As usual this quotient can be related to homomorphisms.
The results carry over from commutative module theory as can be found in \cite{AdLo94}.
Recall that for two right $\f$-modules ${\cal M}$ and ${\cal N}$, a function
 $\phi: {\cal M} \myr{\cal N}$ is a right $\f$-module homomorphism if
 $$\phi({\bf f } + {\bf g}) = \phi({\bf f}) + \phi({\bf g}) \mbox{ for all } {\bf f,g} \in {\cal M}$$
 and
 $$\phi({\bf f}) \rmult g = \phi ({\bf f} \rmult g) \mbox{ for all } {\bf f} \in {\cal M}, g \in \f.$$ 
The homomorphism is called an \betonen{isomorphism} if $\phi$ is one to one and we then write
 ${\cal M}\cong {\cal N}$.
Let ${\cal S} = {\rm ker}(\phi) = \{ {\bf f} \in {\cal M} \mid \phi({\bf f}) = {\bf 0} \}$.
Then ${\cal S}$ is a right submodule of ${\cal M}$ and $\phi({\cal M})$ is a right submodule
 of ${\cal N}$.
Since all are Abelian groups we know ${\cal M}/{\cal S} \cong \phi({\cal M})$
 under the mapping  ${\cal M}/{\cal S} \myr \phi({\cal M})$ with
 ${\bf f} + {\cal S} \mapsto \phi({\bf f})$ which is in fact an isomorphism.
All right submodules of the quotient ${\cal M}/{\cal S}$ are of the form 
 ${\cal L}/{\cal S}$ where ${\cal L}$ is a right submodule of ${\cal M}$
 containing ${\cal S}$.

We can even show that every finitely generated right $\f$-module is of a special form.
\begin{lemma}~\\
{\sl
Every finitely generated right $\f$-module ${\cal M}$ is isomorphic to $\f^k/{\cal N}$
 for some $k \in \n$ and some right submodule ${\cal N}$ of $\f^k$.
}
\end{lemma}
\Ba{}~\\
Let ${\cal M}$ be a  finitely generated right $\f$-module with
 generating set  ${\bf f}_1, \ldots {\bf f}_k \in {\cal M}$.
Consider the mapping $\phi: \f^k \myr {\cal M}$ defined by 
 $\phi(g_1, \ldots, g_k)=\sum_{i=1}^k {\bf f}_i \rmult g_i$.
Then $\phi$ is an $\f$-module homomorphism with image ${\cal M}$.
Let ${\cal N}$ be the kernel of $\phi$, then the First Isomorphism Theorem for modules yields our claim.
Note that
$\phi$ is uniquely defined by specifying the image of each unit vector ${\bf e}_1, \ldots, {\bf e}_k$,
 namely by $\phi({\bf e}_i) = {\bf f}_i$.
\\
\qed

Now, there are two ways to give a finitely generated right $\f$-module ${\cal M} \subset \f^k$.
One is to be given explicit ${\bf f}_1, \ldots {\bf f}_t \in  \f^k$ such that
 ${\cal M} = \langle \{ {\bf f}_1, \ldots {\bf f}_s \} \rangle_r$.
The other way is to give a right submodule ${\cal N}= \langle \{ {\bf g}_1, \ldots {\bf g}_s \} \rangle_r$
 for explicit ${\bf g}_1, \ldots {\bf g}_s \in \f^k$ such that ${\cal M} \cong \f^k/{\cal N}$.
This is called  a \betonen{presentation} of ${\cal M}$.

Presentations are chosen when studying right ideals of $\f$ as right $\f$-modules.
To see how this is done let $\mathfrak{i}$ be the right ideal generated by
 $\{ f_1, \ldots, f_k \}$ in $\f$.
Let us consider the right $\f$-module homomorphism defined as a mapping
 $\phi: \f^k \myr \mathfrak{i}$ with $\phi(g_1, \ldots, g_k) = \sum_{i=1}^k f_i \rmult g_i$.
Then $\mathfrak{i} \cong \f^k/{\rm ker}(\phi)$ as $\f$-modules.
${\rm ker}(\phi)$ is called the \betonen{right syzygy} of $\{ f_1, \ldots, f_k \}$ denoted by
 ${\rm Syz}(f_1, \ldots, f_k)$.
In fact ${\rm Syz}(f_1, \ldots, f_k)$ is the set of all solutions of the linear equation
 $f_1 X_1 + \ldots + f_k X_k = \zero$ in $\f$.
Syzygies play an important role in Gr\"obner basis theory for ordinary polynomial rings.
%The importance of these sets stems from the fact that they can be used to characterize right Gr\"obner
% bases.
%\begin{theorem}~\\
%{\sl
%Let $F = \{ f_1, \ldots, f_k \}$ be a subset of $\f$.
%Then $F$ is a right Gr\"obner basis if and only if all $(g_1, \ldots, g_k) \in {\rm Syz}(f_1, \ldots, f_k)$
% have a right standard representation with respect to $F$.
%\theoend
%}
%\end{theorem}
%Of course this theorem becomes more important if we are able to restrict the test to special
% generating subsets of ${\rm Syz}(f_1, \ldots, f_k)$ and this has been studied by Apel in his habilitation on Gr\"obner bases in graded structures over fields.
%Such a restriction is normally related to the introduction of reduction
% and stable right
% standard representations as introduced in the previous section.

%%% Local Variables: 
%%% mode: latex
%%% TeX-master: "testlauf"
%%% End: 

%% file: ideals.tex
\section{Ideals and Standard Representations}\label{section.ideal.standard}
A subset $\mathfrak{i} \subseteq \f$ is called a (two-sided)
 \betonen{ideal}, if
\begin{enumerate}
\item $\zero \in \mathfrak{i}$,
\item for $f,g \in \mathfrak{i}$ we have $f \radd g \in \mathfrak{i}$, and
\item for $f \in \mathfrak{i}$, $g, h \in \f$ we have
          $g \rmult f \rmult h \in \mathfrak{i}$.
\end{enumerate}
Ideals can also be specified in terms of a
 generating set.
For $F \subseteq \f \backslash \{ \zero\}$ let
 $\ideal{}{}(F) = \{ \sum_{i=1}^n g_i \rmult f_i \rmult h_i \mid
  f_i \in F, g_i, h_i \in \f, n \in \n \} = \{ \sum_{i=1}^m m_i \rmult f_i \rmult l_i \mid
 f_i \in F, m_i, l_i \in \monoms(\f),  n \in \n \}$.
These generated sets are in fact subsets of $\f$ since for
 $f,g \in \f$ we have that $f \rmult g$
 as well as $f \radd g$ are again elements of $\f$, and it is easily checked
 that they are in fact ideals:
 \begin{enumerate}
 \item $\zero \in \ideal{}{}(F)$ since $\zero$ can be written as
        the empty sum.
 \item For two elements $\sum_{i=1}^n g_i \rmult f_i \rmult h_i$ and 
        $\sum_{i=1}^m \tilde{g}_i \rmult \tilde{f}_i \rmult \tilde{h}_i$ in $\ideal{}{}(F)$,
        the sum $\sum_{i=1}^n g_i \rmult f_i \rmult h_i \radd
        \sum_{i=1}^m \tilde{g}_i \rmult \tilde{f}_i \rmult \tilde{h}_i $ is again
        an element in $\ideal{}{}(F)$.
 \item For an element $\sum_{i=1}^n g_i \rmult f_i \rmult h_i$ in $\ideal{r}{}(F)$
        and two polynomials $g,h$ in $\f$, the product
        $g \rmult (\sum_{i=1}^n g_i \rmult f_i \rmult h_i) \rmult h =
          \sum_{i=1}^n (g \rmult g_i) \rmult  f_i \rmult (h_i \rmult h)$  is again
        an element in $\ideal{}{}(F)$.
 \end{enumerate}
Given an ideal 
 $\mathfrak{i} \subseteq \f$ we
 call a set $F \subseteq \f \backslash \{ \zero\}$ a
 basis of $\mathfrak{i}$ if
 $\mathfrak{i} = \ideal{}{}(F)$.
Then every element $g \in \ideal{}{}(F) \backslash \{ \zero\}$
 can have different representations
 of the form
 $$g = \sum_{i=1}^n g_i \rmult f_i \rmult h_i,  
   f_i \in F, g_i, h_i \in \f, n \in \n.$$
Notice that the $f_i$ occurring in this sum are not necessarily different.
The distributivity law in $\f$ allows to convert such a representation
 into one 
 of the form
 $$g = \sum_{j=1}^m m_i \rmult f_i \rmult l_i,  
   f_i \in F, m_i, l_i \in \monoms(\f), n \in \n.$$
Again special representations can be distinguished in order to characterize special ideal bases.
An ordering on $\f$ is used to define appropriate standard representations.
As in the case of right ideals we will first look at generalizations of
 standard representations for the case of function rings over fields.

\subsection{The Special Case of Function Rings over Fields}

Let $\f_{\myk}$ be a function ring over a field $\myk$.
We first look at an analogon to Definition \ref{def.standard.rep}
\begin{definition}\label{def.two-sided.standard.rep}~\\
{\rm
Let $F$ be a set of polynomials in $\f_{\myk}$ and $g$ a non-zero polynomial in $\ideal{}{}(F)$.
A representations of the form
 $$g = \sum_{i=1}^n m_i \rmult f_i \rmult l_i, 
   f_i \in F, m_i, l_i \in \monoms(\f_{\myk}), n \in \n$$
 where additionally $\hterm(g) \succeq \hterm(m_i \rmult f_i \rmult l_i)$ holds for
 $1 \leq i \leq n$ is called a \betonen{standard representation} of $g$ in terms of $F$.
If every $g \in \ideal{}{}(F) \backslash \{ \zero\}$ has such a representation in terms of $F$,
 then $F$ is called a \betonen{standard basis} of $\ideal{}{}(F)$.
\dend
}
\end{definition}
Notice that since we assume  $f \skm \alpha = \alpha \skm f$, we can
 also substitute the monomials $l_i$ by terms $w_i \in \myt$, i.e.~study
 representations 
 of the form
 $$g = \sum_{i=1}^n m_i \rmult f_i \rmult w_i,  
   f_i \in F, m_i \in \monoms(\f), w_i \in \myt, n \in \n.$$
We will use this additional information in some proofs later on.

As with right standard representations, in order to change an arbitrary representation of an 
 ideal element into a standard representation we have to deal with special sums of polynomials.
We get the following analogon to Definition \ref{def.critical.situations}.

\begin{definition}\label{def.two-sided.critical.situations}~\\
{\rm
Let $F$ be a set of polynomials in $\f_{\myk}$ and $t$ an element in $\myt$.
Then we define a set ${\cal C}(F,t)$ to contain all tuples of the form
 $(t, f_1, \ldots, f_k, m_1, \ldots, m_k, l_1, \ldots, l_k)$, $k \in \n$, $f_1, \ldots, f_k \in F$,
 $m_1, \ldots, m_k, l_1, \ldots, l_k \in \monoms(\f_{\myk})$ such that
 \begin{enumerate}
 \item $\hterm(m_i \rmult f_i \rmult l_i) = t$, $1 \leq i \leq k$, and
 \item $\sum_{i=1}^k \hm(m_i \rmult f_i \rmult l_i) = 0$.
 \end{enumerate}
We set ${\cal C}(F) = \bigcup_{t \in \myt} {\cal C}(F, t)$.
\dend
}
\end{definition}
Notice that this definition is motivated by the definition of syzygies of head monomials
 in commutative polynomial rings over rings.
We can characterize standard bases using this concept (compare Theorem \ref{theo.standard.basis}).
\begin{theorem}~\\
{\sl
Let $F$ be a set of polynomials in $\f_{\myk} \backslash \{ \zero\}$.
Then $F$ is a standard basis of $\ideal{}{}(F)$ if and only if
 for every tuple 
 $(t, f_1, \ldots, f_k, m_1, \ldots, m_k, l_1, \ldots, l_k)$ in ${\cal C}(F)$
 the polynomial $\sum_{i=1}^k m_i \rmult f_i \rmult l_i$ (i.e.~the element in $\f_{\myk}$
 corresponding to this sum) has a standard representation
 with respect to $F$.
\theoend
}
\end{theorem}
\Ba{}~\\
In case $F$ is a standard basis
 since the polynomials related to the tuples 
 are all elements of $\ideal{}{}(F)$ they must
 have standard representations with respect to $F$.
\\
To prove the converse, it remains to show that every element in
 $\ideal{}{}(F)$ has a standard representation
 with respect to $F$.
Hence, let $g = \sum_{j=1}^m m_j \rmult f_{j} \rmult l_{j}$ be an arbitrary
 representation of a non-zero  polynomial $g\in \ideal{}{}(F)$ such that
 $f_j \in F$, $m_{j}, l_j \in \monoms(\f_{\myk})$, $m \in \n$.
Depending on this  representation of $g$ and the
 well-founded total ordering $\succeq$ on $\myt$ we define
 $t = \max_{\succeq} \{ \hterm(m_j \rmult f_{j} \rmult l_{j}) \mid 1\leq j \leq m \}$ and
 $K$ as the number of polynomials $m_j \rmult f_j \rmult l_j$ with head term $t$.
%\\
Then $t \succeq \hterm(g)$ and 
 in case $\hterm(g) = t$ this immediately implies that this representation is
 already a  standard  one. 
%\\
Else we proceed by induction
 on $t$.
%\\
Without loss of generality let $f_1, \ldots, f_K$ be the polynomials
 in the corresponding representation
 such that  $t=\hterm(m_j \rmult f_j \rmult l_j)$, $1 \leq j \leq K$.
Then the tuple $(t, f_1, \ldots, f_K, m_1, \ldots, m_K, l_1, \ldots, l_K)$
 is in ${\cal C}(F)$ and let $h = \sum_{j=1}^K m_j \rmult f_j \rmult l_j$.
%\\
We will now change our representation of $g$ in such a way that for the new
 representation of $g$ we have a smaller maximal term.
%\\
Let us assume $h$ is not $\zero$\footnote{In case  $h =\zero$,
 just substitute the empty sum for the representation of $h$
 in the equations below.}. 
%\\
By our assumption, $h$ has a   standard representation
 with respect to $F$, say $\sum_{i=1}^n \tilde{m}_i \rmult \tilde{f}_i \rmult \tilde{l}_i$, 
 where $\tilde{f}_i \in F$, and $\tilde{m}_i,\tilde{l}_i \in \monoms(\f_{\myk})$ 
 and  all terms occurring in the sum are bounded by
 $t \succ \hterm(h)$.
%\\
This gives us: 
\begin{eqnarray}
  g   & = & \sum_{j=1}^K m_j \rmult f_j \rmult l_j + \sum_{j=K+1}^m m_j \rmult f_j \rmult l_j
             \nonumber\\                                                           
  & = & \sum_{i=1}^n \tilde{m}_i \rmult \tilde{f}_i\rmult \tilde{l}_i + \sum_{j=K+1}^m m_j \rmult  f_j \rmult l_j
              \nonumber
\end{eqnarray}
which is a representation of $g$ where the maximal term of the involved monomial
 multiples is decreased.
\\
\qed
Weak Gr\"obner bases can be defined as in Definition \ref{def.weak.gb}.
Since the ordering $\succeq$ and the multiplication $\rmult$ in general are not compatible,
instead of considering multiples of head terms of the generating set $F$
 we look at head terms of monomial multiples of polynomials in $F$.

\begin{definition}\label{def.two-sided.weak.gb}~\\
{\rm
A subset $F$ of $\f_{\myk} \backslash \{ \zero\}$ is called a 
 \betonen{weak Gr\"obner basis} of $\ideal{}{}(F)$ if
$\hterm(\ideal{}{}(F) \backslash \{ \zero \}) =
 \hterm( \{ m \rmult f \rmult l \mid f \in F, m,l \in \monoms(\f_{\myk}) \}
 \backslash \{ \zero \} )$.
\dend
}
\end{definition}

In the next lemma we show that in fact both characterizations of
 special bases, standard bases and weak Gr\"obner bases,
 coincide as in the case of polynomial
 rings over fields (compare Lemma \ref{lem.rsb=gb}).
\begin{lemma}\label{lem.two-sided.rsb=gb}~\\
{\sl
Let $F$ be a subset of $\f_{\myk} \backslash \{ \zero\}$.
Then $F$ is a standard basis if and only
 if it is a weak Gr\"obner basis.
\lemend
}
\end{lemma}
\Ba{}~\\
Let us first assume that $F$ is a standard basis, i.e., every polynomial
 $g$ in $\ideal{}{}(F)$ has a standard  representation with respect
 to $F$.
In case $g \neq \zero$ this implies the existence of a polynomial $f \in F$
 and monomials
 $m,l \in \monoms(\f_{\myk})$ such that $\hterm(g) = \hterm(m \rmult f \rmult l)$.
Hence $\hterm(g) \in \hterm( \{  m \rmult f \rmult l \mid m,l \in \monoms(\f_{\myk}),
 f \in F \}\backslash \{ \zero\})$.
As the converse, namely $\hterm( \{  m \rmult f \rmult l \mid m,l \in \monoms(\f_{\myk}),
 f \in F \}\backslash \{ \zero\})
 \subseteq \hterm(\ideal{}{}(F) \backslash \{ \zero \})$ trivially holds,
 $F$ then is a weak Gr\"obner basis.
\\
Now suppose that $F$ is a weak Gr\"obner basis and again
 let $g \in \ideal{}{}(F)$.
We have to show that $g$ has a standard  representation with respect to $F$.
This will be done by induction on $\hterm(g)$.
In case $g = \zero$ the empty sum is our required standard  representation.
Hence let us assume $g \neq \zero$.
Since then $\hterm(g) \in \hterm(\ideal{}{}(F)\backslash \{ \zero\})$
 by the definition of weak
 Gr\"obner bases we know there exists
 a polynomial $f \in F$ and monomials $m,l \in \monoms(\f_{\myk})$ such that
 $\hterm(g) = \hterm(m \rmult f \rmult l)$.
Then there exists a monomial $\tilde{m} \in \monoms(\f_{\myk})$ such that
 $\hm(g)=\hm(\tilde{m} \rmult f \rmult l )$,
 namely\footnote{Notice that this step requires that we can view $\f_{\myk}$
 as a vector space. In order to get a similar result without introducing
 vector spaces we would have to use a different definition of weak
 Gr\"obner bases. E.g.~requiring that $\hm(\ideal{}{}(F) \backslash \{ \zero \})
 = \hm(\{ m \rmult f \rmult l \mid f \in F, m,l \in \monoms(\f_{\myk})\}\backslash \{ \zero \}\})$
 would be a possibility. However, then no localization of critical situations
 to head terms is possible, which is {\em the} advantage of having a field as 
 coefficient domain.} $\tilde{m}=(\hc(g) \skm \hc(m \rmult f \rmult l)^{-1}) \skm m)$.
Let $g_1 = g - \tilde{m} \rmult f \rmult l$. 
Then $\hterm(g) \succ \hterm(g_1)$ implies the
 existence of a standard  representation for $g_1$ which can be
 added to the multiple $\tilde{m} \rmult f \rmult l$ to give
 the desired standard  representation of $g$.
\\ 
\qed
Inspecting this proof closer we get the following corollary (compare
 Corollary \ref{cor.right.rep}).
\begin{corollary}\label{cor.two-sided.rep}~\\
{\sl
Let a subset $F$  of $\f_{\myk} \backslash \{ \zero\}$ be a 
 weak Gr\"obner basis.
Then every $g \in \ideal{}{}(F)$ has a standard representation
       in terms of $F$ of the form
        $g = \sum_{i=1}^n m_i \rmult f_i \rmult l_i,
         f_i \in F, m_i, l_i \in \monoms(\f_{\myk}), n \in \n$
         such that
         $\hm(g) = \hm(m_1 \rmult f_1 \rmult l_1)$ and
         $\hterm(m_1 \rmult f_1 \rmult l_1) \succ \hterm(m_2 \rmult f_2 \rmult l_2) \succ \ldots \succ
          \hterm(m_n \rmult f_n \rmult l_n)$.
\corend
}
\end{corollary}
Notice that 
 we hence get stronger representations as
 specified in Definition \ref{def.two-sided.standard.rep} for the case that the set $F$
 is a weak Gr\"obner basis or a standard basis.

In order to proceed as before in the case of one-sided ideals we have to extend our
 restriction of the ordering $\succeq$ on $\f$
 to cope with two-sided multiplication similar
 to Definition \ref{def.refined.ordering}.
\begin{definition}\label{def.two-sided.refined.ordering}~\\
{\rm
We will call an ordering $\geq$ on $\myt$ a 
 \betonen{reductive restriction}
 of the ordering $\succeq$ or simply \betonen{reductive}, 
 if the following hold:
\begin{enumerate}
\item $t \geq s$ implies $t \succeq s$ for $t, s \in \myt$.
\item $\geq$ is a partial well-founded ordering on $\myt$ which is
 compatible with multiplication $\rmult$ in the following sense:
if for $t, t_1, t_2, w_1, w_2 \in \myt$
 $t_2 \geq t_1$, $t_1 \succ t$ and
 $t_2 = \hterm(w_1 \rmult t_1 \rmult w_2)$ hold, then $t_2 \succ \hterm(w_1 \rmult t \rmult w_2)$.
\dend
\end{enumerate}
}
\end{definition}
Again we can distiguish special ``divisors'' of
 monomials:
For $m_1, m_2 \in \monoms(\f_{\myk})$ we call $m_1$ a \betonen{(stable) divisor} of $m_2$ if and only if 
 $\hterm(m_2) \geq \hterm(m_1)$
 and there exist $l_1, l_2 \in \monoms(\f_{\myk})$ such that $m_2 = \hm(l_1 \rmult m_1 \rmult l_2)$.
We then call $l_1, l_2$  \betonen{stable multipliers} of $m_1$.
The intention is that for all terms $t$ with $\hterm(m_1) \succ t$ we
 then can conclude $\hterm(m_2) \succ \hterm(l_1 \rmult t \rmult l_2)$.
Reduction relations based on this divisibility of terms will again have the stability
 properties we desire.
In the commutative polynomial ring 
 we can state a reductive restriction of any term ordering by
 $t \geq s$ for two terms $t$ and $s$ if and only if $s$ divides $t$ as 
 a term. 
In the non-commutative polynomial ring we can state a reductive restriction of any term ordering by
 $t \geq s$ for two terms $t$ and $s$ if and only if $s$ is a subword of $t$. 
Let us continue with an algebraic consequence related to
 this reductive ordering by distinguishing special 
 standard representations as we have done in Definition \ref{def.right_reductive}.
\begin{definition}\label{def.two-sided.right_reductive}~\\
{\rm
Let $F$ be a set of polynomials in $\f_{\myk}$
 and $g$ a non-zero polynomial in $\ideal{}{}(F)$.
A representation of the form 
$$g = \sum_{i=1}^n m_i \rmult f_i \rmult l_i,
 f_i \in F, m_i, l_i \in \monoms(\f_{\myk}), n \in \n$$
 such that
 $\hterm(g) = \hterm(m_i \rmult f_i \rmult l_i) =
 \hterm(m_i \rmult \hterm(f_i) \rmult l_i) \geq \hterm(f_i)$, $1 \leq i \leq k$ for some
 $k \geq 1$, and
 $\hterm(g) \succ \hterm(m_i \rmult \hterm(f_i) \rmult l_i)$ 
 for $k < i \leq n$ is called a
 \betonen{reductive standard representation} in terms of $F$.
\dend
}
\end{definition}
Again the empty sum is taken as reductive standard
 representation of $\zero$.

In case we have $\rmult : \myt \times \myt \myr \myt$ the condition
 can be rephrased as
 $\hterm(g) = m_i \rmult f_i \rmult l_i =
 \hterm(m_i \rmult \hterm(f_i) \rmult l_i) \geq \hterm(f_i)$, $1 \leq i \leq k$.
\begin{definition}~\\
{\rm
A set $F \subseteq\f_{\myk} \backslash \{ \zero\}$ is called a
 \betonen{reductive standard basis} (with respect to the reductive ordering $\geq$)
 of $\ideal{}{}(F)$ if every polynomial $f \in \ideal{}{}(F)$
 has a reductive standard representation in terms of $F$. 
\dend
}
\end{definition}
Again, in order to change an arbitrary representation into one fulfilling our
 additional condition of Definition \ref{def.two-sided.right_reductive}
 we have to deal with special sums of polynomials.
\begin{definition}\label{def.two-sided.reductive.critical.situations}~\\
{\rm
Let $F$ be a set of polynomials in $\f_{\myk}$ and $t$ an element in $\myt$.
Then we define the \betonen{critical set}
 ${\cal C}_{r}(t,F)$ to contain all tuples of the form
 $(t, f_1, \ldots, f_k, m_1, \ldots, m_k, l_1, \ldots, l_k)$, $k \in \n$, $f_1, \ldots, f_k \in F$\footnote{As in 
 the case of commutative polynomials, $f_1, \ldots, f_k$ are not
 necessarily different polynomials from $F$.},
 $m_1, \ldots, m_k, l_1, \ldots, l_k \in \monoms(\f)$ such that
 \begin{enumerate}
 \item $\hterm(m_i \rmult f_i \rmult l_i) =\hterm(m_i \rmult \hterm(f_i) \rmult l_i)= t$,
       $1 \leq i \leq k$, 
 \item $\hterm(m_i \rmult f_i \rmult l_i) \geq \hterm(f_i)$, $1 \leq i \leq k$, and
 \item $\sum_{i=1}^k \hm(m_i \rmult f_i \rmult l_i) = \zero$.
 \end{enumerate}
We set ${\cal C}_{r}(F) = \bigcup_{t \in \myt} {\cal C}_{r}(t,F)$.
\dend
}
\end{definition}
Unfortunately, as in the case of right reductive standard bases,
 these critical situations
 will not be sufficient to characterize reductive standard bases (compare
 again Example \ref{exa.free.group}).
But we can give an analogon to Theorem \ref{theo.right_reductive.standard.basis}.
\begin{theorem}\label{theo.two-sided.reductive.standard.basis}~\\
{\sl
Let $F$ be a set of polynomials in $\f_{\myk} \backslash \{ \zero\}$.
Then $F$ is a reductive standard basis of $\ideal{}{}(F)$ if and only if
\begin{enumerate}
 \item for every $f \in F$ and every $m,l \in \monoms(\f_{\myk})$ the multiple
 $m \rmult f \rmult l$ has a reductive standard representation in terms of $F$,
 \item
 for every tuple 
 $(t, f_1, \ldots, f_k, m_1, \ldots, m_k, l_1, \ldots, l_k)$ in ${\cal C}_{r} (F)$
 the polynomial $\sum_{i=1}^k m_i \rmult f_i \rmult l_i$ (i.e., the element in $\f$
 corresponding to this sum) has a reductive standard representation
 with respect to $F$.
\end{enumerate}
\theoend
}
\end{theorem}
\Ba{}~\\
In case $F$ is a reductive standard basis,
 since these polynomials are all elements of $\ideal{}{}(F)$, they must
 have reductive standard representations with respect to $F$.
\\
To prove the converse, it remains to show that every element in
 $\ideal{}{}(F)$ has a reductive standard representation
 with respect to $F$.
Hence, let $g = \sum_{j=1}^m m_j \rmult f_{j} \rmult l_{j}$ be an arbitrary
 representation of a non-zero  polynomial $g\in \ideal{}{}(F)$ such that
 $f_j \in F$, $m_{j}, l_j \in \monoms(\f_{\myk})$, $m\in \n$.
By our first statement every such monomial multiple $m_j \rmult f_j \rmult l_j$ has
 a reductive standard representation in terms of $F$ and we can
 assume that all multiples are replaced by them.
Depending on this  representation of $g$ and the
 well-founded total ordering $\succeq$ on $\myt$ we define
 $t = \max_{\succeq} \{ \hterm(m_j \rmult f_{j} \rmult l_{j}) \mid 1\leq j \leq m \}$ and
 $K$ as the number of polynomials $m_j \rmult f_j \rmult l_j$ with head term $t$.
Then for each monomial multiple $m_j \rmult f_j \rmult l_j$ with
 $\hterm(m_j \rmult f_j \rmult l_j) = t$ we know that 
 $\hterm(m_j \rmult f_j \rmult l_j) = \hterm(m_j \rmult \hterm(f_j) \rmult l_j)
 \geq \hterm(f_j)$ holds.
%\\
Then $t \succeq \hterm(g)$ and 
 in case $\hterm(g) = t$ this immediately implies that this representation is
 already a reductive standard  one. 
%\\
Else we proceed by induction
 on $t$.
%\\
Without loss of generality let $f_1, \ldots, f_K$ be the polynomials
 in the corresponding representation
 such that  $t=\hterm(m_i \rmult f_i \rmult l_i)$, $1 \leq i \leq K$.
Then the tuple $(t, f_1, \ldots, f_K, m_1, \ldots, m_K, l_1, \ldots, l_K)$
 is in ${\cal C}_{r}(F)$ and let $h = \sum_{i=1}^K m_i \rmult f_i \rmult l_i$.
%\\
We will now change our representation of $g$ in such a way that for the new
 representation of $g$ we have a smaller maximal term.
%\\
Let us assume $h$ is not $\zero$\footnote{In case  $h =\zero$,
 just substitute the empty sum for the representation of $h$
 in the equations below.}. 
%\\
By our assumption, $h$ has a reductive standard representation
 with respect to $F$, say $\sum_{j=1}^n \tilde{m}_j \rmult h_j \rmult \tilde{l}_j$, 
 where $h_j \in F$, and $\tilde{m}_j, \tilde{l}_j \in \monoms(\f_{\myk})$ 
 and  all terms occurring in the sum are bounded by
 $t \succ \hterm(h)$ as $\sum_{i=1}^K \hm(m_i \rmult f_i \rmult l_i) = \zero$.
%\\
This gives us: 
\begin{eqnarray}
  g   & = & \sum_{i=1}^K m_i \rmult f_i \rmult l_i + \sum_{i=K+1}^m m_i \rmult f_i \rmult l_i
             \nonumber\\                                                           
  & = & \sum_{j=1}^n \tilde{m}_j \rmult h_j \rmult \tilde{l}_j + \sum_{i=K+1}^m m_i \rmult f_i \rmult l_i
              \nonumber
\end{eqnarray}
which is a representation of $g$ where the maximal term is smaller than $t$.
\\
\qed
An algebraic characterization of weak Gr\"obner bases again can be given by a property of head monomials
 based on stable divisors of terms (compare Definition \ref{def.gb.fr}).
\begin{definition}\label{def.two-sided.gb}~\\
{\rm
A set $F \subseteq \f_{\myk} \backslash \{ \zero\}$ is called a 
 \betonen{weak reductive Gr\"obner basis} of $\ideal{}{}(F)$ (with respect
 to the reductive ordering $\geq$) if
$\hterm(\ideal{}{}(F) \backslash \{ \zero \}) =
 \hterm( \{ m \rmult f \rmult l \mid f \in F, m, l \in \monoms(\f_{\myk}),
 \hterm(m \rmult f \rmult l) = \hterm(m \rmult \hterm(f) \rmult l) \geq \hterm(f) \}
 \backslash \{ \zero \} )$.
\dend
}
\end{definition}
We will later on see that an analogon of the Translation Lemma holds for the reduction relation
 related to reductive standard representations.
Hence weak reductive  Gr\"obner bases and Gr\"obner bases coincide.
This is again due to the fact that the coefficient domain is a field and will not carry over
 for reduction rings as coefficient domains.

The next lemma states that in fact both characterizations of special bases provided so far coincide.
\begin{lemma}\label{lem.two-sided.sb=gb}~\\
{\sl
Let $F$ be a subset of $\f_{\myk} \backslash \{ \zero\}$.
Then $F$ is a reductive standard basis if and only if it is a weak reductive  Gr\"obner basis.
\lemend
}
\end{lemma}
\Ba{}~\\
Let us first assume that $F$ is a reductive standard basis, i.e., every polynomial
 $g$ in $\ideal{}{}(F)$ has a reductive standard representation with respect
 to $F$.
In case $g \neq \zero$ this implies the existence of a polynomial $f \in F$ and monomials
 $m, l \in \monoms(\f_{\myk})$ such that $\hterm(g) = \hterm(m \rmult f \rmult l)
 = \hterm(m \rmult \hterm(f) \rmult l) \geq \hterm(f)$.
Hence $\hterm(g) \in \hterm( \{ m \rmult  f \rmult l \mid m,l \in \monoms(\f_{\myk}),
 f \in F, \hterm(m \rmult f \rmult l) = \hterm(m \rmult \hterm(f) \rmult l) \geq \hterm(f) \}\backslash \{ \zero\})$.
As the converse, namely $\hterm( \{ m \rmult  f \rmult l \mid m,l \in \monoms(\f_{\myk}),
 f \in F, \hterm(m \rmult f \rmult l) = \hterm(m \rmult \hterm(f) \rmult l) \geq \hterm(f) \}\backslash \{ \zero\})
 \subseteq \hterm(\ideal{}{}(F) \backslash \{ \zero \})$ trivially holds,
 $F$ is a weak reductive  Gr\"obner basis.
\\
Now suppose that $F$ is a weak reductive  Gr\"obner basis and again let $g \in \ideal{}{}(F)$.
We have to show that $g$ has a reductive standard representation with respect to $F$.
This will be done by induction on $\hterm(g)$.
In case $g = \zero$ the empty sum is our required  reductive standard representation.
Hence let us assume $g \neq \zero$.
Since then $\hterm(g) \in \hterm(\ideal{}{}(F)\backslash \{ \zero\})$
 by the definition of weak reductive 
 Gr\"obner bases we know there exists
 a polynomial $f \in F$ and monomials $m,l \in \monoms(\f_{\myk})$ such that
 $\hterm(m \rmult f \rmult l) = \hterm(m \rmult \hterm(f) \rmult l) \geq \hterm(f)$ and there exists $\alpha \in \myk$ such that
 $\hc(g) = \hc(m \rmult f \rmult l) \skm \alpha$, i.e., $\hm(g)=\hm(m \rmult f \rmult l \skm \alpha)$.
Let $g_1 = g - m \rmult f \rmult l \skm \alpha$. 
Then $\hterm(g) \succ \hterm(g_1)$ implies the
 existence of a reductive standard representation for $g_1$ which can be
 added to the multiple $m \rmult f \rmult l \skm \alpha$ to give
 the desired reductive standard representation of $g$.
\\ 
\qed
A close inspection of this proof reveals that in fact we can provide
 a stronger condition for standard representations in terms of weak reductive Gr\"obner bases.
\begin{corollary}\label{cor.two-sided.reductive_rr}~\\
{\sl
Let a subset $F$  of $\f_{\myk} \backslash \{ \zero\}$ be a 
 weak reductive  Gr\"obner basis.
Every $g \in \ideal{}{}(F)$ has a reductive standard representation
       in terms of $F$ of the form
        $g = \sum_{i=1}^n m_i \rmult f_i \rmult l_i,
         f_i \in F, m_i, l_i \in \monoms(\f_{\myk}), n \in \n$
         such that
         $\hterm(g) = \hterm(m_1 \rmult f_1 \rmult l_1) \succ \hterm(m_2 \rmult f_2 \rmult l_2) \succ \ldots \succ
          \hterm(m_n \rmult f_n \rmult l_n)$ and
$\hterm(m_i \rmult f_i \rmult l_i) =
         \hterm(m_i \rmult \hterm(f_i) \rmult l_i) \geq \hterm(f_i)$ for all $1 \leq i \leq n$.
}
\end{corollary}
The importance of Gr\"obner bases in commutative polynomial
 rings stems from the fact that they can be characterized by
 special polynomials, the so-called s-polynomials.
This characterization can be combined with a reduction relation to an algorithm
 which computes finite Gr\"obner bases.

We provide a first characterization for our function ring over the field $\myk$.
Here critical situations  lead to s-polynomials as in the original case
 and can be identified by studying term multiples of polynomials.
Let $p$ and $q$ be two non-zero polynomials in $\f_{\myk}$.
We are interested in terms $t, u_1, u_2, v_1, v_2$ such that
 $ \hterm(u_1 \rmult p \rmult v_1) = \hterm(u_1 \rmult \hterm(p) \rmult v_1) = t =
   \hterm(u_2 \rmult q \rmult v_2) = \hterm(u_2 \rmult \hterm(q) \rmult v_2)$
 and $\hterm(p) \leq t$, $\hterm(q) \leq t$.
Let ${\cal C}_s(p,q)$ (this is a specialization of Definition \ref{def.two-sided.critical.situations})
 be the set containing all such tuples $(t,u_1,u_2, v_1, v_2)$ (as a short hand for $(t, p, q, u_1, u_2, v_1, v_2)$.
We call the polynomial $\hc(u_1 \rmult p \rmult v_1)^{-1} \skm u_1 \rmult p \rmult v_1 -
 \hc(u_2 \rmult q \rmult v_2)^{-1} \skm u_2 \rmult q \rmult v_2
 = \spol{}{}(p,q,t,u_1,u_2,v_1,v_2)$ the \betonen{s-polynomial}
 of $p$ and $q$ related to the tuple $(t,u_1,u_2, v_1, v_2)$.

Again these critical situations
 are not sufficient to characterize weak Gr\"obner bases (compare Example \ref{exa.free.group})
 and additionally we have to test monomial multiples of polynomials now from both sides.
\begin{theorem}\label{theo.two-sided.gb_rr}~\\
{\sl
Let $F$ be a set of polynomials in $\f_{\myk} \backslash \{ \zero\}$.
Then $F$ is a weak Gr\"obner basis of $\ideal{}{}(F)$ if and only if
\begin{enumerate}
\item for all $f$ in $F$ and for all $m, l$ in $\monoms(\f_{\myk})$ the multiple
      $m \rmult f \rmult l$ has a reductive standard representation in terms of $F$, and
\item for all $p$ and $q$ in $F$ and every tuple $(t,u_1,u_2, v_1, v_2)$
      in ${\cal C}_s(p,q)$ the respective s-polynomial $\spol{}{}(p,q,t,u_1,u_2,v_1,v_2)$
      has a reductive standard
      representation in terms of $F$.
\end{enumerate}
}
\end{theorem}
\Ba{}~\\
In case $F$ is a weak Gr\"obner basis it is also a reductive standard basis, and
 since the multiples $m \rmult f \rmult l$ as well as
 the respective s-polynomials are all elements of $\ideal{}{}(F)$ they must
 have reductive standard representations in terms of $F$.
\\
The converse will be proven by showing that every element in
 $\ideal{}{}(F)$ has a reductive standard representation in terms of $F$.
Now, let $g = \sum_{j=1}^m \alpha_{j} \skm v_j \rmult f_{j} \rmult w_{j}$ be an arbitrary
 representation of a non-zero  polynomial $g\in \ideal{}{}(F)$ such that
 $\alpha_{j} \in \myk^*, f_j \in F$, and $v_j, w_{j} \in \myt$.
Since by our first assumption every multiple $v_j \rmult f_j \rmult w_j$ in this
 sum has a reductive standard representation we can assume that
 $\hterm(v_j \rmult \hterm(f_j) \rmult w_j) = \hterm(v_j \rmult f_j \rmult w_j) \geq
 \hterm(f_j)$ holds.
\\
Depending on this  representation of $g$ and the
 well-founded total ordering $\succeq$ on $\myt$ we define
 $t = \max_{\succeq} \{ \hterm(v_j \rmult f_{j} \rmult w_{j}) \mid 1\leq j \leq m \}$ and
 $K$ as the number of polynomials $v_j \rmult f_j \rmult w_j$ with head term $t$.
%\\
Without loss of generality we can assume that the polynomial multiples
 with head term $t$ are just $v_1 \rmult f_1  \rmult w_1, \ldots , v_K \rmult f_K \rmult w_K$.
%\\
We proceed by induction
 on $(t,K)$, where
 $(t',K')<(t,K)$ if and only if $t' \prec t$ or $(t'=t$ and
 $K'<K)$\footnote{Note that this ordering is well-founded since $\succ$
                  is well-founded on $\myt$ and $K \in\n$.}.
\\
Obviously, $t \succeq \hterm(g)$ must hold. 
If $K = 1$ this gives us $t = \hterm(g)$ and by our assumptions our
 representation is already of the required form.
Hence let us assume $K > 1$.
Then for the two polynomials $f_1,f_2$
 in the corresponding 
      representation\footnote{Not necessarily $f_2 \neq f_1$.}
      such that  $t=\hterm(v_1 \rmult \hterm(f_1) \rmult w_1) =\hterm(v_1 \rmult f_1 \rmult w_1) =
       \hterm(v_2 \rmult f_2 \rmult w_2) =  \hterm(v_2 \rmult \hterm(f_2) \rmult w_2)$ and
       $t \geq \hterm(f_1)$, $t \geq \hterm(f_2)$.
     Then the tuple $(t, v_1, v_2, w_1,w_2)$ is in ${\cal C}_s(f_1,f_2)$ and
      we have an s-polynomial
      $h= \hc(v_1 \rmult f_1 \rmult w_1)^{-1} \skm  v_1 \rmult f_1
      \rmult w_1 -\hc(v_2 \rmult f_2 \rmult w_2)^{-1}\skm v_2 \rmult f_2 \rmult w_2$ corresponding
      to this tuple.
%\\
We will now change our representation of $g$ by using the additional
 information on this s-polynomial in such a way that for the new
 representation of $g$ we either have a smaller maximal term or
 the occurrences of the term $t$
 are decreased by at least 1.
%\\
Let us assume the s-polynomial is not $\zero$\footnote{In case 
               $h =\zero$,
               just substitute the empty sum
               for the  reductive representation of $h$
               in the equations below.}. 
%\\
By our assumption, $h$ has a reductive standard representation in terms of $F$, say
  $\sum_{i=1}^n \tilde{\alpha}_i \skm \tilde{v}_i \rmult \tilde{f}_i \rmult \tilde{w}_i$, 
  where $\tilde{\alpha}_i \in \myk^*,\tilde{f}_i \in F$, and $\tilde{v}_i,\tilde{w}_i \in \myt$ 
  and  all terms occurring in this sum are bounded by
  $t \succ \hterm(h)$.
%\\
This gives us: 
     \begin{eqnarray}
       &   & \alpha_1 \skm v_1 \rmult  f_1 \rmult w_1 + \alpha_2 \skm v_2 \rmult  f_2 \rmult w_2
             \nonumber\\  
       &   &  \nonumber\\                                                         
       & = &  \alpha_1 \skm v_1 \rmult f_1 \rmult w_1 +
              \underbrace{ \alpha'_2 \skm \beta_1 \skm v_1 \rmult f_1 \rmult w_1
                   - \alpha'_2 \skm \beta_1 \skm v_1 \rmult f_1 \rmult w_1}_{=\, 0}
              \nonumber\\ 
       &   &  + \underbrace{\alpha'_2\skm \beta_2}_{\alpha_2}  \skm v_2 \rmult f_2 \rmult w_2 \nonumber\\
       &   & \nonumber\\ 
       & = & (\alpha_1 + \alpha'_2 \skm \beta_1) \skm v_1 \rmult f_1 \rmult w_1 - \alpha'_2 \skm 
               \underbrace{(\beta_1 \skm v_1 \rmult f_1 \rmult w_1
             -  \beta_2 \skm v_2 \rmult f_2 \rmult w_2)}_{=\,
             h} 
             \nonumber\\
       &   & \nonumber\\ 
       & = & (\alpha_1 + \alpha'_2 \skm \beta_1) \skm v_1 \rmult f_1 \rmult w_1 - \alpha'_2 \skm
                   (\sum_{i=1}^n \tilde{\alpha}_i \skm \tilde{v}_i \rmult \tilde{f}_i \rmult \tilde{w}_i) \label{s1.i}
     \end{eqnarray}
     where $\beta_1=\hc(v_1 \rmult f_1 \rmult w_1)^{-1}$, $\beta_2=\hc(v_2 \rmult f_2 \rmult w_2)^{-1}$
      and  $\alpha'_2 \skm \beta_2 = \alpha_2$.
     Substituting (\ref{s1.i}) in the representation of $g$ gives rise to a smaller one.
\\
\qed
Notice that both test sets in this characterization in general cannot be described in a finitary manner,
 i.e., provide no finite test for the property of being a Gr\"obner basis.

A problem which is related to the fact that the ordering $\succeq$
 and the multiplication $\rmult$ in general are not compatible
 is that an important property fulfilled for representations of 
 polynomials in commutative polynomial rings no longer holds:
As in the case of right ideals the existence of a standard
 representation for some polynomial $f \in \f_{\myk}$ no longer implies the existence of one for a 
 multiple $m \rmult f \rmult l$ where $m,l \in \monoms(\f_{\myk})$.
However there are restrictions where this implication will hold (compare
 Lemma \ref{lem.red.reps}).
\begin{lemma}\label{lem.two-sided.red.reps}~\\
{\sl
Let  $F$ be a subset of $\f_{\myk} \backslash \{ \zero\}$ and
 $p$ a non-zero polynomial in $\f_{\myk}$.
If $p$ has a reductive standard representation with respect to $F$ and
 $m,l$ are monomials such that $ \hterm(m \rmult p \rmult l) =
 \hterm(m \rmult \hterm(p) \rmult l) \geq \hterm(p)$, then the multiple $m \rmult p \rmult l$ again has a
 reductive standard representation with respect to $F$. 
\lemend
}
\end{lemma}
\Ba{}~\\
Let $p = \sum_{i=1}^n m_i \rmult f_i \rmult l_i$ with $n \in \n$, $f_i \in F$,
 $m_i, l_i \in \monoms(\f_{\myk})$
 be a reductive  standard representation of $p$ in terms of $F$, i.e.,
 $\hterm(p) = \hterm(m_i \rmult \hterm(f_i) \rmult l_i) =
 \hterm(m_i \rmult f_i \rmult l_i) \geq \hterm(f_i)$, $1 \leq i \leq k$ and
 $\hterm(p) \succeq  \hterm(m_i \rmult f_i \rmult l_i)$
 for all $k+1 \leq i \leq n$. 
\\
Let us first analyze the multiple $m \rmult m_j \rmult f_j \rmult l_j \rmult l$.\\
Let $\terms(m_j \rmult f_j \rmult l_j) = \{ s_1, \ldots, s_k \}$ with $s_1 \succ s_i$, $2 \leq i \leq l$,
 i.e.~$s_1 = \hterm(m_j \rmult f_j \rmult l_j) = \hterm(m_j \rmult \hterm(f_j) \rmult l_j) = \hterm(p)$.
Hence $\hterm(m \rmult \hterm(p) \rmult l) = \hterm(m \rmult s_1 \rmult l) \geq \hterm(p) = s_1$
 and as $s_1 \succ s_i$, $2 \leq i \leq l$,
 by Definition \ref{def.two-sided.refined.ordering}
 we can conclude $\hterm(m \rmult \hterm(p) \rmult l) 
 = \hterm(m \rmult s_1 \rmult l) \succ m \rmult s_i \rmult l \succeq
 \hterm(m \rmult s_i \rmult l)$
 for $2 \leq i \leq l$.
This implies $\hterm(m \rmult \hterm(m_j \rmult f_j \rmult l_j) \rmult l) = \hterm(m \rmult m_j \rmult f_j \rmult l_j \rmult l)$.
Hence we get
\begin{eqnarray*}
 \hterm(p \rmult m) & = & \hterm(m \rmult \hterm(p) \rmult l) \\
                    & = & \hterm(m \rmult \hterm(m_j \rmult f_j \rmult l_j) \rmult l), \mbox{ as } 
                          \hterm(p) = \hterm(m_j \rmult f_j \rmult l_j) \\
                    & = & \hterm(m \rmult m_j \rmult f_j \rmult l_j \rmult l)
\end{eqnarray*}
and
since $\hterm(m \rmult p \rmult l) \geq \hterm(p) \geq \hterm(f_j)$ we can conclude
 $\hterm(m \rmult m_j \rmult f_j \rmult l_j \rmult l) \geq \hterm(f_j)$.
It remains to show that $m \rmult m_j \rmult f_j \rmult l_j \rmult l$ has
 a reductive standard representation in terms of $F$.
First we show that $\hterm(m \rmult m_j \rmult \hterm(f_j) \rmult l_j \rmult l)
 \geq \hterm(f_j)$.
We know $m_j \rmult \hterm(f_j) \rmult l_j \succeq
  \hterm(m_j \rmult \hterm(f_j) \rmult l_j) = \hterm(m_j \rmult f_j \rmult l_j)$\footnote{Notice that $m_j \rmult \hterm(f_j) \rmult l_j$ can be a polynomial and
 hence we cannot conclude $m_j \rmult \hterm(f_j) \rmult l_j =
 \hterm(m_j \rmult \hterm(f_j) \rmult l_j)$.}
 and hence $\hterm(m \rmult m_j \rmult \hterm(f_j) \rmult l_j \rmult l) =
 \hterm(m \rmult \hterm(m_j \rmult f_j \rmult l_j) \rmult l) =
 \hterm(m \rmult m_j \rmult f_j \rmult l_j \rmult l) \geq \hterm(f_j)$.
Now in case $m \rmult m_j, l_j \rmult l \in \monoms(\f_{\myk})$ we are done as then $(m_j \rmult m) \rmult f_j \rmult ( l_j \rmult l)$
 is a  reductive standard representation in terms of $F$.
\\
%Since $\hterm(m \rmult \hterm(p) \rmult l) =
% \hterm(m \rmult p \rmult l) \geq \hterm(p)$ and for all terms
% $s \neq \hterm(m_1 \rmult f_1 \rmult l_1)=\hterm(p)$
% in the polynomial $m_1 \rmult f_1 \rmult l_1$
% we have $\hterm(p) \succ s$  we get
% $\hterm(m \rmult \hterm(p) \rmult l) \succ \hterm(m \rmult s \rmult l)$ and in particular
% $\hterm(m \rmult p \rmult l) = \hterm(m \rmult \hterm(p) \rmult l) \succ \hterm(m \rmult s \rmult l)$.
%This now implies
% $\hterm(m \rmult \hterm(m_1 \rmult \hterm(f_1) \rmult l_1) \rmult l) 
%  = \hterm(m \rmult \hterm(m_1 \rmult f_1 \rmult l_1) \rmult l)
%  = \hterm(m \rmult m_1 \rmult f_1 \rmult l_1 \rmult l)$.
%Hence we get
%\begin{eqnarray*}
% \hterm(m \rmult p \rmult l) & = & \hterm(m \rmult \hterm(p) \rmult l) \\
%                    & = & \hterm(m \rmult \hterm(m_1 \rmult f_1 \rmult l_1) \rmult l), \mbox{ as } 
%                          \hterm (p) = \hterm(m_1 \rmult f_1 \rmult l_1) \\
%                    & = & \hterm(m \rmult m_1 \rmult f_1 \rmult l_1 \rmult l)
%\end{eqnarray*}
%and
%since $\hterm(m \rmult p \rmult l) \geq \hterm(f_1)$ we can conclude
% $\hterm(m \rmult m_1 \rmult f_1 \rmult l_1 \rmult l) \geq \hterm(f_1)$.
%\\
%It remains to show that $m \rmult (m_1 \rmult f_1 \rmult l_1) \rmult l = 
% (m_1 \rmult m) \rmult f_1 \rmult (l_1 \rmult l)$ has
% a reductive standard representation in terms of $F$.
%\\
%In case $m_1 \rmult m, l_1 \rmult l \in \monoms(\f_{\myk})$
%  we are done as then $(m_1 \rmult m) \rmult f_1 \rmult (l_1 \rmult l)$
% is one.
%\\
Hence let us assume
 $m \rmult m_j = \sum_{i=1}^{k_1} \tilde{m}_i$,
 $l_j \rmult l = \sum_{i'=1}^{k_1'} \tilde{l}_{i'}$,
 $\tilde{m}_i, \tilde{l}_{i'} \in \monoms(\f_{\myk})$.
Let $\terms(f_j) = \{ t_1, \ldots, t_w \}$
 with $t_1 \succ t_i$,  $2 \leq i \leq w$,
 i.e.~$t_1 = \hterm(f_j)$.
As $\hterm(m_j \rmult  \hterm(f_j) \rmult l_j) \geq \hterm(f_j) \succ t_p$,
 $2 \leq p \leq w$, again by Definition \ref{def.two-sided.refined.ordering}
 we can conclude $\hterm(m_j \rmult  \hterm(f_j) \rmult l_j) \succ
 m_j \rmult t_p \rmult l_j \succeq \hterm(m_j \rmult t_p \rmult l_j)$,
 and $m_j \rmult  \hterm(f_j) \rmult l_j \succ
 \sum_{p=2}^w m_j \rmult t_p \rmult l_j $.
Then for each $s_i$, $2 \leq i \leq l$
 there exists $t_q \in \terms(f_1)$ such that
 $s_i \in \supp(m_j \rmult t_q \rmult l_j)$.
Since $\hterm(p) \succ s_i$ and even
 $\hterm(p) \succeq m_j \rmult t_q \rmult l_j$ we find that
 either
 $\hterm(m \rmult p \rmult l) \succeq \hterm(m \rmult (m_j \rmult t_q \rmult l_j) \rmult l) = \hterm((m \rmult m_j) \rmult t_q \rmult (l_j \rmult l))$
 in case $\hterm(m_j \rmult t_q \rmult l_j ) = \hterm(m_j \rmult f_j \rmult l_j)$
 or 
 $\hterm(m \rmult p \rmult l) \succ \hterm(m \rmult (m_j \rmult t_q \rmult l_j) \rmult l) =
 \hterm((m \rmult m_j) \rmult  t_q \rmult (l_j \rmult l))$.
Hence we can conlude $\tilde{m}_i \rmult f_j \rmult \tilde{l}_{i'} \predeq
  \hterm(m \rmult p \rmult l)$, $1 \leq i \leq k_1$, $1 \leq i' \leq k_1'$ and
 for at least one such multiple we get
 $\hterm(\tilde{m}_i \rmult f_1 \rmult \tilde{l}_{i'}) = \hterm(m \rmult m_j \rmult f_j \rmult l_j \rmult l) \geq \hterm(f_j)$.
\\
It remains to analyze the situation for the function
 $(\sum_{i=k+1}^n m \rmult (m_i \rmult f_i \rmult l_i) \rmult l$.
Again we find that for all terms $s$ in the $m_i \rmult f_i \rmult l_i$,
 $k+1 \leq i \leq n$,
 we have $\hterm(p) \succ s$ and we get 
 $\hterm(m \rmult p \rmult l) \succeq \hterm(m \rmult s \rmult l)$.
Hence all polynomial multiples of the $f_i$  in the representation 
 $\sum_{i=k+1}^n ((\sum_{j=1}^{k_i} \tilde{m}^i_j) \rmult f_i \rmult (\sum_{j=1}^{k'_i} \tilde{l}^i_j))$, where
 $m \rmult m_i = \sum_{j=1}^{k_i}\tilde{m}^i_j$,
 $l_i \rmult l = \sum_{j=1}^{k'_i}\tilde{l}^i_j$, are bounded by 
 $\hterm(m \rmult p \rmult l)$.
\\
\qed
Notice that this lemma  no longer holds in case we only require
 $\hterm(m \rmult \hterm(p) \rmult l) =
 \hterm(m \rmult p \rmult l) \succeq \hterm(p)$, as then $\hterm(p) \succ s$
 no longer implies 
 $\hterm(m \rmult p \rmult l) \succ \hterm(m \rmult s \rmult l)$.

Our standard representations from Definition \ref{def.two-sided.right_reductive}
 are closely related to a reduction relation
 based on the divisibility of terms as defined in the context of 
  reductive restrictions of orderings on page \pageref{def.two-sided.refined.ordering}.
\begin{definition}\label{def.two-sided.red}~\\
{\rm
Let $f,p$ be two non-zero polynomials in $\f_{\myk}$.
We say $f$ \betonen{reduces} $p$ \betonen{to} $q$ \betonen{at a monomial}
 $\alpha \skm t$ \betonen{in one step}, denoted by $p \red{}{\myr}{}{f} q$, if
 there exist $m,l \in \monoms(\f_{\myk})$ such that
\begin{enumerate}
\item $t \in \supp(p)$ and $p(t) = \alpha$,
\item $\hterm(m \rmult \hterm(f)\rmult l) = \hterm(m \rmult f \rmult l) = t \geq \hterm(f)$,
\item $\hm(m \rmult f \rmult l) = \alpha \skm t$, and
\item $q = p - m \rmult f \rmult l$.
\end{enumerate}
We write $p \red{}{\myr}{}{f}$ if there is a polynomial $q$ as defined
above and $p$ is then called  reducible by $f$. 
%\\
Further, we can define $\red{*}{\myr}{}{}, \red{+}{\myr}{}{}$ and
 $\red{n}{\myr}{}{}$ as usual.
%\\
 Reduction by a set $F \subseteq \f_{\myk} \backslash \{ \zero \}$ is denoted by
 $p \red{}{\myr}{}{F} q$ and abbreviates $p \red{}{\myr}{}{f} q$
 for some $f \in F$.
\dend
}
\end{definition}
Due to the fact that the coefficients lie in a field, again if for some terms $w_1,w_2 \in \myt$ we
 have $\hterm(w_1 \rmult f \rmult w_2) = \hterm(w_1 \rmult \hterm(f) \rmult w_2) = t \geq \hterm(f)$
 this implies reducibility at the monomial $\alpha \skm t$.
\begin{lemma}~\\ 
{\sl 
Let $F$ be a set of polynomials in $\f_{\myk} \backslash \{ \zero\}$.
\begin{enumerate}
\item For $p,q \in \f_{\myk}$ we have that 
       $p \red{}{\myr}{}{F} q$ implies $p \succ q$, in particular $\hterm(p)
       \succeq \hterm(q)$.
\item $\red{}{\myr}{}{F}$ is Noetherian.
\lemend
\end{enumerate}
}
\end{lemma}
\Ba{}
\begin{enumerate}
\item Assuming that the reduction step takes place at a monomial $\alpha \skm t$,
      by Definition \ref{def.two-sided.red} we know $\hm(m_1 \rmult f \rmult m_2) = \alpha \skm t$ which yields
      $p \succ p  -  m_1 \rmult f \rmult m_2$
      since $\hm( m_1 \rmult f \rmult m_2) \succ \reductum( m_1 \rmult f \rmult m_2)$.
\item This follows directly from 1. as the ordering $\succeq$ on $\myt$ is well-founded
 (compare Lemma \ref{lem.wellfounded}).
\end{enumerate}\renewcommand{\baselinestretch}{1}\small\normalsize
\qed
The next lemma shows how reduction sequences and  reductive standard
 representations are related.
\begin{lemma}\label{lem.two-sided.red.rep}~\\
{\sl
Let $F$ be a set of polynomials in $\f_{\myk}$ and $p$
 a non-zero polynomial in $\f_{\myk}$.
Then $p \red{*}{\myr}{}{F} \zero$ implies that $p$ has a   reductive standard
 representation in terms of $F$.
\lemend
}
\end{lemma}
\Ba{}~\\
This follows directly by adding up the polynomials used in the
 reduction steps occurring in the reduction sequence $p \red{*}{\myr}{}{F} \zero$.
\\
\qed
If $p \red{*}{\myr}{}{F} q$, then $p$ has a  reductive standard representation in 
 terms of $F \cup \{ q \}$, especially $p-q$ has one in terms of $F$.

As stated before an analogon to the Translation Lemma holds.
\begin{lemma}\label{lem.two-sided.trans}~\\
{\sl
Let $F$ be a set of polynomials in $\f_{\myk}$ and $p,q,h$ polynomials in $\f_{\myk}$.
\begin{enumerate}
\item
Let $p-q \red{}{\myr}{}{F} h$.
Then there exist $p', q' \in \f_{\myk}$ such that $p \red{*}{\myr}{}{F} p'$
 and $q \red{*}{\myr}{}{F} q'$ and $h = p'-q'$.
\item
Let $\zero$ be a normal form of $p-q$ with respect to $F$.
Then there exists $g \in \f_{\myk}$ such that $p \red{*}{\myr}{}{F} g$ and
 $q \red{*}{\myr}{}{F} g$.
\end{enumerate}
\lemend}
\end{lemma}
\pagebreak
\Ba{}~\\
\begin{enumerate}
\item
Let $p-q \red{}{\myr}{}{F} h$ at the monomial $\alpha \skm t$, i.e.,
 $h = p-q - m \rmult f \rmult l$ for some $m, l \in \monoms(\f_{\myk})$ such that
 $\hterm(m \rmult \hterm(f)\rmult l) = \hterm(m \rmult f \rmult l) = t \geq \hterm(f)$
 and $\hm(m \rmult f \rmult l) = \alpha \skm t$.
We have to distinguish three cases:
\begin{enumerate}
\item $t \in \supp(p)$ and $t \in \supp(q)$:
      Then we can eliminate the occurence of $t$ in the respective polynomials
        by reduction and get
         $p \red{}{\myr}{}{f} p - \alpha_1 \skm (m \rmult f \rmult l)
        = p'$, $q \red{}{\myr}{}{f} q - \alpha_2 \skm (m \rmult f \rmult l) = q'$, 
        where $\alpha_1 \skm \hc(m \rmult f \rmult l)$
        and $\alpha_2 \skm \hc(m \rmult f \rmult l)$ are the coefficients of $t$ in
        $p$ respectively $q$.
      Moreover, $\alpha_1 \skm \hc(m \rmult f \rmult l) - \alpha_2 \skm \hc(m \rmult f \rmult l) = 
       \alpha$ and hence $\alpha_1 - \alpha_2 = 1$, as $\hc(m \rmult f \rmult l) = \alpha$.
      This gives us $p' - q' = p - \alpha_1 \skm (m \rmult f \rmult l) - 
        q + \alpha_2 \skm (m \rmult f \rmult l) =
        p - q - (\alpha_1 - \alpha_2) \skm (m \rmult f \rmult l) =
        p-q - m \rmult f \rmult l = h$.
\item $t \in \supp(p)$ and $t \not\in \supp(q)$:
      Then we can eliminate the term $t$ in the polynomial $p$ by right reduction
        and get $p \red{}{\myr}{}{f} p - m \rmult f \rmult l = p'$, $q=q'$,
        and, therefore, $p' - q' = p - m \rmult f \rmult l - q =h$.
\item $t \in \supp(q)$ and $t \not\in \supp(p)$:
      Then we can eliminate the term $t$ in the polynomial $q$ by right reduction
        and get $q \red{}{\myr}{}{f} q + m \rmult f \rmult l = q'$, $p=p'$, and, therefore,
        $p' - q' = p - (q + m \rmult f \rmult l) =h$.
\end{enumerate}
\item We show our claim by induction on $k$, where $p-q \red{k}{\myr}{}{F} \zero$.
In the base case $k=0$ there is nothing to show as then $p=q$.
Hence, let $p-q \red{}{\myr}{}{F} h \red{k}{\myr}{}{F} \zero$.
Then by 1.~there are polynomials $p',q' \in \f_{\myk}$ such that 
 $p \red{*}{\myr}{}{F} p'$
 and $q \red{*}{\myr}{}{F} q'$ and $h = p'-q'$.
Now the induction hypothesis for $p'-q' \red{k}{\myr}{}{F} \zero$ yields the
 existence of a polynomial $g \in \f_{\myk}$ such that
 $p \red{*}{\myr}{}{F} g$ and $q \red{*}{\myr}{}{F} g$.
\end{enumerate}
\qed 
The essential part of the proof is that reducibility as defined in Definition \ref{def.two-sided.red}
 is connected to stable divisors of terms and not to coefficients.
We will later see that for function rings over reduction rings, when the coefficient
 is also involved in the reduction step, this lemma no longer holds.

Next we state the definition of Gr\"obner bases based on the reduction relation.
\begin{definition}\label{def.two-sided.gb.reduction}~\\
{\rm
A subset $G$ of  $\f_{\myk}$ is called a 
 \betonen{Gr\"obner basis} (with respect to the reduction relation
 $\red{}{\myr}{}{}$) of the
 ideal ${\mathfrak i}= \ideal{}{}(G)$, if
 $\red{*}{\longleftrightarrow}{}{G} = \;\;\equiv_{{\mathfrak i}}$
 and $\red{}{\myr}{}{G}$ is confluent.
}
\end{definition}

Remember the free group ring in Example \ref{exa.free.group} where the polynomial $b + \lambda$ lies
 in the ideal generated by the polynomial $a + \lambda$.
Then of course $b+\lambda$  also lies in the ideal generated by $a+\lambda$.
Unlike in the case of polynomial rings over fields where for any set
 of polynomials $F$ we have
 $\red{*}{\longleftrightarrow}{b}{F} = \;\;\equiv_{\ideal{}{}(F)}$,
 here we have
 $b + \lambda \equiv_{\ideal{}{}(\{a + \lambda \})} 0$ but
 $b + \lambda \nred{*}{\longleftrightarrow}{}{a + \lambda} 0$.
Hence the first condition of Definition \ref{def.two-sided.gb.reduction} is again
 neccessary.

Now by Lemma \ref{lem.two-sided.trans} and Theorem \ref{theo.weak+translation}
 weak Gr\"obner bases are Gr\"obner bases
 and can be characterized 
 as follows:

\begin{corollary}\label{lem.two-sided.gb=red.to.zero}~\\
{\sl
Let $G$ be a set of polynomials in $\f_{\myk} \backslash \{ \zero\}$.
$G$ is a (weak) Gr\"obner basis of $\ideal{}{}(G)$
 if and only if for every $g \in \ideal{}{}(G)$
 we have $g\red{*}{\myr}{}{G} \zero$.  
}
\end{corollary}

Finally we can characterize Gr\"obner bases similar to Theorem 
 \ref{theo.buchberger.completion}.
\begin{theorem}\label{theo.two-sided.s-pol}~\\
{\sl
Let $F$ be a set of polynomials in $\f_{\myk} \backslash \{ \zero \}$.
Then $F$ is a Gr\"obner basis of $\ideal{}{}(G)$ if and only if
\begin{enumerate}
\item for all $f$ in $F$ and for all $m, l$ in $\monoms(\f_{\myk})$ we have
       $m \rmult f \rmult l \red{*}{\myr}{}{F} \zero$, and
\item for all  $p$ and $q$ in $F$ and every tuple $(t,u_1,u_2,v_1,v_2)$
      in ${\cal C}(p,q)$ and the respective s-polynomial $\spol{}(p,q,t,u_1,u_2,v_1,v_2)$
       we have $\spol{}(p,q,t,u_1,u_2,v_1,v_2) \red{*}{\myr}{}{F} \zero$.
\end{enumerate}
\theoend
}
\end{theorem}
We will later on prove a stronger version of this theorem.

The importance of Gr\"obner bases in the classical case stems
 from the fact that we only have to check a finite set of s-polynomials for $F$
 in order to decide, whether $F$ is a Gr\"obner basis.
Hence, we are interested in localizing the test sets in Theorem \ref{theo.two-sided.s-pol} --
 if possible to finite ones.
\begin{definition}\label{def.two-sided.weak.sat}~\\
{\rm
A set of polynomials $F \subseteq \f_{\myk} \backslash \{ \zero\}$ is called
 \betonen{weakly saturated}, if for all monomials $m,l$ in $\monoms(\f_{\myk})$
 and every polynomial
 $f \in F$ we have
 $m\rmult f \rmult l \red{*}{\myr}{}{F} \zero$.
\dend
}
\end{definition}
This of course implies that for a weakly saturated set $F$ and any
 $m,l \in\monoms(\f_{\myk})$, $f \in F$ the multiple
 $m \rmult f \rmult l$ has a reductive standard representation in terms of $F$.

Notice that since the coefficient domain is a field we could restrict
 ourselves to multiples with elements of $\myt$.
However, as we will later on allow reduction rings as coefficient
 domains, we present this more general definition.
\begin{definition}\label{def.two-sided.saturator}~\\
{\rm
Let $F$ be a set of polynomials in $\f_{\myk} \backslash \{ 0 \}$.
A set
 $\sat(F) \subseteq \{ m \rmult f \rmult l \mid f \in F, m, l \in \monoms(\f_{\myk}) \}$
 is called a \betonen{stable saturator} for $F$ if for any $f \in F$, $m, l \in \monoms(\f_{\myk})$
 there exist $s \in \sat(F)$, $m', l' \in\monoms(\f_{\myk})$ such that
 $m \rmult f \rmult l = m' \rmult s \rmult l'$,
 $\hterm(m \rmult f \rmult l) = \hterm(m' \rmult \hterm(s) \rmult l')
 \geq \hterm(s)$. 
}
\end{definition}
\begin{corollary}~\\
{\sl
Let $\sat(F)$ be a stable saturator of a set $F \subseteq \f_{\myk}$.
Then for any $f \in F$, $m, l \in \monoms(\f_{\myk})$ there exists $s \in \sat(F)$
 such that $m \rmult f \rmult l \red{}{\myr}{}{s} \zero$.
}
\end{corollary}
\begin{lemma}~\\
{\sl
Let $F$ be a set of polynomials in $\f_{\myk} \backslash \{ 0 \}$.
If for all $s\in \sat(F)$ we have $ s \red{*}{\myr}{}{F} \zero$, then
 for every $m$, $l$ in $\monoms(\f_{\myk})$
 and every polynomial $f$ in $F$ the multiple
 $m \rmult f \rmult l$ has a reductive standard  representation in terms of $F$.
}
\end{lemma}
\Ba{}~\\
This follows immediately from Lemma \ref{lem.two-sided.red.reps} and
 Lemma \ref{lem.two-sided.red.rep}.
\\
\qed
\begin{definition}~\\
{\rm
Let $p$ and $q$ be two non-zero polynomials in $\f_{\myk}$.
Then a subset $C \subseteq \{ \spol{}(p,q,t,u_1,u_2, v_1,v_2) \mid
 (t,u_1,u_2, v_1, v_2) \in {\cal C}_s(p,q) \}$ is called a \betonen{stable
 localization} for the critical situations if for every s-polynomial
 $\spol{}(p,q,t,u_1,u_2,v_1,v_2)$ related to a tuple $(t,u_1,u_2, v_1, v_2)$
 in ${\cal C}_s(p,q)$ there exists a polynomial $h \in C$ and
 monomials $\alpha \skm w_1, 1 \skm w_2 \in \monoms(\f_{\myk})$ such that
 \begin{enumerate}
\item $\hterm(h) \leq \hterm(\spol{}(p,q,t,u_1,u_2,v_1,v_2))$,
 \item $\hterm(w_1 \rmult h \rmult w_2) = \hterm(w_1 \rmult \hterm(h) \rmult w_2) =
       \hterm(\spol{}(p,q,t,u_1,u_2,v_1,v_2))$,
 \item $\spol{}(p,q,t,u_1,u_2,v_1,v_2) = (\alpha \skm w_1) \rmult h \rmult  w_2$.
\dend
 \end{enumerate}
}
\end{definition}
The idea behind this definition is to reduce the number of s-polynomials, which
 have to be considered when checking for the Gr\"obner basis property.
\begin{corollary}~\\
{\sl
Let $C \subseteq \{ \spol{}(p,q,t,u_1,u_2,v_1,v_2) \mid
 (t,u_1,u_2,v_1,v_2) \in {\cal C}_s(p,q) \}$ be a stable localization for
 two polynomials $p,q \in \f_{\myk}$.
Then for any s-polynomial $\spol{}(p,q,t,u_1,u_2,v_1,v_2)$ there exists $h \in C$
 such that $\spol{}(p,q,t,u_1,u_2,v_1,v_2) \red{}{\myr}{}{h} \zero$.
}
\end{corollary}
\begin{lemma}~\\
{\sl
Let $F$ be a set of polynomials in $\f_{\myk} \backslash \{ 0 \}$.
If for all $h$ in 
      a stable localization  $C \subseteq \{ \spol{}(p,q,t,u_1,u_2,v_1,v_2) \mid
    (t,u_1,u_2,v_1,v_2) \in {\cal C}_s(p,q) \}$,
  we have $h \red{*}{\myr}{}{F} \zero$, then
 for every $(t,u_1,u_2,v_1,v_2)$
 in ${\cal C}_s(p,q)$ the s-polynomial  $\spol{}(p,q,t,u_1,u_2,v_1,v_2)$
 has a reductive standard  representation in terms of $F$.
}
\end{lemma}
\Ba{}~\\
This follows immediately from Lemma \ref{lem.two-sided.red.reps} and
 Lemma \ref{lem.two-sided.red.rep}.
\\
\qed
\begin{theorem}\label{theo.two-sided.s-pol.2}~\\
{\sl
Let $F$ be a set of polynomials in $\f_{\myk} \backslash \{ 0 \}$.
Then $F$ is a Gr\"obner basis if and only if
\begin{enumerate}
\item for all $s$ in $\sat(F)$  we have
       $s \red{*}{\myr}{}{F} \zero$, and
\item for all  $p$ and $q$ in $F$, and every polynomial $h$ in 
      a stable localization  $C \subseteq \{ \spol{}(p,q,t,u_1,u_2,v_1,v_2) \mid
      (t,u_1,u_2,v_1,v_2) \in {\cal C}(p,q) \}$,
      we have $h \red{*}{\myr}{}{F} \zero$.
\end{enumerate}
\theoend
}
\end{theorem}
\Ba{}~\\
In case $F$ is a Gr\"obner basis by Lemma \ref{lem.two-sided.gb=red.to.zero}
  all elements of $\ideal{}{}(F)$ must reduce to zero by $F$.
Since the polynomials in the saturator and the respective
 localizations of the s-polynomials all belong to the ideal generated
 by $F$ we are done.
\\
The converse will be proven by showing that every element in
 $\ideal{}{}(F)$ has a  reductive standard representation in terms of $F$.
Now, let $g = \sum_{j=1}^n (\alpha_{j} \skm w_j) \rmult f_{j} \rmult z_{j}$ be an arbitrary
 representation of a non-zero  polynomial $g\in \ideal{}{}(F)$ such that
 $\alpha_{j} \in \myk^*, f_j \in F$, and $w_{j}, z_j \in \myt$.
\\
By the definition of the stable saturator 
 for every multiple $w_j \rmult f_j \rmult z_j$ in this
 sum we have some $s \in \sat(F)$, $m, l \in \monoms(\f_{\myk})$ such that $w_j \rmult f_{j} \rmult z_{j} =
 m \rmult s \rmult l$
 and $\hterm(w_j \rmult f_{j} \rmult z_{j}) = \hterm(m \rmult s \rmult l) = 
 \hterm(m \rmult \hterm(s) \rmult l) \geq \hterm(s)$.
Since we have $s \red{*}{\myr}{}{F} \zero$, by Lemma
 \ref{lem.two-sided.red.reps} we can conclude that each $w_j \rmult f_{j} \rmult z_{j}$ has a reductive
 standard representation in terms of $F$.
Therefore, we can assume that
 $\hterm(w_j \rmult \hterm(f_j) \rmult z_j) = \hterm(w_j \rmult f_j \rmult z_j) \geq
 \hterm(f_j)$ holds.
\\
Depending on this  representation of $g$ and the
 well-founded total ordering $\succeq$ on $\myt$ we define
 $t = \max_{\succeq} \{ \hterm(w_j \rmult f_{j} \rmult z_{j}) \mid 1\leq j \leq n \}$ and
 $K$ as the number of polynomials $w_j \rmult f_j \rmult z_j$ with head term $t$.
\\
Without loss of generality we can assume that the polynomial multiples
 with head term $t$ are just $(\alpha_1 \skm w_1) \rmult  f_1  \rmult z_1, \ldots , (\alpha_K \skm w_K) \rmult  f_K \rmult z_K$.
We proceed by induction
 on $(t,K)$, where
 $(t',K')<(t,K)$ if and only if $t' \prec t$ or $(t'=t$ and
 $K'<K)$\footnote{Note that this ordering is well-founded since $\succ$
                  is well-founded on $\myt$ and $K \in\n$.}.
\\
Obviously, $t \succeq \hterm(g)$ must hold. 
If $K = 1$ this gives us $t = \hterm(g)$ and by our assumption our
 representation is already of the required form.
\\
Hence let us assume $K > 1$, then for the  two
 not necessarily different polynomials $f_1,f_2$
 in the corresponding   representation we have
       $t=\hterm(w_1 \rmult \hterm(f_1) \rmult z_1) =\hterm(w_1 \rmult f_1 \rmult z_1) =
       \hterm(w_2 \rmult f_2 \rmult z_2) =  \hterm(w_2 \rmult \hterm(f_2) \rmult z_2)$ and
       $t \geq \hterm(f_1)$, $t \geq \hterm(f_2)$.
Then the tuple $(t, w_1,w_2,z_1,z_2)$ is in ${\cal C}(f_1,f_2)$ and
      we have a polynomial $h$ in a stable localization 
      $C \subseteq \{ \spol{}(f_1,f_2,t,w_1,w_2,z_1,z_2)
       \mid (t, w_1, w_2,z_1,z_2) \in {\cal C}(f_1,f_2) \}$
      and $\alpha \skm w, 1 \skm z \in \monoms(\f_{\myk})$ such that 
      $\spol{}(f_1,f_2,t,w_1,w_2,z_1,z_2)= \hc(w_1 \rmult f_1 \rmult z_1)^{-1} \skm  w_1 \rmult f_1
      \rmult z_1 -\hc(w_2 \rmult f_2 \rmult z_2)^{-1}\skm w_2 \rmult f_2 \rmult z_2
      = (\alpha \skm w) \rmult h \rmult z$ and 
       $\hterm(\spol{}(f_1,f_2,t,w_1,w_2,z_1,z_2) =
        \hterm(w \rmult h\rmult z) = \hterm(w \rmult \hterm(h)\rmult z) \geq \hterm(h)$.
\\
We will now change our representation of $g$ by using the additional
 information on this situation in such a way that for the new
 representation of $g$ we either have a smaller maximal term or
 the occurrences of the term $t$
 are decreased by at least 1.
%\\
Let us assume the s-polynomial is not $\zero$\footnote{In case 
               $h =\zero$,
               just substitute the empty sum
               for the right reductive representation of $h$
               in the equations below.}. 
%\\
By our assumption, $h \red{*}{\myr}{}{F} \zero$ 
 and by Lemma \ref{lem.two-sided.red.rep} $h$ 
 has a reductive standard  representation in terms of $F$. 
Then by Lemma \ref{lem.two-sided.red.reps} the multiple
 $(\alpha \skm w) \rmult h \rmult z$ again
 has a right reductive standard representation in terms of $F$, say
  $\sum_{i=1}^n m_i \rmult h_i \rmult l_i$, 
  where $h_i \in F$, and $m_i,l_i \in \monoms(\f_{\myk})$ and all terms occurring in this sum are bounded by
  $t \succ \hterm((\alpha \skm w) \rmult h \rmult z)$.
%\\
This gives us: 
     \begin{eqnarray}
       &   & (\alpha_1 \skm w_1) \rmult f_1 \rmult z_1 + (\alpha_2 \skm w_2) \rmult f_2 \rmult z_2
             \nonumber\\  
       &   &  \nonumber\\                                                         
       & = &  (\alpha_1 \skm w_1) \rmult f_1 \rmult z_1 +
              \underbrace{ (\alpha'_2 \skm \beta_1 \skm w_1) \rmult f_1 \rmult z_1
                   - (\alpha'_2 \skm  \beta_1 \skm w_1) \rmult f_1 \rmult z_1}_{=\, 0} \nonumber\\  
       &   &  + \underbrace{(\alpha'_2\skm \beta_2 }_{= \alpha_2} \skm w_2) \rmult f_2 \rmult z_2 \nonumber\\
       &   & \nonumber\\ 
       & = & ((\alpha_1 + \alpha'_2 \skm \beta_1) \skm w_1) \rmult f_1 \rmult z_1 - \alpha'_2 \skm 
               \underbrace{((\beta_1 \skm w_1) \rmult f_1 \rmult z_1
             -  (\beta_2 \skm w_2) \rmult f_2 \rmult z_2)}_{=\,
             (\alpha \skm w) \rmult h \rmult z} 
             \nonumber\\
       &   & \nonumber\\ 
       & = & ((\alpha_1 + \alpha'_2 \skm \beta_1) \skm w_1) \rmult f_1 \rmult z_1 - \alpha'_2 \skm
                   (\sum_{i=1}^n m_i \rmult h_{i} \rmult l_{i}) \label{s2.i}
     \end{eqnarray}
     where $\beta_1=\hc(w_1 \rmult f_1 \rmult z_1)^{-1}$, $\beta_2=\hc(w_2 \rmult f_2 \rmult z_2)^{-1}$
      and  $\alpha'_2 \skm \beta_2 = \alpha_2$.
     By substituting (\ref{s2.i}) in our representation of $g$ the representation
      becomes smaller.
\\
\qed
Obviously this theorem states a criterion for when a set is a  Gr\"obner basis.
As in the case of completion procedures such as the Knuth-Bendix procedure or Buchberger's algorithm,
 elements from these test sets which do not reduce to zero can be added to the set being tested,
 to gradually describe a not necessarily finite  Gr\"obner basis.
Of course in order to get a computable completion procedure certain assumptions on the test
 sets have to be made, e.g.~they should themselves be recursively enumerable, and normal forms with respect 
 to finite sets have to be computable.
The examples from page \pageref{section.polynomials} can also be studied with
 respect to two-sided ideals.
For polynomial rings, skew-polynomial rings commutative monoid rings and commutative
 respectively poly-cyclic group rings finite Gr\"obner bases can be computed in the
 respective setting.
%%%%%%%%%%%%%%%%%%%%%%%%%%%%%%%%%%%%%%%%%%%%%%%%%%%%%%%%%%%%%%%%%%%%%%%%%%%%%%%%%%%%%%%%%%%%%%%%%
\subsection{Function Rings over Reduction Rings}\label{section.ideal.rr}

The situation becomes more complicated if $\rr$ is not a field.

Let $\rr$ be a non-commutative ring with a reduction relation $\Longrightarrow_B$
associated with  subsets $B \subseteq\rr$ as described in Section \ref{section.reductionrings}.

When following the path of linking special standard representations and reduction relations
 we get the same results as in Section \ref{section.right.rr}, i.e., 
 such
 representations naturally arise from the respective reduction relations.
Hence we proceed by studying a special reduction relation which subsumes the
 two reduction relations presented for one-sided ideals in function rings
 over reduction rings.
%Now the question arises how we can give a completion procedure for
% finitely generated  ideals in function rings.
%At least an enumeration can be achieved by using the following idea
% (of course always assuming that our function ring is computable):
%Assuming that $\rr$ is equipped with a reduction relation $\Longrightarrow$ fulfilling
% (A1) -- (A3), the question of whether a polynomial $p$ has a  reductive
% standard representation in terms of
% a set of polynomials $F$ can be answered using a reduction relation on $\f$.
As before for our ordering $>_{\rr}$ on $\rr$ we require:
 for $\alpha, \beta \in \rr$, $\alpha >_{\rr} \beta$ if and only if
 there exists a finite set $B \subseteq \rr$ such that $\alpha
 \red{+}{\Longrightarrow}{}{B} \beta$.
This ordering will ensure that the reduction relation on $\f$ is terminating.
The reduction relation on $\rr$ can be used to define various reduction relations on the function ring.
Here we want to present a reduction relation which in some sense is based on the ``divisibility''
 of the term to be reduced by the head term of the polynomial used for reduction.
\begin{definition}\label{def.two-sided.rred}~\\
{\rm
Let $f,p$ be two non-zero polynomials in $\f$.
We say $f$ \betonen{reduces} $p$ \betonen{to} $q$ \betonen{at a monomial}
 $\alpha \skm t$ \betonen{in one step}, denoted by $p \red{}{\myr}{}{f} q$, if
 there exist monomials $m,l \in \monoms(\f)$ such that
\begin{enumerate}
\item $t \in \supp(p)$ and $p(t) = \alpha$,
\item $\hterm(m \rmult \hterm(f)\rmult l) =
       \hterm(m \rmult f \rmult l) = t \geq \hterm(f)$, 
\item $\alpha \Longrightarrow_{\hc(m \rmult f \rmult l)} \beta$,
       with\footnote{Remember that by Axiom (A2) for reduction rings
       $\alpha \R_{\gamma} \beta$ implies $\alpha - \beta \in \ideal{}{}(\gamma)$
       and hence $\alpha = \sum_{i=1}^k \gamma_i \skm \gamma \skm \delta_i + \beta$,
       $\gamma_i, \delta_i \in \rr$.}  
      $\alpha = \sum_{i=1}^k \gamma_i \skm \hc(m \rmult f \rmult l)
         \skm \delta_i + \beta$
      for some $\beta, \gamma_i, \delta_i \in \rr$, $1 \leq i \leq k$, and
\item $q = p - \sum_{i=1}^k \gamma_i \skm m \rmult f \rmult l \skm \delta_i$.
\end{enumerate}
We write $p \red{}{\myr}{}{f}$ if there is a polynomial $q$ as defined
above and $p$ is then called  reducible by $f$. 
%\\
Further, we can define $\red{*}{\myr}{}{}, \red{+}{\myr}{}{}$ and
 $\red{n}{\myr}{}{}$ as usual.
%\\
Reduction by a set $F \subseteq \f \backslash \{ \zero \}$ is denoted by
 $p \red{}{\myr}{}{F} q$ and abbreviates $p \red{}{\myr}{}{f} q$
 for some $f \in F$.
\dend
}
\end{definition}
By specializing item 3.~ of this definition to
 $$3.~\alpha \Longrightarrow_{\hc(m \rmult f \rmult l)} \mbox{ such that }
 \alpha = \hc(m \rmult f \rmult l)$$
 we get an analogon to Definition \ref{def.rred_rr}.
Similarly, specializing 3.~to
 $$3.~\alpha \Longrightarrow_{\hc(m \rmult f \rmult l)} \beta \mbox{ such that }
 \hc(m \rmult f \rmult l) + \beta$$
 gives us an analogon to Definition \ref{def.rred.rr}.
 
Reviewing Example \ref{exa.reduction.not.terminating}
 we find that the reduction relation is not terminating
 when using infinite sets of polynomials for reduction.
But for finite sets we get the following analogon of Lemma
 \ref{lem.sred.rr}.
\begin{lemma}~\\ 
{\sl 
Let $F$ be a finite set of polynomials in $\f \backslash \{ \zero\}$.
\begin{enumerate}
\item For $p,q \in \f$, $p \red{}{\myr}{}{F} q$ implies $p \succ q$, in particular $\hterm(p)
  \succeq \hterm(q)$.
\item $\red{}{\myr}{}{F}$ is Noetherian.
\lemend
\end{enumerate}
}
\end{lemma}
\Ba{}
\begin{enumerate}
\item Assuming that the reduction step takes place at a monomial $\alpha \skm t$,
      by Definition \ref{def.two-sided.rred} we know
      $\hm(p - \sum_{i=1}^k \gamma_i \skm m_1 \rmult f \rmult m_2 \skm \delta_i) = \beta \skm t$ 
      which yields
      $p \succ p  - \sum_{i=1}^k \gamma_i \skm m_1 \rmult f \rmult m_2 \skm \delta_i$
      since $\alpha >_{\rr} \beta$.
\item This follows from 1. and Axiom (A1) as long as only finite sets of polynomials are involved.
\end{enumerate}\renewcommand{\baselinestretch}{1}\small\normalsize
\qed
As for the one-sided case a Translation Lemma does not hold for this
 reduction relation.
Hence we have to distinguish between weak Gr\"obner bases and Gr\"obner
 bases.
\begin{definition}~\\
{\rm
A set $F \subseteq \f \backslash \{ \zero \}$ is called a weak
 Gr\"obner basis of $\ideal{}{}(F)$ if for all
 $g \in \ideal{}{}(F)$ we have $g \red{*}{\myr}{}{F} \zero$. 
\dend
}
\end{definition}
Now as for one-sided weak Gr\"obner bases, weak Gr\"obner bases
 allow special representations of the polynomials in the ideal
 they generate.
\begin{corollary}\label{cor.two-sided.representation.rr}~\\
{\sl
Let $F$ be a set of polynomials in $\f$ 
 and $g$ a non-zero polynomial in $\ideal{}{}(F)$
 such that $g \red{*}{\myr}{}{F} \zero$.
Then $g$ has a representation of the form 
$$g = \sum_{i=1}^n m_i \rmult f_i \rmult l_i,
 f_i \in F, m_i, l_i \in \monoms(\f), n \in \n$$
 such that
 $\hterm(g) = \hterm(m_i \rmult \hterm(f_i) \rmult l_i) =
 \hterm(m_i \rmult f_i \rmult l_i) \geq \hterm(f_i)$ for $1 \leq i \leq k$, and
 $\hterm(g) \succ \hterm(m_i \rmult f_i \rmult l_i)$ 
 for all $k+1 \leq i \leq n$.
}
\end{corollary}
\Ba{}~\\
We show our claim by induction on $n$ where $g \red{n}{\myr}{}{F} \zero$.
If $n=0$ we are done.
Else let $g \red{1}{\myr}{}{F} g_1 \red{n}{\myr}{}{F} \zero$.
In case the reduction step takes place at the head monomial,
 there exists a polynomial $f \in F$
 and  monomial $m,l \in \monoms(\f)$
 such that $\hterm(m \rmult \hterm(f) \rmult l) = \hterm(m \rmult f \rmult l) = \hterm(g)
 \geq \hterm(f)$ and $\hc(g) \R_{\hc(m \rmult f \rmult l)} \beta$ with
 $\hc(g) \R_{\hc(m \rmult f \rmult l)} \beta$ with
 $\hc(g) = \beta + \sum_{i=1}^k \gamma_i \skm \hc(m \rmult f\rmult l) \skm \delta_i$
 for some $\gamma_i, \delta_i \in \rr$, $1 \leq i \leq k$.
Moreover the induction hypothesis then is applied to  $g_1 = 
 g - \sum_{i=1}^k \gamma_i \skm m \rmult f \rmult l \skm \delta_i$.
If the reduction step takes place at a monomial with term smaller $\hterm(g)$
 for the respective monomial multiple $m \rmult f \rmult l$ we immediately get
 $\hterm(g) \succ \hterm(m \rmult f \rmult l)$ and we can apply our induction hypothesis
 to the resulting polynomial $g_1$.
In both cases we can arrange the monomial multiples $m \rmult f \rmult l$ arising from
 the reduction steps in such a way that gives us th desired representation.
\\
\qed 
As in Theorem \ref{theo.two-sided.gb_rr} we can characterize weak Gr\"obner bases using g- and m-polynomials instead of s-polynomials.
\begin{definition}\label{def.two-sided.gpol}~\\
{\rm
Let $P = \{ p_1, \ldots, p_k \}$ be a multiset of   (not necessarily different)
 polynomials in $\f$ and
 $t$ an element in $\myt$
 such that there are 
 $u_1, \ldots, u_k, v_1, \ldots, v_k \in \myt$ with
 $\hterm(u_i \rmult p_i \rmult v_i) = \hterm(u_i \rmult \hterm(p_i) \rmult v_i) = t$, for all $1 \leq i \leq k$.
Further let $\gamma_i = \hc(u_i \rmult p_i \rmult v_i)$ for  $1 \leq i \leq k$.
\\
Let $G$ be a (weak)
  Gr\"obner basis of $\{ \gamma_1, \ldots, \gamma_k \}$ with respect to $\R$ 
 in $\rr$ and
$$\alpha = \sum_{i=1}^k \sum_{j=1}^{n_i} \beta_{i,j} \skm \gamma_i \skm \delta_{i,j}$$
 for $\alpha \in G$, $\beta_{i,j}, \delta_{i,j} \in \rr$,  $1 \leq i \leq k$, and $1 \leq j \leq n_i$.
Then we define the  \betonen{g-polynomials (Gr\"obner polynomials)}
 corresponding to $p_1, \ldots, p_k$ and $t$ by setting
$$ g_{\alpha} = \sum_{i=1}^k  \sum_{j=1}^{n_i}  \beta_{i,j} \skm u_i \rmult p_i \rmult  v_i \skm  \delta_{i,j}.$$
Notice that $\hm(g_{\alpha})= \alpha \skm t$.
\\
We define the 
 \betonen{m-polynomials (module polynomials)}
 corresponding to $P$ and $t$ as those
$$h =  \sum_{i=1}^k \sum_{j = 1}^{n_i}  \beta_{i,j} \skm  u_i \rmult p_i \rmult v_i  \skm  \delta_{i,j}$$
where $\sum_{i=1}^k \sum_{j = 1}^{n_i}  \beta_{i,j} \skm \gamma_i  \skm  \delta_{i,j}=0$.
Notice that $\hterm(h) \prec t$.
\dend
}
\end{definition}
Notice that while we allow the multiplication of two terms to have influence on
 the coefficients of the result\footnote{Skew-polynomial rings are a classical example,
 see Section \ref{section.skewpolynomials}.} we require that $t \skm \alpha = \alpha \skm t$. 

Given a set of polynomials $F$, the set of corresponding
 g- and m-polynomials is again defined for all possible multisets of
 polynomials in $F$ and appropriate terms $t$ as specified by
 Definition \ref{def.two-sided.gpol}.
Notice that given a finite set of polynomials the corresponding sets of
 g- and m-polynomials in general will be infinite.

We can use g- and m-polynomials to characterize special bases
 in function rings over reduction rings satisfying Axiom (A4)
 in case we add an additional condition as before.
\begin{theorem}\label{theo.cp.i}~\\
{\sl
Let $F$ be a finite set of polynomials in $\f \backslash \{ \zero \}$
 where the reduction ring satisfies (A4).
Then $F$ is a weak Gr\"obner basis if and only if
\begin{enumerate}
\item for all $f$ in $F$ and for all $m,l$ in $\monoms(\f)$ we have
       $m \rmult f \rmult l \red{*}{\myr}{}{F} \zero$, and
\item all g- and all m-polynomials corresponding to $F$ as specified in
       Definition \ref{def.two-sided.gpol}
       reduce to zero using $F$.
\end{enumerate}
\theoend
}
\end{theorem}
\Ba{}~\\
In case $F$ is a weak Gr\"obner basis,
 since the multiples $m \rmult f \rmult l$ as well as the respective
 g- and m-polynomials are all elements of $\ideal{}{}(F)$ they must
 reduce to zero using $F$.
\\
The converse will be proven by showing that every element in
 $\ideal{}{}(F)$ is reducible by $F$.
Then as $g \in \ideal{}{}(F)$ and $g \red{}{\myr}{}{F} g'$ implies
 $g' \in \ideal{}{}(F)$ we have $g\red{*}{\myr}{}{F} \zero$.
Notice that this only holds in case the reduction relation $\red{}{\myr}{}{F}$
 is Noetherian.
This follows as by our assumption $F$ is finite.
\\
Let $g \in \ideal{}{}(F)$ have a representation in terms of $F$ of the following form:
 $g = \sum_{j=1}^n m_j \rmult f_{j} \rmult l_{j}$,
 $f_j \in F$ and $m_j,l_{j} \in \monoms(\f)$.
Depending on this  representation of $g$ and the
 well-founded total ordering $\succeq$ on $\myt$ we define
 $t = \max_{\succeq} \{ \hterm(m_j \rmult f_{j} \rmult l_{j}) \mid 1\leq j \leq n \}$ and
 $K$ as the number of polynomials $m_j \rmult f_j \rmult l_j$ with head term $t$.
We show our claim by induction on $(t,K)$, where
 $(t',K')<(t,K)$ if and only if $t' \prec t$ or $(t'=t$ and
 $K'<K)$.
\\
Since by our first assumption every multiple $m_j \rmult f_j \rmult l_j$ in this
 sum reduces to zero using $F$ and hence
 has a representation as described in Corollary \ref{cor.two-sided.representation.rr} we can assume that
 $\hterm(m_j \rmult \hterm(f_j) \rmult l_j) = \hterm(m_j \rmult f_j \rmult l_j) \geq
 \hterm(f_j)$ holds.
Without loss of generality we can assume that the polynomial multiples
 with head term $t$ are just $m_1 \rmult f_1  \rmult l_1, \ldots , m_K \rmult f_K \rmult l_K$.
\\
Obviously, $t \succeq \hterm(g) = \hterm(m_1 \rmult \hterm(f_1) \rmult l_1) \geq \hterm(f_1)$ must hold. 
If $K=1$ this gives us $t = \hterm(g)$ and even
 $\hm(g) = \hm(m_1 \rmult f_1 \rmult l_1)$, implying that $g$ is
 reducible at $\hm(g)$ by $f_1$.
\\
Hence let us assume $K>1$.
\\
First let $\sum_{j=1}^K \hm(m_j \rmult f_j \rmult l_j) = \zero$.
Then there
 is a m-polynomial $h$, corresponding to the
 polynomials $f_1, \ldots, f_K$ and the term $t$ such that 
 $\sum_{j=1}^K l_j \rmult f_j \rmult m_j = h$.
%\\
We will now change our representation of $g$ by using the additional
 information on this m-polynomial in such a way that for the new
 representation of $g$ we  have a smaller maximal term.
%\\
Let us assume the m-polynomial is not $\zero$\footnote{In case 
               $h =\zero$,
               just substitute the empty sum
               for the  reductive representation of $h$
               in the equations below.}. 
%\\
By our assumption, $h$ is reducible to
 zero using $F$ and hence has a representation with respect to $F$ as described in Corollary \ref{cor.two-sided.representation.rr}, say
  $\sum_{i=1}^n  \tilde{m}_i \rmult \tilde{f}_i \rmult \tilde{l}_i$, 
  where $\tilde{f}_i \in F$, $\tilde{m}_i, \tilde{l}_i \in \monoms(\f)$ 
  and  all terms occurring in the sum are bounded by
  $t \succ \hterm(h)$.
Hence replacing the sum $\sum_{j=1}^K m_j \rmult f_j \rmult l_j $ by
 $\sum_{i=1}^n \tilde{m}_i \rmult \tilde{f}_i \rmult \tilde{l}_i$
 gives us a smaller representation of $g$.
\\
Hence let us assume
 $\sum_{j=1}^K \hm(m_j \rmult f_j \rmult l_j) \neq 0$.
Then
 we have $\hterm(m_1 \rmult f_1 \rmult l_1 + \ldots + m_K \rmult f_K \rmult l_K) = t = \hterm(g)$,
 $\hc(g) = \hc(m_1 \rmult f_1 \rmult l_1 + \ldots + m_K \rmult f_K \rmult l_K) \in \ideal{r}{}(\{ \hc(m_1 \rmult f_1 \rmult l_1), \ldots, 
 \hc(m_K \rmult f_K \rmult l_K) \})$ and 
 even $\hm(m_1 \rmult f_1 \rmult l_1 + \ldots + m_K \rmult f_K \rmult l_K) = \hm(g)$.
Hence $\hc(g)$ is $\R$-reducible  by $\alpha$, $\alpha \in G$, $G$<a
 (weak) right Gr\"obner basis of $\ideal{r}{}(\{ \hc(m_1 \rmult f_1 \rmult l_1), \ldots, 
 \hc(m_K \rmult f_K \rmult l_K) \})$ in $\rr$ with respect to the reduction relation $\R$.
Let $g_{\alpha}$ be the respective g-polynomial corresponding to $\alpha$.
Then we know that $g_{\alpha}\red{*}{\myr}{}{F} \zero$.
Moreover, we know that the head monomial of $g_{\alpha}$ is reducible
 by some polynomial $f \in F$ and we assume
 $\hterm(g_{\alpha}) = \hterm(m \rmult \hterm(f) \rmult l) = \hterm(m \rmult f \rmult l) \geq \hterm(f)$
 and $\hc(g_{\alpha}) \R_{\hc(m \rmult f \rmult l)}$.
Then, as $\hc(g)$ is $\R$-reducible by $\hc(g_{\alpha})$, 
 $\hc(g_{\alpha})$ is $\R$-reducible to zero and (A4) holds,
 the head monomial of $g$ is also reducible by some $f' \in F$ and we are done. 
\\
\qed
Of course this theorem is still true for infinite $F$ if we can show
 that for the respective function ring the reduction relation is terminating.

Now the question arises when the critical situations in this characterization
 can be localized to subsets of the respective sets.
Reviewing the Proof of Theorem \ref{theo.two-sided.s-pol.2} we find that 
 Lemma \ref{lem.two-sided.red.reps} is central as it describes when multiples of
 polynomials which have a reductive standard representation in terms
 of some set $F$ again have such a representation.
As before, this does not hold for function rings over reduction
 rings in general.
We have stated that it is not natural to link right reduction as defined
 in Definition \ref{def.two-sided.rred} to special standard representations.
Hence, to give localizations of Theorem \ref{theo.cp.i} another
 property for $\f$ is sufficient:
\begin{definition}\label{def.two-sided.stable.loc}~\\
{\rm
A set $C \subset S \subseteq \f$ is called a \betonen{stable localization} of
 $S$ if for every $g \in S$ there exists $f \in C$
 such that $g \red{}{\myr}{}{f} \zero$.
% and if $\hm(f)$ is reducible by $F$, so is $\hm(g)$.
\dend
}
\end{definition}
In case $\f$ and $\red{}{\myr}{}{}$ allow such stable localizations,
 we can rephrase Theorem \ref{theo.cp.i} as follows:
\begin{theorem}\label{theo.two-sided.loc}~\\
{\sl
Let $F$ be a finite set of polynomials in $\f \backslash \{ \zero \}$
 where the reduction ring satisfies (A4).
Then $F$ is a weak Gr\"obner basis of $\ideal{}{}(F)$ if and only if
\begin{enumerate}
\item for all $s$ in a stable localization of 
       $\{ m \rmult f \rmult l \mid f \in \f, m,l \in \monoms(\f) \}$ we have 
       $s \red{*}{\myr}{}{F} \zero$, and
\item for all $h$ in a stable localization of the g- and m-polynomials corresponding to $F$ as specified in
       Definition \ref{def.two-sided.gpol} we have 
       $h \red{*}{\myr}{}{F} \zero$.
\end{enumerate}
\theoend
}
\end{theorem} 
We have stated that for arbitrary reduction relations in $\f$ it is not natural to
 link them to special standard representations.
Still, when proving Theorem \ref{theo.two-sided.loc}, we will find that in order to change the
 representation of an arbitrary ideal element,
 Definition \ref{def.two-sided.stable.loc} is not enough to ensure reducibility.
However, we can substitute the critical situation using an analogon of
 Lemma \ref{lem.two-sided.red.reps}, which while not related to reducibility in this case will
 still be sufficient to make the representation smaller.
\begin{lemma}\label{lem.two-sided.red.reps_rr}~\\
{\sl
Let  $F \subseteq \f \backslash \{ \zero\}$ and
 $f$, $p$  non-zero polynomials in $\f$.
If $p \red{}{\myr}{}{f} \zero$ and $f \red{*}{\myr}{}{F} \zero$,
 then $p$ has a standard representation of the form
$$p = \sum_{i=1}^n m_i \rmult f_i \rmult l_i,
 f_i \in F, m_i,l_i \in \monoms(\f), n \in \n$$
 such that
 $\hterm(p) = \hterm(m_i \rmult \hterm(f_i) \rmult l_i) =
 \hterm(m_i \rmult f_i \rmult l_i) \geq \hterm(f_i)$ for $1 \leq i \leq k$ and
 $\hterm(p) \succ \hterm(m_i \rmult f_i \rmult l_i)$ 
 for all $k+1 \leq i \leq n$.
\lemend
}
\end{lemma}
\Ba{}~\\
If $p \red{}{\myr}{}{f} \zero$ then $p = \sum_{j=1}^s \gamma_j \skm m' \rmult f \rmult l' \skm \delta_j$ with $m',l' \in \monoms(\f)$, $\gamma_j, \delta_j \in \rr$,
 and $\hterm(p) = \hterm(m \rmult \hterm(f) \rmult l) =
 \hterm(m \rmult f \rmult l) \geq \hterm(f)$.
Similarly $f \red{*}{\myr}{}{F} \zero$ implies\footnote{Notice that in this
 representation we write the products of the form 
 $\gamma \skm m$ respectively $l \skm \delta$ arising in the reduction steps
 as simple monomials.}
 $f = \sum_{i=1}^n m_i \rmult f_i \rmult l_i,
 f_i \in F, m_i, l_i \in \monoms(\f), n \in \n$
 such that
 $\hterm(f) = \hterm(m_i \rmult \hterm(f_i) \rmult l_i) =
 \hterm(m_i \rmult f_i \rmult l_i) \geq \hterm(f_i)$, $1 \leq i \leq k$, and
 $\hterm(f) \succ \hterm(m_i \rmult f_i \rmult l_i)$ 
 for all $k+1 \leq i \leq n$. 
\\
We show our claim for all multiples with $\gamma_j \skm m'$ and $l' \skm \delta_j$
 for $1 \leq j \leq s$.
Let $m = \gamma_j \rmult m'$ and $l = l' \skm \delta_j$ and
let us analyze $m \rmult m_i \rmult f_i \rmult l_i \rmult l$ with $\hterm(m_i \rmult f_i \rmult l_i) = \hterm(f)$,
 $1 \leq i \leq k$.
Let $\terms(m_i \rmult f_i \rmult l_i) = \{ s_1^i, \ldots, s_{w_i}^i \}$
 with $s_1^i \succ s_j^i$, $2 \leq j \leq w_i$, i.e.~$s_1^i = \hterm(m_i \rmult f_i \rmult l_i)= \hterm(p)$.
Hence $m \rmult \hterm(p) \rmult l = m \rmult s_1^i \rmult l \geq \hterm(p) = s_1^i$
 and as $s_1^i \succ s_j^i$, $2 \leq j \leq w_i$, by Definition
 \ref{def.two-sided.refined.ordering} we can conclude
 $\hterm(m \rmult \hterm(p) \rmult l) = \hterm(m \rmult s_1^i \rmult l)
 \succ m \rmult s_j^i \rmult l \succeq \hterm(m \rmult s_j^i \rmult l)$ for
 $2 \leq j \leq w_i$.
This implies $\hterm(m \rmult \hterm(m_i \rmult f_i \rmult l_i) \rmult l) = 
 \hterm(m \rmult m_i \rmult f_i \rmult l_i \rmult l)$
Hence we get
\begin{eqnarray*}
 &  & \hterm(m \rmult f \rmult l) \\
                    & = & \hterm(m \rmult \hterm(f) \rmult l) \\
                    & = & \hterm(m \rmult \hterm(m_i \rmult f_i \rmult l_i) \rmult l), \mbox{ as } 
                          \hterm(p) = \hterm(m_i \rmult f_i \rmult l_i) \\
                    & = & \hterm(m \rmult m_i \rmult f_i \rmult l_i \rmult l)
\end{eqnarray*}
and
since $\hterm(m \rmult f \rmult l) \geq \hterm(f) \geq \hterm(f_i)$ we can conclude
 $\hterm(m \rmult m_i \rmult f_i \rmult l_i \rmult l) \geq \hterm(f_i)$.
It remains to show that $m \rmult (m_i \rmult f_i \rmult l_i) \rmult l = (m \rmult m_i) \rmult f_i \rmult (l_i \rmult l)$ has
  representations of the desired form in terms of $F$.
First we show that
 $\hterm((m \rmult m_i \rmult \hterm(f_i) \rmult l_i \rmult l) \geq \hterm(f_i)$.
We know
 $m_i \rmult \hterm(f_i) \rmult l_i \succeq \hterm(m_i \rmult \hterm(f_i) \rmult l_i)
 = \hterm(m_i \rmult f_i \rmult l_i)$ and hence 
 $\hterm(m \rmult m_i \rmult \hterm(f_i) \rmult l_i \rmult l)=
 \hterm(m \rmult \hterm(m_i \rmult f_i \rmult l_i) \rmult l) =
 \hterm(m \rmult m_i \rmult f_i \rmult l_i \rmult l) \geq \hterm(f_i)$.
\\
Now in case $m \rmult m_i, l_i \rmult l \in \monoms(\f)$
 we are done as then $(m \rmult m_i) \rmult f_i \rmult (l_i \rmult l)$
 is a representation of the desired form.
\\
Hence let us assume
 $m \rmult m_i = \sum_{j=1}^{k_i} \tilde{m}^i_j$,$l_i \rmult l = \sum_{j'=1}^{k'_i} \tilde{l}^i_j$ $\tilde{m}^i_j,\tilde{l}^i_{j'} \in \monoms(\f)$.
Let $\terms(f_i) = \{ t^i_1, \ldots, t^i_{w} \}$ with $t^i_1 \succ t^i_j$, 
 $2 \leq j \leq w$,
 i.e.~$t^i_1 = \hterm(f_i)$.
As $\hterm(m_i \rmult  \hterm(f_i) \rmult l_i) \geq \hterm(f_i) \succ t_j$,
 $ 2 \leq j \leq w$, again by Definition \ref{def.two-sided.refined.ordering}
 we can conclude that
 $\hterm(m_i \rmult \hterm(f_i) \rmult l_i) \succ m_i \rmult t^i_j \rmult l_i
 \succeq \hterm(m_i \rmult t^i_j \rmult l_i)$, $2 \leq j \leq l$,
 and $m_i \rmult \hterm(f_i) \rmult l_i \succ \sum_{j=2}^w
  m_i \rmult t^i_j \rmult l_i$.
Then for each $s^i_j$, $2 \leq j \leq w_i$, there exists $t^i_{j'} \in
 \terms(f_i)$ such that $s^i_j \in \supp(m_i \rmult t^i_{j'} \rmult l_i)$.
Since $\hterm(f) \succ s^i_j$ and even
 $\hterm(f) \succ m_i \rmult t^i_{j'} \rmult l_i$ we find that
 either
 $\hterm(m \rmult f \rmult l) \succeq \hterm(m \rmult (m_i \rmult t^i_{j'} \rmult l_i) \rmult l) = \hterm((m \rmult m_i) \rmult t^i_{j'} \rmult (l_i \rmult l))$
 in case $\hterm(m_i \rmult t^i_{j'} \rmult l_i) = \hterm(m_i \rmult f_1 \rmult l_i)$
 or 
 $\hterm(m \rmult f \rmult l) \succ m \rmult (m_i \rmult t^i_{j'} \rmult l_i) \rmult l = (m \rmult m_i) \rmult t^i_{j'} \rmult (l_i \rmult l)$.
Hence we can conclude $\tilde{m}^i_j \rmult f_i \rmult \tilde{l}^i_{j'} \predeq \hterm(m \rmult f \rmult l)$, $1 \leq j \leq k_i$, $1 \leq j' \leq K_i$ and
 for at least one $\tilde{m}^i_j$, $\tilde{l}^i_{j'}$ we get
 $\hterm(\tilde{m}^i_{j} \rmult f_i \rmult \tilde{l}^i_{j'}) = \hterm(m \rmult m_i \rmult f_i \rmult l_i \rmult l) \geq \hterm(f_i)$.
\\
It remains to analyze the situation for the functions
 $(\sum_{i=k+1}^n m \rmult (m_i \rmult f_i \rmult l_i) \rmult l$.
Again we find that for all terms $s$ in the $m_i \rmult f_i \rmult l_i$, $k+1 \leq i \leq n$,
 we have $\hterm(f) \succeq s$ and we get 
 $\hterm(m \rmult f \rmult l) \succeq \hterm(m \rmult s \rmult l)$.
Hence all polynomial multiples of the $f_i$  in the representation 
 $\sum_{i=k+1}^n (\sum_{j=1}^{k_i} \tilde{m}^i_j) \rmult  f_i \rmult (\sum_{j=1}^{K_i} \tilde{l}^i_{j'})$, where
 $m \rmult m_i = \sum_{j=1}^{k_i}\tilde{m}^i_j$, $l_i \rmult l = \sum_{j=1}^{K_i} \tilde{l}^i_{j'}$, are bounded by 
 $\hterm(m \rmult f \rmult l)$.
\\
\qed

\Ba{Theorem \ref{theo.two-sided.loc}}~\\
The proof is basically the same as for Theorem \ref{theo.cp.i}.
Due to Lemma \ref{lem.two-sided.red.reps_rr} we can substitute the multiples $m_j \rmult f_j \rmult l_j$
 by appropriate representations without changing $(t,K)$.
Hence, we only have to ensure that despite testing less polynomials we are able to
 apply our induction hypothesis.
Taking the notations from the proof of Theorem \ref{theo.cp.i},
let us first check the situation for m-polynomials.
\\
Let $\sum_{j=1}^K \hm(m_j \rmult f_j \rmult l_j) = \zero$.
Then by Definition \ref{def.two-sided.gpol} there exists
 a module polynomial $h =\sum_{j=1}^K m_j \rmult f_j \rmult l_j$
 and by our assumption there is a polynomial $h'$
 in the stable localization such that $h \red{}{\myr}{}{h'} \zero$.
Moreover, $h' \red{*}{\myr}{}{F} \zero$.
Then by Lemma \ref{lem.two-sided.red.reps_rr} the m-polynomial $h$ has a
 standard representations
 bounded by $t$.
Hence we can change the representation of $g$ by substituting $h$ by its
 representation giving us a smaller representation
 and by our induction hypothesis $g$ is reducible by $F$ and we are done.
\\
It remains to study the case where $\sum_{j=1}^K \hm(m_j \rmult f_j \rmult l_j) \neq 0$.
Then
 we have $\hterm(\sum_{j=1}^K m_j  \rmult f_j \rmult l_j) = t = \hterm(g)$,
 $\hc(g) = \hc(\sum_{j=1}^K m_j  \rmult f_j \rmult l_j) \in \ideal{}{}(\{ \hc(m_1 \rmult f_1 \rmult l_1), \ldots, 
 \hc(m_K \rmult f_K \rmult l_K) \})$ and 
 even $\hm(\sum_{j=1}^K m_j  \rmult f_j \rmult l_j) = \hm(g)$.
Hence $\hc(g)$ is $\R$-reducible  by $\alpha$, $\alpha \in G$, $G$ a
 (weak) Gr\"obner basis of $\ideal{}{}(\{ \hc(m_1 \rmult f_1 \rmult l_1), \ldots, 
 \hc(m_K \rmult f_K \rmult l_K) \})$ in $\rr$ with respect to the reduction relation $\R$.
Let $g_{\alpha}$ be the respective g-polynomial corresponding to $\alpha$.
Then we know that $g_{\alpha}\red{}{\myr}{}{g_{\alpha}'} \zero$ for some
 $g_{\alpha}'$ in the stable localization and $g_{\alpha}'\red{*}{\myr}{}{F} \zero$.
Moreover, we know that the head monomial of $g_{\alpha}'$ is reducible
 by some polynomial $f \in F$ and we assume
 $\hterm(g_{\alpha}) = \hterm(m \rmult \hterm(f) \rmult l) = \hterm(m \rmult f \rmult l) \geq \hterm(f)$
 and $\hc(g_{\alpha}) \R_{\hc(m \rmult f \rmult l)}$.
Then, as $\hc(g)$ is $\R$-reducible by $\hc(g_{\alpha})$, 
 $\hc(g_{\alpha})$ is $\R$-reducible by $\hc(g_{\alpha}')$, 
 $\hc(g_{\alpha}')$ is $\R$-reducible to zero and (A4) holds,
 the head monomial of $g$ is also reducible by some $f' \in F$ and we are done. 
\\
\qed
Again, if for infinite $F$ we can assure that the reduction relation is Noetherian,
 the proof still holds.
%%%%%%%%%%%%%%%%%%%%%%%%%%%%%%%%%%%%%%%%%%%%%%%%%%%%%%%%%
\subsection{Function Rings over the Integers}\label{section.ideal.integers}

In the previous section we have seen that for the reduction relation for
 $\f$  based on the abstract
 notion of the reduction relation $\R_{\rr}$ there is not
 enough information on the reduction step involving the coefficient and hence we cannot
 prove an analogon of the Translation Lemma.

As in the case of studying one-sided ideals,
 when studying special reduction rings where we have more information
 on the specific reduction relation $\R_{\rr}$ the situation often can
 be improved.
Again we go into the details for the case that
 $\rr$ is the ring of the integers $\z$.
The reduction relation presented in Definition \ref{def.two-sided.rred}
 then can be reformulated for this special case as follows:

\begin{definition}\label{def.two-sided.rred.z}~\\
{\rm
Let $p$, $f$ be two non-zero polynomials in $\f_{\z}$. 
We say $f$ \betonen{reduces} $p$ \betonen{to} $q$ \betonen{at} $\alpha \skm t$
 \betonen{in one step},
 i.e. $p \red{}{\myr}{}{f} q$, if there exist $ u,v \in \terms(\f_{\z})$ such that
\begin{enumerate}
\item $t \in \supp(p)$ and $p(t) = \alpha$,
\item $\hterm(u \rmult \hterm(f)\rmult v) = \hterm(u \rmult f \rmult v) = t \geq \hterm(f)$, 
\item $\alpha \geq_{\z} \hc(u \rmult f \rmult v) > 0$ and 
         $\alpha \R_{\hc(u \rmult f \rmult v)} \delta$ where
         $\alpha = \hc(u \rmult f \rmult v) \skm \beta + \delta$ with 
         $\beta, \delta  \in \z$, 
         $0 \leq \delta < \hc(u \rmult f \rmult v)$, and
\item $q = p -  u \rmult f \rmult v \skm \beta$.
\end{enumerate}
We write $p \red{}{\myr}{}{f}$ if there is a polynomial $q$ as defined
above and $p$ is then called reducible by $f$. 
%\\
Further, we can define $\red{*}{\myr}{}{}, \red{+}{\myr}{}{}$ and
 $\red{n}{\myr}{}{}$ as usual.
%\\
Reduction by a set $F \subseteq \f \backslash \{ \zero \}$ is denoted by
 $p \red{}{\myr}{}{F} q$ and abbreviates $p \red{}{\myr}{}{f} q$
 for some $f \in F$.
\dend
}
\end{definition}
As before, for this reduction relation we can still have $t \in \supp(q)$.
The important part in showing termination now is that if we still have $t \in \supp(q)$
 then its coefficient will be smaller according to our ordering chosen
 for $\z$ (compare Section \ref{section.right.integers}) and since this ordering is well-founded we are done. 
Due to the additional information on the coefficents, again we do not have to restrict
 ourselves to finite sets of polynomials in order to ensure termination.
\begin{corollary}~\\ 
{\sl 
Let $F$ be a set of polynomials in $\f_{\z} \backslash \{ \zero\}$.
\begin{enumerate}
\item For $p,q \in \f_{\z}$, $p \red{}{\myr}{}{F} q$ implies $p \succ q$, in particular $\hterm(p)
  \succeq \hterm(q)$.
\item $\red{}{\myr}{}{F}$ is Noetherian.
\lemend
\end{enumerate}
}
\end{corollary}
Similarly, the additional information we have on the coefficients before
 and after the reduction step now enables us to prove
 an analogon of the Translation Lemma for function rings over the integers.
The first and second part of the lemma are only needed to prove the essential
 third part.

\begin{lemma}\label{lem.two-sided.trans.integers}~\\
{\sl
Let $F$ be a set of polynomials in $\f_{\z}$ and $p,q,h$ polynomials in $\f_{\z}$.
\begin{enumerate}
\item
Let $p-q \red{}{\myr}{}{F} h$ such that the reduction step takes place
 at the monomial $\alpha \skm t$ and we additionally have $t \not\in \supp(h)$.
Then there exist $p', q' \in \f_{\z}$ such that $p \red{*}{\myr}{}{F} p'$
 and $q \red{*}{\myr}{}{F} q'$ and $h = p'-q'$.
\item 
Let $\zero$ be the unique normal form of $p$ with respect to $F$ and $t = \hterm(p)$.
Then there exists a polynomial $f \in F$ such that $p \red{}{\myr}{}{f} p'$
 and $t \not \in \supp(p')$.
\item
Let $\zero$ be the unique normal form of $p-q$ with respect to $F$.
Then there exists $g \in \f_{\z}$ such that $p \red{*}{\myr}{}{F} g$ and
 $q \red{*}{\myr}{}{F} g$.
\end{enumerate}
\lemend}
\end{lemma}
\Ba{}~\\
\begin{enumerate}
\item
Let $p-q \red{}{\myr}{}{F} h$ at the monomial $\alpha \skm t$, i.e.,
 $h = p-q - u \rmult f \rmult v \skm \beta$ for some $u,v \in \terms(\f_{\z}),
 \beta \in \z$ such that
 $\hterm(u \rmult \hterm(f)\rmult v) = \hterm(u \rmult f \rmult v) = t \geq \hterm(f)$
 and $\hc(u \rmult f \rmult v) > 0$.
Remember that $\alpha$ is the coefficient
 of $t$ in $p-q$.
Then as $t \not\in \supp(h)$ we know $\alpha = \hc(u \rmult f \rmult v)\skm \beta$.
Let $\alpha_1$ respectively $\alpha_2$ be the coefficients of $t$ in
 $p$ respectively $q$ and $\alpha_1 = (\hc(u \rmult f \rmult v)\skm \beta) \skm \beta_1 + \gamma_1$
 respectively $\alpha_2 = (\hc(u \rmult f \rmult v)\skm \beta) \skm \beta_2 + \gamma_2$ for some
 $\beta_1, \beta_2, \gamma_1, \gamma_2 \in \z$ where $0 \leq \gamma_1, \gamma_2
 < \hc(u \rmult f \rmult v)\skm \beta$.
Then $\alpha = \hc(u \rmult f \rmult v)\skm \beta = \alpha_1 - \alpha_2 =
 (\hc(u \rmult f \rmult v)\skm \beta) \skm (\beta_1 - \beta_2) + (\gamma_1 - \gamma_2)$, and
 as $\gamma_1 - \gamma_2$ is no multiple of $\hc(u \rmult f \rmult v)\skm \beta$ we have
 $\gamma_1 - \gamma_2 = 0$ and hence $\beta_1 - \beta_2 = 1$.
We have to distinguish two cases:
\begin{enumerate}
\item $\beta_1 \neq 0$ and $\beta_2 \neq 0$:
      Then $p \red{}{\myr}{}{F} p -  u \rmult f \rmult v \skm \beta \skm \beta_1 = p'$,
        $q \red{}{\myr}{}{F} q -  u \rmult f \rmult v \skm \beta \skm \beta_2 = q'$ and
        $p'-q' = p -  u \rmult f \rmult v \skm \beta \skm \beta_1 -
         q + u \rmult f \rmult v \skm \beta \skm \beta_2=
        p-q - u \rmult f \rmult v \skm \beta \skm \beta = h$.
\item $\beta_1 = 0$ and $\beta_2 = - 1$ (the case $\beta_2 = 0$
       and $\beta_1 = 1$ being symmetric):
      Then $p' = p$, 
        $q \red{}{\myr}{}{F} q - u \rmult f \rmult v \skm \beta \skm \beta_2 = 
         q + u \rmult f \rmult v \skm \beta = q'$ and
        $p'-q' = p-q - u \rmult f \rmult v \skm \beta = h$.
\end{enumerate}
\item Since $p \red{*}{\myr}{}{F} \zero$, $\hm(p) = \alpha \skm t$
 must be $F$-reducible.
Let $f_i \in F$, $i \in I$ be a series of all not necessarily different polynomials in $F$
 such that $\alpha \skm t$ is reducible by them involving terms $u_i,v_i$.
Then $\hc(u_i \rmult f_i \rmult v_i) > 0$.
Moreover, let $\gamma = \min_{\leq} \{ \hc(u_i \rmult f_i \rmult v_i) \mid i \in I \}$
 and without loss of generality $\hm(u \rmult f \rmult v) = \gamma \skm t$ for some
 $f \in F$, $\hterm(u \rmult \hterm(f) \rmult v) = \hterm(u \rmult f \rmult v) \geq \hterm(f)$.
We claim that for $p \red{}{\myr}{}{f} p - \beta \skm u \rmult f \rmult v = p'$ where $\alpha =
 \beta \skm \gamma + \delta$, $\beta, \delta \in \z$, $0 \leq \delta < \gamma$,
 we have $t \not \in \supp(p')$.
Suppose $\hterm (p') = t$.
Then by our definition of reduction we must have $0 < \hc(p') < \hc(u \rmult f \rmult v)$.
But then $p'$ would no longer be $F$-reducible contradicting our assumption
 that $\zero$ is the unique normal form of $p$.
\item Since $\zero$ is the unique normal form of $p-q$ by 2.~there exists
 a reduction sequence $p-q \red{}{\myr}{}{f_{i_1}} h_1 \red{}{\myr}{}{f_{i_2}} h_2 \red{}{\myr}{}{f_{i_3}}
 \ldots \red{}{\myr}{}{f_{i_k}} \zero$ such that $\hterm(p-q) \succ \hterm(h_1)
 \succ \hterm(h_2) \succ \ldots$.
We show our claim by induction on $k$, where $p-q \red{k}{\myr}{}{F} \zero$
 is such a reduction sequence.
In the base case $k=0$ there is nothing to show as then $p=q$.
Hence, let $p-q \red{}{\myr}{}{F} h \red{k}{\myr}{}{F} \zero$.
Then by 1.~there are polynomials $p',q' \in \f_{\z}$ such that 
 $p \red{*}{\myr}{}{F} p'$
 and $q \red{*}{\myr}{}{F} q'$ and $h = p'-q'$.
Now the induction hypothesis for $p'-q' \red{k}{\myr}{}{F} \zero$ yields the
 existence of a polynomial $g \in \f_{\z}$ such that
 $p \red{*}{\myr}{}{F} g$ and $q \red{*}{\myr}{}{F} g$.
\end{enumerate}
\qed 

Hence weak Gr\"obner bases are in fact Gr\"obner bases and
 can hence be characterized as
 follows (compare Definition \ref{def.weak.gb}):

\begin{definition}~\\
{\rm
A set $F \subseteq \f_{\z} \backslash \{ \zero \}$ is called a (weak)
 Gr\"obner basis of $\ideal{}{}(F)$ if for all
 $g \in \ideal{}{}(F)$ we have $g \red{*}{\myr}{}{F} \zero$. 
\dend
}
\end{definition}
\begin{corollary}\label{cor.two-sided.standard.rep}~\\
{\sl
Let $F$ be a set of polynomials in $\f_{\z}$ 
 and $g$ a non-zero polynomial in $\ideal{}{}(F)$
 such that $g \red{*}{\myr}{}{F} \zero$.
Then $g$ has a representation of the form 
$$g = \sum_{i=1}^n m_i \rmult f_i \rmult l_i,
 f_i \in F, m_i, l_i \in \monoms(\f_{\z}), n \in \n$$
 such that
 $\hterm(g) = \hterm(m_i \rmult \hterm(f_i) \rmult l_i) =
 \hterm(m_i \rmult f_i \rmult l_i) \geq \hterm(f_i)$, $1 \leq i \leq k$, and
 $\hterm(g) \succ \hterm(m_i \rmult f_i \rmult l_i) = 
 \hterm(m_i \rmult \hterm(f_i) \rmult l_i)$ 
 for all $k+1 \leq i \leq n$.
\\
In case $\zero$ is the unique normal form of $g$ with respect to $F$ we even can find a
 representation where additionally $\hterm(m_1 \rmult f_1 \rmult l_1) \succ
 \hterm(m_2 \rmult f_2 \rmult l_2) \succ \ldots \succ \hterm(m_n \rmult
 f_n \rmult l_n)$.
}
\end{corollary}
\Ba{}~\\
We show our claim by induction on $n$ where $g \red{n}{\myr}{}{F} \zero$.
If $n=0$ we are done.
Else let $g \red{1}{\myr}{}{F} g_1 \red{n}{\myr}{}{F} \zero$.
In case the reduction step takes place at the head monomial,
 there exists a polynomial $f \in F$
 and $u, v \in \terms(\f_{\z}), \beta \in \z$
 such that $\hterm(u \rmult\hterm(f) \rmult v) = \hterm(u \rmult f \rmult v) = \hterm(g)
 \geq \hterm(f)$ and $\hc(g) \R_{\hc(u \rmult f \rmult v)} \delta$ with
 $\hc(g) = \hc(u \rmult f\rmult v) \skm \beta + \delta$ for some $\beta, \delta \in \z$,
 $0 \leq \delta < \hc(u \rmult f \rmult v)$.
Moreover the induction hypothesis then is applied to  $g_1 = 
 g - u \rmult f \rmult v \skm \beta$.
If the reduction step takes place at a monomial with term smaller $\hterm(g)$
 for the respective monomial multiple $u \rmult f \rmult v \skm \beta$ we immediately get
 $\hterm(g) \succ u \rmult f \rmult v \skm \beta$ and we can apply our induction hypothesis
 to the resulting polynomial $g_1$.
In both cases we can arrange the monomial multiples $u \rmult f \rmult v \skm \beta$ arising from
 the reduction steps in such a way that gives us the desired representation.
\\
\qed 
Now Gr\"obner bases
 can be characterized using the concept of 
 s-polynomials combined with the technique of saturation
 which is neccessary in order to describe the
 whole ideal congruence by the reduction relation.
\begin{definition}\label{def.two-sided.s-poly.z}~\\
{\rm
Let $p_{1}, p_{2}$ be polynomials in $\f_{\z}$.
If there are respective terms
 $t,u_1, u_2, v_1, v_2 \in \myt$ such that
 $\hterm(u_i \rmult \hterm(p_i) \rmult v_i) = \hterm(u_i \rmult p_i \rmult v_i)=t
 \geq \hterm(p_i)$ let $HC(u_i \rmult p_{i} \rmult v_i) = \gamma_i$.\\
Assuming $\gamma_1 \geq \gamma_2 > 0$\footnote{Notice that $\gamma_i > 0$
 can always be achieved by studying the situation for $- p_i$ in case 
 we have $HC(u_i \rmult p_{i} \rmult v_i) < 0$.},
 there are $\beta, \delta \in \z$ such that
 $\gamma_1 = \gamma_2 \skm \beta + \delta$ and $0 \leq \delta < \gamma_2$
 and we get the following s-polynomial
      $$\spol{}(p_{1}, p_{2},t,u_1, u_2,v_1, v_2) =  u_2 \rmult p_2 \rmult v_2 \skm \beta -  u_1 \rmult p_1 \rmult v_1.$$
The set $\spols(\{p_1,p_2\})$ then is the set of all such
 s-polynomials corresponding to $p_1$ and $p_2$.
\dend
}
\end{definition}
Again these sets in general are not finite.
\begin{theorem}\label{theo.two-sided.gb.reduction.z}~\\
{\sl
Let $F$ be a  set of polynomials in $\f_{\z} \backslash \{ \zero \}$.
Then $F$ is a Gr\"obner basis if and only if
\begin{enumerate}
\item for all $f$ in $F$ and for all $m,l$ in $\monoms(\f_{\z})$ we have 
       $m \rmult f \rmult l \red{*}{\myr}{}{F} \zero$, and
\item all s-polynomials corresponding to $F$ as specified in
       Definition \ref{def.two-sided.s-poly.z}
       reduce to $\zero$ using $F$.
\end{enumerate}
\theoend
}
\end{theorem}
\Ba{}~\\
In case $F$ is a Gr\"obner basis,
 since these polynomials are all elements of $\ideal{}{}(F)$ they must
 reduce to zero using $F$.
\\
The converse will be proven by showing that every element in
 $\ideal{}{}(F)$ is reducible  by $F$.
Then as $g \in \ideal{}{}(F)$ and $g \red{}{\myr}{}{F} g'$ implies
 $g' \in \ideal{}{}(F)$ we have $g\red{*}{\myr}{}{F} \zero$.
Notice that this is sufficient as the reduction relation $\red{}{\myr}{}{F}$
 is Noetherian.
\\
Let $g \in \ideal{}{}(F)$ have a representation in terms of $F$ of
 the following form:
$g = \sum_{j=1}^n v_j \rmult f_{j} \rmult w_{j} \skm \alpha_j$ such that
 $f_j \in F$, $v_j,w_{j} \in \myt$ and $\alpha_j \in \z$.
Depending on this  representation of $g$ and the
 well-founded total ordering $\succeq$ on $\myt$ we define
 $t = \max_{\succeq} \{ \hterm(v_j \rmult f_{j} \rmult w_{j}) \mid 1\leq j \leq m \}$, 
 $K$  as the number of polynomials $f_j \rmult w_j$ with head term $t$, and
 $M = \{\{ \hc(v_j \rmult f_j \rmult w_j) \mid \hterm(v_j \rmult f_j \rmult w_j) = t \}\}$ a multiset
 in $\z$.
We show our claim by induction on $(t,M)$, where
 $(t',M')<(t,M)$ if and only if $t' \prec t$ or $(t'=t$ and $M' \ll M)$.
\\
Since by our first assumption every multiple $v_j \rmult f_j \rmult w_j$ in this
 sum reduces to zero using $F$ and hence
 has a  representation as specified
 in Corollary \ref{cor.two-sided.standard.rep}, we can assume that
 $\hterm(v_j \rmult \hterm(f_j) \rmult w_j) = \hterm(v_j \rmult f_j \rmult w_j) \geq
 \hterm(f_j)$ holds.
Moreover, without loss of generality we can assume that the polynomial multiples
 with head term $t$ are just $v_1 \rmult f_1  \rmult w_1, \ldots , v_K \rmult f_K \rmult w_K$ and
 additionally we can assume $\hc(v_j \rmult f_j \rmult w_j) > 0$\footnote{This can easily be achieved
 by adding $-f$ to $F$ for all $f \in F$ and using $v_j \rmult (-f_j) \rmult w_j \skm (- \alpha_j)$
 in case $\hc(v_j \rmult f_j \rmult w_j) < 0$.}.
\\
Obviously, $t \succeq \hterm(g)$ must hold. 
If $K = 1$ this gives us $t = \hterm(g)$ and even
 $\hm(g) = \hm(v_1 \rmult f_1 \rmult w_1 \skm \alpha_1)$, implying that $g$ is right
 reducible at $\hm(g)$ by $f_1$.
\\
Hence let us assume $K>1$.
\\
Without loss of generality we can assume that $\hc(v_1 \rmult f_1 \rmult w_1) \geq
 \hc(v_2 \rmult f_2 \rmult w_2) > 0$ and there are $\alpha, \beta \in \z$ such that
 $\hc(v_2 \rmult f_2 \rmult w_2) \skm \alpha + \beta = \hc(v_1 \rmult f_1 \rmult w_1)$ and
 $\hc(v_2 \rmult f_2 \rmult w_2) > \beta \geq 0$.
Since $t = \hterm(v_1 \rmult f_1 \rmult w_1) =  \hterm(v_2 \rmult f_2 \rmult w_2)$ by Definition
 \ref{def.two-sided.s-poly.z} we have an s-polynomial $\spol{}(f_1,f_2,t,v_1, v_2,w_1,w_2)
 =v_2 \rmult f_2 \rmult w_2 \skm \alpha - v_1 \rmult f_1 \rmult w_1$.
If $\spol{}(f_1,f_2,t,v_1, v_2,w_1,w_2) \neq \zero$\footnote{In case $\spol{}(f_1,f_2,t,v_1, v_2,w_1,w_2) = \zero$
 the proof is similar. We just have to subsitute $\zero$ in the equations below which
 immediately gives us a smaller representation of $g$.}
 then $\spol{}(f_1,f_2,t,v_1,v_2,w_1,w_2)\red{*}{\myr}{}{F} \zero$ implies
 $\spol{}(f_1,f_2,t,v_1,v_2,w_1,w_2) = \sum_{i=1}^k m_i \rmult h_i \rmult l_i$,
 $h_i \in F$, $m_i, l_i \in \monoms(\f_{\z})$ where this
 sum is a representation in the sense of Corollary \ref{cor.two-sided.standard.rep}
 with terms bounded by $\hterm(\spol{}(f_1,f_2,t,v_1,v_2,w_1,w_2)) \leq t$.
This gives us
\begin{eqnarray}
 &  & v_1 \rmult f_1 \rmult w_1 \skm \alpha_1 + v_2 \rmult f_2 \rmult w_2 \skm \alpha_2  \\ \nonumber
  & = & v_1 \rmult f_1 \rmult w_1 \skm \alpha_1 +
         \underbrace{v_2 \rmult f_2 \rmult w_2 \skm \alpha_1 \skm \alpha - v_2 \rmult f_2 \rmult w_2 \skm \alpha_1 \skm \alpha}_{=\zero} +
        v_2 \rmult f_2 \rmult w_2 \skm \alpha_2  \\ \nonumber
  & = &  v_2 \rmult f_2 \rmult w_2 \skm (\alpha_1 \skm \alpha + \alpha_2)
          -  \underbrace{(v_2 \rmult f_2 \rmult  w_2 \skm \alpha - v_1 \rmult f_1 \rmult w_1)}_{=\spol{}(f_1,f_2,t,v_1,v_2,w_1,w_2)}\skm \alpha_1  \\ \nonumber
  & = & v_2 \rmult f_2 \rmult w_2 \skm (\alpha_1 \skm \alpha + \alpha_2) - (\sum_{i=1}^k m_i \rmult  h_i \rmult l_i) \skm \alpha_1 \\ \nonumber
\end{eqnarray}
and substituting this in the representation of $g$ we get a new representation with
 $t' = \max_{\succeq} \{ \hterm(v_j \rmult f_{j} \rmult w_{j}),
  \hterm(m_j \rmult h_{j} \rmult l_{j}) \mid f_j, h_j \mbox{ appearing in the new representation } \}$, and
 $M' = \{\{ \hc(v_j \rmult f_j \rmult w_j), \hc(m_j \rmult h_{j} \rmult l_{j}) \mid \hterm(v_j \rmult f_j \rmult w_j) = \hterm(m_j \rmult h_{j} \rmult l_{j}) = t' \}\}$
 and either $t' \pred t$ and we have a smaller representation for $g$
 or in case $t'=t$ we have to distinguish two cases:
 \begin{enumerate}
 \item $\alpha_1 \skm \alpha  + \alpha_2 = 0$. \\
       Then $M' = (M - \{\{ \hc(v_1 \rmult f_1 \rmult w_1), \hc(v_2 \rmult f_2 \rmult w_2) \}\}) \cup 
     \{ \{ \hc(m_j \rmult h_j \rmult l_j) \mid \hterm(m_j \rmult h_j \rmult l_j) = t \}\}$.
       As those polynomials $h_j$ with $\hterm(m_j \rmult h_j \rmult l_j) = t$ are used to reduce the
       monomial $\beta \skm t = \hm(\spol{}(f_1,f_2,t,v_1,v_2,w_1,w_2))$ we know that for them we have
       $0 < \hc(m_j \rmult h_j \rmult l_j) \leq \beta < \hc(v_2 \rmult f_2 \rmult w_2) \leq 
     \hc(v_1 \rmult f_1 \rmult w_1)$ and hence $M' \ll M$ and we have a smaller representation for $g$.
 \item $\alpha_1 \skm \alpha  + \alpha_2\neq 0$. \\
       Then $M' = (M - \{\{ \hc(v_1 \rmult f_1 \rmult w_1) \}\}) \cup 
     \{ \{ \hc(m_j \rmult h_j \rmult l_j) \mid \hterm(m_j \rmult h_j \rmult l_j) = t \}\}$.
       Again  $M' \ll M$ and we have a smaller representation for $g$.
 \end{enumerate}
Notice that the case $t'=t$ and $M'\ll M$ cannot occur infinitely often but has to result
 in either $t' <t$ or will lead to $t'=t$ and $K=1$ and hence to reducibility by $\red{}{\myr}{}{F}$.
\\
\qed
Now the question arises when the critical situations in this characterization
 can be localized to subsets of the respective sets
 as in Theorem \ref{theo.two-sided.s-pol.2}.
Reviewing the Proof of Theorem \ref{theo.two-sided.s-pol.2} we find that 
 Lemma \ref{lem.two-sided.red.reps} is central as it describes when multiples of
 polynomials which have a reductive standard representation in terms
 of some set $F$ again have such a representation.
As we have seen before, this will not hold for function rings over reduction
 rings in general.
As in Section \ref{section.ideal.rr}, to give localizations of Theorem \ref{theo.two-sided.gb.reduction.z} the
 concept of stable subsets is sufficient:
\begin{definition}\label{def.two-sided.stable.loc.z}~\\
{\rm
A set $C \subset S \subseteq \f_{\z}$ is called a \betonen{stable localization} of
 $S$ if for every $g \in S$ there exists $f \in C$
 such that $g \red{}{\myr}{}{f} \zero$.
% and if $\hm(f)$ is reducible by $F$, so is $\hm(g)$.
\dend
}
\end{definition}
In case $\f_{\z}$ and $\red{}{\myr}{}{}$ allow such stable localizations,
 we can rephrase Theorem \ref{theo.two-sided.gb.reduction.z} as follows:
\begin{theorem}\label{theo.two-sided.loc.z}~\\
{\sl
Let $F$ be a set of polynomials in $\f_{\z} \backslash \{ \zero \}$.
Then $F$ is a Gr\"obner basis of $\ideal{}{}(F)$ if and only if
\begin{enumerate}
\item for all $s$ in a stable localization of 
       $\{ m \rmult f \rmult l \mid f \in \f_{\z}, m,l \in \monoms(\f_{\z}) \}$ we have 
       $s \red{*}{\myr}{}{F} \zero$, and
\item for all $h$ in a stable localization of the s-polynomials corresponding to $F$ as specified in
       Definition \ref{def.two-sided.s-poly.z} we have 
       $h \red{*}{\myr}{}{F} \zero$.
\end{enumerate}
\theoend
}
\end{theorem} 
When proving Theorem \ref{theo.two-sided.loc.z}, we can substitute the critical situation using an analogon of
 Lemma \ref{lem.two-sided.red.reps}, which  will be sufficient to make the representation used in the proof smaller.
It is a direct consequence of Lemma \ref{lem.two-sided.red.reps_rr}.
\begin{corollary}\label{cor.two-sided.red.reps.z}~\\
{\sl
Let  $F \subseteq \f_{\z} \backslash \{ \zero\}$ and
 $f$, $p$  non-zero polynomials in $\f_{\z}$.
If $p \red{}{\myr}{}{f} \zero$ and $f \red{*}{\myr}{}{F} \zero$,
 then $p$ has a representation of the form
$$p = \sum_{i=1}^n m_i \rmult f_i \rmult l_i,
 f_i \in F, m_i,l_i \in \monoms(\f_{\z}), n \in \n$$
 such that
 $\hterm(p) = \hterm(m_i \rmult \hterm(f_i) \rmult l_i) =
 \hterm(m_i \rmult f_i \rmult l_i) \geq \hterm(f_i)$ for $1 \leq i \leq k$ and
 $\hterm(p) \succ \hterm(m_i \rmult f_i \rmult l_i)$ 
 for all $k+1 \leq i \leq n$.
\corend
}
\end{corollary}

\Ba{Theorem \ref{theo.two-sided.loc.z}}~\\
The proof is basically the same as for Theorem \ref{theo.two-sided.gb.reduction.z}.
Due to Corollary \ref{cor.two-sided.red.reps.z} we can substitute the multiples $v_j \rmult f_j \rmult w_j$
 by appropriate representations.
Hence, we only have to ensure that despite testing less polynomials we are able to
 apply our induction hypothesis.
Taking the notations from the proof of Theorem \ref{theo.two-sided.gb.reduction.z},
let us  check the situation for $K>1$.
\\
Without loss of generality we can assume that $\hc(v_1 \rmult f_1 \rmult w_1) \geq
 \hc(v_2 \rmult f_2 \rmult w_2) > 0$ and there are $\alpha, \beta \in \z$ such that
 $\hc(v_2 \rmult f_2 \rmult w_2) \skm \alpha + \beta = \hc(v_1 \rmult f_1 \rmult w_1)$ and
 $\hc(v_2 \rmult f_2 \rmult w_2) > \beta \geq 0$.
Since $t = \hterm(v_1 \rmult f_1 \rmult w_1) =  \hterm(v_2 \rmult f_2 \rmult w_2)$ by Definition
 \ref{def.two-sided.s-poly.z} we have an s-polynomial $h$ in the stable localization of
 $\spols(f_1,f_2)$ such that
 $v_2 \rmult f_2 \rmult w_2 \skm \alpha - v_1 \rmult f_1 \rmult w_1 \red{}{\myr}{}{h} \zero$.
If $h \neq \zero$\footnote{In case $h = \zero$
 the proof is similar.
 We just have to subsitute $\zero$ in the equations below which
 immediately gives us a smaller representation of $g$.}
 then by Corollary \ref{cor.two-sided.red.reps.z}
 $h\red{*}{\myr}{}{F} \zero$ implies
 $v_2 \rmult f_2 \rmult w_2 \skm \alpha - v_1 \rmult f_1 \rmult w_1 = \sum_{i=1}^k m_i \rmult h_i \rmult l_i$,
 $h_i \in F$, $m_i,l_i \in \monoms(\f_{\z})$ where this
 sum is a representation in the sense of Corollary \ref{cor.two-sided.standard.rep}
 with terms bounded by $\hterm(m \rmult h \rmult l) \leq t$.
This gives us
\begin{eqnarray}
 &  & v_1 \rmult f_1 \rmult w_1 \skm \alpha_1 + v_2 \rmult f_2 \rmult w_2 \skm \alpha_2 \\ \nonumber
  & = & v_1 \rmult f_1 \rmult w_1 \skm \alpha_1 + \underbrace{v_2 \rmult f_2 \rmult w_2 \skm \alpha_1 \skm \alpha - v_2 \rmult f_2 \rmult w_2 \skm \alpha_1 \skm \alpha}_{=\zero} +
        v_2 \rmult f_2 \rmult w_2 \skm \alpha_2 \\ \nonumber
  & = & v_2 \rmult f_2 \rmult w_2 \skm (\alpha_1 \skm \alpha + \alpha_2) -  
       (v_2 \rmult f_2 \rmult  w_2 \skm \alpha - v_1 \rmult f_1 \rmult w_1)
              \skm \alpha_1 \\ \nonumber
  & = & v_2 \rmult f_2 \rmult w_2 \skm (\alpha_1 \skm \alpha + \alpha_2) - (\sum_{i=1}^k m_i \rmult h_i \rmult l_i) \skm \alpha_1 \\ \nonumber
\end{eqnarray}
and substituting this in the representation of $g$ we get a new representation with
 $t' = \max_{\succeq} \{ \hterm(v_j \rmult f_{j} \rmult w_{j}),
  \hterm(m_j \rmult h_{j} \rmult l_{j}) \mid f_j, h_j \mbox{ appearing in the new representation } \}$, and
 $M' = \{\{ \hc(v_j \rmult f_j \rmult w_j), \hc(m_j \rmult h_{j} \rmult l_{j}) \mid \hterm(v_j \rmult f_j \rmult w_j) = \hterm(m_j \rmult h_{j} \rmult l_{j}) = t' \}\}$
 and  either $t' \pred t$ 
 or  $(t'=t$ and $M' \ll M)$ and in both cases we have a smaller representation for $g$.
Notice that the case $t'=t$ and $M'\ll M$ cannot occur infinitely often but has to result
 in either $t' <t$ or will lead to $t'=t$ and $K=1$ and hence to reducibility by $\red{}{\myr}{}{F}$.
\\
\qed
%%%%%%%%%%%%%%%%%%%%%%%%%%%%%%%%%%%%%%%%%%%%%%%%%%%%%%%%%%%%%%%%%%%%%%%
\section{Two-sided Modules}\label{section.two-sided.module}
Given a function ring $\f$ with unit ${\bf 1}$ and a natural number $k$, 
 let $\f^k = \{ (f_1, \ldots, f_k) \mid f_i \in \f \}$
 be the set of all vectors of length $k$ with coordinates in $\f$.
Obviously $\f^k$ is an additive commutative group with respect
 to ordinary vector addition.
Moreover,
$\f^k$ is such an \betonen{$\f$-module} with respect to
 the scalar multiplication $f \rmult (f_1, \ldots, f_k) = 
 (f \rmult f_1, \ldots, f \rmult f_k)$ and $ (f_1, \ldots, f_k) \rmult f = 
(f_1 \rmult f , \ldots, f_k \rmult f)$.
Additionally $\f^k$ is called \betonen{free} as it has a basis\footnote{Here the term {\em basis} is
 used in the meaning of being a linearly independent set of generating vectors.}.
One such basis is the set of unit vectors ${\bf e}_1 = ({\bf 1}, \zero, \ldots, \zero), {\bf e}_2 = (\zero, {\bf 1}, \zero, \ldots, \zero),
 \ldots, {\bf e}_k = (\zero, \ldots, \zero, {\bf 1})$.
Using this basis the elements of $\f^k$ can be written uniquely as ${\bf f} = \sum_{i=1}^k f_i \rmult {\bf e}_i$
 where ${\bf f} =  (f_1, \ldots, f_k)$.

\begin{definition}~\\
{\rm
A subset of $\f^k$ which is again an $\f$-module is called a \betonen{submodule} of $\f^k$.
}
\end{definition}
As before any ideal of $\f$ is an $\f$-module and even a submodule of the $\f$-module $\f^1$.
Provided a set of vectors $S = \{ {\bf f}_1, \ldots, {\bf f}_s \}$ the set
 $\{ \sum_{i=1}^s\sum_{j=1}^{n_i}  g_{ij} \rmult {\bf f}_i \rmult {h_{ij}} \mid g_{ij}, {h_{ij}} \in \f \}$ is a submodule of $\f^k$.
This set is denoted as $\langle S \rangle$ and $S$ is called a generating set.

\begin{theorem}~\\
{\sl
Let $\f$ be Noetherian.
Then every submodule of $\f^k$ is finitely generated.
}
\end{theorem}
\Ba{}~\\
Let ${\cal S}$ be a submodule of $\f^k$.
Again we show our claim by induction on $k$.
For $k=1$ we find that ${\cal S}$ is in fact an ideal in $\f$ and hence by our hypothesis finitely
 generated.
For $k>1$ let us look at the set $I = \{ f_1 \mid (f_1, \ldots, f_k) \in {\cal S} \}$.
Then again $I$ is an ideal in $\f$ and hence finitely generated.
Let $\{g_1, \ldots, g_s \mid g_i \in \f \}$ be a generating set of $I$.
Choose ${\bf g}_1, \ldots, {\bf g}_s \in {\cal S}$ such that the first coordinate of ${\bf g}_i$ is $g_i$.
Note that the set $\{ (f_2, \ldots, f_k) \mid (\zero,f_2, \ldots, f_k) \in  {\cal S}\}$ is a submodule
 of $\f^{k-1}$ and hence finitely generated by some set $\{ (n_2^i, \ldots, n_k^i), 1 \leq i \leq w \}$.
Then the set $\{{\bf g}_1, \ldots, {\bf g}_s \} \cup \{ {\bf n}_i = (\zero, n_2^i, \ldots, n_k^i) \mid 
 1 \leq i \leq w \}$ is a generating set for ${\cal S}$.
To see this assume ${\bf m} = (m_1, \ldots, m_k)\in {\cal S}$.   
Then $m_1 = \sum_{i=1}^s \sum_{j=1}^{n_i} h_{ij} \rmult g_i \rmult {h_{ij}}'$ for some $h_{ij}, {h_{ij}}' \in \f$ and
 ${\bf m'} = {\bf m} - \sum_{i=1}^s \sum_{j=1}^{n_i} h_{ij} \rmult {\bf g}_i \rmult {h_{ij}}' \in {\cal S}$ with first coordinate $\zero$.
Hence ${\bf m'} = \sum_{i=1}^w \sum_{j=1}^{m_i} l_{ij} \rmult {\bf n}_i \rmult {l_{ij}}'$ for some $l_{ij}, {l_{ij}}' \in \f$ giving rise to
$${\bf m} = {\bf m'} + \sum_{i=1}^s \sum_{j=1}^{n_i} h_{ij} \rmult {\bf g}_i \rmult {h_{ij}}'  = 
\sum_{i=1}^w \sum_{j=1}^{m_i} l_{ij} \rmult {\bf n}_i \rmult {l_{ij}}' + 
\sum_{i=1}^s \sum_{j=1}^{n_i} h_{ij} \rmult {\bf g}_i \rmult {h_{ij}}' .$$
\qed
$\f^k$ is called Noetherian if and only if all its submodules are finitely generated.

If $\f$ is a reduction ring Section \ref{section.two-sided.module} outlines how the existence of 
 Gr\"obner bases for submodules can be shown.

Now given a submodule ${\cal S}$ of $\f^k$, we can define 
 $\f^k/{\cal S} = \{ {\bf f} + {\cal S} \mid {\bf f} \in \f^k \}$.
Then with addition defined as $({\bf f} + {\cal S}) + ({\bf g} + {\cal S}) = ({\bf f}+ {\bf g}) + {\cal S}$
 the set $\f^k/{\cal S}$ is an abelian group and can be turned into
 an $\f$-module by the action $g \rmult ({\bf f} + {\cal S}) \rmult h = g \rmult {\bf f} \rmult h + {\cal S}$
 for $g,h \in \f$.
$\f^k/{\cal S}$ is called the \betonen{quotient module} of $\f^k$ by ${\cal S}$.

As usual this quotient can be related to homomorphisms.
The results carry over from commutative module theory as can be found in \cite{AdLo94}.
Recall that for two $\f$-modules ${\cal M}$ and ${\cal N}$, a function
 $\phi: {\cal M} \myr{\cal N}$ is an $\f$-module homomorphism if
 $$\phi({\bf f } + {\bf g}) = \phi({\bf f}) + \phi({\bf g}) \mbox{ for all } {\bf f,g} \in {\cal M}$$
 and
 $$\phi(g \rmult {\bf f} \rmult h) = g \rmult \phi ({\bf f}) \rmult h \mbox{ for all } {\bf f} \in {\cal M}, g,h \in \f.$$ 
The homomorphism is called an \betonen{isomorphism} if $\phi$ is one to one and we write
 ${\cal M}\cong {\cal N}$.
Let ${\cal S} = {\rm ker}(\phi) = \{ {\bf f} \in {\cal M} \mid \phi({\bf f}) = {\bf 0} \}$.
Then ${\cal S}$ is a submodule of ${\cal M}$ and $\phi({\cal M})$ is a submodule
 of ${\cal N}$.
Since all are abelian groups we know ${\cal M}/{\cal S} \cong \phi({\cal M})$
 under the mapping  ${\cal M}/{\cal S} \myr \phi({\cal M})$ with
 ${\bf f} + {\cal S} \mapsto \phi({\bf f})$ which is in fact an isomorphism.
All submodules of the quotient ${\cal M}/{\cal S}$ are of the form 
 ${\cal L}/{\cal S}$ where ${\cal L}$ is a submodule of ${\cal M}$
 containing ${\cal S}$.

Unfortunately, contrary to the one-sided case we can no longer show that 
 every finitely generated $\f$-module ${\cal M}$ is isomorphic to some
 quotient of $\f^k$.
Let ${\cal M}$ be a  finitely generated $\f$-module with
 generating set  ${\bf f}_1, \ldots {\bf f}_k \in {\cal M}$.
Consider the mapping $\phi: \f^k \myr {\cal M}$ defined by 
 $\phi(g_1, \ldots, g_k)=\sum_{i=1}^k g_i \rmult {\bf f}_i $ for ${\cal M}$.
The image of the
 $\f$-module homomorphis is no longer ${\cal M}$.

%%% Local Variables: 
%%% mode: latex
%%% TeX-master: "testlauf"
%%% End: 

%% file: applications.tex
In this chapter we outline how the concept of Gr\"obner bases can be used
 to describe algebraic questions and when  solutions can be achieved.
We will describe the problems in the following manner
\\
\\
\problem{Problem}
{A description of the algebraic setting of the problem.}
{A description of the problem itself.}
{&& A description of how the problem can be analyzed using Gr\"obner bases.}

In a first step we do not require finiteness or computability of the operations, especially
 of a Gr\"obner basis.
Since an ideal itself is always a Gr\"obner basis itself, the assumption ``Let G be a respective
 Gr\"obner basis'' always holds and means a Gr\"obner basis of the ideal generated by $G$.

In case a Gr\"obner basis is computable (though not necessarily finite) and the normal form
 computation for a polynomial with respect to a finite set is effective, 
 our so-called proceedings give rise to procedures
 which can then be used to treat the problem in a constructive manner.
If additionally the Gr\"obner basis computation terminates, these procedures terminate as
 well and the instance of the problem is decidable.
In case Gr\"obner basis computation always terminates for a chosen setting the whole problem
 is decidable in this setting.

Of course ``termination'' here is meant in a theoretical sense while as we know practical
 ``termination'' is already often not achievable for the Gr\"obner basis computation in the ordinary
 polynomial ring due to complexity issues although finite Gr\"obner bases always exist.

The terminology extends to one-sided ideals and we note those problems, where the one-sided case
 also makes sense.

We will also note when weak Gr\"obner bases are sufficient for the solution of a
 problem.

\section{Natural Applications}
The most obvious problem related to Gr\"obner bases is the ideal membership problem.
Characterizing Gr\"obner bases with respect to a reduction relation uses the important fact that an
 element belonging to the ideal will reduce to zero using the Gr\"obner basis.

\problem{Ideal Membership Problem}
{A set $F \subseteq \f$ and an element $f \in \f$.}
{$f \in \ideal{}{}(F)$?}
{&1.& Let $G$ be a Gr\"obner basis of $\ideal{}{}(F)$.\\
 &2.& If $f \red{*}{\myr}{}{G} \zero$, then $f \in \ideal{}{}(F)$.}

Hence Gr\"obner bases give a semi-answer to this question in case they are computable
 and the normal form computation is effective.
To give a negative answer the Gr\"obner basis computation must either terminate or
 one must explicitly prove, e.g.~using properties of the enumerated Gr\"obner basis,
 that the element will never reduce to zero.

These results carry over to one-sided ideals using the appropriate one-sided Gr\"obner bases.

Moreover, weak Gr\"obner bases are sufficient to solve the problem.

A normal form computation always gives rise to a special representation in terms of the polynomials
 used for reduction and in case the normal form is zero such representations are special standard
 representations.
We give two instances of this problem.

\problem{Representation Problem 1}
{A Gr\"obner basis  $G \subseteq \f$ and an element $f \in \ideal{}{}(G) $.}
{Give a representation of $f$ in terms of $G$.}
{&& Reducing $f$ to $\zero$ using $G$ yields such a representation.
}

In case the normal form computation is effective, we can collect the polynomials and
 multiples used in the reduction process and combine them to the desired representation.
Notice that since we know that the element is in the ideal, it is enough to additionally
 require that the Gr\"obner basis is recursively enumerable as a set.

The result carries over to one-sided ideals using the appropriate one-sided Gr\"obner bases.

Again, weak Gr\"obner bases are sufficient to solve the problem.

Often the ideal is not presented in terms of a Gr\"obner basis.
Then additional information is necessary which in the computational case is related to
 collecting the history of polynomials created during completion.
Notice that the proceedings in this case require some equivalent to
 Lemma \ref{lem.two-sided.red.reps} to hold and hence the problem is restricted 
 to function rings over fields.

\problem{Representation Problem 2}
{A set $F \subseteq \f_{\myk}$ and an element $f \in \ideal{}{}(F) $.}
{Give a representation of $f$ in terms of $F$.}
{&1.& Let $G$ be a Gr\"obner basis of $\ideal{}{}(F)$. \\
 &2.& Let $g = \sum_{i=1}^{k_g} m_i \rmult f_i \rmult \tilde{m}_i$ be representations of the elements $g \in G$ 
 in terms of $F$. \\
 &3.& Let $f = \sum_{j=1}^k n_i \rmult g_i \rmult \tilde{n}_i$ be a representation of $f$ in terms of $G$. \\
 &4.& The sums in 2. and 3. yield a representation of $f$ in terms of $F$.
}

In case the Gr\"obner basis is computable by a completion procedure the procedure has to keep track of the
 history of polynomials by storing their representations in terms of $F$.
If the completion stops we can reduce $f$ to zero and substitute the representations of the
 polynomials used by their ``history representation''.
If the Gr\"obner basis is only recursively enumerable both processes have to be interwoven and to continue until
 the normal form computation for $f$ reaches $\zero$.

The result carries over to one-sided ideals using the appropriate one-sided Gr\"obner bases.

Moreover, weak Gr\"obner bases are sufficient to solve the problem.

Other problems are related to the comparison of ideals.
For example given two ideals one can ask whether one is included in the other.
 
\problem{Ideal Inclusion Problem}
{Two sets $F_1,F_2 \subseteq \f$.}
{$\ideal{}{}(F_1) \subseteq  \ideal{}{}(F_2)$?}
{&1.& Let $G$ be a Gr\"obner basis of $\ideal{}{}(F_2)$ .\\
 &2.& If $F_1 \red{*}{\myr}{}{G} \zero$, then $\ideal{}{}(F_1) \subseteq \ideal{}{}(F_2)$. \\
}

In case the Gr\"obner basis is computable and the normal form computation is effective this
 yields a semi-decision procedure for the problem.
If additionally the Gr\"obner basis computation terminates for $F_1$
 or we can prove that some 
 element of the set $F_1$ does not belong to $\ideal{}{}(F_2)$, e.g.~by deriving knowledge from the
 enumerated Gr\"obner basis, we can also give a negative answer.

The result carries over to one-sided ideals using the appropriate one-sided Gr\"obner bases.

Weak Gr\"obner bases are sufficient to solve the problem.

Applying the inclusion problem in both directions we get a characterization for equality of ideals.

\problem{Ideal Equality Problem}
{Two sets $F_1,F_2 \subseteq \f$.}
{$\ideal{}{}(F_1) = \ideal{}{}(F_2)$?}
{&1.& Let $G_1$, $G_2$ be Gr\"obner bases of $\ideal{}{}(F_1)$ respectively $\ideal{}{}(F_2)$.\\
 &2.& If $F_1 \red{*}{\myr}{}{G_2} \zero$, then $\ideal{}{}(F_1) \subseteq \ideal{}{}(F_2)$. \\
 &3.& If $F_2 \red{*}{\myr}{}{G_1} \zero$, then $\ideal{}{}(F_2) \subseteq \ideal{}{}(F_1)$. \\
 &4.& If 2. and 3. both hold, then $\ideal{}{}(F_1) = \ideal{}{}(F_2)$.
}

Again, Gr\"obner bases at least give a semi-answer in case they are computable and the normal form
 procedure is effective.
We can confirm whether two generating sets are bases of one ideal.
Of course, in case the computed Gr\"obner bases are finite, we can also give a negative answer.
However, if the Gr\"obner bases are not finite, a negative answer is only possible, if
we can prove either $F_1 \not\subseteq \ideal{}{}(F_2)$ or $F_2 \not\subseteq \ideal{}{}(F_1)$.

The result carries over to one-sided ideals using the appropriate one-sided Gr\"obner bases.

Again, weak Gr\"obner bases are sufficient to solve the problem.

In case $\f$ contains a unit say ${\bf 1}$, we can ask whether an ideal is equal to the trivial ideal
 in $\f$ generated by the unit.

\problem{Ideal Triviality Problem 1}
{A set $F \subseteq \f$.}
{$\ideal{}{}(F) = \ideal{}{}(\{ {\bf 1} \})$?}
{&1.& Let $G$ be a respective Gr\"obner basis. \\
 &2.& If ${\bf 1} \red{*}{\myr}{}{G} \zero$, then $\ideal{}{}(F) = \ideal{}{}(\{ {\bf 1} \})$.
}

Again Gr\"obner bases give a semi-answer in case they can be computed.
If the Gr\"obner basis is additionally finite or we can prove that ${\bf 1} \not\in\ideal{}{}(F)$,
 then we can also confirm $\ideal{}{}(F) \neq \ideal{}{}(\{ {\bf 1} \})$.

Since $\ideal{}{}(\{ 1 \}) = \f$ one can also rephrase the question for rings without a unit.

\problem{Ideal Triviality Problem 2}
{A set $F \subseteq \f$.}
{$\ideal{}{}(F) = \f$?}
{&1.& Let $G$ be a Gr\"obner basis of $\ideal{}{}(F)$. \\
 &2.& If for every $t \in \myt$, $t \red{*}{\myr}{}{G} \zero$, then $\ideal{}{}(F) = \f$.
}

Of course now we have the problem that the test set $\myt$ in general will not be finite.
Hence a Gr\"obner basis can give a semi-answer in case we can restrict this test set to a finite
 subset.
If the Gr\"obner basis is additionally finite or we can prove that $t \nred{*}{\myr}{}{G} \zero$
 for some $t$ in the finite sub test set of $\myt$,
 then we can also confirm $\ideal{}{}(F) \neq \f$.

Both of these result carry over to one-sided ideals using the appropriate one-sided Gr\"obner bases.

As before, weak Gr\"obner bases are sufficient to solve the problem.

\problem{Ideal Union Problem}
{Two sets $F_1, F_2 \subseteq \f$ and an element $f \in \f$.}
{$f \in \ideal{}{}(F_1) \cup \ideal{}{}(F_2)$?}
{&1.& Let $G_1$, $G_2$ be Gr\"obner bases of $\ideal{}{}(F_1)$ respectively $\ideal{}{}(F_2)$. \\
 &2.& If $f \red{*}{\myr}{}{G_1} \zero$, then $f \in \ideal{}{}(F_1) \cup \ideal{}{}(F_2)$.\\
 &3.& If $f \red{*}{\myr}{}{G_2} \zero$, then $f \in \ideal{}{}(F_1) \cup \ideal{}{}(F_2)$.
}

Notice that $\ideal{}{}(F_1) \cup \ideal{}{}(F_2) \neq \ideal{}{}(F_1 \cup F_2)$.
Moreover $G_1 \cup G_2$
 is neither a Gr\"obner basis of $\ideal{}{}(F_1) \cup \ideal{}{}(F_2)$, which
 in general is no ideal itself,  nor of
 $\ideal{}{}(F_1 \cup F_2)$. 

Again, weak Gr\"obner bases are sufficient to solve the problem.

The ideal generated by the set $F_1 \cup F_2$ is called the sum of the two ideals.
\begin{definition}~\\
{\rm
For two ideals $\mathfrak{i}, \mathfrak{j} \subseteq \f$ the \betonen{sum} is defined as the set
$$\mathfrak{i}+\mathfrak{j} = \{ f \radd g \mid f \in \mathfrak{i}, g \in \mathfrak{j} \}.$$
}
\end{definition}
As in the case of commutative polynomials one can show the following theorem.
\begin{theorem}~\\
{\sl
For two ideals $\mathfrak{i}, \mathfrak{j} \subseteq \f$ the sum $\mathfrak{i}+\mathfrak{j}$ is again
 an ideal.
In fact, it is the smallest ideal containing both, $\mathfrak{i}$ and $\mathfrak{j}$.
If $F$ and $G$ are the respective generating sets for $\mathfrak{i}$ and $\mathfrak{j}$, then
 $F \cup G$ is a generating set for $\mathfrak{i}+\mathfrak{j}$.
}
\end{theorem}
\Ba{}
First we check that the sum is indeed an ideal:
\begin{enumerate}
\item as $\zero \radd \zero = \zero$ we get $\zero \in \mathfrak{i}+\mathfrak{j}$,
\item for $h_1, h_2 \in \mathfrak{i}+\mathfrak{j}$ we have that there are $f_1, f_2 \in \mathfrak{i}$ and
          $g_1, g_2 \in \mathfrak{j}$ such that $h_1 = f_1 \radd g_1$ and $h_2=f_2  \radd g_2$.
      Then $h_1 \radd h_2 = (f_1 \radd g_1) \radd (f_2 \radd g_2) = (f_1 \radd  f_2) \radd (g_1 \radd g_2)
           \in \mathfrak{i}+\mathfrak{j}$, and
\item for $h_1 \in \mathfrak{i}+\mathfrak{j}$, $h_2 \in \f$ we have that there are $f \in \mathfrak{i}$
       and $g \in \mathfrak{j}$ such that $h_1 = f \radd g$.
       Then $h_1 \rmult h_2 = (f \radd g) \rmult h_2 = f \rmult h_2 \radd g \rmult h_2 \in \mathfrak{i}+\mathfrak{j}$
       as well as $h_2 \rmult h_1 = h_2 \rmult (f \radd g) = h_2 \rmult f \radd h_2 \rmult g 
       \in \mathfrak{i}+\mathfrak{j}$.
\end{enumerate}
Since any ideal containing $\mathfrak{i}$ and $\mathfrak{j}$ contains $\mathfrak{i}+\mathfrak{j}$, this is
 the smallest ideal containing them.
It is easy to see that $F \cup G$ is a generating set for the sum.
\qed

Of course $F \cup G$ in general will not be a Gr\"obner basis.
This becomes immediately clear when looking at the following corollary.

\begin{corollary}~\\
{\sl
For $F \subseteq \f$ we have
$$\ideal{}{}(F) = \bigcup_{f\in F} \ideal{}{}(f).$$
}
\end{corollary}
But we have already seen that for function rings a polynomial in general is no Gr\"obner basis
 of the ideal or one-sided ideal it generates.

\problem{Ideal Sum Problem}
{Two sets $F_1, F_2 \subseteq \f$ and an element $f \in \f$.}
{$f \in \ideal{}{}(F_1) + \ideal{}{}(F_2)$?}
{&1.& Let $G$  be a  Gr\"obner basis of $\ideal{}{}(F_1 \cup F_2)$. \\
 &2.& If $f \red{*}{\myr}{}{G} \zero$, then $f \in \ideal{}{}(F_1) + \ideal{}{}(F_2)$.
}

Both of these result carry over to one-sided ideals using the appropriate one-sided Gr\"obner bases.

As before, weak Gr\"obner bases are sufficient to solve the problem.

Similar to sums for commutative function rings
 we can define products of ideals.
\begin{definition}~\\
{\rm
For two ideals $\mathfrak{i}, \mathfrak{j}$ in a commutative function
 ring $\f$ the \betonen{product}
 is defined as the set
$$\langle\mathfrak{i}\rmult \mathfrak{j}\rangle = \ideal{}{}(\{  f_i \rmult g_i \mid f_i \in \mathfrak{i}, g_i \in \mathfrak{j} \}).$$
\dend
}
\end{definition}

\begin{theorem}~\\
{\sl
For two ideals $\mathfrak{i}, \mathfrak{j}$ in a commutative function ring $\f$
 the product $\langle\mathfrak{i}\rmult\mathfrak{j}\rangle$ is again
 an ideal.
If $F$ and $G$ are the respective generating sets for $\mathfrak{i}$ and $\mathfrak{j}$, then
 $F \rmult G = \{ f \rmult g \mid  f \in F, g \in G \}$ 
 is a generating set for $\mathfrak{i}\rmult\mathfrak{j}$.
}
\end{theorem}
\Ba{}
First we check that the product is indeed an ideal:
\begin{enumerate}
\item as $\zero \in \mathfrak{i}$ and $\zero \in \mathfrak{j}$ we get $\zero \in \mathfrak{i} \rmult \mathfrak{j}$,
\item for $f,g \in \mathfrak{i}\rmult\mathfrak{j}$ we have $f \radd g  \in \mathfrak{i}\rmult \mathfrak{j}$
       by our definition, and
\item for $f \in \mathfrak{i}\rmult\mathfrak{j}$, $h \in \f$ we have that there are $f_i \in \mathfrak{i}$
       and $g_i \in \mathfrak{j}$ such that $f = \sum_{i=1}^k f_i \rmult g_i$
       and then $f \rmult h = (\sum_{i=1}^k f_i \rmult g_i) \rmult h = \sum_{i=1}^k f_i \rmult (g_i \rmult h)
       \in \mathfrak{i}\rmult\mathfrak{j}$.
\end{enumerate}
It is obvious that $\ideal{}{}(F \rmult G) \subseteq \langle\mathfrak{i}\rmult\mathfrak{j}\rangle$ as
 $F \rmult G \subseteq \mathfrak{i}\rmult\mathfrak{j}$.
On the other hand every polynomial in $\langle\mathfrak{i}\rmult\mathfrak{j}\rangle$ can be written as a sum of products
 $\tilde{f} \rmult \tilde{g}$ where $\tilde{f} = \sum_{i=1}^n h_i \rmult f_i \in \mathfrak{i}$, $f_i \in F$,
 $h_i \in \f$ and $\tilde{g} = \sum_{j=1}^m g_j \rmult \tilde{h}_j$, $g_j \in G$, $\tilde{h}_j \in \f$.
Hence every such product  $\tilde{f} \rmult \tilde{g}$ is again of the desired form.
\\
\qed

\problem{Ideal Product Problem}
{Two subsets $F_1, F_2$ of a commutative function ring $\f$ and an element $f \in \f$.}
{$f \in \langle\ideal{}{}(F_1) \rmult \ideal{}{}(F_2)\rangle$?}
{&1.& Let $G$  be a  Gr\"obner basis of $\ideal{}{}(F_1 \rmult F_2)$. \\
 &2.& If $f \red{*}{\myr}{}{G} \zero$, then $f \in \langle\ideal{}{}(F_1) \rmult \ideal{}{}(F_2)\rangle$.
}

Again, weak Gr\"obner bases are sufficient to solve the problem.

%Notice that $\mathfrak{i}\rmult\mathfrak{j}$ is still an ideal in a non-commutative function ring:
%We only have to show that right and left multiples by elements in $\f$ are again elements of the product.
%Remember that for $f \in \mathfrak{i}\rmult\mathfrak{j}$ there are $f_i \in \mathfrak{i}$
% and $g_i \in \mathfrak{j}$ such that $f = \sum_{i=1}^k f_i \rmult g_i$.
%Then for $g \in \f$ we have  
% $f \rmult g = (\sum_{i=1}^k f_i \rmult g_i) \rmult g = \sum_{i=1}^k f_i \rmult (g_i \rmult g)
%     \in \mathfrak{i}\rmult\mathfrak{j}$ since $g_i \rmult g \in \mathfrak{j}$, as well as
% $ g \rmult f = g \rmult (\sum_{i=1}^k f_i \rmult g_i) = \sum_{i=1}^k (g \rmult f_i) \rmult g_i
% \in \mathfrak{i}\rmult\mathfrak{j}$ since $g \rmult f_i \in \mathfrak{i}$.
%But the generating set of the product can no longer be described in this easy way.
%One possibility would be to add special multiples as in $\{ f \rmult g \mid f \in F \rmult \monoms(\f), g \in G \}$
% which of course is no nice characterization when hoping for a finitary solution.

We close this section by showing how Gr\"obner bases can help to detect the existence of inverse
 elements in $\f$ in case $\f$ has a unit say ${\bf 1}$.

\begin{definition}~\\
{\rm
Let $\f$ be a function ring with unit ${\bf 1}$ and $f \in \f$.
An element $g \in \f$ is called a \betonen{ right inverse} of $f$ in $\f$ if  $f \rmult g = {\bf 1}$.
Similarly $g$ is called a \betonen{left inverse} of $f$ in $\f$ if $g \rmult f = {\bf 1}$.
\dend
}
\end{definition}

\problem{Inverse Element Problem}
{An element $f \in \f$.}
{Does $f$ have a right or left inverse in $\f$?}
{&1.& Let $G_r$ be a respective right Gr\"obner basis of $\ideal{r}{}(f)$. \\
 &2. & If ${\bf 1} \red{*}{\myr}{r}{G_r} \zero$, then $f$ has a right inverse. \\
 &1'. & Let $G_{\ell}$ be a respective left Gr\"obner basis of $\ideal{\ell}{}(f)$. \\
 &2'. & If ${\bf 1} \red{*}{\myr}{r}{G_{\ell}} \zero$, then $f$ has a left inverse. \\
}

To see that this is correct we give the following argument for the right inverse case:
It is clear that $f$ has a right inverse in $\f$ if and only if 
 $\ideal{r}{}(\{ f \}) = \f$ since $f \rmult g -{\bf 1}= \zero$ for some
 $g \in \f$ if and only if
 ${\bf 1} \in \ideal{r}{}(\{ f\})$.
So, in order to decide whether $f$ has a right inverse in $\f$ one
 has to distinguish the following two cases provided we have a right Gr\"obner basis $G_r$ of $\ideal{r}{}(\{ f \})$:
If ${\bf 1} \nred{*}{\myr}{r}{G_r} \zero$ then $f$ has no right inverse.
If ${\bf 1} \red{*}{\myr}{r}{G_r} \zero$ then we know ${\bf 1} \in \ideal{r}{}(\{ f \})$,
 i.e.~there exist $h  \in \f$ such that
 ${\bf 1} =  f \rmult h$ and hence $h$ is a right inverse of $f$ in $\f$.

A symmetric argument holds for the case of left inverses.

Of course in case $\f$ is commutative, left inverses and right inverses coincide in case they exist
 and we can use the fact that $f \rmult g -{\bf 1} = g \rmult f -{\bf 1} = \zero$ if and only if ${\bf 1} \in \ideal{}{}(\{ f \})$.

Again, weak Gr\"obner bases are sufficient to solve the problem.

It is also possible to ask for the existence of left and right inverses for elements of the quotient rings described
 in the next section.

\section{Quotient Rings}
Let $F$ be a subset of $\f$ generating an ideal $\mathfrak{i} = \ideal{}{}(F)$.
The canonical homomorphism from $\f$ onto $\f/\mathfrak{i}$ is defined as
$$ f \longmapsto [f]_{\mathfrak{i}} $$
with $[f]_{\mathfrak{i}} = f + \mathfrak{i}$ denoting the congruence class of $f$ modulo $\mathfrak{i}$.
The ring operations are given by
$$ [f]_{\mathfrak{i}} + [g]_{\mathfrak{i}} = [f+g]_{\mathfrak{i}},$$
$$ [f]_{\mathfrak{i}} \mrm [g]_{\mathfrak{i}} = [ f \rmult g]_{\mathfrak{i}}.$$
A natural question now is whether two elements of $\f$ are in fact in the same congruence class modulo
 $\mathfrak{i}$.

\pagebreak
\problem{Congruence Problem}
{A set $F \subseteq \f$ and two elements $f,g \in \f$.}
{$f = g$ in $\f/\ideal{}{}(F)$?}
{&1.& Let $G$ be a Gr\"obner basis of $\ideal{}{}(F)$. \\
 &2.& If $f-g \red{*}{\myr}{}{G} \zero$, then $f = g$ in $\f/\ideal{}{}(F)$.
}

Hence if $G$ is a Gr\"obner basis for which normal form computation is effective, the congruence problem
 is solvable.

Usually one element of the congruence class is identified as its representative and since normal forms
 with respect to Gr\"obner bases are unique, they can be chosen as such representatives.

Notice that for weak Gr\"obner bases unique representations for the quotient
 can no longer be determined by reduction (review Example \ref{exa.Z}).

\problem{Unique Representatives Problem}
{A set $F \subseteq \f$ and an element $f \in \f$.}
{Determine a unique representative for $f$ in $\f/\ideal{}{}(F)$.}
{&1.& Let $G$ be a respective Gr\"obner basis. \\
 &2.& The normal form of $f$ with respect to $G$ is a unique representative.
}

Provided a Gr\"obner basis of $\mathfrak{i}$ together with an effective normal form algorithm
we can specify unique representatives by
$$ [f]_{\mathfrak{i}}:= {\rm normal\_form}(f, G),$$
and define addition and multiplication in the quotient by
$$ [f]_{\mathfrak{i}} + [g]_{\mathfrak{i}} := {\rm normal\_form}(f+g, G),$$
$$ [f]_{\mathfrak{i}} \mrm [g]_{\mathfrak{i}} := {\rm normal\_form}(f \rmult g, G).$$ 
Similar to the case of polynomial rings for a function ring over a field $\myk$ we can
 show that this structure is a $\myk$-vector space with a special basis.

\begin{lemma}~\\
{\sl
For any ideal $\mathfrak{i} \subseteq \f_{\myk}$ the following hold:
\begin{enumerate}
\item $\f_{\myk}/\mathfrak{i}$ is a $\myk$-vector space.
\item The set $B = \{ [t]_{\mathfrak{i}} \mid t \in \myt \}$ is a vector space basis
       and we can chose $[t]_{\mathfrak{i}} = {\rm monic}({\rm normal\_form}(t,G))$ for $G$ being a Gr\"obner basis
       of $\mathfrak{i}$.
\end{enumerate}
\lemend}
\end{lemma}
\Ba{}
\begin{enumerate}
\item We have to show that the following properties hold for $V = \f_{\myk}/\mathfrak{i}$:
      \begin{enumerate}
      \item There exists a mapping $K \times V \longrightarrow  V$, 
            $(\alpha , [f]_{\mathfrak{i}}) \longmapsto  \alpha \skm [f]_{\mathfrak{i}}$
            called multiplication with scalars.
      \item $(\alpha \skm \beta) \skm [f]_{\mathfrak{i}} = 
            \alpha \skm ( \beta \skm [f]_{\mathfrak{i}})$ for all
            $\alpha, \beta \in \myk$, $[f]_{\mathfrak{i}} \in V$.
      \item $\alpha \skm ([f]_{\mathfrak{i}} + [g]_{\mathfrak{i}}) = 
            \alpha \skm [f]_{\mathfrak{i}} + \alpha \skm [g]_{\mathfrak{i}}$
            for all $\alpha \in \myk$, $[f]_{\mathfrak{i}},
            [g]_{\mathfrak{i}} \in V$.
      \item $(\alpha + \beta) \skm [f]_{\mathfrak{i}} =
            \alpha \skm [f]_{\mathfrak{i}} + \beta \skm [f]_{\mathfrak{i}}$
            for all $\alpha, \beta \in \myk$, $[f]_{\mathfrak{i}} \in V$.
      \item ${\bf 1} \skm [f]_{\mathfrak{i}} = [f]_{\mathfrak{i}}$ for all
            $[f]_{\mathfrak{i}} \in V$.
      \end{enumerate}
      It is easy to show that this follows from the natural definition
       $$\alpha \skm [f]_{\mathfrak{i}} := [ \alpha \skm f]_{\mathfrak{i}}$$
      for $\alpha \in \myk$, $[f]_{\mathfrak{i}} \in V$.
\item It follows immediately that $B$ generates the quotient $\f_{\myk}/\mathfrak{i}$. 
 So it remains to show that this basis is free in the sense that $\zero$ cannot be
 represented as a non-trivial linear combination of elements in $B$.
 Let $G$ be a Gr\"obner basis of $\mathfrak{i}$.
 Then we can choose the elements of $B$ as the normal forms of the elements in $\myt$ with respect
 to $G$. 
 Since for a polynomial in normal form all its terms are also in normal form we can conclude that these
 normal forms are elements of $\monoms(\f_{\myk})$ and since $\myk$ is a field we can make them monic.
 This leaves us with a basis 
 $\{ \tilde{t} = {\rm monic}({\rm normal\_form}(t, G)) \mid t \in \myt \}$ .
 Now let us assume that $B$ is not free, i.e.~there exists $k \in \n$ minimal 
 with $\alpha_i \in \myk \backslash \{ 0 \}$ and $[t_i]_{\mathfrak{i}} \in B$,
 $1 \leq i \leq k$
 such that $\sum_{i}^k \alpha_i \skm  [t_i]_{\mathfrak{i}} = \zero$.
 Since then we also get ${\rm normal\_form}(\sum_{i}^k \alpha_i \skm  \tilde{t}_i, G) = \zero$
 and all $\tilde{t_i}$ are different and in normal form, all $\alpha_i$ must equal $0$
 contradicting our assumption.
\end{enumerate}
\qed

If we can compute normal forms for the quotient elements, we can give a multiplication table for
 the quotient in terms of the vector space basis by 
$$[t_i]_{\mathfrak{i}} \mrm [t_j]_{\mathfrak{i}} = [t_i \rmult t_j]_{\mathfrak{i}} =
 {\rm normal\_form}(t_i \mm t_j, G).$$

Notice that for a function ring over a reduction ring the set $B = \{ [t]_{\mathfrak{i}} \mid t \in \myt \}$
 also is a generating set where we can chose $[t]_{\mathfrak{i}} = {\rm normal\_form}(t,G)$.
But we can no longer choose the representatives to be a subset of $\myt$.
This is due to the fact that if a monomial $\alpha \skm t$ is reducible by some polynomial $g$ this
 does not imply that some other  monomial $\beta \skm t$ or even the term $t$ is reducible by $g$.
For example let $\rr = \z$, $\myt = \{ a, \lambda \}$ and $a \rmult a = 2 \skm a$, $\lambda \rmult \lambda = \lambda$,
 $a \rmult \lambda = \lambda \rmult a = a$.
Then $2 \skm a$ is reducible by $a$ while of course $a$ isn't.  

In case $\f_{\myk}$ contains a unit say ${\bf 1}$ we can ask whether an element of $\f_{\myk} / \mathfrak{i}$ is invertible.

\begin{definition}~\\
{\rm
Let $f \in \f_{\myk}$.
An element $g \in \f_{\myk}$ is called a \betonen{right inverse} of $f$ in 
 $\f_{\myk}/\mathfrak{i}$ if  $f \rmult g = {\bf 1} \mbox{ mod } \mathfrak{i}$.
Similarly $g$ is called a \betonen{left inverse} of $f$ in $\f_{\myk}/\mathfrak{i}$ 
 if $g \rmult f =  {\bf 1} \mbox{ mod }\mathfrak{i}$.
\dend
}
\end{definition}

In case $\f_{\myk}$ is commutative, right and left inverses coincide if they exist and we can
 tackle the problem by using the fact that $f$ has an inverse in $\mathfrak{i}$ if and only if
 $f \rmult g -{\bf 1} \in  \mathfrak{i}$ if and only if ${\bf 1} \in \mathfrak{i} + \ideal{}{}(\{ f \})$.
Hence, if we have a Gr\"obner basis $G$ of the  ideal $\mathfrak{i} + \ideal{}{}(\{ f \})$ 
 the existence of an inverse of $f$
 is equivalent to ${\bf 1} \red{*}{\myr}{}{G} \zero$.

Even,  weak Gr\"obner bases are sufficient to solve the problem.

For the non-commutative case we introduce a new non-commuting tag variable $z$ 
 by lifting the multiplication $z \rmult z = z$, $z \rmult t = zt$ and $t \rmult z = tz$ for $t \in \myt$
 and extending $\myt$ to $z\myt = \{ z^{i}t_1zt_2z\ldots zt_kz^{j} \mid k \in \n, i,j \in \{0,1\}, t_i \in \myt \}$.
The order on this enlarged set of terms is induced by combining a syllable ordering with respect to $z$
 with the original ordering on $\myt$. 
By $\f_{\myk}^{z\myt}$ we denote the function ring over $z\myt$.

This technique of using a tag variable  now allows to study the right ideal generated by $f$
 in $\f_{\myk}/\mathfrak{i}$, where $\mathfrak{i} = \ideal{}{}(F)$ for some set $F \subseteq \f_{\myk}$, by studying
 the ideal generated by $F \cup \{ z \rmult f \}$ in $\f_{\myk}^{z\myt}$ because of the following fact:

\begin{lemma}~\\
{\sl
Let $F \subseteq \f_{\myk}$ and $f \in \f_{\myk}$. Then $\ideal{}{\f_{\myk}^{z\myt}}(F \cup \{ z \rmult f \})$ has a Gr\"obner
 basis of the form $G \cup \{ z \rmult p_i \mid i \in I, p_i \in \f_{\myk} \}$ with $G \subseteq \f_{\myk}$.
In fact the set $\{ p_i \mid i \in I \}$ then is a right Gr\"obner basis of $\ideal{r}{\f_{\myk}/\mathfrak{i}}(\{f\})$.
\lemend
}
\end{lemma}
\Ba{}~\\
Let $G \subseteq \f_{\myk}$ be a Gr\"obner basis of $\ideal{}{\f_{\myk}}(F)$.
Then obviously $\ideal{}{\f_{\myk}^{z\myt}}(F \cup \{ z \rmult f \}) = \ideal{}{\f_{\myk}^{z\myt}}(G \cup \{ z \rmult f \})$.
Theorem \ref{theo.two-sided.s-pol.2}
 specifies a criterion to check whether a set is a Gr\"obner basis and gives rise to test sets for 
 a completion procedure.
Notice that due to the ordering on $z\myt$ which uses the tag variable to induce syllables, we can state
 the following important result:
\begin{quote}
If for a polynomial $q \in \f_{\myk}$ the multiple $z \rmult q$
 has a standard representation, then so has every
 multiple $u \rmult ( z \rmult q ) \rmult z \rmult v$ for $u,v \in z\myt$.
\end{quote}
Moreover, since $G$ is already a Gr\"obner basis, no critical situation for polynomials in $G$ have to
 be considered.
\\
Then a completion of $G \cup \{ z \rmult f \}$ can be obtained as follows:
\\
In a first step only three kinds of critical situations have to be considered:
 \begin{enumerate}
 \item s-polynomials of the form $zu \rmult g \rmult v -z \rmult f \rmult w$ where $u,v,w \in \myt$ such that
        $\hterm(zu \rmult g \rmult v) = \hterm(z \rmult f \rmult w)$,
 \item s-polynomials of the form $z \rmult f \rmult u - z \rmult f \rmult v$ where $u,v \in \myt$ such that
        $\hterm(z \rmult f \rmult u) = \hterm(z \rmult f \rmult v)$,  and
 \item polynomials of the form $z \rmult f \rmult u$ where $u \in \myt$ such that $\hterm(f \rmult u) \neq
       \hterm(f) \rmult u$.
 \end{enumerate}
Since normal forms of polynomials of the form $z \rmult p$, $p \in \f_{\myk}$, with respect to subsets
 of $\f_{\myk} \cup z \rmult \f_{\myk}$ are again elements of $z \rmult \f_{\myk} \cup \{ \zero \}$, we can assume that from then on we are
 completing a set $G \cup \{ z \rmult q_i \mid q_i \in \f_{\myk} \}$ and again three kinds of critical
 situations have to be considered:
 \begin{enumerate}
 \item s-polynomials of the form $zu \rmult g \rmult v -z \rmult q_i \rmult w$ where $u,v,w \in \myt$ such that
        $\hterm(zu \rmult g \rmult v) = \hterm(z \rmult q_i \rmult w)$,
 \item s-polynomials of the form $z \rmult q_i \rmult u - z \rmult q_j \rmult v$ where $u,v \in \myt$ such that
        $\hterm(z \rmult q_i \rmult u) = \hterm(z \rmult q_j \rmult v)$,  and
 \item polynomials of the form $z \rmult p_i \rmult u$ where $u \in \myt$ such that $\hterm(p_i \rmult u) \neq
       \hterm(p_i) \rmult u$.
 \end{enumerate} 
Normal forms again are elements of $z \rmult \f_{\myk} \cup \{ \zero \}$.
Hence a Gr\"obner basis of the form $G \cup \{ z \rmult p_i \mid i \in I, p_i \in \f_{\myk} \}$ with $G \subseteq \f_{\myk}$
 must exist.
\\
It remains to show that the set $\{ p_i \mid i \in I \}$ is in fact a right Gr\"obner basis of 
 $\ideal{r}{\f_{\myk}/\mathfrak{i}}(\{f\})$.
This follows immediately if we recall the history of the polynomials $p_i$.
In the first step they arise as a normal form with respect to $G \cup \{ z \rmult f \}$ of a polynomial
 either of the form $zu \rmult g \rmult v -z \rmult f \rmult w$, $z \rmult f \rmult u - z \rmult f \rmult v$
 or $z \rmult f \rmult u$, hence belonging to $\ideal{r}{\f_{\myk}/\mathfrak{i}}(\{f\})$.
In the iteration step, the new $p_n$ arises as a normal form with respect to $G \cup \{ z \rmult p_i \mid i \in I_{old} \}$ of a polynomial
 either of the form $zu \rmult g \rmult v -z \rmult p_i \rmult w$, $z \rmult p_i \rmult u - z \rmult p_j \rmult v$
 or $z \rmult p_i \rmult u$, hence belonging to $\ideal{r}{\f_{\myk}/\mathfrak{i}}(\{ p_i \mid i \in I_{old} \}) =
 \ideal{r}{\f_{\myk}/\mathfrak{i}}(\{f\})$.
\qed

Since we require $\f_{\myk}$ to have a unit (otherwise looking for inverse elements makes no sense), $\f_{\myk}^{z\myt}$ then
 will contain $z$.
\newpage
\problem{Inverse Element Problem}
{An element $f \in \f_{\myk}$ and a generating set $F$ for $\mathfrak{i}$.}
{Does $f$ have a right or left inverse in $\f_{\myk}/\mathfrak{i}$?}
{&1. & Let $G$ be a Gr\"obner basis of $\ideal{}{\f_{\myk}^{z\myt}}(F \cup \{ z \rmult f\})$. \\
 &2. & If $z \red{*}{\myr}{}{G} \zero$, then $f$ has a right inverse. \\
 &1'. & Let $G$ be a Gr\"obner basis of $\ideal{}{\f_{\myk}^{z\myt}}(F \cup \{ f \rmult z\})$. \\
 &2'. & If $z \red{*}{\myr}{}{G} \zero$, then $f$ has a left inverse. \\
}

To see that this is correct we give the following argument for the case of right inverses:
It is clear that $f$ has a right inverse in $\f_{\myk}/\mathfrak{i}$ if and only if 
 $f \rmult g -{\bf 1} \in \mathfrak{i}$ for some $g \in \f_{\myk}$.
On the other hand we get $f \rmult g -{\bf 1} \in \mathfrak{i}$ if and only if 
 $z \rmult f \rmult g - z \in \ideal{}{\f_{\myk}^{z\myt}}(F) \cap z \rmult \f_{\myk}$:
$f \rmult g -{\bf 1} \in \mathfrak{i}$ immediately
 implies $z \rmult (f \rmult g -{\bf 1}) \in \ideal{}{\f_{\myk}^{z\myt}}(F) \cap z \rmult \f_{\myk}$ as 
 $\mathfrak{i} \subseteq \ideal{}{\f_{\myk}^{z\myt}}(F)$, $z \in z\myt \subseteq \f_{\myk}^{z\myt}$ and 
 $z \rmult (f \rmult g -{\bf 1}) \in z \rmult \f_{\myk}$.
On the other hand, if $z \rmult f \rmult g - z \in \ideal{}{\f_{\myk}^{z\myt}}(F) \cap z \rmult \f_{\myk} \subseteq
 \ideal{}{\f_{\myk}^{z\myt}}(F)$, then we have a representation
 $z \rmult f \rmult g - z = \sum_{i=1}^k h_i \rmult f_i \rmult \tilde{h}_i$, $h_i, \tilde{h}_i \in \f_{\myk}^{z\myt}$,
 $f_i \in F \subseteq \f_{\myk}$.
For a polynomial $p \in \f_{\myk}^{z\myt}$ and some element $\alpha \in \myk$
 let $p [z=\alpha ]$ be the polynomial which arises from $p$ by substituting
 $\alpha$ for the variable $z$.
Then by substituting $z={\bf 1}$ we get $f \rmult g - {\bf 1} = \sum_{i=1}^k h_i[z={\bf 1}] \rmult f_i \rmult \tilde{h}_i[z={\bf 1}]$
 with $h_i[z={\bf 1}], \tilde{h}_i[z={\bf 1}] \in \f_{\myk}$ and are done.

Now, in order to decide whether $f$ has a right inverse in $\mathfrak{i}$ one
 has to distinguish the following two cases provided we have a Gr\"obner basis $G$ of
 $\ideal{}{\f_{\myk}^{z\myt}}(F \cup \{ z \rmult f\})$:
If $z \nred{*}{\myr}{}{G} \zero$ then there exists no $g \in \f_{\myk}$ such that $f \rmult g - {\bf 1} \in \mathfrak{i}$
 and hence $f$ has no right inverse.
If $z \red{*}{\myr}{}{G} \zero$ then we know $z \in \ideal{}{\f_{\myk}^{z\myt}}(F \cup \{ z \rmult f \})$,
 and even $z \in \ideal{r}{\f_{\myk}/\mathfrak{i}}(\{z \rmult f \})$ 
Hence there exist $m_i,\tilde{m}_i, n_j  \in \monoms(\f_{\myk}^{z\myt})$, $f_i \in F$ such that
 $$z = \sum_{i=1}^k m_i \rmult f_i \rmult \tilde{m}_i + \sum_{j=1}^l z \rmult f \rmult n_j.$$
Now substituting $z={\bf 1}$ gives us that for $h = \sum_{j=1}^l n_j$ we have
 $f \rmult h = {\bf 1} (\mbox{ mod }\mathfrak{i})$ and we are done.

As before,  weak Gr\"obner bases are sufficient to solve the problem.
%%%%%%%%%%%%%%%%%%%%%%%%%%%%%%%%%%%%%%%%%%%%%%%%%%%%%%%%%%%%%%%%%%%%%%%%%
\section{Elimination Theory}\label{section.elimination}
In ordinary polynomial rings special term orderings called elimination orderings can
 be used to produce Gr\"obner bases with useful properties.
Many problems, e.g.~the ideal intersection problem or the subalgebra problem, can be
 solved using tag variables.
The elimination orderings are then  used to separate the ordinary variables from these
 additional tag variables.
Something similar can be achieved for function rings.

Let $Z = \{ z_i \mid i \in I \}$ be a set of new tag variables commuting with terms.
The multiplication $\rmult$ can be extended by $z_i \rmult z_j = z_iz_j$, $z \rmult t = zt$
 and $t \rmult z = zt$ for $z,z_i,z_j \in Z$ and $t \in \myt$.
The ordering $\succeq$ is lifted to $Z^*\myt = \{ wt \mid w \in Z^*, t \in \myt \}$ by
 $w_1t_1 \succeq w_2t_2$ if and only if $w_1 \geq_{\rm lex} w_2$ or $(w_1 = w_2$ and $t_1 \succeq t_2$)
 for all $w_1,w_2 \in Z^*$, $t_1,t_2 \in \myt$.
Moreover, we require $w \succ t$ for all $w \in Z^*$, $t \in \myt$.
This ordering is called an elimination ordering.

Up to now we have studied ideals in $\f^{\myt}$.
Now we can view $\f^{\myt}$ as a subring of $\f^{Z^*\myt}$ and study ideals in both rings.
For a generating set $F \subset \f^{\myt}$ we have $\ideal{}{\f^{\myt}}(F) \subseteq \ideal{}{\f^{Z^*\myt}}(F)$.
This follows immediately since for every $f = \sum_{i=1}^k m_i \rmult f_i \rmult \tilde{m}_i$, $m_i, \tilde{m}_i \in \monoms(\f^{\myt})$
 this immediately implies $m_i, \tilde{m}_i \in \monoms(\f^{Z^*\myt})$.

\begin{lemma}\label{lem.elimination}~\\
{\sl
Let $G$ be a weak Gr\"obner basis of an ideal in $\f^{Z^*\myt}$ with respect to an elimination ordering.
Then the following hold:
\begin{enumerate}
\item $\ideal{}{\f^{Z^*\myt}}(G) \cap \f^{\myt} = \ideal{}{\f^{\myt}}(G \cap \f^{\myt})$.
\item $G \cap \f^{\myt}$ is a weak Gr\"obner basis for $\ideal{}{\f^{\myt}}(G \cap \f^{\myt})$ with respect to $\succeq$.
\item If $G$ is a Gr\"obner basis, then $G \cap \f^{\myt}$ is a Gr\"obner basis for $\ideal{}{\f^{\myt}}(G \cap \f^{\myt})$ with respect to $\succeq$.

\end{enumerate}
}
\end{lemma}
\Ba{}
\begin{enumerate}
\item 
  \begin{itemize}
  \item $\ideal{}{\f^{Z^*\myt}}(G) \cap \f^{\myt} \subseteq \ideal{}{\f^{\myt}}(G \cap \f^{\myt})$:\\
        Let $ f \in \ideal{}{\f^{Z^*\myt}}(G) \cap \f^{\myt}$.
        By the elimination ordering property for $w \in Z^*$ and $t \in \myt$
         we have that $wt \succ w \succ t$ holds and we get that $\hterm(f) \in \myt$
         if and only if $f \in \f^{\myt}$.
        Since $ f \in \ideal{}{\f^{Z^*\myt}}(G)$ we know that $f \red{*}{\myr}{}{G} \zero$ and as all 
         monomials in $f$ are also in $\f^{\myt}$ for each $g \in G$ used in this reduction sequence we know
         $\hterm(g) \in \myt$ and hence $g \in\f^{\myt}$.
        Moreover, the reduction sequence gives us a representation $f = \sum_{i=1}^k m_i \rmult f_i \rmult \tilde{m}_i$
         with $f_i \in G \cap \f^{\myt}$ and $m_i, \tilde{m}_i \in \monoms(\f^{\myt})$, 
         implying $f \in \ideal{}{\f^{\myt}}(G \cap \f^{\myt})$.
  \item $\ideal{}{\f^{\myt}}(G \cap \f^{\myt}) \subseteq \ideal{}{\f^{Z^*\myt}}(G) \cap \f^{\myt}$:\\
        Let $ f \in \ideal{}{\f^{\myt}}(G \cap \f^{\myt})$.
        Then $ f = \sum_{i=1}^k m_i \rmult f_i \rmult \tilde{m}_i$ with 
         $f_i \in G \cap \f^{\myt}$ and $m_i, \tilde{m}_i \in \monoms(\f^{\myt})$.
        Hence $f \in \ideal{}{\f^{\myt}}(G) \subseteq \ideal{}{\f^{Z^*\myt}}(G)$ and $f \in \f^{\myt}$
         imply $f \in \ideal{}{\f^{Z^*\myt}}(G) \cap \f^{\myt}$.
  \end{itemize}
\item We show this by proving that for every
       $f \in \ideal{}{\f^{\myt}}(G \cap \f^{\myt})$ we have
       $f \red{*}{\myr}{}{G \cap \f^{\myt}} \zero$.
      Since $G$ is a weak Gr\"obner basis of $\ideal{}{\f^{Z^*\myt}}(G)$ and 
       $\ideal{}{\f^{\myt}}(G \cap \f^{\myt}) \subseteq \ideal{}{\f^{Z^*\myt}}(G \cap \f^{\myt}) 
       \subseteq \ideal{}{\f^{Z^*\myt}}(G)$ we get $f \red{*}{\myr}{}{G} \zero$.
      On the other hand, as every monomial in $f$ is an element of $\f^{\myt}$, only elements of
       $G \cap \f^{\myt}$ are applicable for reduction.
\item Let $G$ be a Gr\"obner basis with respect to some reduction relation $\myr$.
        To show that $G \cap \f^{\myt}$ is a Gr\"obner basis of $\ideal{}{\f^{\myt}}(G \cap \f^{\myt})$
        we proceed in two steps:
        \begin{enumerate}
        \item $\red{*}{\lr}{}{G \cap \f^{\myt}} = \;\;\equiv_{\ideal{}{\f^{\myt}}(G \cap \f^{\myt})}$:\\
        $\red{*}{\lr}{}{G \cap \f^{\myt}} \subseteq \;\;\equiv_{\ideal{}{\f^{\myt}}(G \cap \f^{\myt})}$ trivially holds as because of Axiom (A2) reduction steps stay within the ideal
        congruence.
        To see the converse let $f \equiv_{\ideal{}{}(G \cap \f^{\myt})} g$ 
        for $f,g \in \f^{\myt}$.
        Then, as $G$ is a Gr\"obner basis and also
        $f \equiv_{\ideal{}{\f^{Z^*\myt}}(G)} g$ holds, we know 
        $f \red{*}{\lr}{}{G} g$ and as $\hterm(f), \hterm(g) \in\f^{\myt}$,
         only elements from $G \cap \f^{\myt}$ can be involved and we are done.  
        \item $\red{}{\myr}{}{G \cap \f^{\myt}}$ is confluent:\\
        Let $g, g_1, g_2 \in\f^{\myt}$ such that
         $g\red{}{\myr}{}{G \cap \f^{\myt}}g_1$ and
         $g \red{}{\myr}{}{G \cap \f^{\myt}} g_2$. 
        Then, as $\red{}{\myr}{}{G}$ is confluent we know that there exists
         $f \in \f^{Z^*\myt}$ such that $g_1 \red{*}{\myr}{}{G} f$ and
         $g_2 \red{*}{\myr}{}{G} f$.
        Now since $\hterm(g) \in \f^{\myt}$ we can conclude that $g_1, g_2, f \in
        \f^{\myt}$ and hence all polynomials used for reduction in the reduction sequences
        lie in $G \cap \f^{\myt}$ proving our claim.
        \end{enumerate}
\end{enumerate}
\qed
Given an ideal $\mathfrak{i} \subseteq \f^{Z^*\myt}$ the set $\mathfrak{i} \cap \f^{\myt}$ is again
 an ideal, now in $\f^{\myt}$.
This follows as 
\begin{enumerate}
\item $\zero \in \mathfrak{i} \cap \f^{\myt}$ since $\zero \in \mathfrak{i}$
      and $\zero \in \f^{\myt}$.
\item For $f,g \in \mathfrak{i} \cap \f^{\myt}$ we have $f+g \in \mathfrak{i}$ as $f,g \in \mathfrak{i}$
       and $f+g \in \f^{\myt}$ as $f,g \in \f^{\myt}$ yielding $f+g \in \mathfrak{i} \cap \f^{\myt}$.
\item For $f \in \mathfrak{i} \cap \f^{\myt}$ and $h \in \f^{\myt}$ we have that $f \rmult h, h \rmult f \in 
       \mathfrak{i}$ as $f \in \mathfrak{i}$ and $f \rmult h, h \rmult f \in \f^{\myt}$ as $f, h \in \f^{\myt}$
       yielding $f \rmult h, h \rmult f \in \mathfrak{i} \cap \f^{\myt}$.
\end{enumerate}
The ideal $\mathfrak{i} \cap \f^{\myt}$ is called the elimination ideal of $\mathfrak{i}$
 with respect to $Z$ since the occurrences of the tag variables $Z$ are eliminated.

\begin{definition}~\\
{\rm
For an ideal $\mathfrak{i}$ in $\f$ the set
$$\surd\mathfrak{i} = \{ f \in \f \mid \mbox{ there exists } m \in \n \mbox{ with } f^m \in \mathfrak{i} \}$$
is called the \betonen{radical} of $\mathfrak{i}$.
\dend
}
\end{definition}
Obviously we always have $\mathfrak{i} \subseteq \surd\mathfrak{i}$.
Moreover, if $\f$ is commutative the radical of an ideal is again an ideal.
This follows as
\begin{enumerate}
\item $\zero \in \surd\mathfrak{i}$ since $\zero \in \mathfrak{i}$,
\item For $f,g \in \surd\mathfrak{i}$ we know $f^m, g^n \in \mathfrak{i}$ for some $m,n \in \n$.
      Now $f+g \in \surd\mathfrak{i}$ if we can show that $(f+g)^q \in \mathfrak{i}$ for some $q \in \n$.
      Remember that for $q = m+n-1$ every term in the binomial expansion of $(f+g)^q$ has a factor
       of the form $f^i \rmult g^j$ with $i+j = m+n-1$.
      As either $i \geq m$ or $j \geq n$ we find $f^i \rmult g^j \in \mathfrak{i}$ yielding
       $(f+g)^q \in \mathfrak{i}$ and hence $f+g \in \surd\mathfrak{i}$.
      Notice that commutativity is essential in this setting.
\item For $f \in \surd\mathfrak{i}$ we know $f^m \in \mathfrak{i}$ for some $m \in \n$.
      Hence for $h \in \f^{\myt}$ we get $(f \rmult h)^m = f^m \rmult h^m \in \mathfrak{i}$
       yielding $f \rmult h \in \surd\mathfrak{i}$.
      Again commutativity is essential in the proof.
\end{enumerate}
Unfortunately this no longer holds for non-commutative function rings.
For example take $\myt = \{ a,b \}^*$ with concatenation as multiplication.
Then for $\mathfrak{i} = \ideal{}{}(\{a^2 \}) = \{ \sum _{i=1}^n \alpha_i \skm u_ia^2v_i \mid n \in \n,
 \alpha_i \in \q, u_i,v_i \in \myt \}$ we get $a \in \surd\mathfrak{i}$.
But for $b \in \f$ there exists no $m \in \n$ such that $(ab)^m \in \mathfrak{i}$
 and hence $\surd\mathfrak{i}$ is no ideal.

In the commutative polynomial ring the question whether some polynomial $f$ lies in the radical of
 some ideal generated by a set $F$ can be answered by introducing a tag variable $z$ and computing a 
 Gr\"obner basis of the ideal generated by the set $F \cup \{ fz -{\bf 1} \}$.
It can be shown that if a commutative function ring $\f$ contains a unit ${\bf 1}$ we get a similar result. 

\begin{theorem}\label{theo.radicalmembership}~\\
{\sl
Let $F \subseteq \f$ and $f \in \f$ where $\f$ is a commutative function ring containing a unit ${\bf 1}$.
Then $f \in\surd\ideal{}{\f^{\myt}}(F)$ if and only if ${\bf 1} \in \ideal{}{\f^{\{z\}^*\myt}}(F \cup \{ z \rmult f -{\bf 1} \})$
 for some new tag variable $z$.
}
\end{theorem}
\Ba{}~\\
If $f \in\surd\ideal{}{\f^{\myt}}(F)$, then 
 $f^m \in  \ideal{}{\f^{\myt}}(F)\subseteq \ideal{}{\f^{\{z\}^*\myt}}(F \cup \{ z \rmult f -{\bf 1} \})$ for some $m \in \n$.
But we also have that $z \rmult f -{\bf 1} \in \ideal{}{\f^{\{z\}^*\myt}}(F \cup \{ z \rmult f -{\bf 1} \})$.
Remember that for the tag variable we have $t \rmult z = zt$ for all $t \in \myt$ and hence
 $f \rmult z = z \rmult f$ yielding
\begin{eqnarray*}
{\bf 1} & = & z^m \rmult f^m - (z^m \rmult f^m -{\bf 1}) \\
  & = &  \underbrace{z^m \rmult f^m}_{\in  \ideal{}{\f^{\myt}}(F)}
         - \underbrace{(z \rmult f -{\bf 1}) \rmult (\sum_{i=0}^{m-1}z^i \rmult f^i)}_{\ideal{}{\f^{\{z\}^*\myt}}(F \cup \{ z \rmult f -{\bf 1} \})} 
\end{eqnarray*} 
and hence ${\bf 1} \in  \ideal{}{\f^{\{z\}^*\myt}}(F \cup \{ z \rmult f -{\bf 1} \})$ and we are done.
\\
On the other hand, ${\bf 1} \in \ideal{}{\f^{\{z\}^*\myt}}(F \cup \{ z \rmult f -{\bf 1} \})$ implies
 ${\bf 1} = \sum_{i=1}^k m_i \rmult f_i \rmult \tilde{m}_i + \sum_{j=1}^l n_j \rmult (z \rmult f -{\bf 1}) \rmult \tilde{n}_j$ with $m_i, \tilde{m}_i, n_j, \tilde{n}_j \in \monoms(\f^{\{z\}^*\myt})$.
Moreover, since for the tag variable we have $z \rmult t = t \rmult z = zt$ for all $t \in \myt$
 all terms occurring in $\sum_{i=1}^k g_i \rmult f_i \rmult h_i$ are of the form $z^jt$ for some $t \in \myt$,
 $j \in \n$.
Now, since $z \rmult f -{\bf 1} \in \ideal{}{\f^{\{z\}^*\myt}}(F \cup \{ z \rmult f -{\bf 1} \})$, we have
 $z^jt \rmult f^j = t \rmult z^j \rmult f^j = t$ as well as $f^j \rmult z^j \rmult t = z^j \rmult f^j \rmult t = t$.
Hence, the occurrences of $z$ in a term $z^jt$ with $t \in \myt$ can be ``cancelled'' by
 multiplication with $f^m$, $m \geq j$.
Therefore, by choosing $m \in \n$ sufficiently large to cancel {\em all} occurrences of $z$ in the
 terms of $\sum_{i=1}^k m_i \rmult f_i \rmult \tilde{m}_i$, multiplying the equation with
 $f^m$ from both sides yields
$$f^{2m} = \sum_{i=1}^k (f^m \rmult m_i) \rmult f_i \rmult (\tilde{m}_i \rmult f^m)$$
 and $f^m \rmult m_i, \tilde{m}_i \rmult f^m \in \f^{\myt}$.
Hence $f^{2m} \in \ideal{}{\f^{\myt}}(F)$ and therefore $f \in \surd\ideal{}{\f^{\myt}}(F)$.
\\
\qed

This theorem now enables us to describe the membership problem for radicals of ideals in terms of
 Gr\"obner bases.

\problem{Radical Membership Problem}
{A set $F \subseteq \f$ and an element $f \in \f$, $\f$ containing a unit ${\bf 1}$.}
{$f \in \surd\ideal{}{}(F)$?}
{&1.& Let $G$ be a respective Gr\"obner basis of $\ideal{}{\f^{\{z\}^*\myt}}(F \cup \{ z \rmult f -{\bf 1} \})$
 for some new tag variable $z$.\\
 &2.& If ${\bf 1} \red{*}{\myr}{}{G} \zero$, then $f \in \surd\ideal{}{}(F)$.
}

If additionally the function ring is commutative, remember that then $\surd\mathfrak{i}$ is an ideal
 and we then describe the equality problem for radicals of ideals.

Notice that weak Gr\"obner bases are sufficient to solve the problem.

\problem{Radical Equality Problem}
{Two sets $F_1,F_2 \subseteq \f$, $\f$ commutative containing a unit.}
{$\surd\ideal{}{}(F_1) = \surd\ideal{}{}(F_2)$?}
{&1.& If for all $f \in F_1$ we have $f \in \surd\ideal{}{}(F_2)$, 
       then $\surd\ideal{}{}(F_1) \subseteq \surd\ideal{}{}(F_2)$. \\
 &2.& If for all $f \in F_2$ we have $f \in \surd\ideal{}{}(F_1)$, 
       then $\surd\ideal{}{}(F_2) \subseteq \surd\ideal{}{}(F_1)$. \\
 &3.& If 1. and 2. both hold, then $\surd\ideal{}{}(F_1) = \surd\ideal{}{}(F_2)$.
} 

Correctness can be shown as follows:
Let us assume that for all $f \in F_1$ we have $f \in \surd\ideal{}{}(F_2)$.
Then, as $\f$ is commutative $\ideal{}{}(F_1) \subseteq \surd\ideal{}{}(F_2)$ holds.
Now let $f \in \surd\ideal{}{}(F_1)$.
Then for some $m \in \n$ we have $f^m \in \ideal{}{}(F_1) \subseteq \surd\ideal{}{}(F_2)$
 and hence $\surd\ideal{}{}(F_1) \subseteq \surd\ideal{}{}(F_2)$.

If $\f$ is not commutative, $\ideal{}{}(F_1) \subseteq \surd\ideal{}{}(F_2)$ need not hold.
Remember the function ring with $\myt = \{ a,b\}^*$.
Take $F_1 = \{ a \}$ and $F_2 = \{ a^2 \}$.
Then $a \in \surd\ideal{}{}(F_2)$ since $a^2 \in \ideal{}{}(F_2)$.
But while $ab \in \ideal{}{}(F_1)$ we have $ab \not\in \surd\ideal{}{}(F_2)$.

Radicals of one-sided ideals can be defined as well and Theorem \ref{theo.radicalmembership} is also valid in
 this setting and can be used to state the radical membership problem for one-sided ideals.

Another problem which can be handled using tag variables and elimination orderings in the
 commutative polynomial ring is that of ideal intersections.
Something similar can be done for function rings containing a unit.

\begin{theorem}\label{theo.intersection}~\\
{\sl
Let $\mathfrak{i}$ and $\mathfrak{j}$ be two ideals in $\f$ and $z$ a new tag variable.
Then 
$$ \mathfrak{i} \cap \mathfrak{j} = \ideal{}{\f^{\{z\}^*\myt}}(z \rmult \mathfrak{i} \cup (z-{\bf 1}) \rmult \mathfrak{j}) \cap \f $$
where $z \rmult\mathfrak{i} = \{ z \rmult f \mid f \in \mathfrak{i} \}$ and $(z-{\bf 1}) \rmult\mathfrak{j} = \{ (z-{\bf 1}) \rmult f \mid f \in \mathfrak{j} \}$.
}
\end{theorem}
\Ba{}~\\
Every polynomial $f \in \mathfrak{i} \cap \mathfrak{j}$ can be written as $f = z \rmult f - (z-{\bf 1}) \rmult f$
 and hence $f \in \ideal{}{\f^{\{z\}^*\myt}}(z \rmult \mathfrak{i} \cup (z-{\bf 1}) \rmult \mathfrak{j}) \cap \f$.
On the other hand, $f \in \ideal{}{\f^{\{z\}^*\myt}}(z \rmult \mathfrak{i} \cup (z-{\bf 1}) \rmult \mathfrak{j}) \cap \f$
 implies $f = \sum_{i=1}^k m_i \rmult z \rmult f_i \rmult \tilde{m}_i + \sum_{j=1}^l n_j \rmult (z-{\bf 1}) \rmult\tilde{f}_j  \rmult \tilde{n}_j$
 with $f_i \in \mathfrak{i}$, $\tilde{f}_j \in \mathfrak{j}$ and $m_i, \tilde{m}_i, n_j, \tilde{n}_j \in\monoms(\f^{\{z\}^*\myt})$.
Since $f \in \f^{\myt}$, substituting $z= {\bf 1}$ gives us $f \in \mathfrak{i}$ and $z=0$ gives us $f \in \mathfrak{j}$
 and hence $f \in \mathfrak{i} \cap \mathfrak{j}$.
\\
\qed

Moreover, combining this result with Lemma \ref{lem.elimination} gives us the means to
 characterize a Gr\"obner basis of the intersection ideal.

\problem{Intersection Problem}
{Two  sets $F_1,F_2 \subseteq \f$.}
{Determine a basis of $\ideal{}{}(F_1) \cap \ideal{}{}(F_2)$.}
{&1.& Let $G$ be a Gr\"obner basis of $\ideal{}{\f^{\{z\}^*\myt}}(z \rmult \mathfrak{i} \cup (z-{\bf 1}) \rmult \mathfrak{j})$ with respect to an elimination ordering with $z > \myt$. \\
&2.& Then $G \cap \f$ is a Gr\"obner basis of $\ideal{}{}(F_1) \cap \ideal{}{}(F_2)$.}

These ideas extend to one-sided ideals as well.

Again, weak Gr\"obner bases are sufficient to solve the problem.

Of course Theorem \ref{theo.intersection} can be generalized to intersections of more than two ideals.

The techniques can also be applied to treat quotients of ideals in case $\f_{\myk}$ is commutative.

\begin{definition}~\\
{\rm
For two ideals $\mathfrak{i}$ and $\mathfrak{j}$ in a commutative
 function ring $\f_{\myk}$ we define the quotient to be the set
 $$\mathfrak{i}/\mathfrak{j} = \{ g \mid g \in \f_{\myk} \mbox{ with } g \rmult \mathfrak{j} \subseteq \mathfrak{i}\}$$
where $g \rmult \mathfrak{j} = \{ g \rmult f \mid f \in \mathfrak{j} \}$.
\dend
}
\end{definition}

\begin{lemma}~\\
{\sl
Let $\f_{\myk}$ be a commutative function ring.
Let $\mathfrak{i}$ and $\mathfrak{j} = \ideal{}{}(F)$ be two ideals in $\f_{\myk}$.
Then
$$\mathfrak{i}/\mathfrak{j} = \bigcap_{f \in F} ( \mathfrak{i} / \ideal{}{}(\{ f \}).$$
}
\end{lemma}
\Ba{}
First let $g \in \mathfrak{i}/\mathfrak{j}$. Then $g \rmult \mathfrak{j} \subseteq \mathfrak{i}$.
Since $\mathfrak{j} = \ideal{}{}(F)$ we get $g \rmult f \in \mathfrak{i}$ for all $f \in F$.
As $\f_{\myk}$ is commutative we can conclude $g \rmult \ideal{}{}(\{f\}) \subseteq \mathfrak{i}$ for all $f \in F$
 and hence $g \in \mathfrak{i}/\ideal{}{}(\{ f \})$ for all $f \in F$ yielding
 $g \in \bigcap_{f \in F} ( \mathfrak{i} / \ideal{}{}(\{ f \})$.
\\
On the other hand, $g \in \bigcap_{f \in F} ( \mathfrak{i} / \ideal{}{}(\{ f \})$ implies
 $g \in \mathfrak{i}/\ideal{}{}(\{ f \})$ for all $f \in F$ and hence
 $g \rmult \ideal{}{}(\{ f \}) \subseteq \mathfrak{i}$  for all $f \in F$.
Since $\mathfrak{j} = \ideal{}{}(F)$ then $g \rmult \mathfrak{j} \subseteq \mathfrak{i}$ and
 hence $g \in \mathfrak{i}/\mathfrak{j}$.
\qed

Hence we can describe quotients of ideals in terms of quotients of the special
 form $\mathfrak{i} / \ideal{}{}(\{ f \})$.
These special quotients now can be described using ideal intersection in case $\f_{\myk}$ contains a unit
 element ${\bf 1}$.

\begin{lemma}~\\
{\sl
Let $\f_{\myk}$ be a commutative function ring.
Let $\mathfrak{i}$ be an ideal and $f \neq \zero$ a polynomial in $\f_{\myk}$.
Then 
$$\mathfrak{i} / \ideal{}{}(\{f\}) = (\mathfrak{i} \cap \ideal{}{}(\{ f \})) \rmult f^{-1}$$
where $f^{-1}$ is an element in $\f_{\myk}$ such that $f \rmult f^{-1} = {\bf 1}$.
}
\end{lemma}
\Ba{}~\\
First let $g \in \mathfrak{i}/\ideal{}{}(\{ f \})$.
Then $g \rmult \ideal{}{}(\{ f \}) \subseteq \mathfrak{i}$ and $g \rmult f \in \mathfrak{i}$, even
 $g \rmult f \in \mathfrak{i} \cap \ideal{}{}(\{ f \})$.
Hence $g \in (\mathfrak{i} \cap \ideal{}{}(\{ f \})) \rmult f^{-1}$.
\\
On the other hand let $g \in (\mathfrak{i} \cap \ideal{}{}(\{ f \})) \rmult f^{-1}$.
Then $g \rmult f \in \mathfrak{i} \cap \ideal{}{}(\{ f \}) \subseteq \mathfrak{i}$.
Since $\f_{\myk}$ is commutative, this implies $g \rmult \ideal{}{}(\{ f \}) \subseteq \mathfrak{i}$
 and hence $g \in \mathfrak{i}/\ideal{}{}(\{ f \})$.
\\
\qed

Hence we can study the quotient of $\mathfrak{i}$ and $\mathfrak{j} = \ideal{}{}(F)$ by studying
 $(\mathfrak{i} \cap \ideal{}{}(\{ f \})) \rmult f^{-1}$ for all $f \in F$.

%%%%%%%%%%%%%%%%%%%%%%%%%%%%%%%%%%%%%%%%%%%%%%%%%%%%%%%%%%%%%%%%
\section{Polynomial Mappings}
In this section we are interested in $\myk$-algebra homomorphisms between the non-commutative
 polynomial ring $\myk[Z^*]$ where $Z = \{ z_1,\ldots, z_n \}$, and $\f_{\myk}^{\myt}$.
Let
$$\phi: \myk[Z^*] \myr \f_{\myk}^{\myt}$$
be a ring homomorphism which is determined by a linear mapping
$$\phi: z_i \longmapsto f_i$$
with $f_i \in \f_{\myk}^{\myt}$, $1 \leq i \leq n$.
Then for a non-commutative polynomial $g \in \myk[Z^*]$ with $g = \sum_{j=1}^m \alpha_j \skm w_j$,
 $w_j \in Z^*$ we get $\phi(g) = \sum_{j=1}^m \alpha_j \skm \phi(w_j)$ where
 $\phi(w_j) = w_j[z_1 \longmapsto f_1, \ldots, z_n \longmapsto f_n]$.
The kernel of such a mapping is defined as 
 $$\kernel(\phi) = \{ g \in \myk[Z^*] \mid \phi(g) = \zero \}$$
and the image is defined as
 $$\image(\phi) = \{ f \in \f_{\myk}^{\myt} \mid \mbox{ there exists } g \in \myk[Z^*] \mbox{ such that } \phi(g) = f \}.$$
Note that $\image(\phi)$ is a subalgebra of $\f_{\myk}^{\myt}$.
\begin{lemma}~\\
{\sl 
Let $\phi: \myk[Z^*] \myr \f_{\myk}^{\myt}$ be a ring homomorphism.
Then $\myk[Z^*] / \kernel(\phi) \cong \image(\phi)$.
\lemend
}
\end{lemma}
\Ba{}~\\
To see this inspect the mapping $\psi:\myk[Z^*] / \kernel(\phi) \myr \image(\phi)$
 defined by $g + \kernel(\phi) \mapsto \phi(g)$.
Then $\psi$ is an isomorphism. 
\begin{enumerate}
\item $\psi(g + \kernel(\phi)) = \zero$ for $g \in \kernel(\phi)$ by the definition of $\kernel(\phi)$.
\item $\psi((g_1 + \kernel(\phi)) + (g_2 + \kernel(\phi))) = \phi(g_1 + g_2) = \psi(g_1 + \kernel(\phi)) + \psi(g_2 + \kernel(\phi))$.
\item $\psi((g_1 + \kernel(\phi)) \rmult (g_2 + \kernel(\phi))) = \phi(g_1 \rmult g_2) =
       \psi(g_1 + \kernel(\phi)) \rmult \psi(g_2 + \kernel(\phi))$, as for $g \in\myk[Z^*]$ and $h \in\kernel(\phi)$
       we have $\psi(g \rmult h) = \psi(h \rmult g) = \zero$.
\item $\psi$ is onto as its image is the image of $\phi$ and by the definition of
       the latter for each $f \in \image(\phi)=\image(\psi)$ there exists $g \in \myk[Z^*]$ such that
       $\phi(g) = f$.  
      Since for all $h \in \kernel(\phi)$ we have $\phi(h) = \zero$ then $\psi(g + \kernel(\phi))=\psi(g)=\phi(g)$.  
\item Assume that for $g_1,g_2 \in \myk[Z^*]$ we have $\psi(g_1 + \kernel(\phi)) = \psi(g_2 + \kernel(\phi))$.
      Then $\phi(g_1) = \phi(g_2)$ and this immediately implies that $g_1 - g_2 \in \kernel(\phi)$ and hence $\psi$
       is also a monomorphism.
\end{enumerate} 
\qed
Now the theory of elimination described in the previous section 
 can be used to provide a Gr\"obner basis for $\kernel(\phi)$.
Remember that the tag variables commute with the elements on $\myt$.
Again we use the function ring $\f_{\myk}^{Z^*\myt}$ and the
 fact that $\myk[Z^*] \subseteq \f_{\myk}^{Z^*\myt}$ by mapping the polynomials
 to the respective functions in $\f_{\myk}^{Z^*} \subseteq \f_{\myk}^{Z^*\myt}$.
\begin{theorem}\label{theo.kernel}~\\
{\sl
Let ${\mathfrak i} = \ideal{}{}(\{ z_1 - f_1, \ldots, z_n - f_n \}) \subseteq
 \f_{\myk}^{Z^*\myt}$.
Then $\kernel(\phi) = {\mathfrak i} \cap \myk[Z^*]$.
\theoend
}
\end{theorem}
\Ba{}~\\
Let $g \in {\mathfrak i} \cap \myk[Z^*]$ .
Then $g = \sum_{j=1}^n h_j \rmult s_j \rmult h_j'$ with $s_j \in \{ z_1 - f_1, \ldots, z_n - f_n \}$,
 $h_j, h_j' \in \f_{\myk}^{Z^*\myt}$.
As $\phi(z_j - f_j) = \zero$ for all $ 1 \leq j \leq n$ we get $\phi(g) = \zero$
 and hence $g \in \kernel(\phi)$. \\
To see the converse let $g \in \kernel(\phi)$.
Then $g \in \myk[Z^*]$ and hence $g = \sum_{j=1}^m \alpha_j \skm w_j$ where $w_j \in Z^*$, $1 \leq j \leq m$.
On the other hand we know $\phi(g) = \zero$.
Then 
\begin{eqnarray}
g & = & g - \phi(g) \nonumber\\ 
  & = & \sum_{j=1}^m \alpha_j \skm w_j - \sum_{j=1}^m \alpha_j \skm \phi(w_j) \nonumber\\ 
  & = & \sum_{j=1}^m \alpha_j \skm (w_j - \phi(w_j)) \nonumber
\end{eqnarray}
It remains to show that $w - \phi(w) \in {\mathfrak i}$ for all $w \in Z^*$
 as this implies $g \in{\mathfrak i} \cap \myk[Z^*]$.
This will be done by induction on $k = |w|$.
For $k = 1$ we get $w = z_i$ for some $1 \leq i \leq n$ and $w - \phi(w) = z_i - f_i \in {\mathfrak i}$.
In the induction step let $w \id a_1 \ldots a_k$, $a_i \in Z$.
Then we get 
\begin{eqnarray}
&&a_1(a_2 \ldots a_k - \phi(a_2 \ldots a_k)) + (a_1 - \phi(a_1))\phi(a_2 \ldots a_k)\nonumber\\ 
&=&a_1a_2 \ldots a_k - a_1\phi(a_2 \ldots a_k) + a_1\phi(a_2 \ldots a_k)-\phi(a_1)\phi(a_2 \ldots a_k)\nonumber\\ 
&=&a_1a_2 \ldots a_k -\phi(a_1 \ldots a_k)\nonumber 
\end{eqnarray}
Then, as $|a_2 \ldots a_k| = k-1$ the induction hypothesis yields $a_2 \ldots a_k - \phi(a_2 \ldots a_k)
 \in {\mathfrak i}$ and as of course $a_1 - \phi(a_1) \in {\mathfrak i}$
 we find that $a_1a_2 \ldots a_k -\phi(a_1 \ldots a_k) \in {\mathfrak i}$.
\\
\qed
Now  if $G$ is a (weak) Gr\"obner basis of
 ${\mathfrak i}$ in $\f_{\myk}^{Z^*\myt}$ with respect to
 an elimination ordering where the elements in $Z^*$ are made
 smaller than those in $\myt$,
 then $G \cap \myk[Z^*]$ is a (weak) Gr\"obner basis of the kernel of $\phi$.
Hence, in case finite such bases exist or bases allowing to solve the membership problem,
 they can be used to treat the following question.

\problem{Kernel of a Polynomial Mapping}
{A set $F = \{ z_1 - f_1, \ldots, z_n - f_n \} \subseteq \f_{\myk}^{Z^*\myt}$
 encoding a mapping $\phi :  \myk[Z^*] \myr \f_{\myk}^{\myt}$ and an element $f \in \myk[Z^*]$.}
{$f \in \kernel(\phi)$?}
{&1.& Let $G$ be a (weak) Gr\"obner basis of $\ideal{}{}(\{ z_1 - f_1, \ldots, z_n - f_n \})$
      with respect to an elimination ordering.\\
 &2.& Let $G' = G \cap \myk[Z^*]$. \\
 &3.& If $f \red{*}{\myr}{}{G'} \zero$, then $f \in \kernel(\phi)$.}

A similar question can be asked for the image of a polynomial mapping.

\problem{Image of a Polynomial Mapping}
{A set $F = \{ z_1 - f_1, \ldots, z_n - f_n \} \subseteq \f_{\myk}^{Z^*\myt}$
 encoding a mapping $\phi :  \myk[Z^*] \myr \f_{\myk}^{\myt}$ and an element $f \in \f_{\myk}^{\myt}$.}
{$f \in \image(\phi)$?}
{&1.& Let $G$ be a Gr\"obner basis of $\ideal{}{}(\{ z_1 - f_1, \ldots, z_n - f_n \})$
      with respect to an elimination ordering.\\
 &2.& If $f \red{*}{\myr}{}{G} h$, with $h \in \myk[Z^*]$, 
       then $f \in \image(\phi)$.}

The basis for this solution is the following theorem.
\begin{theorem}\label{theo.image}~\\
{\sl
Let ${\mathfrak i} = \ideal{}{}(\{ z_1 - f_1, \ldots, z_n - f_n \}) \subseteq
 \f_{\myk}^{Z^*\myt}$ and let $G$ be a Gr\"obner basis of ${\mathfrak i}$ with respect to an 
 elimination ordering where the elements in $Z^*$ are smaller than those in $\myt$.
Then $f \in \f_{\myk}^{\myt}$ lies in the image of $\phi$ if and only if there exists 
 $h \in \myk[Z^*]$ such that $f \red{*}{\myr}{}{G} h$.
Moreover, $f = \phi(h)$.
\theoend
}
\end{theorem}
\Ba{}~\\
Let $f \in \image(\phi)$, i.e., $f \in \f_{\myk}^{\myt}$.
Then $f = \phi(g)$ for some $g \in \myk[Z^*]$.
Moreover, $f - g = \phi(g) - g$, and similar to the proof of Theorem \ref{theo.kernel} we can show
 $f -g \in {\mathfrak i}$.
Hence, $f$ and $g$ must reduce to the same normal form $h$ with respect to $G$.
As $g \in \myk[Z^*]$ this implies $h \in \myk[Z^*]$ and we are done. \\
To see the converse, for $f \in \f_{\myk}^{\myt}$ let $f \red{*}{\myr}{}{G} h$ with $h \in \myk[Z^*]$.
Then $f -h \in {\mathfrak i}$ and hence $f-h = \sum_{j=1}^k g_j \rmult s_j \rmult g_j'$ with
 $s_j \in \{ z_1 - f_1, \ldots, z_n - f_n \}$, $g_j, g_j' \in  \f_{\myk}^{\myt}$.
As $\phi(s_j) = \zero$ we get $f - \phi(h) = \zero$ and hence $f = \phi(h)$ is in the image of $\phi$. 
\\
\qed
Obviously the question of whether an element lies in the image of $\phi$ then can be answered
 in case we can compute a unique normal form of the element with respect to the Gr\"obner basis
 of ${\mathfrak i} = \ideal{}{}(\{ z_1 - f_1, \ldots, z_n - f_n \})$.

Another question is whether the mapping $\phi:  \myk[Z^*] \myr \f_{\myk}^{\myt}$ is onto.
This is the case if for every $t \in \myt$ we have $t \in \image(\phi)$.
A simpler solution can be found in case $\myt \subseteq \Sigma^*$ for some finite set of letters
 $\Sigma = \{ a_1, \ldots, a_k \}$ and additionally $\myt$ is subword closed as a subset of 
 $\Sigma^*$. 
\begin{theorem}~\\
{\sl
Let ${\mathfrak i} = \ideal{}{}(\{ z_1 - f_1, \ldots, z_n - f_n \}) \subseteq
 \f_{\myk}^{Z^*\myt}$ and let $G$ be a Gr\"obner basis of ${\mathfrak i}$ with respect to an 
 elimination ordering where the elements in $Z^*$ are smaller than those in $\myt$.
Then $f \in \f_{\myk}^{\myt}$ is onto if and only if for each $a_j \in \Sigma$,
 we have $a_j \red{*}{\myr}{}{G} h_j$ where $h_j \in \myk[Z^*]$.
Moreover, $a_j = \phi(h_j)$.
\theoend
}
\end{theorem}
\Ba{}~\\
Remember that $\phi$ is onto if and only if $a_j \in \image(\phi)$ for
 $1 \leq j \leq k$.\\
Let us first assume that $\phi$ is onto, i.e., $a_1, \ldots, a_k \in \image(\phi)$.
Then by Theorem \ref{theo.image} there exist $h_j \in \myk[Z^*]$ such that
 $a_j \red{*}{\myr}{}{G} h_j$, $1 \leq j \leq k$. \\
To see the converse,
 again, by Theorem \ref{theo.image} the existence of $h_j \in \myk[Z^*]$ such that
 $a_j \red{*}{\myr}{}{G} h_j$, $1 \leq j \leq k$ now implies $a_1, \ldots, a_k \in \image(\phi)$  
 and we are done.
\\
\qed

%%%%%%%%%%%%%%%%%%%%%%%%%%%%%%%%%%%%%%%%%%%%%%%%%%%%%%%%%%%%%%%%
\section{Systems of One-sided Linear Equations in Function Rings over the Integers}
Let $\f_{\z}$ be the function ring over the integers $\z$ as specified in Section
 \ref{section.right.integers}.
Additionally we require that multiplying terms by terms results in terms,
 i.e.,~$\rmult : \myt \times \myt \myr \myt$.
Then a reduction relation can be defined for $\f_{\z}$ as follows:
\begin{definition}~\\
{\rm
Let $p$, $f$ be two non-zero polynomials in $\f_{\z}$. 
We say $f$ \betonen{reduces} $p$ \betonen{to} $q$ \betonen{at} $\alpha \skm t$ \betonen{in one step},
 i.e. $p \red{}{\myr}{}{g} q$, if
\begin{enumerate}
\item[(a)] $t = \hterm(f \rmult u) = \hterm(f) \rmult u$ for some $u \in \myt$.
\item[(b)] $\hc(f) > 0$ and $\alpha = \hc(f) \skm \beta + \delta$ with $\beta,\delta \in \z$, $\beta \neq 0$, and
           $0 \leq \delta < \hc(f)$.
\item[(c)] $q = p - f \rmult (\beta \skm u)$.
\end{enumerate}
}
\end{definition}
The definition of s-polynomials can be derived from Definition \ref{def.s-poly.z}.
\begin{definition}~\\
{\rm
Let $p_{1}, p_{2}$ be  two polynomials in $\f_{\z}$.
If there are respective terms
 $t,u_1, u_2 \in \myt$ such that
 $\hterm(p_i) \rmult u_i = \hterm(p_i \rmult u_i)=t
 \geq \hterm(p_i)$ let $HC(p_{i}) = \gamma_i$.\\
Assuming $\gamma_1 \geq \gamma_2 > 0$\footnote{Notice that $\gamma_i > 0$
 can always be achieved by studying the situation for $- p_i$ in case 
 we have $HC(p_{i}) < 0$.},
 there are $\beta, \delta \in \z$ such that
 $\gamma_1 = \gamma_2 \skm \beta + \delta$ and $0 \leq \delta < \gamma_2$
 and we get the following s-polynomial
      $$\spol{}(p_{1}, p_{2},t,u_1, u_2) =  \beta \skm p_2 \rmult u_2 -  p_1 \rmult u_1.$$
The set $\spols(\{p_1,p_2\})$ then is the set of all such
 s-polynomials corresponding to $p_1$ and $p_2$.
\dend
}
\end{definition}
Notice that two polynomials can give rise to infinitely many s-polynomials.
A subset $C$ of these possible s-polynomials $\spols(p_1, p_2)$ is called a
 stable localization
 if for any possible s-polynomial $p \in \spols(p_1, p_2)$ there exists a special
 s-polynomial $h \in C$ such that $p \red{}{\myr}{}{h} \zero$.
% \begin{enumerate}
% \item $p$ arises from the overlap created by $t = \hterm(g_1 \rmult t_1) = \hterm(g_1) \rmult t_1
%        = \hterm(g_2) \rmult t_2 = \hterm(g_2 \rmult t_2)$,
%        $t \geq \hterm(g_1)$, $t \geq \hterm(g_2)$,
% \item $s$ arises from the overlap $u = \hterm(g_1 \rmult u_1) = \hterm(g_1) \rmult u_1
%        = \hterm(g_2) \rmult u_2 = \hterm(g_2 \rmult u_2)$,
%        $u \geq \hterm(g_1)$, $u \geq \hterm(g_2)$,
% \item $t \geq u$,
% \item $t_1 = u_1 \rmult v$, $t_2 = u_2 \rmult v$, $v \in \myt$,
% \item $p = s \rmult v$ for some $v \in \myt$.
% \end{enumerate}

In the following let $f_1, \ldots, f_m \in \f_{\z}$.
We describe a generating set of solutions for the linear one-sided inhomogeneous equation 
 $f_1 \rmult X_1 + \ldots + f_m \rmult X_m = f_0$
 in the variables $X_1, \ldots, X_m$ provided a finite computable right Gr\"obner basis of the right ideal 
 generated by $\{ f_1, \ldots, f_m \}$ in $\f_{\z}$ exists.
%Additionally we assume that we can verify the Gr\"obner basis property by checking two
% kinds of test sets
% \begin{itemize}
% \item ${\sf SPOL}(F)$ corresponding to critical polynomials of the generating set $F$
%   called s-, g-, or m-polynomials depending on the context,
% \item ${\sf SAT}(F)$ corresponding to critical polynomial multiples of the generating set $F$ (compare
%   the saturation techniques described e.g.~in \cite{MaRe95}). 
% \end{itemize}
%Both sets are dependent on the reduction relation chosen for the function ring and we require the
% polynomials in this set to be restrictions of all candidates in the sense that
% for any s-polynomial orespectively saturating polynomial $h$ there exists a polynomial $s$ in the
% respective set ${\sf SPOL}(F)$ or ${\sf SAT}(F)$ such that $h = s \rmult m$ for some $m \in \monoms$,
% $\hm(h) = \hm(\hm(s) \rmult m)$, for any term $t \in \term$ with $t \pred \hterm(s)$ we have 
% $\hterm(t \rmult m) \pred \hterm(t)$.

In order to find a generating set of solutions we have to find {\em one} solution of
\begin{eqnarray}
f_1 \rmult X_1 + \ldots + f_m \rmult X_m & = & f_0 \label{equation.inhomogeneous}
\end{eqnarray}
and if possible a finite set of generators for the solutions of the homogeneous equation
\begin{eqnarray}
f_1 \rmult X_1 + \ldots + f_m \rmult X_m & = & \zero. \label{equation.homogeneous}
\end{eqnarray}
We proceed as follows assuming that we have a {\em finite} right Gr\"obner basis of the right
 ideal generated by $\{ f_1, \ldots, f_m \}$:
\begin{enumerate}
\item Let $G = \{ g_1, \ldots, g_n \}$ be a right Gr\"obner basis of the  right ideal
 generated by $\{ f_1, \ldots, f_m \}$ in $\f_{\z}$, and ${\bf f} = (f_1, \ldots, f_m)$, 
 ${\bf g} = (g_1, \ldots, g_n)$ the corresponding vectors. 
There are two linear mappings given by matrices $P \in {\sf M}_{m \times n}(\f_{\z})$, 
 $Q \in {\sf M}_{n \times m}(\f_{\z})$ such that ${\bf f} \skm P = {\bf g}$ and 
 ${\bf g} \skm Q = {\bf f}$.
\item Equation \ref{equation.inhomogeneous} is solvable if and only if 
 $f_0 \in \ideal{r}{}(\{ f_1, \ldots, f_m \})$. 
 This is equivalent to $f_0 \red{*}{\myr}{r}{G} 0$ and the reduction sequence gives
 rise to a representation $f_0 = \sum_{i=1}^n g_i \rmult h_i = {\bf g} \skm {\bf h}$ where
 ${\bf h} = (h_1, \ldots, h_n)$.
 Then, as ${\bf f} \skm P = {\bf g}$, we get ${\bf g} \skm {\bf h} = ({\bf f} \skm P) \skm {\bf h}$
 and $P \skm {\bf h}$ is such a solution of equation \ref{equation.inhomogeneous}.
\item Let $\{ {\bf z}_1, \ldots, {\bf z}_r \}$ be a generating set for the solutions of the  
 homogeneous equation
\begin{eqnarray}
g_1 \rmult X_1 + \ldots + g_n \rmult X_n & = & 0 \label{equation.homogeneous.G}
\end{eqnarray}
 and let $I_m$ be the $m \times m$ identity matrix. 
Further let ${\bf w}_1, \ldots, {\bf w}_m$ be the columns of the matrix $P \skm Q - I_m$. 
Since ${\bf f} \skm (P \skm Q - I_m) = {\bf f} \skm P \skm Q - {\bf f} \skm I_m = {\bf g} \skm Q - {\bf f} = 0$ 
 these are solutions of equation \ref{equation.homogeneous}. 
We can even show that the set $\{ P \skm {\bf z}_1, \ldots, P \skm {\bf z}_r, {\bf w}_1, \ldots, {\bf w}_m \}$ 
 generates all solutions of equation \ref{equation.homogeneous}:\\
Let ${\bf q} = (q_1, \ldots, q_m)$ be an arbitrary solution of equation \ref{equation.homogeneous}. 
Then $Q \skm {\bf q}$ is a solution of equation \ref{equation.homogeneous.G} as ${\bf f} = {\bf g} \skm Q$. 
Hence there are $h_1, \ldots, h_r \in \f_{\z}$ such that
 $Q \skm {\bf q} = {\bf z}_1 \skm h_1 + \ldots {\bf z}_r \skm h_r$. 
Further we find
$$ {\bf q} = P \skm Q \skm {\bf q} - (P \skm Q - I_m) \skm {\bf q} =
             P \skm  {\bf z}_1 \skm h_1 + \ldots P \skm {\bf z}_r \skm h_r + {\bf w}_1 \skm q_1 + \ldots + {\bf w}_m \skm q_m$$
and hence ${\bf q}$ is a right linear combination of elements in 
 $\{ P \skm {\bf z}_1, \ldots, P \skm {\bf z}_r, {\bf w}_1, \ldots, {\bf w}_m \}$. 
\end{enumerate}
Now the important part is to find a generating set for the solutions of the homogeneous
 equation \ref{equation.homogeneous.G}.
In commutative polynomial rings is was sufficient to look at special vectors arising from
 those situations causing s-polynomials.
These situations are again important in our setting:

For every $g_i,g_j \in G$ not necessarily different such that the stable  
 localization $C_{i,j} \subseteq \spols(g_i,g_j)$ for the s-polynomials is not empty  
 and additionally we require these sets to be finite, we compute
 vectors ${\bf a}^{\ell}_{ij}$, $1 \leq \ell \leq |C|$ as follows:  
\\
Let $t = \hterm(g_i \rmult u) = \hterm(g_i) \rmult u 
 = \hterm(g_j) \rmult v = \hterm(g_j \rmult v)$, 
 $t \geq \hterm(g_i)$, $t \geq \hterm(g_j)$, be the overlapping term corresponding
 to  $h_{\ell} \in C_{i,j}$.
Further let $\hc(g_i) \geq \hc(g_j) > 0$ and $\hc(g_i) = \alpha \skm \hc(g_j) + \beta$ for some
 $\alpha, \beta \in \z$, $0 \leq \beta < \hc(g_j)$. 
Then
$$h_{\ell} = g_i \rmult u - g_j \rmult (\alpha \skm v) = \sum_{l=1}^n g_l \rmult h_l,$$
       where the polynomials $h_l \in \f_{\z}$ are due to the reduction sequence
    $h_{\ell}  \red{*}{\myr}{r}{G} 0$.
       \\
       Then ${\bf a}_{ij}^{\ell} = ( a_1, \ldots, a_n)$, where
       \begin{eqnarray}
         a_i & = & h_i - u, \nonumber\\ 
         a_j & = & h_j + \alpha \skm v, \nonumber\\
         a_l & = & h_l,      \nonumber
       \end{eqnarray}
       $l \neq i,j$, is a solution of \ref{equation.homogeneous.G} as
       $\sum_{l=1}^n g_l \rmult h_l - g_i \rmult u + g_j \rmult \alpha \skm v = 0$.

If all sets $\spols(g_i,g_j)$ are empty for $g_i, g_j \in G$, in the case of ordinary Gr\"obner bases
 in polynomial rings one could conclude that the homogeneous equation \ref{equation.homogeneous.G}
 had no solution.
This is no longer true for arbitrary function rings.

\begin{example}\label{exa.a+1}~\\
{\rm 
Let $\z[\m]$ be a monoid ring where $\m$ is presented by the complete string rewriting system
 $\Sigma = \{a,b\}$, $T = \{ ab \myr \lambda \}$.
Then for the homogeneous equation
$$(a+1)\rmult X_1 + (b+1) \rmult X_2 = 0$$
we find that the set $\{a+1, b+1\}$ is a prefix  Gr\"obner basis of the right ideal it
 generates.
Moreover neither of the head terms of the polynomials in this basis is prefix of the other and
 hence no s-polynomials with respect to prefix reduction exist.
Still the equation can be solved: $(b, -1)$ is a solution since $(a+1) \rmult b - (b+1) = b+1 - (b+1) = 0$.
}
\end{example}

Hence inspecting s-polynomials is not sufficient to describe all solutions.
This phenomenon is due to the fact that as seen before in most function rings s-polynomials are not sufficient
 for a Gr\"obner basis test.
Additionally the concept of saturation has to be incorporated.
In Example \ref{exa.a+1} we know that $(a+1) \rmult b = 1+b$, i.e.~$b+1 \in \sat(a+1)$.
Of course $(a+1) \rmult b \red{}{\myr}{}{b+1} 0$ and hence $(a+1) \rmult b = b+1$ gives rise to
 a solution $(b, -1)$ as required above.

More general we can express these additional solutions as follows:
For every $g_i \in G$ with $\sat(g_i)$ a stable saturator for $\{ g_i \}$
 and again we additionally require it to be finte,
 we define vectors
 ${\bf b}_{i,\ell} = ( b_1, \ldots, b_n)$
 $1 \leq \ell \leq |\sat(g_i)|$ as follows:
For  $g_i \rmult w_{\ell} \in \sat(g_i)$ we know 
 $g_i \rmult w_\ell= \sum_{l=1}^n g_l \rmult h_l$ as $G$ is a  Gr\"obner basis.
 Then ${\bf b}_{i,\ell} = ( b_1, \ldots, b_n)$, where
 \begin{eqnarray}
 b_i & = & h_i - w_{\ell},\nonumber\\
 b_l & = & h_l,     \nonumber
 \end{eqnarray}
 $l \neq i$, is a solution of equation \ref{equation.homogeneous.G}
 as $\sum_{l=1}^n g_l \rmult h_l - g_i \rmult w_\ell = 0$.
\begin{lemma}~\\
{\sl
Let $\{ g_1, \ldots, g_n \}$ be a finite right Gr\"obner basis.
For $g_i,g_j$ let $C_{i,j}$ be a stable localization of $\spols(g_i,g_j)$.
The finitely many vectors ${\bf a}^{\ell_1}_{i,j}, {\bf b}_{i,\ell_2}$, $1 \leq i,j \leq n$, $1 \leq \ell_1 \leq |C_{i,j}|$, $1 \leq \ell_2 \leq |\sat(g_i)|$
 form a right generating set for all solutions of equation \ref{equation.homogeneous.G}.
}
\end{lemma}
\Ba{}~\\
Let ${\bf p} = (p_1, \ldots, p_n)$ be an arbitrary (non-trivial) solution of equation \ref{equation.homogeneous.G},
 i.e., $\sum_{i=1}^n g_i \rmult p_i = 0$.
Let $T_p = \max \{ \hterm(g_i \rmult t_j^{p_i}) \mid 1 \leq i \leq n, p_i = \sum_{j=1}^{n_i} \alpha_j^{p_i} \skm t_j^{p_i} \}$, 
 $K_p$  the number of multiples $g_i \rmult t_j^{p_i}$ with
 $T_p = \hterm(g_i \rmult t_j^{p_i}) \neq \hterm(g_i) \rmult t_j^{p_i}$, and
 $M_p = \{ \{ \hc(g_i) \mid \hterm(g_i \rmult t_j^{p_i}) = T_p \} \}$ a multiset in $\z$.
A solution ${\bf q}$ is called smaller than ${\bf p}$ if either $T_q \pred T_p$ or ($T_q = T_p$ and
 $K_q < K_p$) or ($T_q = T_p$ and $K_q = K_p$ and $M_q \ll M_p$).
We will prove our claim by induction on $T_p$, $K_p$ and $M_p$ and have to distinguish two cases:
\begin{enumerate}
\item  If there is $1 \leq i \leq n$, $1 \leq j \leq n_i$ such that 
  $T_p = \hterm(g_i \rmult t_j^{p_i})\neq \hterm(g_i) \rmult t_j^{p_i}$, 
 then there exists $s_{\ell} \in \sat(g_i)$ such that $g_i \rmult t_j^{p_i} = s_{\ell} \rmult v$ for
 some $v \in \myt$, $\hterm(s_{\ell} \rmult v) = \hterm(s_{\ell}) \rmult v$
 and $s_{\ell} = g_i \rmult w_{\ell}$, $w_{\ell} \in \myt$.
 Then we can set ${\bf q} = {\bf p} + \alpha_j^{p_i} \skm {\bf b}_{i,\ell} \rmult v$ with
 \begin{eqnarray}
 q_i & = & p_i + \alpha_j^{p_i} \skm (h_i - w_{\ell}) \rmult v \nonumber \\
 q_l & = & p_l + \alpha_j^{p_i} \skm h_l \rmult v \mbox{ for }l \neq i\nonumber
 \end{eqnarray}
 which is again a solution of equation \ref{equation.homogeneous.G}.
 It remains to show that it is a smaller one.
 To see this we have to examine the multiples $g_l \rmult t_j^{q_l}$ for all $1 \leq l \leq n$, $1 \leq j \leq m_l$
 where $q_l = \sum_{j=1}^{m_l} \alpha_j^{q_l} \skm t_j^{q_l}$.
 Remember that
        $\hterm(s_{\ell}) \leq \hterm(s_{\ell} \rmult v) = \hterm(s_{\ell}) \rmult v = T_p$.
 Moreover, for all terms $w_j^{h_l}$ in $h_l = \sum_{j=1}^{m_l} \beta_j^{h_l} \skm w_j^{h_l}$ we know
  $w_j^{h_l} \predeq \hterm(s_{\ell})$, as the $h_l$ arise from
  the reduction sequence $g_i \rmult w_{\ell} \red{*}{\myr}{p}{G} 0$,
  and hence $\hterm(w_j^{h_l} \rmult v) \predeq \hterm(s_{\ell} \rmult v) = T_p$.
 \begin{enumerate}
 \item For $l = i$ we get $g_i \rmult q_i = g_i \rmult (p_i + \alpha_j^{p_i} \skm (h_i - w_{\ell}) \rmult v)
        = g_i \rmult p_i + \alpha_j^{p_i} \skm g_i \rmult h_i \rmult v - \alpha_j^{p_i} \skm g_i \rmult w_{\ell} \rmult v$
        and as $\hterm(g_i \rmult t_j^{p_i}) = \hterm(g_i \rmult w_{\ell} \rmult v)$ and the resulting monomials
        add up to zero we get $\max \{ \hterm(g_i \rmult w_j^{h_i}) \mid 1 \leq j \leq m_i \} \leq T_p$.
 \item For $l \neq i$ we get $g_l \rmult q_l = g_l \rmult (p_l + \alpha_j^{p_i} \skm h_l \rmult v) 
        = g_l \rmult p_l  + \alpha_j^{p_i} \skm g_l \rmult h_l \rmult v$ and 
        $\max \{ \hterm(g_i \rmult w_j^{h_l}) \mid 1 \leq j \leq m_l \} \predeq T_p$ as well as
        $\max \{ \hterm(g_i \rmult w_j^{h_l}) \mid 1 \leq j \leq m_l \} \predeq T_p$.
 \end{enumerate}
 Hence while still in one of the cases we must have $T_q = T_p$, the element $g_i \rmult t_j^{p_i}$ is replaced by the sum
  $\sum_{l=1}^n g_l \rmult h_l \rmult v$ where the $h_l$ arise from the reduction sequence 
  $s_{\ell} \red{*}{\myr}{}{G} 0$.
 Let $h_l = \sum_{j=1}^{k_l} \alpha_j^{h_l} \skm t_j^{h_l}$. 
 Since $s_{\ell}$ is stable, for all elements $g_l \rmult t_j^{h_l}$ involved in the reduction of the head term
  of $s_{\ell}$ we know $\hterm(g_l \rmult t_j^{h_l} \rmult v) = \hterm(g_l) \rmult t_j^{h_l} \rmult v = T_p$
  and no other elements result in this term.
 Hence $K_q < K_p$ and ${\bf q}$ is smaller than ${\bf p}$.
\item Let us now assume there are $1 \leq i_1,i_2 \leq n$, $1 \leq j_1 \leq n_{i_1}$,
  $1 \leq j_2 \leq n_{i_2}$ such that $\hterm(g_{i_1} \rmult t_{j_1}^{p_{i_1}})  
  = \hterm(g_{i_1}) \rmult t_{j_1}^{p_{i_1}} =   T_p = \hterm(g_{i_2}) \rmult t_{j_2}^{p_{i_2}}
  = \hterm(g_{i_2} \rmult t_{j_2}^{p_{i_2}})$.
 Moreover, we assume $\hc(g_{i_1}) \geq \hc(g_{i_2}) > 0$ and $\hc(g_{i_1}) = \alpha \skm \hc(g_{i_2}) + \beta$,
  $\alpha, \beta \in \z$, $0 \leq \beta < \hc(g_{i_2})$.
 Let $h_{\ell_2} \in C_{i_1,i_2}$ such that for
  the corresponding s-polynomial $p = g_{i_1} \rmult t_{j_1}^{p_{i_1}} - \alpha \skm g_{i_2} \rmult t_{j_2}^{p_{i_2}}$
  we have $p = h_{\ell_2} \rmult v$ and
  $h_{\ell_2} = g_{i_1} \rmult u_1 - g_{i_2} \rmult (\alpha \skm u_2)$.
 Since we have  a vector ${\bf a}_{{i_1},{i_2}}^{\ell_2}$ corresponding to  $h_{\ell_2}$,
  we can define a new solution ${\bf q} = {\bf p} + \alpha_{j_1}^{p_{i_1}}  \skm {\bf a}_{{i_1},{i_2}} \rmult v$ with
 \begin{eqnarray}
 q_{i_1} & = & p_{i_1} + \alpha_{j_1}^{p_{i_1}} \skm (h_{i_1} - u_1) \rmult v \nonumber \\
 q_{i_2} & = & p_{i_2} + \alpha_{j_1}^{p_{i_1}} \skm (h_{i_2} + \alpha \skm u_2) \rmult v \nonumber \\
 q_l & = & p_l + \alpha_{j_1}^{p_{i_1}} \skm h_l \rmult v \mbox{ for }l \neq i,j.\nonumber
 \end{eqnarray}
 It remains to show that this solution indeed is smaller.
 To do this we examine the  multiples $g_l \rmult t_j^{q_l}$ for all $1 \leq l \leq n$, $1 \leq j \leq m_l$
  where $q_l = \sum_{j=1}^{m_l} \alpha_j^{q_l} \skm t_j^{q_l}$.
 Let $h_l = \sum_{j=1}^{k_l} \alpha_j^{h_l} \skm t_j^{h_l}$.
 Since the elements $g_l \rmult t_j^{h_l}$ arise from the reduction sequence
 $h_{\ell_2} \red{*}{\myr}{}{G} 0$ and the s-polynomial is stable we have
 additional information on how these elements affect the size of the new solution ${\bf q}$.
 Since $\hterm(g_l \rmult t_j^{h_l}) = \hterm(g_l) \rmult t_j^{h_l} \leq \hterm(h_{\ell_2})$
 we can conclude $\hterm(g_l \rmult t_j^{q_l}) \leq \hterm(h_{\ell_2}) \rmult v \predeq T_p$
 and we get the following boundaries:
 \begin{enumerate}
 \item For $l \neq i_1,i_2$ we get
       $g_l \rmult q_l = g_l \rmult p_l + \alpha_{j_1}^{p_{i_1}} \skm g_l \rmult h_l \rmult v$.
       This implies $\max \{\hterm(g_{l} \rmult t_j^{q_l}) \mid 1 \leq j \leq m_l\} \predeq T_p$.
 \item For $l = i_1$ we get
       $g_{i_1} \rmult q_{i_1} = g_{i_1} \rmult p_{i_1} + \alpha_{j_1}^{p_{i_1}} \skm g_{i_1}
       \rmult h_{i_1} \rmult v -  
       \alpha_{j_1}^{p_{i_1}} \skm g_{i_1} \rmult u_1 \rmult v$. 
       Since $\alpha_{j_1}^{p_{i_1}} \skm \hm(g_{i_1}) \rmult t_{j_1}^{p_{i_1}} = 
       \alpha_{j_1}^{p_{i_1}} \skm \hm(g_{i_1}) \rmult u_1 \rmult v$ we get 
       $\max \{\{\hterm(g_{i_1} \rmult t_j^{q_{i_1}}) \mid 1 \leq j \leq m_{i_1} \}
       \backslash \{ \hterm(g_{i_1}) \rmult t_{j_1}^{p_{i_1}}, \hterm(g_{i_1}) \rmult u_1 \rmult v \} \} \predeq T_p$.
 \item For $l = i_2$ we get
       $g_{i_2} \rmult q_{i_2} = g_{i_2} \rmult p_{i_2} + \alpha_{j_1}^{p_{i_1}}
       \skm g_{i_2} \rmult h_{i_2} \rmult v 
       + \alpha_{j_1}^{p_{i_1}} \skm g_{i_2} \rmult  \alpha \rmult u_2 \rmult v$.
       Again $\max \{\hterm(g_{i_1} \rmult t_j^{q_{i_1}}) \mid 1 \leq j \leq m_{i_1} \} \predeq T_p$.
 \end{enumerate}
Now in case $\beta =0$ we know that the equations are strict as then
 $\hterm(h_{\ell_2}) \rmult v \pred T_p$ holds.
Then either $T_q \pred T_p$ or $(T_q = T_p$ and $K_q < K_p)$.
If $\beta \neq 0$ we have to be more carefull and have to show that then $M_q \ll M_p$.
For the elements $g_l \rmult t_j^{h_l}$ arising from reducing the head of the s-polynomial we know
 that $g_l \rmult t_j^{h_l} \rmult v$ again has the same head coefficient as $g_l \rmult t_j^{h_l}$.
Now as $\hc(h_{\ell_2}) = \beta$, by the definition of our reduction relation
 we know that only $g_l$ with $\hc(g_l) \leq \beta$ are applicable.
Hence while two elements $\hc(g_{i_1}), \hc(g_{i_2})$ are removed from the multiset $M_p$
 only ones less equal to  $\beta < \hc(g_{i_2}) \leq \hc(g_{i_1})$ are added and hence the multiset
 becomes smaller.

Hence we find that in all cases above either 
$T_q \pred T_p$ or ( $T_q = T_p$ and $K_q < K_p$) or ($T_q = T_p$, $K_q = K_p$ and $M_q \ll M_p$).
Therefore, in all cases, we can reach a smaller solution and since our ordering on solutions is
 well-founded, or claim holds.
\end{enumerate}
\qed
\begin{corollary}~\\
{\sl
Let $\{ g_1, \ldots, g_n \}$ be a finite right Gr\"obner basis.
For not necessarily finite localizations
 $C_{i,j} \subseteq \spols(g_i,g_j)$ and $\sat(g_i)$ the not necessarily finite set
 of vectors ${\bf a}^{\ell_1}_{i,j}, {\bf b}_{i,\ell_2}$, $1 \leq i,j \leq n$, 
 $h_{\ell_1} \in C_{i,j}$, $s_{\ell_2} \in \sat(g_i)$
 forms a right generating set for all solutions of equation \ref{equation.homogeneous.G}.
\corend
}
\end{corollary}

The approach extends to systems of linear equations by using Gr\"obner bases in right modules.
A study of the situation for one-sided equations in integer monid and group rings can be found in \cite{Re00}.

%%% Local Variables: 
%%% mode: latex
%%% TeX-master: "testlauf.tex"
%%% TeX-master: "testlauf"
%%% End: 

%% file: conclusions.tex
\chapter{Conclusions}
The aim of this work was to give a guide for introducing reduction relations and Gr\"obner basis theory to algebraic structures.
We chose function rings as they allow a representation of their elements
 by formal sums.
This gives a natural link to those algebraic structures known in the literature where the Gr\"obner basis method works.
At the same time function rings provide enough flexibility to subsume these algebraic structures.

In the general setting of function rings we introduced the algebraic terms which are vital in Gr\"obner basis theory: head monomials, head terms, standard representations, standard bases, reduction relations and of course (weak) Gr\"obner bases.
Incorporating the technique of saturation we could give characterizations of Gr\"obner bases in terms of critical situations similar to the original approach.

We have established the theory first for right ideals in function rings 
over fields as this is the easiest setting.
This has been generalized to function rings over reduction rings - a very general setting.
Then in order to show how more knowledge on the reduction relation can be used to get deeper results on characterizing Gr\"obner bases, we have 
studied the special reduction ring $\z$, which is of interest in the literature.
The same approach has been applied to two-sided ideals in function rings
 with of course weaker results but still providing characterizations of 
Gr\"obner bases.

Important algebraic structures where the Gr\"obner basis method has 
been successfully applied in the literature have been outlined in the setting of function rings.
It has also been shown how special applications from Gr\"obner basis theory in polynomial rings can be lifted to function rings.

What remains to be done is to find out if this approach can be extended 
to function rings allowing infinite formal sums as elements.
Such an extension would allow to subsume the work of Mora et. al. on power series which resulted in the tangent cone algorithm.
These rings are covered by graded structures as defined by Apel in his habilitation (\cite{Ap98}), by monomial structures as defined by Pesch in
 his PhD Thesis (\cite{Pe97}) and by Mora in ``The Eigth variation'' (on
 Gr\"obner bases).
However, these approaches require admissible orderings and hence do not 
cover general monoid rings.
%%% Local Variables: 
%%% mode: latex
%%% TeX-master: "testlauf"
%%% End: 